\pdfoutput=1
\RequirePackage{silence}
\documentclass[a4paper,svgnames]{amsart}
\usepackage[hmarginratio={1:1},vmarginratio={1:1},lmargin=60.0pt,tmargin=60.0pt]{geometry}

\usepackage{booktabs}
\allowdisplaybreaks

\usepackage[T1]{fontenc}
\usepackage{latexsym,exscale,amsfonts,amssymb,mathtools}
\usepackage{amsmath,amsthm,amsfonts,amssymb,amscd,textcomp,bbm}
\usepackage{mathrsfs,stackrel}
\usepackage{etoolbox}
\usepackage{tikz,tikz-cd}
\usetikzlibrary{matrix,arrows.meta,decorations.markings,shapes.multipart,decorations.pathreplacing,decorations.pathmorphing}
\usepackage{aliascnt}
\usepackage{ytableau}
\usepackage{xparse}
\usepackage{graphicx}
\usepackage[export]{adjustbox}

\usepackage{cite}

\usepackage[normalem]{ulem}
\usepackage{listings}

\usepackage{array}
\newcolumntype{C}{>{$}c<{$}}

\usepackage{xcolor,colortbl}

\definecolor{mygray}{gray}{0.6}
\definecolor{mygraydark}{gray}{0.4}
\definecolor{mygraylight}{gray}{0.85}
\definecolor{spinach}{RGB}{46,139,87}
\definecolor{tomato}{RGB}{255,99,71}
\definecolor{orchid}{RGB}{143,40,194}
\definecolor{neonblue}{RGB}{0,30,255}
\definecolor{pumpkin}{RGB}{224,180,80}
\definecolor{citron}{RGB}{190,180,90}

\definecolor{macaroniandcheese}{RGB}{255,185,123}
\definecolor{cybergrape}{RGB}{88,66,124}
\definecolor{cottoncandy}{RGB}{255,166,201}

\definecolor{lava}{RGB}{207,16,32}
\definecolor{cream}{RGB}{255,253,208}
\definecolor{verdigris}{RGB}{67,179,174}
\definecolor{Black}{RGB}{0,0,0}
\definecolor{mydarkblue}{RGB}{10,10,170}
\definecolor{darkspinach}{RGB}{20,70,20}
\definecolor{darktomato}{RGB}{155,40,30}
\definecolor{darkorchid}{RGB}{50,10,100}
\definecolor{darklava}{RGB}{150,8,16}

\usepackage{todonotes}

\usepackage{enumitem}
\setlist[enumerate]{itemsep=0.15cm,label=\emph{\upshape(\alph*)}}
\setlist[enumerate,2]{itemsep=0.15cm,label=\emph{\upshape(\roman*)}}
\setlist[enumerate,3]{itemsep=0.15cm,label=\emph{\upshape(\Alph*)}}

\let\emph\relax
\DeclareTextFontCommand{\emph}{\em}

\renewcommand{\dots}{\text{...}}
\renewcommand{\vdots}{\rotatebox{90}{\text{...}}}

\newcommand{\placeholder}{{}_{-}}

\usepackage[all]{xy}
\SelectTips{cm}{}

\usepackage{tikz}
\tikzstyle{densely dotted}=[dash pattern=on \pgflinewidth off .5pt]
\tikzset{
anchorbase/.style={baseline={([yshift=-0.5ex]current bounding box.center)}},
tinynodes/.style={font=\tiny,text height=0.25ex,text depth=0.05ex},
smallnodes/.style={font=\scriptsize, text height=0.75ex, text depth=0.15ex},
usual/.style={line width=1.0,color=black},
usual2/.style={line width=1.0,color=brown},
usuald/.style={line width=1.0,color=magenta,densely dotted},
soergelone/.style={line width=1.0,color=neonblue},
soergeltwo/.style={line width=1.0,color=tomato,densely dotted},
soergelthree/.style={line width=1.0,color=spinach,densely dashed},
soergelfour/.style={line width=1.0,color=orchid,densely dashdotted},
soergelfive/.style={cyan, decoration={snake, amplitude=0.2mm, segment length=3mm}, decorate, thick},
mor/.style={line width=0.75,color=darkorchid,fill=cream},
rex/.style={line width=0.75,color=darkorchid,fill=orchid},
}

\newcommand{\drawtwomvalentvertex}[5]{
\begin{scope}[shift={(#4,#5)}]
\pgfmathsetmacro{\xend}{0.1*#1}
\pgfmathsetmacro{\yend}{0.05*(#1-1)}
\foreach \n in {1,...,#1}{
\ifodd\n
\pgfmathsetmacro\yoneStart{0} 
\pgfmathsetmacro\ytwoStart{0.1*(#1-1)} 
\else
\pgfmathsetmacro\yoneStart{0.1*(#1-1)} 
\pgfmathsetmacro\ytwoStart{0} 
\fi
\draw[#2] (0.2*\n-0.1,\yoneStart) to (\xend,\yend);   
\draw[#3] (0.2*\n-0.1,\ytwoStart) to (\xend,\yend);   
}
\end{scope}
}

\usepackage{xparse}
\NewDocumentCommand{\diagrammaticmorphism}{ m m m +m }{%
\begin{scope}[shift={(#1,#2)}]
\foreach \var [count=\i] in {#4} {
\ifnum\var=0
\draw[usual] (0.2*\i-0.1,0) -- ++(0,0.1*#3);
\else
\ifnum\var=1
\draw[soergelone] (0.2*\i-0.1,0) -- ++(0,0.1*#3);
\else
\ifnum\var=2
\draw[soergeltwo] (0.2*\i-0.1,0) -- ++(0,0.1*#3);
\else
\ifnum\var=3
\draw[soergelthree] (0.2*\i-0.1,0) -- ++(0,0.1*#3);
\else
\ifnum\var=4
\draw[soergelfour] (0.2*\i-0.1,0) -- ++(0,0.1*#3);
\fi
\fi
\fi
\fi
\fi
}
\end{scope}
}

\newcommand{\rotcircleCCW}[3]{%
\draw[thick, decoration={markings,
mark=at position 0.17 with {\arrow[scale=#3,rotate=-10]{>}},
mark=at position 0.5 with {\arrow[scale=#3,rotate=-10]{>}},
mark=at position 0.83 with {\arrow[scale=#3,rotate=-10]{>}}},
postaction=decorate] (#1,#2) circle [radius=0.2*#3];
}

\newcommand{\varrotcircleCW}[3]{%
\draw[thick, decoration={markings,
mark=at position 0.33 with {\arrow[scale=#3,rotate=190]{>}},
mark=at position 0.66 with {\arrow[scale=#3,rotate=190]{>}},
mark=at position 1 with {\arrow[scale=#3,rotate=190]{>}}},
postaction=decorate] (#1,#2) circle [radius=0.2*#3];
}

\newcommand{\varrotcircleCCWforflip}[3]{%
\draw[thick, decoration={markings,
mark=at position 0.166 with {\arrow[scale=#3,rotate=10]{>}},
mark=at position 0.5 with {\arrow[scale=#3,rotate=10]{>}},
mark=at position 0.833 with {\arrow[scale=#3,rotate=10]{>}}},
postaction=decorate] (#1,#2) circle [radius=0.2*#3];
}

\newcommand{\varrotcircleCWforflip}[3]{%
\draw[thick, decoration={markings,
mark=at position 0.33 with {\arrow[scale=#3,rotate=-190]{>}},
mark=at position 0.66 with {\arrow[scale=#3,rotate=-190]{>}},
mark=at position 1 with {\arrow[scale=#3,rotate=-190]{>}}},
postaction=decorate] (#1,#2) circle [radius=0.2*#3];
}

\newcommand{\trimergeCW}[4]{%
\begin{tikzpicture}[anchorbase,scale=1]
\draw[mor] (0,0) rectangle (1,0.5) node[black,pos=0.5]{$#1$};
\draw[mor] (2,0) rectangle (3,0.5) node[black,pos=0.5]{$#2$};
\draw[mor] (1,1.5) rectangle (2,2) node[black,pos=0.5]{$#3$};
\draw[usual] (0.5,0.5) to (1.5,1) to (1.5,1.5);
\draw[usual] (2.5,0.5) to (1.5,1);
\varrotcircleCW{1.5}{1}{#4}
\end{tikzpicture}
}
\newcommand{\trimergeCCW}[4]{%
\begin{tikzpicture}[anchorbase,scale=1]
\draw[mor] (0,0) rectangle (1,0.5) node[black,pos=0.5]{$#1$};
\draw[mor] (2,0) rectangle (3,0.5) node[black,pos=0.5]{$#2$};
\draw[mor] (1,1.5) rectangle (2,2) node[black,pos=0.5]{$#3$};
\draw[usual] (0.5,0.5) to (1.5,1) to (1.5,1.5);
\draw[usual] (2.5,0.5) to (1.5,1);
\rotcircleCCW{1.5}{1}{#4}
\end{tikzpicture}
}

\newcommand{\trisplitCW}[4]{%
\begin{tikzpicture}[anchorbase,scale=1]
\begin{scope}[yscale=-1]
\draw[mor] (0,0) rectangle (1,0.5) node[black,pos=0.5]{$#1$};
\draw[mor] (2,0) rectangle (3,0.5) node[black,pos=0.5]{$#2$};
\draw[mor] (1,1.5) rectangle (2,2) node[black,pos=0.5]{$#3$};
\draw[usual] (0.5,0.5) to (1.5,1) to (1.5,1.5);
\draw[usual] (2.5,0.5) to (1.5,1);
\varrotcircleCCWforflip{1.5}{1}{#4}
\end{scope}
\end{tikzpicture}
}

\newcommand{\trisplitCCW}[4]{%
\begin{tikzpicture}[anchorbase,scale=1]
\begin{scope}[yscale=-1]
\draw[mor] (0,0) rectangle (1,0.5) node[black,pos=0.5]{$#1$};
\draw[mor] (2,0) rectangle (3,0.5) node[black,pos=0.5]{$#2$};
\draw[mor] (1,1.5) rectangle (2,2) node[black,pos=0.5]{$#3$};
\draw[usual] (0.5,0.5) to (1.5,1) to (1.5,1.5);
\draw[usual] (2.5,0.5) to (1.5,1);
\varrotcircleCWforflip{1.5}{1}{#4}
\end{scope}
\end{tikzpicture}
}

\usetikzlibrary{cd}
\usetikzlibrary{decorations}
\usetikzlibrary{decorations.markings}
\usetikzlibrary{decorations.pathreplacing}
\usetikzlibrary{decorations.pathmorphing}
\usetikzlibrary{arrows.meta,shapes,positioning,matrix,calc}
\usetikzlibrary{shapes.callouts}
\tikzstyle directed=[postaction={decorate,decoration={markings,mark=at position #1 with {\arrow[line width=0.3mm, black]{>}}}}]
\tikzstyle rdirected=[postaction={decorate,decoration={markings,mark=at position #1 with {\arrow[line width=0.3mm, black]{<}}}}]
\tikzstyle marked=[postaction={decorate,decoration={markings,
mark=at position #1 with {\fill[black] (0,0) circle (.055cm);}}}]
\tikzstyle markedone=[postaction={decorate,decoration={markings,
mark=at position #1 with {\fill[neonblue] (0,0) circle (.055cm);}}}]
\tikzstyle markedtwo=[postaction={decorate,decoration={markings,
mark=at position #1 with {\fill[tomato] (0,0) circle (.055cm);}}}]
\tikzstyle markedthree=[postaction={decorate,decoration={markings,
mark=at position #1 with {\fill[spinach] (0,0) circle (.055cm);}}}]
\tikzstyle markedfour=[postaction={decorate,decoration={markings,
mark=at position #1 with {\fill[orchid] (0,0) circle (.055cm);}}}]

\let\emph\relax
\DeclareTextFontCommand{\emph}{\bfseries\em}

\allowdisplaybreaks

\renewcommand{\dots}{\text{...}}
\renewcommand{\vdots}{\raisebox{-0.05cm}{\rotatebox{90}{\text{...}}}}

\usepackage{stmaryrd}
\newcommand{\dsum}{\inplus}

\newcommand{\munit}{\mathbbm{1}}

\newcommand{\catstuff}[1]{\mathbf{#1}}

\newcommand{\setstuff}[1]{#1}

\newcommand{\obstuff}[1]{#1}
\newcommand{\morstuff}[1]{#1}

\DeclareMathOperator{\End}{End}
\DeclareMathOperator{\Hom}{Hom}

\newcommand{\Sym}{\setstuff{Sym}}
\DeclareMathOperator{\len}{len}
\DeclareMathOperator{\Rex}{Rex}
\DeclareMathOperator{\id}{id}
\DeclareMathOperator{\rank}{rank}
\DeclareMathOperator{\proj}{proj}
\DeclareMathOperator{\incl}{incl}
\DeclareMathOperator{\Indec}{Indec}

\newcommand{\C}{\mathbb{C}}
\newcommand{\Z}{\mathbb{Z}}
\newcommand{\R}{\mathbb{R}}
\newcommand{\F}{\mathbb{F}}

\newcommand{\K}{\mathbb{K}}

\newcommand{\N}{\mathbb{N}}

\newcommand{\rdual}[1]{{#1}{}^{\star}}

\newcommand{\idcoeff}[1]{\kappa(#1)}
\newcommand{\un}{\underline}
\newcommand{\onetensor}{1^{\otimes}}
\newcommand{\ot}{\otimes}
\DeclareMathOperator{\LIF}{LIF}
\DeclareMathOperator{\ptr}{tr}
\DeclareMathOperator{\something}{SS}
\newcommand{\DD}{\mathbb{D}}
\DeclareMathOperator{\Kar}{Kar}
\newcommand{\ep}{\varepsilon}

\makeatletter
\newcommand{\uset}[3][0ex]{%
\mathrel{\mathop{#3}\limits_{
\vbox to#1{\kern-2\ex@
\hbox{$\scriptscriptstyle#2$}\vss}}}}
\makeatother
\newcommand{\muless}{\uset[1ex]{1}{<}}

\newcommand{\notrightnow}[1]{}

\numberwithin{equation}{section}
\usepackage[hypertexnames=false]{hyperref}
\hypersetup{
pdftoolbar=true,        
pdfmenubar=true,        
pdffitwindow=false,     
pdfstartview={FitH},    
pdftitle={Idempotents, traces, and dimensions in Hecke categories},    
pdfauthor={Ben Elias, Liam Rogel and Daniel Tubbenhauer},     
pdfsubject={},   
pdfcreator={Ben Elias, Liam Rogel and Daniel Tubbenhauer},   
pdfproducer={Ben Elias, Liam Rogel and Daniel Tubbenhauer}, 
pdfkeywords={}, 
pdfnewwindow=true,      
colorlinks=true,       
linkcolor=mydarkblue,          
citecolor=teal,        
filecolor=magenta,      
urlcolor=orchid,          
linkbordercolor=lava,
citebordercolor=teal,
urlbordercolor=orchid,  
linktocpage=true
}

\usepackage{aliascnt}

\def\makeautorefname#1#2{\csdef{#1autorefname}{#2}}
\makeautorefname{section}{Section}%
\makeautorefname{subsection}{Section}%
\makeautorefname{subsubsection}{Section}%

\newcommand\NewTheorem[2][\relax]{%
\newaliascnt{#2}{equation}%
\newtheorem{#2}[#2]{#2}%
\aliascntresetthe{#2}%
\expandafter\def\csname#2autorefname\endcsname{#2}%
\ifx#1\relax\else
\AtEndEnvironment{#2}{\null\hfill$#1$}%
\newaliascnt{#2*}{equation}%
\newtheorem{#2*}[#2*]{#2}%
\aliascntresetthe{#2*}%
\expandafter\def\csname#2*autorefname\endcsname{#2}%
\fi
}

\def\equationautorefname~#1\null{(#1)\null}

\numberwithin{equation}{subsection}

\NewTheorem{Proposition}
\NewTheorem{Theorem}
\NewTheorem{Corollary}
\NewTheorem{Lemma}
\theoremstyle{definition}
\NewTheorem[\Diamond]{Definition}
\NewTheorem[\Diamond]{Notation}
\NewTheorem[\Diamond]{Example}
\NewTheorem[\Diamond]{Examples}
\theoremstyle{remark}
\NewTheorem[\Diamond]{Remark}
\NewTheorem[\Diamond]{Assumption}

\begin{document}
\vbadness=10001
\hbadness=10001
\overfullrule=1mm
\title[Idempotents, traces, and dimensions in Hecke categories]{Idempotents, traces, and dimensions in Hecke categories}
\author[Ben Elias, Liam Rogel and Daniel Tubbenhauer]{Ben Elias, Liam Rogel and Daniel Tubbenhauer}

\address{B.E.: Department of Mathematics, Fenton Hall, Room 210, 
University of Oregon, Eugene, OR 97403, United States, \href{https://pages.uoregon.edu/belias/}{pages.uoregon.edu/belias/}}
\email{belias@uoregon.edu}

\address{L.R.: Department of Mathematics, Office 48-420, University of Kaiserslautern-Landau (RPTU), Gottlieb-Daimler-Straße 47, 67663 Kaiserslautern, Germany}
\email{rogel@mathematik.uni-kl.de}

\address{D.T.: The University of Sydney, School of Mathematics and Statistics F07, Office Carslaw 827, NSW 2006, Australia, \href{http://www.dtubbenhauer.com}{www.dtubbenhauer.com}, https://orcid.org/0000-0001-7265-5047}
\email{daniel.tubbenhauer@sydney.edu.au}

\begin{abstract}
We explain how to compute idempotents that correspond to the indecomposable objects in the Hecke category. Closed formulas are provided for some common coefficients that appear in these idempotents. We also explain how to compute categorical dimensions in the asymptotic Hecke category. In many cases, we reduce this to a computation of a partial trace and give recursive formulas for some common partial traces. In the sequel, we apply this technology and perform additional (computer) calculations to complete the description of the asymptotic Hecke category for finite Coxeter groups in all but three cells.
\end{abstract}


\maketitle
%



\section{Introduction}\label{section:intro}

The primary accomplishment of this paper is to develop the theory of idempotents and partial traces in Hecke categories, and to explain common situations when such computations have
closed formulas or can be found by a simple recursive algorithm. This is a topic of its own intrinsic merit, see \autoref{intro:idempandtrace}, but our original motivation was to compute
categorical dimensions in asymptotic Hecke categories, and this is where we begin the introduction. In the sequel \cite{ElRoTu-ah-2} we use these techniques and additional computer
calculations to describe the asymptotic Hecke category for most of the cells which were not yet known (all of which are in types $H_3$ and $H_4$).

\subsection{Asymptotic Hecke algebras and categories}\label{subsection:introcells}

Let $W$ be a Coxeter group. Kazhdan and Lusztig \cite{KaLu-reps-coxeter-groups} introduced the study of \emph{cell theory}, which can be concisely defined as follows. Associated to $W$ is its Hecke algebra, a
$\Z[v,v^{-1}]$-algebra which deforms the group algebra $\Z[W]$, with a distinguished basis $\{b_w\}_{w \in W}$ now called the \emph{Kazhdan--Lusztig basis}. Cell theory studies
\emph{based ideals}, ideals which are spanned by subsets of the basis. We say that $x \le_{LR} y$ if every based two-sided ideal containing $b_y$ also contains $b_x$. One defines
$\le_R$ and $\le_{L}$ similarly, using right or left ideals. The equivalence classes under $\le_{LR}$ are \emph{two-sided cells}, and similarly for \emph{right cells} and
\emph{left cells} using $\le_R$ and $\le_{L}$. The Hecke algebra has a natural filtration by based two-sided ideals, indexed by the set of two-sided cells.

A motivation for studying based ideals and cells comes from categorification. If an algebra $A$ is categorified by a monoidal additive Krull--Schmidt category, the indecomposable
objects descend to a basis for $A$, and the tensor ideals descend to based ideals in $A$. The Hecke algebra is famously categorified by the \emph{Hecke category}, which arises in
many different guises (e.g. perverse sheaves on flag varieties, translation functors acting on category $\mathcal{O}$, Soergel bimodules, diagrammatic algebra) and is ubiquitous in
representation theory.

The Hecke category is graded, with grading shift categorifying multiplication by the parameter $v$. (We work throughout with the $\R$-linear
version of the Hecke category.)
In the associated graded algebra (for the two-sided cell filtration), the structure coefficients in each two-sided cell $\mathcal{J}$ are Laurent polynomials with exponents between $v^{-a(\mathcal{J})}$ and $v^{+a(\mathcal{J})}$, where $a(\mathcal{J})$ is Lusztig's a-function. We can renormalize the basis by $v^{-a(\mathcal{J})}$ so that structure coefficients live in $\Z[v^{-1}]$, and then specialize $v^{-1} \mapsto 0$. The result is the \emph{asymptotic Hecke algebra} or \emph{$J$-ring} $A_{\mathcal{J}}$ associated to this two-sided cell, an associative unital $\Z$-algebra which Lusztig proves is semisimple after specialization to any field; see e.g. \cite{Lu-cells-affine-weyl-2}. One of the
advantages of $A_{\mathcal{J}}$ is that it is usually much simpler than the Hecke algebra, yet one can classify the simple representations of the Hecke algebra in terms of simple
representations of various $A_{\mathcal{J}}$, see e.g. \cite{Lu-leading-coeff-hecke}.

There are also categorifications of the asymptotic Hecke algebras $A_{\mathcal{J}}$ by semisimple monoidal categories $\catstuff{A}_{\mathcal{J}}$, the \emph{asymptotic Hecke
categories}. They were first constructed for Weyl groups using perverse sheaves in \cite{Lu-cells-tensor-cats}. After the Soergel conjecture was proven in \cite{ElWi-hodge-sbim}, the same
constructions were ported to Soergel bimodules and adapted to arbitrary Coxeter groups \cite[Section 18.15]{Lu-hecke-book} (this material can be found in the updated version of
\cite{Lu-hecke-book} on the arXiv, but not in the published version of the book). In \cite{ElWi-relative-lefschetz} it was proven that the asymptotic Hecke categories are pivotal,
hence are fusion categories.

The intersection of a left cell and a right cell (in the same two-sided cell) is often called an \emph{$H$-cell}. These are typically neither left nor right nor two-sided cells.
There is a bijection between left and right cells induced by taking inverses, called \emph{adjunction}. (In the literature, adjunction is often called duality, though we
reserve the name ``duality'' for something else.) Of particular note is the $H$-cell given by the intersection of a left cell and its adjoint right cell, which is called a
\emph{diagonal ($H$-)cell}. Associated to a diagonal cell $\mathcal{H} \subset \mathcal{J}$ we have a ``subalgebra'' $A_{\mathcal{H}} \subset A_{\mathcal{J}}$, and we can think of
$A_{\mathcal{H}}$ as being even simpler than $A_{\mathcal{J}}$. The inclusion map sends the unit of $A_{\mathcal{H}}$ to an idempotent in $A_{\mathcal{J}}$, so technically $A_{\mathcal{H}}$ is not a subalgebra but an idempotent truncation of $A_{\mathcal{J}}$.

\begin{Example} \label{Ex:typeAintroMorita} Let $W$ be the symmetric group of type $A_{n}$, and let $\mathcal{J}$ be any two-sided cell. Then $A_{\mathcal{J}}$ is isomorphic to a
matrix algebra, where each element of $\mathcal{J}$ corresponds to a matrix entry. Each row is a right cell, and each column is a left cell, and each $H$-cell is a
singleton. A diagonal cell $\mathcal{H}$ is a single matrix entry on the diagonal, and $A_{\mathcal{H}}$ is the ground field $\K$. By the usual Morita equivalence between a matrix
algebra and its ground field, the representation theory of $A_{\mathcal{J}}$ is determined by that of $A_{\mathcal{H}}$. \end{Example}

Henceforth the only $H$-cells we discuss in this paper are diagonal cells. Outside of type $A$ and dihedral groups, it is common for different diagonal $H$-cells in the same two-sided cell to have different sizes, and hence non-isomorphic asymptotic Hecke algebras.

\begin{Example}\label{E:H4} In type $H_4$ there is a \emph{middle cell} $\mathcal{J}_{6}$, a two-sided cell which contains 24 diagonal cells. Of these, eight have size $14$, ten have size $18$, and six have size $24$. This leads to three different isomorphism classes of algebras $A_{\mathcal{H}}$. Using \cite{Al-j-h4} and a computer calculation one can determine the dimension of the center of these semisimple algebras and verify that the algebras of dimension $18$ and $24$ are Morita equivalent, but are not Morita equivalent to the $14$-dimensional algebra.
\end{Example}

When $A_{\mathcal{J}}$ is Morita equivalent to $A_{\mathcal{H}}$ for a diagonal cell $\mathcal{H} \subset \mathcal{J}$, we call this phenomenon \emph{(algebraic) $H$-reduction}. 
Outside of type $A$ (and dihedral groups of type $I_2(m)$ for $m$ even) $H$-reduction typically fails, see e.g. \cite[Theorem 2A.17, Theorem 3A.5]{Tu-sandwich-cellular}. Sometimes the diagonal cells are not even Morita equivalent to each other, as in \autoref{E:H4}, and this phenomenon is not well-understood.

There is also a fusion category $\catstuff{A}_{\mathcal{H}}$ associated to each diagonal cell $\mathcal{H} \subset \mathcal{J}$. Now we examine the categorical analog of
$H$-reduction. A monoidal category can act (by endofunctors) on another category, making that category a \emph{$2$-representation} of the monoidal category. A good analog of a
simple module in this context is a simple transitive 2-representation, or 2-simple for short \cite{MaMi-cell-2-reps}. In \cite{MaMaMiTuZh-soergel-2reps} it was proven for finite
Coxeter groups that asymptotic Hecke categories give enough information to classify all 2-simples of the Hecke category. In a major improvement over the decategorified setting, the
category $\catstuff{A}_{\mathcal{H}}$ is always 2-Morita equivalent to $\catstuff{A}_{\mathcal{J}}$, for any diagonal cell $\mathcal{H} \subset \mathcal{J}$ (see \cite[Theorem
4.32]{MaMaMiTuZh-bireps}). That is, unlike the decategorified setting, categorical $H$-reduction always holds.

Also remarkable is that the Drinfeld centers of $\catstuff{A}_{\mathcal{J}}$ (and hence also $\catstuff{A}_{\mathcal{H}}$) are modular categories, by the above mentioned properties of these categories and standard results (as, for example, in \cite[Corollary 8.20.14]{EtGeNiOs-tensor-categories}). Moreover, the corresponding S-matrices are believed to coincide with the so-called Fourier matrices of \cite{Lu-exotic-fourier} (see also the appendix of
this paper by Malle) and \cite{BrMa-hecke}.

\begin{Example} Continuing \autoref{Ex:typeAintroMorita}, for any $\mathcal{H} \subset \mathcal{J}$, the category $\catstuff{A}_{\mathcal{H}}$ is the category of vector spaces, and $\catstuff{A}_{\mathcal{J}}$ is a matrix category with entries in the category of vector spaces. In this case \cite[Theorem 4.32]{MaMaMiTuZh-bireps} is a direct categorification of Morita equivalence. \end{Example}

To summarize, the categorical representation theory of Hecke categories can be reduced in some sense to the categorical representation theory of asymptotic Hecke categories for diagonal $H$-cells, which motivates the study of these categories $\catstuff{A}_{\mathcal{H}}$. For finite Coxeter groups, up to 2-Morita equivalence, we need only study one diagonal cell in each two-sided cell.

\subsection{Expliciting asymptotic Hecke algebras and categories}\label{subsection:explicitasymp}

Having discussed the general algebraic properties of asymptotic Hecke algebras and categories, let us focus on what is known about them. What algebras or categories are they? Note
that the rings $A_{\mathcal{J}}$ are often much less accessible than their smaller counterparts $A_{\mathcal{H}}$.

It turns out that $A_{\mathcal{H}}$ has a pleasant description for all diagonal cells in all finite Coxeter types: for finite Weyl types see \cite{Lu-characters-reductive-groups} and
\cite{Lu-leading-coeff-hecke}; for the remaining cases see \cite{Al-left-cells-h4,duCl-positivity-finite-hecke,Al-j-h4}. (There are even results beyond the finite cases, see e.g. \cite{MR4034792,MR4633004}.) Except for the case where $\mathcal{J} = \mathcal{J}_{6}$ is the infamous middle cell
in type $H_{4}$, the algebras $A_{\mathcal{H}}$ are isomorphic to well-known and uncomplicated algebras, such as group rings of elementary abelian $2$-groups. See \cite{ElRoTu-ah-2} for details. The middle cell $\mathcal{J}_{6}$ in type $H_{4}$ has three algebraically distinct diagonal cells, and their asymptotic Hecke algebras are not well-known, and
are famous examples of ``exotic'' semisimple algebras. The ``best thing one can hope for'' is an explicit description of their structure constants and their nontrivial subalgebras.
This was carried out in \cite{Al-j-h4}.

In addition to having a nice answer, the calculation of $A_{\mathcal{H}}$ can be fully automated. For example, in rank $\leq 8$, including all exceptional Coxeter types $F_{4}$, $E_{6}$, $E_{7}$, $E_{8}$, $H_{3}$, and $H_{4}$, the computation of $A_{\mathcal{H}}$ is completely computer verified by \cite{duCl-positivity-finite-hecke}. Between explicit identification and computer verification, altogether we have a quite complete picture of the structure of the asymptotic Hecke algebra. 

In contrast, our understanding of asymptotic Hecke categories is far less complete, especially outside of classical types.

When $A_{\mathcal{H}}$ is a well-known and uncomplicated algebra one can also hope that $\catstuff{A}_{\mathcal{H}}$ is a well-known and uncomplicated fusion category, such as the
category $\catstuff{Vec}(G)$ of vector spaces graded by a finite group $G$, or a variant thereupon. For Weyl groups this is indeed the case. Bezrukavnikov--Finkelberg--Ostrik
\cite{BeFiOs-tensor-cats-III} used a geometric approach to compute $\catstuff{A}_{\mathcal{H}}$ for all cells in Weyl groups, with the exception of the cells $\mathcal{J}_{17}$ in
type $E_{7}$ and $\mathcal{J}_{13}$, $\mathcal{J}_{13}^{\prime}$ in type $E_{8}$ (in the notation of \cite{ElRoTu-ah-2}). Ostrik \cite{Os-tensor-exceptional-cells} went on to
compute the remaining $\catstuff{A}_{\mathcal{H}}$ by Lie theoretical methods. For classical groups $\catstuff{A}_{\mathcal{H}}\cong\catstuff{Vec}(G)$, where $G$ is an elementary
abelian $2$-group. In exceptional types one can also get $\catstuff{Vec}(G)$ for other groups $G$, or twists thereof in the aforementioned $\mathcal{J}_{17}$, $\mathcal{J}_{13}$ and
$\mathcal{J}_{13}^{\prime}$ cases.


Let us pause to note that while this classification is complete for Weyl groups, the approach is very involved and far from being hands-on. In particular, the approaches in the
literature do not seem to be amenable to computer calculations, and to the best of our knowledge, no computer verification has been done.

For non-crystallographic Coxeter groups one cannot use a geometric or Lie theoretic approach. Dihedral groups can be attacked directly, see \cite{El-two-color-soergel,El-q-satake,MaTu-soergel,RoTh-ah-cat}, but very little was known in types $H_3$ and $H_4$. In particular, the fusion categories associated to diagonal cells in the middle cell $\mathcal{J}_6$ in type $H_4$ have long been sought out as examples of interesting exotic fusion categories, but there were few techniques in the literature to compute their categorical structure.



The goal of this paper series is to rectify this situation. We provide general techniques for translating questions about asymptotic Hecke categories to computations in the
diagrammatic Hecke category, and techniques for performing those computations, or translating them to a setting which is accessible to computer assisted computation. This paper is focused on computational tools in the diagrammatic Hecke category itself. In the sequel, we use these diagrammatic methods to verify known results in various types, and compute new results in type $H_3$ and $H_4$. More specifically:
\begin{enumerate}[label=$\blacktriangleright$]

\item For the middle cell $\mathcal{J}_{6}$ in type $H_{4}$, for one diagonal cell $\mathcal{H}$ in each isomorphism class (e.g. one each of size 14, 18, and 24), we compute all the categorical dimensions in the category $\catstuff{A}_{\mathcal{H}}$, and identify all its sub-fusion categories. This is ``the best one can hope for.''

\item Aside from $\mathcal{J}_6$ and three other two-sided cells, we identify $\catstuff{A}_{\mathcal{H}}$ explicitly for at least one diagonal cell in all other two-sided cells in type $H_{4}$. In particular, this completes the classification of 2-simples of the Hecke category and the identification of the S-matrices (which do coincide with the Fourier matrices) in all but these cases.

\item We explain how these methods can be implemented on a machine. This implementation is partially complete (see \cite{ElRoTu-ah-2}), and most of our calculations are computer verified.
\end{enumerate}

\begin{Remark} Let $W$ be a finite Coxeter group with longest element $w_0$. If $\mathcal{J}$ is a two-sided cell, then the elements of the form $w_0 x$ for $x \in \mathcal{J}$ also form a two-sided cell, denoted $w_0 \mathcal{J}$ and called the \emph{$w_0$-dual} of $\mathcal{J}$. The three two-sided cells which we are not able to compute have distinguished involutions of high length (almost $\len(w_0)$). Consequently their $w_0$-duals have short distinguished involutions and are computationally accessible. Because it is true in all other known cases, one could conjecture that $w_0$-dual cells (resp. diagonal cells) have isomorphic asymptotic Hecke categories. If this were proven, it would obviate the need to study these three missing cells (and would reduce the work dramatically in other types as well). \end{Remark}

\subsection{Elaboration on categorical dimension}\label{subsection:introcatdim}

Existing information about the Grothendieck group and known classification results for fusion categories together imply (for most diagonal cells) that each particular fusion category
$\catstuff{A}_{\mathcal{H}}$ is one of a small list of possible options, which can be distinguished by the categorical dimensions of certain objects. Thus we focus our attention on the computation of categorical dimensions.

The Hecke category has a contravariant duality functor $\DD$, and (unrelatedly) every object $B$ has a biadjoint object $\rdual{B}$. The indecomposable objects in the Hecke category
are, up to isomorphism and grading shift, parametrized by $W$. We write $B_w$ for the self-dual indecomposable object associated to $w \in W$, and $B_w(k)$ for its grading shift (for
$k \in \Z$). We have $\DD(B_w) = B_w$ and $\rdual{B_w} \cong B_{w^{-1}}$. Morphism spaces in the category are $R$-bimodules, where $R$ is the polynomial ring associated to the
reflection representation $V$ of $W$, and the endomorphism ring of the monoidal identity $\End(\munit)$ is $R$.

The loose construction of the asymptotic Hecke category $\catstuff{A}_{\mathcal{J}}$ goes as follows. First, kill the tensor ideal associated to all cells lower than $\mathcal{J}$,
which amounts to killing the objects $B_w$ for $w$ in those cells (and killing morphisms which factor through them). Let $a = a(\mathcal{J})$ be the value of Lusztig's $a$-function.
Given elements $x$, $y$, and $z$ in $\mathcal{J}$, the tensor product $B_x \ot B_y$ decomposes into summands, and $B_z(k)$ only appears as a summand for $-a \le k \le a$. One defines a
new tensor product $\bullet$ where $B_x \bullet B_y$ only consists of those summands $B_w(a)$ with shift exactly $a$, which we call the \emph{$a$-degree submodule}. These summands are
then shifted by $(-a)$ so that they are self-dual, and by restricting to morphisms of degree zero, one obtains the asymptotic Hecke category. For clarity, we let $A_w$ denote the
object in the asymptotic Hecke category associated to $B_w$, and write $A_x \star A_y$ for the monoidal structure.

The main subtlety in this construction is that the $a$-degree submodule is canonically a subobject of $B_x \ot B_y$ but not canonically a summand; it has no canonical complement;
it has a canonical inclusion map but not a canonical projection map. This situation is resolved in \cite{ElWi-relative-lefschetz}, which uses relative Hodge theory to construct
projection maps given the inclusion maps and a choice of regular dominant weight $\rho \in V$. In this way the
$a$-degree submodule is made into an explicit summand of $B_x \ot B_y$, and one can compute with it diagrammatically (e.g. in the Karoubi envelope of the diagrammatic Hecke
category). To elaborate, for $f \in R$, let $M_f$ be the endomorphism of $B_x \ot B_y$ which multiplies in
the middle (at the tensor product) by $f$. The element $\rho^a \in R$ has degree $2a$, and one uses $M_{\rho^a}$ to identify projection maps in degree $-a$ with projection maps in degree $+a$.

In each diagonal cell $\mathcal{H}$ there is a unique \emph{Duflo involution} $d$, and $A_d$ is the monoidal identity of $\catstuff{A}_{\mathcal{H}}$. For any $x \in \mathcal{H}$
there is a one-dimensional space of morphisms $B_d \to B_x \ot B_{x^{-1}}$ in degree $-a$ (modulo lower cells), and any element $\incl = \incl_{d}^{x,x^{-1}}$ of this space is a valid
inclusion map. Let $p = p_{x,x^{-1}}^d := \DD(\incl) \circ M_{\rho^a}$, a morphism $B_x \ot B_{x^{-1}} \to B_d$ of degree $+a$. The composition $p \circ \incl$ is a scalar multiple of the identity map (modulo lower cells), and \cite{ElWi-relative-lefschetz} shows
that this scalar is a positive real number $\mu_x$, so $\proj = \proj_{x,x^{-1}}^d := \frac{1}{\mu_x} p$ is a valid projection map paired with $\incl_d^{x,x^{-1}}$.

By adjunction, we also get a morphism $\incl_x^{d,x}$ of degree $-a$. The composition $\DD(\incl_x^{d,x}) \circ M_{\rho^a} \circ \incl_x^{d_x}$ is a scalar multiple of $\id_{B_x}$
(modulo lower cells), a different positive real number $\lambda_x$. In the sequel we prove that
\begin{equation} \dim(A_x) = \mu_x / \lambda_x. \end{equation}

If we represent these morphisms of degree $-a$ as trivalent vertices, the statements above correspond to
\begin{gather} \label{introdimpicture}
\begin{tikzpicture}[anchorbase,scale=1]
\draw[mor] (1,-1) rectangle (2,-0.5) node[black,pos=0.5]{$x$};
\draw[mor] (0,0) rectangle (1,0.5) node[black,pos=0.5]{$d$};
\draw[mor] (2,0) rectangle (3,0.5) node[black,pos=0.5]{$x$};
\draw[mor] (1,1) rectangle (2,1.5) node[black,pos=0.5]{$x$};
\node[scale=0.7] at (1.5,0.25) {$\rho^{a}$};
\draw[usual] (0.75,0.5) to (1.5,0.75) to (1.5,1);
\draw[usual] (2.25,0.5) to (1.5,0.75);
\draw[usual] (0.75,0) to (1.5,-0.25) to (1.5,-0.5);
\draw[usual] (2.25,0) to (1.5,-0.25);
\end{tikzpicture}
=
\lambda_x \cdot
\begin{tikzpicture}[anchorbase,scale=1]
\draw[mor] (0,0) rectangle (1,0.5) node[black,pos=0.5]{$x$};
\end{tikzpicture},
\qquad
\begin{tikzpicture}[anchorbase,scale=1]
\draw[mor] (1,-1) rectangle (2,-0.5) node[black,pos=0.5]{$d$};
\draw[mor] (0,0) rectangle (1,0.5) node[black,pos=0.5]{$x$};
\draw[mor] (2,0) rectangle (3,0.5) node[black,pos=0.5]{$x^{-1}$};
\draw[mor] (1,1) rectangle (2,1.5) node[black,pos=0.5]{$d$};
\node[scale=0.7] at (1.5,0.25) {$\rho^{a}$};
\draw[usual] (0.75,0.5) to (1.5,0.75) to (1.5,1);
\draw[usual] (2.25,0.5) to (1.5,0.75);
\draw[usual] (0.75,0) to (1.5,-0.25) to (1.5,-0.5);
\draw[usual] (2.25,0) to (1.5,-0.25);
\end{tikzpicture}
=
\mu_x \cdot
\begin{tikzpicture}[anchorbase,scale=1]
\draw[mor] (0,0) rectangle (1,0.5) node[black,pos=0.5]{$d$};
\end{tikzpicture}
, \qquad
\dim(A_x) = \frac{\mu_x}{\lambda_x}.
\end{gather}

In the first picture of \eqref{introdimpicture} we enclose the polynomial $\rho^a$ within a certain bigon. When the polynomial $\rho^a$ is replaced with a polynomial $f$ of the same degree as $\rho^a$, the result is also a scalar multiple of the identity map, and we refer to the coefficient of the identity as $\lambda_x(f)$.

\subsection{What happens in this paper}\label{subsection:whathappens}

In order to actually compute these categorical dimensions, there are several intermediate steps:
\begin{enumerate} \item Compute the idempotents which describe the objects $B_d$ and $B_x$ as objects in the Hecke category.
\item Find a suitable morphism $\incl_d^{x,x^{-1}}$.
\item Compute the compositions above to find the numbers $\lambda_x$ and $\mu_x$.
\end{enumerate}

Idempotents have a recursive construction as linear combinations of morphisms involving previously-compu\-ted
idempotents, and the coefficients in this recursive formula are related to so-called \emph{local intersection forms}. For example, each reduced expression for an element $x$ gives a
recursive construction of an idempotent endomorphism (of a Bott-Samelson bimodule), if only one knew the local intersection forms ``along'' this expression. We elaborate further in \autoref{intro:idempandtrace}.

We give closed formulas for the most common family of local intersection forms, see \autoref{Thm:LIFrecursionWORKS}. It turns out that this family of local intersection forms is
determined by the position of each element within a coset for a dihedral group. This dramatically reduces the number of local intersection forms which a (human or electronic) computer
would need to determine. A reduced expression is called \emph{recursible} (see \autoref{def:recursible}) if every local intersection form along the expression is known from our closed
formula (there is a precise criterion). We are thus able to automate the computation of idempotents for elements with a recursible reduced expression.

For the remaining steps we focus on two special cases. In \autoref{subsection:simplifications} we assume that $x$ is larger than $d$ in the weak right Bruhat order (i.e. a reduced expression for $d$ can be extended to one for $x$). In this case we, loosely speaking, reduce computations involving $x$ to computations involving $d$. First we give an explicit construction for $\incl_d^{x,x^{-1}}$ in terms of $\incl_d^{d,d}$ (it is straightforward). Then in \autoref{subsection:simplifications} we show that $\mu_x = \lambda_d(f)$ where
$f$ is a so-called \emph{partial trace} of $\rho^a$ (see \autoref{D:partialtrace}). This reduces the computation of categorical dimension to the calculation of partial traces, and the calculation of $\lambda_d(f)$
for various polynomials $f$.

In \autoref{subsection:simplifications2} we focus on diagonal cells for which the distinguished involution $d$ is the longest element of a parabolic subgroup of $W$. We call these
\emph{parabolic cells}. For any element $x$ in a parabolic cell, $x$ is larger than $d$ in the weak right Bruhat order, so by the previous paragraph we can focus on the object $d$. In
this case we give an explicit construction of $\incl_d^{d,d}$, and we prove that $\lambda_d(f) = \partial_d(f)$, where $\partial_d$ is the Demazure operator of $d$.

\begin{Remark} It is fairly common in finite Coxeter groups for two-sided cells to contain at least one parabolic cell. This is true for all two-sided cells in type $A$. If it is false for a given two-sided cell, it is typically true for the $w_0$-dual two-sided cell (to the best of our knowledge, this is not in the literature). The
middle cell $\mathcal{J}_6$ in type $H_4$ has 24 diagonal cells, coming in three isomorphism classes, and two of the three have a parabolic cell. \end{Remark}

Another accomplishment of this paper is to automate the computation of partial traces in many cases. In \autoref{thm:linearpartialtracerecursion} we provide a recursive formula which
computes partial traces under the assumption that there is a \emph{linear} extension of the reduced expression for $d$ to one for $x$, see (see \autoref{def:recursible}). As for local
intersection forms, even when the computation cannot be fully automated, our techniques should still be helpful.

In the sequel we use computers to compute local intersection forms and categorical dimensions for various non-parabolic cells. We do use some of the technology developed here for
finding idempotents, but admittedly, we do not use much of the work done here to compute partial traces, and perhaps there is room for improvement as a result. One should view these
two papers as being different forks of the same drive to understand categorical dimensions, this paper focusing on abstract advances, and that paper focusing on brute force computation.

\subsection{More about idempotents}\label{intro:idempandtrace}

Let $B$ be an indecomposable object whose degree zero endomorphism space consists only of scalar multiples of the identity map. Let $X$ be an arbitrary object, and $k \in \Z$. The \emph{local intersection form} of $B(k)$ at $X$ is the pairing
\[ \Hom^k(X,B) \times \Hom^{-k}(B,X) \to \End^0(B) = \R \]
given by composition. The rank of the local intersection form determines the multiplicity of $B(k)$ as a summand of $X$, and dual bases give rise to (orthogonal) inclusion and projection maps to these summands. Using those inclusion and projection maps one can construct orthogonal idempotents in $\End(X)$. For our recursive construction of idempotents, we are most interested in the case when $X = B_w B_s$, where $s$ is a simple reflection such that $ws>w$.

Whenever we speak of ``computing'' a local intersection form, implicitly we have chosen a particular basis for these hom spaces. For example, if $X$ and $B$ are Bott--Samelson
bimodules then one could choose a basis of diagrams. When this happens the local intersection form might be valued in $\Z$ instead of $\R$, because the diagrammatic category has an
integral form where composition is defined. Then they could be specialized to $\F_p$. By \cite{ElWi-hodge-sbim}, the ranks of local intersection forms over $\R$ are determined purely in the Grothendieck group (they are
``as non-degenerate as possible''). Meanwhile, over $\F_p$ the local intersection forms are known to have additional degeneracies, and computing their rank is difficult. A
great deal of work has gone into computing local intersection forms over $\Z$ or in finite characteristic, with applications to modular representation theory, see e.g. \cite{Wi-torsion-explosion,MR3839376,JeWi-p-canonical,GiJeWi} and many others. In both cases these local intersection forms are named after and agree with certain local intersection forms in the geometry of flag varieties, either for perverse sheaves
(in characteristic zero) or parity sheaves (in finite characteristic), see for example \cite{MR3966750}.

In characteristic zero, the Soergel conjecture holds \cite{ElWi-hodge-sbim}, giving us a great deal of control on the size of morphism spaces and the decomposition of objects. In
\autoref{section:idempotentsinhecke} we exploit the Soergel conjecture and the Soergel hom formula to build produce a nice basis for particular spaces of projections and inclusions. Our
projections and inclusions are constructed using certain morphisms called $(w,i)$-trivalent vertices, introduced and developed in \autoref{section:idempotentsinhecke}. These morphisms
were known to experts in the field, but we are not aware of them appearing in the literature or their properties being developed outside of the special cases which appear in
\cite{El-thick-soergel-typea}. Using the properties we develop, we are able to prove a closed formula for many common local intersection forms which we believe is new and quite
powerful, see \autoref{Thm:LIFrecursionWORKS}. Our formula agrees with one found in \cite{ElLi-universal-soergel} for universal Coxeter groups.

The literature on local intersection forms in the Hecke category is growing, but it tends to fall into the following dichotomy: abstract results in characteristic zero, or
computational results in finite characteristic. There has not been a significant effort, before now, to compute local intersection forms in characteristic zero, taking advantage of the
Soergel conjecture. Perhaps this is because the ranks were known without computation. Nonetheless, there are numerous reasons one might wish to compute them, even ignoring their
connections to geometry. As motivated above, they play a major role in the computation of categorical dimensions. Knowing which primes divide determinants of local intersection forms
is the ``first line of defense'' against modular representation theory. A first question is: for which elements of $W$ does characteristic $p$ differ from characteristic $0$, i.e. does
$p$-canonical basis \cite{JeWi-p-canonical} disagree with the Kazhdan--Lusztig basis. Knowing local intersection forms will allow one to definitively find elements where they differ
(with major consequences, see \cite{Wi-torsion-explosion}), and can precisely find the minimal elements in the Bruhat order where they differ. In addition, knowing
the local intersection forms in characteristic zero can often let one bootstrap the local intersection forms in characteristic $p$: this was done in the parallel world of webs and the
Temperley--Lieb algebra in 
\cite{BuLiSe-tl-char-p,TuWe-quiver-tilting,MaSp-pl-jones-wenzl}.

%

Having discussed local intersection forms, let us discuss how they are used to build idempotents. For example, here is one way to describe the indecomposable Soergel bimodule $B_w$.
Pick a reduced expression $\un{w}$ for $w$, and consider the corresponding Bott--Samelson bimodule $B_{\un{w}}$. Then $B_w$ appears as a direct summand of $B_{\un{w}}$ with
multiplicity one, and all other summands have the form $B_y(k)$ for $y<w$ and $k \in \Z$. Because of this we say that $B_w$ is a \emph{top summand} of $B_{\un{w}}$, and an idempotent
in $B_{\un{w}}$ which projects to it will be called a \emph{top idempotent}. Inductively we have already constructed $B_y$ as the image of top idempotent in some other object, so we
may study maps to and from $B_y$, and we can compute various local intersection forms. Take the identity element of $B_{\un{w}}$ and subtract all the idempotents projecting to various
$B_y(k)$; what remains is a top idempotent. We can perform this same kind of construction not just for Bott--Samelson bimodules $B_{\un{w}}$ but for any object with $B_w$ as a top
summand. For example, if $w = xs > x$, then $B_w$ is a top summand in $B_x B_s$. The recursive formula that builds the idempotent projecting to $B_w$, knowing already the idempotent
projecting to $B_x$, can be found in \autoref{Thm:inductiveclasp}.

The idea to construct top idempotents recursively in this fashion is well-known in the community. In the parallel context of webs, this idea dates back to Wenzl's recursive formulas
for Jones--Wenzl idempotents (as nicely summarized and masterfully used in \cite{KaLi-TL-recoupling}); see \cite{El-ladders-clasps} and references therein for further developments and discussion. For the Hecke category, this is the idea used in
computer code in finite characteristic in \cite{GiJeWi}. The fact that top idempotents need not be unique is also known, though this was mostly discussed in finite characteristic,
see e.g. \cite{LiWi-non-perverse-soergel-bimodule} or \cite[Example 3.13]{SuTuWeZh-mixed-tilting}.

However, to the best of our knowledge, the literature lacks any methodical study or theory of top idempotents, which we provide in \autoref{section:idempotents}. This theory should apply to numerous other ``mixed
categories,'' including the Hecke category in finite characteristic, though we stick to the Hecke category in characteristic zero for pedagogical reasons. In
\autoref{thm:claspvsunique} we give a precise result stating when a top idempotent in the Hecke category (in any characteristic) is unique: this happens if and only if that idempotent
satisfies certain orthogonality properties, making it what we call a \emph{clasp idempotent}. In characteristic zero, whether a clasp idempotent exists is determined purely in the
Grothendieck group by \autoref{Thm:criterionforclasp}. These results did not seem to be in the literature before now. Only in
\autoref{Thm:criterionforclasp} and \autoref{section:idempotentsinhecke} do we seriously use special features of the Hecke category, as well as the Soergel conjecture.

\begin{Remark} The Hecke category differs from, for example, Temperley--Lieb and web categories in that the former is graded, leading to numerous subtleties in the definition of clasp idempotents. We warn the reader that some other graded diagrammatic categories, like categorified quantum groups, do not have top idempotents and are not amenable to our theory. \end{Remark}

\subsection{Organization of the paper}

In \autoref{section:hecke} we recall the basic properties of Hecke algebras and Hecke categories. This section is mostly review, but we include an expanded version of the Soergel hom formula which is more precise than most statements in the literature. We do not introduce or need cell theory until later in \autoref{section:asymptotichecke}.

The literature does not contain a good exposition on the general theory of idempotents in mixed categories, like the Hecke category. In \autoref{section:idempotents} we develop the
theory of top idempotents and clasp idempotents.

In \autoref{section:idempotentsinhecke} we shift our focus to properties of idempotents and related morphisms which are specific to the Hecke category, and explain how idempotents can
be recursively constructed. We give explicit formulas for the most frequent local intersection forms. We also explain how to compute some partial traces. Most of this material is new.

In \autoref{section:asymptotichecke} we recall the asymptotic Hecke category, and reformulate it so that it can interface with the most computable construction of the Hecke category (the Karoubi envelope of the diagrammatic Hecke category). We outline the construction of the asymptotic Hecke category and the computation of categorical dimensions, leaving some technical details and justifications to the sequel. We prove several simplifications for the computation of categorical dimensions, especially for parabolic cells.

\begin{Remark}
We recommend reading the paper in color, but we have made effort to make the paper readable in black-and-white by using dotted, dashed, dashdotted lines and various shadings.
\end{Remark}

\noindent\textbf{Acknowledgments.}
We like to thank Joel Gibson and Geordie Williamson for help with the Soergel bimodule implementation in Magma. We thank Victor
Ostrik, Leonardo Patimo and Ulrich Thiel for valuable discussions. Part of this project was done while the second author visited the first author at the University of Oregon, and their hospitality is gratefully acknowledged.

B.E. was supported by NSF grant DMS-2201387, and appreciates the support given to his research group by DMS-2039316. L.R. is supported by the SFB-TRR 195 ``Symbolic Tools in Mathematics and their Application'' of the German Research Foundation (DFG). D.T. was supported by many bad habits and the ARC Future Fellowship FT230100489.

\section{Hecke algebras and categories}\label{section:hecke}

We recall salient features of the Hecke algebra and the Hecke category. Most of the material is well-known and can be found in many places, but we additionally present a slightly more refined version of the Soergel hom formula \autoref{thm:SHFplus} than is typically found in the literature.


\subsection{Setup}\label{subsection:setup}

Let $W=(W,S)$ be a Coxeter group with fixed Coxeter system, i.e.
\begin{gather*}
W=
\langle
s\in S\mid
s^{2}=1, \quad
\underbrace{\dots sts}_{m_{s,t}}=\underbrace{\dots tst}_{m_{s,t}}
\rangle.
\end{gather*}
Here $m_{s,t}\in\N_{\geq 2}\cup\{\infty\}$ for $s \ne t \in S$, and the braid relation is vacuous in case $m_{s,t}=\infty$. We use standard conventions when we draw the Coxeter graph of $W$. We write $<$ for the Bruhat order on $W$, and $\len$ for the length function.


\begin{Notation}
Our simple reflections will be indexed by integers, i.e. $S = \{s_1, s_2, \ldots, s_n\}$. When simple reflections are used as indices or in sequences (e.g. when they are used for combinatorial purposes) we write $i$ instead of $s_{i}\in S$ for short. So $(s_1, s_2, s_1)$ would be abbreviated as $(1,2,1)$. When $W$ acts on something we typically do not abbreviate and write $s_i$.
\end{Notation}

An \emph{expression} is a sequence $\un{w} = (i_{1}, \ldots, i_k)$ of simple reflections, where $\len(\un{w}) = k$ is called the \emph{length} of $\un{w}$. The expression $\un{w}$ \emph{expresses} the element $w = i_{1} \cdots i_k \in W$, and it is \emph{reduced} if $\len(w) = \len(\un{w})$. For $w \in W$, let $\Rex(w)$ denote the set of reduced expressions of $w$.

\begin{Notation}
The identity element of $W$ will be denoted $\id$, so that it is not confused with $1$, which is shorthand for $s_1$. However, when elements of $W$ are used as indices (as in the Kazhdan--Lusztig basis $\{b_w\}$), we typically use $\emptyset$ instead of $\id$, referencing the fact that the empty sequence is a reduced expression for $\id$. Thus we use $\id \in W$ only rarely, to reduce ambiguity with various other identities.\end{Notation}

We fix a realization $(V,V^*,\Delta, \Delta^\vee)$ of $(W,S)$ over the real numbers $\R$ (or, if desired, the complex numbers $\C$). That is, $V$ is a finite-dimensional vector space over $\R$, $\Delta = \{\alpha_i \mid i \in S\} \subset V$ is a collection of \emph{simple roots}, and $\Delta^\vee = \{\alpha_i^\vee \mid i \in S\} \subset V^*$ is a collection of \emph{simple coroots}, satisfying two conditions. We set $a_{ij} = \alpha_i^\vee(\alpha_j) \in \R$, and these are the entries in the \emph{Cartan matrix} of the realization. The first condition is that $a_{ii} = 2$. The second condition is that there is an action of $W$ on $V$ via \begin{equation}\label{eq:sacts} s_i(v) := v-\alpha_{i}^\vee(v) \cdot \alpha_{i}. \end{equation} (In \cite{ElWi-soergel-calculus} there is a third condition for realizations, which is automatic in a characteristic zero field.)

In this paper only the Cartan matrix of the realization will play a significant role (e.g. there is no dependence on whether $\Delta$ is a spanning set for $V$ or not, or other more subtle features).

We assume our realization is \emph{balanced}, see \cite[Definition 3.6]{ElWi-soergel-calculus}. For small values of $m_{ij}$, being balanced is equivalent to
\begin{equation} \begin{cases} a_{ij} = 0 & \text{ when } m_{ij} = 2, \\
a_{ij} = -1 & \text{ when } m_{ij} = 3, \\
a_{ij} \cdot a_{ji} = 2 & \text{ when } m_{ij} = 4, \\
a_{ij} \cdot (a_{ij} a_{ji} - 2) = -1 & \text{ when } m_{ij} = 5. \end{cases} \end{equation}
One implication is that $a_{ij} = a_{ji}$ when $m_{ij}$ is odd. We say that a realization is \emph{standard} if 
\begin{equation} a_{ij} = -2 \cos(\pi/m_{ij}), \end{equation}
for $m_{ij}< \infty$, and $a_{ij} = -2$ for $m_{ij} = \infty$.

We use the Soergel conjecture (see below) in a crucial way, and the techniques used to prove this conjecture in \cite{ElWi-hodge-sbim} require the positivity property \cite[Equation (3.1)]{ElWi-hodge-sbim}, which can fail for some non-standard realizations (e.g. when $m_{ij} = 5$ and $a_{ij} = -2 \cos(3 \pi/5)$). It is expected that the Soergel conjecture should also hold for many non-standard realizations (such as the one just mentioned). While our computer calculations in the sequel use a standard realization, the techniques we discuss apply to general realizations for which the Soergel conjecture holds. 

\subsection{The (algebraic) Hecke category}\label{subsection:sbim}

Let $\setstuff{R}=\Sym(V)$ be the $\Z$-graded ring of polynomials (whose linear terms are $V$), graded so that $\deg(v) = 2$ for $v \in V$. Then \eqref{eq:sacts} induces an action of $W$ on $\setstuff{R}$. For $i\in S$ we let $\setstuff{R}^{i}$ be the subring of $i$-invariants. Define the graded $(\setstuff{R},\setstuff{R})$-bimodule
\begin{equation} \obstuff{B}_{i}=v\setstuff{R}\otimes_{\setstuff{R}^{i}}\setstuff{R}. \end{equation}
Here $v$ denotes a grading shift, with
$v(1\otimes 1)$ being in degree $-1$. Tensor products (over $R$) of various $\obstuff{B}_i$ are called \emph{Bott--Samelson bimodules}. More precisely, if $\un{w} = (i_1, \ldots, i_k)$ is an expression, then the correponding Bott--Samelson bimodule is
\[ \obstuff{B}_{\un{w}} := \obstuff{B}_{i_{1}} \ot_R \cdots \ot_R \obstuff{B}_{i_{k}} = \obstuff{B}_{i_{1}}\dots\obstuff{B}_{i_{k}}.\]
For the empty expression we obtain the monoidal identity:
\[ \munit = \obstuff{B}_{\emptyset} = \setstuff{R}. \]
As is common in this field, we omit the tensor product symbol $\otimes = \otimes_R$ when discussing bimodules, though we continue to use it when discussing elements of bimodules. When discussing tensor powers of objects we write e.g. $\obstuff{B}^{\otimes k}$.

For most graded $\setstuff{R}$-bimodules $\obstuff{B}$ we write $\obstuff{B}^k$ to denote the degree $k$ part of $\obstuff{B}$. For example, $\obstuff{B}_i^{-1} =
\R \cdot (1 \otimes 1)$. However, for the ring $\setstuff{R}$ itself we do not use this convention, to avoid conflict with setting $\setstuff{R}^2$ to be the subring of
invariants under $s_2$. We primarily care about the degree zero part of $R$, which is $\R \cdot \id$, and the degree $2$ part of $R$, which is $V$. Henceforth $V$ should always be interpreted as a subspace of $R$.

For $R$-bimodules $\obstuff{B}$ and $\obstuff{B}'$ we write $\Hom
(\obstuff{B},\obstuff{B}')$ for the graded $R$-bimodule of bimodule morphisms, and $\Hom^0
(\obstuff{B},\obstuff{B}')$ for the subspace of degree-preserving bimodule morphisms. We also set $\Hom^k
(\obstuff{B},\obstuff{B}') := \Hom^0
(\obstuff{B},v^k \obstuff{B}')$.

The category of \emph{Soergel bimodules} $\catstuff{S}=\catstuff{S}(W,S)$ is defined as
\begin{gather*}
\catstuff{S}:=\mathrm{cl}^{\otimes,\oplus,\dsum,v}
\big(\{\obstuff{B}_{i}\mid i\in S\}\big).
\end{gather*}
That is, it is the smallest full subcategory of $\Z$-graded $R$-bimodules which contains the bimodules $\obstuff{B}_{i}$, and is closed under tensor products $\otimes$, direct sums $\oplus$, direct summands $\dsum$, and shifts $v$.

\begin{Notation}
As in the previous display, we will use the convention the superscript of $\mathrm{cl}$ indicates the type of closure we take. Closure under tensor products tacitly includes the monoidal identity.
\end{Notation}

\begin{Remark} The meaning of closure under direct summands is somewhat subtle, see \autoref{subsection:karoubibullshit}. \end{Remark}

\begin{Example}\label{example:soergel1}
The baby example to keep in mind is $\catstuff{S}$ for the Coxeter group $(W,S=\{1\})$ of type $A_{1}$. We choose the one-dimensional realization where $V = \R \cdot x$, on which  $W$ acts by the sign representation. The polynomial algebra is $R=\R[x]$ and the invariant polynomials are $R^{1}=\R[x^{2}]$, both $\Z$-graded with $\deg(x)=2$. We have
\begin{equation}
\obstuff{B}_{\emptyset}=R = \munit,
\quad
\obstuff{B}_{1}=vR\otimes_{R^{1}}R.
\end{equation}
One can verify that
\begin{equation} \label{B1B1} \obstuff{B}_{1} \otimes \obstuff{B}_1 \cong
v\obstuff{B}_{1} \oplus v^{-1} \obstuff{B}_1.
\end{equation}
In particular, all tensor products $\obstuff{B}_1^{\otimes k}$ of positive length $k$ will decompose into shifts of $\obstuff{B}_1$, and the only indecomposable objects in $\catstuff{S}$ are shifts of $\obstuff{B}_{\emptyset}$ and $\obstuff{B}_{1}$. That is, $\catstuff{S}=\mathrm{cl}^{\oplus,v}
\big(\{\obstuff{B}_{\emptyset},\obstuff{B}_{1}\}\big)$.
\end{Example}

The category $\catstuff{S}$ is rigid, and left and right duals are isomorphic (see e.g. the isotopy relations in \cite[\S 10.2.4]{ElMaThWi-soergel}). Meanwhile, the category $\catstuff{S}$ also admits a monoidal \emph{(Poincar\'{e})
duality} functor $\DD$ (unrelated to the left and right duals in its rigidity structure). To differentiate between these notions of duality, we refer to $\rdual{-}$ as
\emph{adjunction}, and to $\DD$ as \emph{duality}. Note that the generating objects $B_i$ are both self-adjoint and self-dual: $\rdual{\obstuff{B}_{i}} \cong \obstuff{B}_{i}$ and
$\DD(B_i) \cong B_i$. However, duality is monoidally covariant and adjunction is monoidally contravariant, so $(B_i B_j)^* \cong B_j B_i$ whereas $\DD(B_i B_j) \cong B_i B_j$.

Let $\un{w} \in \Rex(w)$. By \cite[Theorem 1]{So-hcbim}, there is an indecomposable direct summand of $\obstuff{B}_{\un{w}}$ which is not a shift of a direct summand of some shorter expression, and this summand is unique up to isomorphism (though it need not be unique as a submodule). This summand is also independent (up to isomorphism) of the choice of reduced expression. Choose such a summand and denote it $\obstuff{B}_w$. All other direct summands of $\obstuff{B}_{\un{w}}$ have isomorphic to $v^k \obstuff{B}_{y}$ for $y < w$ and $k \in \Z$. That is, 
\begin{equation} \label{eq:topsummand} \obstuff{B}_{\un{w}} \cong \obstuff{B}_{w} \oplus \bigoplus_{y < w} c_{\un{w},y} \obstuff{B}_y \end{equation}
for some graded multiplicities $c_{\un{w},y} \in \N[v,v^{-1}]$. For this reason we call $\obstuff{B}_{w}$ the \emph{top summand} of $\obstuff{B}_{\un{w}}$. Note that $\DD(B_w) \cong B_w$ and $B_w^* \cong B_{w^{-1}}$.

\begin{Remark} Is it better to think of $\obstuff{B}_w$ as a bimodule, or as an isomorphism class of bimodule? It is only uniquely-defined in the latter sense. For the remainder of this section all statements only depend on the isomorphism class of $\obstuff{B}_w$. In the following section, we discuss the issues involved with pinning down $\obstuff{B}_w$ explicitly as a particular summand of $\obstuff{B}_{\un{w}}$. \end{Remark}

By \cite[Theorem 1]{So-hcbim}, every indecomposable Soergel bimodule is a shift of the bimodule $\obstuff{B}_w$ described above, for some $w \in W$. We can choose a set of representatives for the isomorphism classes of indecomposable objects as follows:
\begin{gather}
\Indec(\catstuff{S})=
\{v^{k}\obstuff{B}_{w}\mid w\in W,k\in\Z\}.
\end{gather}
We write $\id_w$ as shorthand for $\id_{B_w}$.

\begin{Remark} \label{Remark:notreduced}
If $\un{w}$ is any expression, not necessarily reduced, then the direct summands of $\obstuff{B}_{\un{w}}$ all have the form $v^k \obstuff{B}_{y}$ where $y$ is expressed by a subexpression of $\un{w}$. \end{Remark}


\subsection{The Hecke algebra and decategorification}\label{subsection:heckealg}

Let $H=H(W,S)$ be the \emph{Hecke algebra} over $\Z[v,v^{-1}]$ attached to $W$, for $v$ a formal variable. For the definition and basic details, see e.g. \cite[Chapter 3]{ElMaThWi-soergel}. It has two well-studied bases, 
the \emph{standard basis} denoted by $\{h_w \mid w \in W\}$, and 
the \emph{Kazhdan--Lusztig basis} or \emph{KL basis} denoted by $\{b_{w}\mid w\in W\}$. The Hecke algebra admits a $\Z$-linear involution (the bar involution, denoted using a bar) for which the elements $b_w$ are self-dual, and for which $\overline{v} = v^{-1}$. We have
\begin{equation} \label{eq:bw} b_w = h_w + \sum_{y < w} p_{w,y} h_y, \end{equation}
where the \emph{Kazhdan--Lusztig polynomials} $p_{w,y}$ satisfy $p_{w,y} \in v \Z[v]$. For $y < w$ let $\mu(w,y) \in \Z$ denote the coefficient of $v^1$ in $p_{w,y}$; these are called \emph{$\mu$-coefficients}. For any $i \in S$ for which $wi > w$ we have
\begin{equation} \label{decatinductivedecomp} b_w b_i = b_{wi} + \sum_{yi < y < w} \mu(w,y) b_y. \end{equation}
Meanwhile, if $wi < w$ we have
\begin{equation} \label{widown} b_w b_i = (v+v^{-1}) b_w. \end{equation}
For a proof see \cite[Lemma 9.10]{Lu-hecke-book}.

The Hecke algebra has a \emph{standard pairing}, denoted $(-,-)$, with values in $\Z[v,v^{-1}]$. It is sesquilinear: $(vb,b') = v^{-1}(b,b')$ and $(b,vb') = v(b,b')$. The standard basis is ``false orthonormal.'' It is not an orthonormal basis, but nonetheless, for \textit{self-dual} elements $b = \sum q_y h_y$ and $b' = \sum r_y h_y$ in $H$ we have
\begin{equation} \label{eq:stdpairing} (b,b') = \sum q_y r_y. \end{equation}

By \cite[Theorem 1]{So-hcbim} there is an algebra isomorphism $[\catstuff{S}] \cong H$ sending $[B_i] \mapsto b_i$, with $v$ being the decategorification of the grading shift. Soergel also proved a formula (see \cite[Theorem 3.6]{ElWi-hodge-sbim}) for the size of morphism spaces between Soergel bimodules, known as \emph{Soergel's hom formula} (or the Soergel hom formula): for $B, B' \in \catstuff{S}$, it states that $\Hom(B,B')$ is free as a right (or left) $R$-module, of graded rank equal to $([B],[B'])$.

Soergel conjectured that the isomorphism $[\catstuff{S}] \cong H^{v}$ sends $[B_w] \mapsto [b_w]$ for all $w \in W$ (at least for realizations in characteristic zero, which is our current setting), and proved this conjecture for Weyl groups. The general case was proven in \cite[Theorem 1.1]{ElWi-hodge-sbim}. A consequence is that sums of KL basis elements in the Hecke algebra lift to direct sums of Soergel bimodules. Thus \eqref{decatinductivedecomp} and \eqref{widown} lift to
\begin{equation} \label{inductivedecomp} \obstuff{B}_w \obstuff{B}_i  \cong \obstuff{B}_{wi} \oplus \bigoplus_{yi < y < w} \obstuff{B}_y^{\oplus \mu(w,y)} \end{equation}
when $wi > w$, and
\begin{equation} \label{downdecomp} \obstuff{B}_w \obstuff{B}_i  \cong v \obstuff{B}_w \oplus v^{-1} \obstuff{B}_w \end{equation}
when $wi < w$.

Combining the Soergel conjecture with the Soergel hom formula is very powerful. For example, the only morphisms of degree $0$ between various $B_w$ are identity maps (up to scalar).

\begin{Lemma} \label{lem:deg012initial} For $y, w \in W$ with $y \ne w$ we have
\begin{equation}\label{eq:nonegmaps} \dim \Hom^k(\obstuff{B}_w,\obstuff{B}_y) = 0 \quad \text{for } k \le 0, \end{equation}
\begin{equation} \label{eq:dimendzeroisone} \dim \End^0(\obstuff{B}_w) = 1, \end{equation}
\begin{equation} \label{eq:degree1soergelhom}\dim \Hom^1(\obstuff{B}_y, \obstuff{B}_w) = \dim \Hom^1(\obstuff{B}_w,\obstuff{B}_y) = \mu(w,y), \text{(see below)} \end{equation}
\begin{equation} \label{eq:deg2maps} \dim \End^2(\obstuff{B}_w) = \dim V + \sum_{x < w} \mu(w,x)^2. \end{equation}
Recall that $\mu(w,y)$ is only defined when $y < w$, so we have abused notation in \eqref{eq:degree1soergelhom}. More accurately, these degree 1 hom spaces have dimension equal to zero unless either $y<w$ or $w<y$, in which case they have dimension $\mu(w,y)$ or $\mu(y,w)$ respectively.
\end{Lemma}

\begin{proof} Most of this is immediate from the Soergel hom formula and \eqref{eq:stdpairing} and \eqref{eq:bw}, given that $p_{w,y} \in v \Z[v]$ for $y < w$. We elaborate only on
\eqref{eq:deg2maps}. One computes that $\End(B_w)$ is free as a right (resp. left) $R$-module with one generator (the identity) in degree zero, $\sum_{x<w} \mu(w,x)^2$ generators
in degree $2$, and additional generators in higher degrees. Thus the dimension in degree $2$ is $\sum_x \mu(w,x)^2$ plus the dimension of $R$ in degree 2, which is $\dim V$. \end{proof}

\begin{Remark} In most of the computations we do later in this paper, $\mu(w,y) \in \{0,1\}$. \end{Remark}

\begin{Remark} \label{rmk:nolowerindeg1} Let $y < w$. One should think of $p_{w,y}$ as counting the degrees of morphisms which go from $B_w$ ``directly'' to $B_y$, and it is a
lower bound for the graded dimension of $\Hom(B_w, B_y)$. There are also morphisms from $B_w$ to $B_y$ which factor through other indecomposable objects $B_z$. If $z \le w,y$
then the graded dimension of $\Hom(B_w,B_z)$ is at least $p_{w,z}$ and the graded dimension of $\Hom(B_z,B_y)$ is at least $p_{y,z}$, so one expects at least $p_{w,z} p_{y,z}$
(linearly independent) morphisms which factor as $B_w \to B_z \to B_y$. Finally, one obtains additional morphisms in $\Hom(B_w,B_y)$ by taking other morphisms and multiplying
by $R$. Both of these additional kinds of morphisms are evident in \eqref{eq:deg2maps}. The coefficient of $v^2$ in $p_{w,w}$ is zero, but nonetheless $\End^2(B_w)$ has a basis coming from elements of $R$ times $\End^0(B_w)$, and coming from morphisms which factor through $B_x$ for $x<w$.

Morphisms of degree $1$ are special, in that the coefficient of $v^1$ in $p_{w,y}$ agrees with the dimension of $\Hom^1(B_w,B_y)$, see \eqref{eq:degree1soergelhom}. As $\Hom^{-1}(B_w, B_y) = 0$, there are no morphisms of degree $+1$ which come from multiplying by positive degree elements of $R$. A degree $a$ map $B_w \to B_z$ composes with a degree $b$ map $B_z \to B_y$ to give a map $B_w \to B_y$ of degree $a+b$. When $a+b=1$ and $a \ge 0$ and $b \ge 0$ (thanks to \eqref{eq:nonegmaps}), then either $a=0$ or $b=0$. This implies (by \eqref{eq:nonegmaps} again) that $z=w$ or $z=y$, and the degree zero moprhism to or from $B_z$ is just a multiple of the identity. The moral of the story is: there is no interesting factorization of degree $+1$ maps. \end{Remark}

%

\subsection{Calculations in the Hecke algebra}\label{subsection:hecke-computer}

Calculations in Hecke algebra can be easily done using a computer, even in large types. The main program is due to du Cloux \cite{duCl-positivity-finite-hecke}. On the website
\cite{ElRoTu-ah-cat-code} we explain how this can be done. Because of this, we do not feel it necessary to justify certain calculations in the Hecke algebra (e.g. that $b_w b_x = b_y + b_z + b_u$ in some particular example) which can be easily verified by computer.

%

\subsection{Filtrations by ideals and Soergel's hom formula}\label{subsection:filtration}

This section does not use the monoidal structure on the Hecke category, but only the structure of composition.

\begin{Notation} For a Krull--Schmidt category $\catstuff{C}$, let $\Indec(\catstuff{C})$ denote the set of isoclasses of indecomposable objects in $\catstuff{C}$, or perhaps a choice of representatives for these isoclasses. \end{Notation}

\begin{Notation} For a Krull--Schmidt category $\catstuff{C}$, let $\setstuff{X}$ be a set of objects in $\catstuff{C}$. Let $\catstuff{I}_{\setstuff{X}}$ be the (ordinary = two-sided, not monoidal) ideal inside $\catstuff{C}$
generated by the identity maps of all objects in $\setstuff{X}$. Equivalently, a morphism $\obstuff{Y} \to \obstuff{Z}$ is in $\catstuff{I}_{\setstuff{X}}(\obstuff{Y},\obstuff{Z})$ if and only if it is a
linear combination of morphisms, each of which factors as a composition $\obstuff{Y} \to \obstuff{X} \to \obstuff{Z}$ for some $\obstuff{X} \in \setstuff{X}$. \end{Notation}

In any additive category we have $\obstuff{X} \cong \obstuff{Y} \oplus \obstuff{Z}$ if and only if there are projection and inclusion maps $p_{\obstuff{Y}}, p_{\obstuff{Z}}, i_{\obstuff{Y}}, i_{\obstuff{Z}}$ satisfying
\begin{equation} p_{\obstuff{Y}} \circ i_{\obstuff{Y}} = \id_{\obstuff{Y}}, \; \; p_{\obstuff{Z}} \circ i_{\obstuff{Z}} = \id_{\obstuff{Z}}, \; \;
p_{\obstuff{Y}} \circ i_{\obstuff{Z}} = 0, \; \; p_{\obstuff{Z}} \circ i_{\obstuff{Y}} = 0, \; \;
\id_{\obstuff{X}} = i_{\obstuff{Y}} \circ p_{\obstuff{Y}} + i_{\obstuff{Z}} \circ p_{\obstuff{Z}}. \end{equation}
As a consequence, $\id_{\obstuff{X}}$ is in the ideal generated by $\id_{\obstuff{Y}}$ and $\id_{\obstuff{Z}}$, and conversely, both $\id_{\obstuff{Y}}$ and $\id_{\obstuff{Z}}$ are in the ideal generated by $\id_{\obstuff{X}}$.

Given a set $\setstuff{X}$ of objects of $\catstuff{C}$, let $\overline{\setstuff{X}}$ denote the set of objects whose identity maps live in $\catstuff{I}_{\setstuff{X}}$. Clearly $\setstuff{X} \subset \overline{\setstuff{X}}$ and $\catstuff{I}_{\setstuff{X}} = \catstuff{I}_{\overline{\setstuff{X}}}$. By the previous paragraph, $\overline{\setstuff{X}}$ is closed under taking direct sums or direct summands. By the Krull--Schmidt property, any set of objects closed under taking direct sums or summands is determined by the indecomposable objects in this set.

\begin{Lemma} \label{lemma:onlyneedtoknowindecomposablesforideal} If $\setstuff{X}$ and $\setstuff{X}'$ are two sets of objects in $\catstuff{C}$, then $\catstuff{I}_{\setstuff{X}} \subset \catstuff{I}_{\setstuff{X}'}$ if and only if $\overline{\setstuff{X}} \cap \Indec(\catstuff{C}) \subset \overline{\setstuff{X}'} \cap \Indec(\catstuff{C})$. \end{Lemma}

\begin{proof} The proof is obvious given the discussion above. \end{proof}

Now we return to the setting of $\catstuff{S}$. Fix $w \in W$. Let $\catstuff{I}_{< w}$ denote the ideal in $\catstuff{S}$ generated by the identity maps of $\obstuff{B}_y$ for all $y < w$. This ideal is closed under the right and left actions of $\setstuff{R}$ on morphism spaces. One has a filtration of $\catstuff{S}$ by ideals, indexed by the poset $W$.

\begin{Lemma} \label{lem:imbored} The ideal $\catstuff{I}_{< w}$ agrees with the ideal generated by the identity maps of $\obstuff{B}_{\underline{y}}$ for all reduced expressions for $y$ with $y < w$. \end{Lemma}

\begin{proof} This is an immediate consequence of \autoref{lemma:onlyneedtoknowindecomposablesforideal} and \eqref{eq:topsummand}. \end{proof}

Similarly one can define ideals $\catstuff{I}_{\le w}$ or $\catstuff{I}_{\ngeq w}$, and one has the analogue of \autoref{lem:imbored}.

One can also include Bott--Samelson bimodules for non-reduced expressions in \autoref{lem:imbored}, so long as all direct summands have the form $\obstuff{B}_{y}$ for $y < w$. It is often helpful to factor morphisms through non-reduced expressions. While not as philosophically as significant as the ideals $\catstuff{I}_{< w}$, it is often convenient to use the ideals below to phrase our discussion, so that we need not distinguish between reduced and non-reduced expressions when factoring.

\begin{Definition} Fix $w \in W$. Let $\catstuff{I}_{< \len(w)}$ denote the ideal in $\catstuff{D}$ generated by the identity maps of $\obstuff{B}_y$ for all $y \in W$ with $\len(y) < \len(w)$. \end{Definition}

\begin{Lemma} The ideal $\catstuff{I}_{< \len(w)}$ agrees with the ideal generated by the identity maps of $\obstuff{B}_{\underline{y}}$ for all expressions (reduced or nonreduced) $\underline{y}$ of length less than $\len(w)$. \end{Lemma}

\begin{proof} This is an immediate consequence of \autoref{lemma:onlyneedtoknowindecomposablesforideal} and \autoref{Remark:notreduced}. \end{proof}

We are now prepared to discuss an expanded version of Soergel's hom formula, which is well-known but not commonly stated aloud in the literature.

\begin{Theorem} \label{thm:SHFplus} Let
$\obstuff{B}$ and $\obstuff{B}'$ be self-dual Soergel bimodules categorifying elements $b = \sum q_y h_y$ and $b' = \sum r_y h_y$. Then $\catstuff{I}_{< w}(\obstuff{B}, \obstuff{B}')$ is a free right $R$-module of graded rank 
\begin{equation} \label{strongerSHF} \rank \catstuff{I}_{< w}(\obstuff{B},\obstuff{B}') = \sum_{y < w} q_y r_y. \end{equation}
Similar statements can be made about other ideals like $\catstuff{I}_{\le w}$ or $\catstuff{I}_{< \len(w)}$. Moreover, the ideal filtration on hom spaces is split as a right $R$-module: the quotient $\catstuff{I}_{\le w}(\obstuff{B},\obstuff{B}')/\catstuff{I}_{< w}(\obstuff{B},\obstuff{B}')$ is free as a right $R$-module of graded rank $q_w r_w$.

Let $y \in W$. By the previous paragraph, $\Hom(\obstuff{B},B_y)/\catstuff{I}_{< y}(\obstuff{B},B_y)$ is free of rank $q_y$. Choose any elements $\{f_1, \ldots, f_{q_y}\} \in \Hom(B,B_y)$ which descend to a basis modulo $\catstuff{I}_{< y}$. Then $\{\DD(f_1), \ldots, \DD(f_{q_y})\}$ also descends to a basis for $\Hom(B_y,B)$. Similarly, choose elements $\{g_1, \ldots, g_{r_y}\}$ descending to a basis for  $\Hom(\obstuff{B}',B_y) / \catstuff{I}_{< y}(\obstuff{B}',B_y)$. Then the compositions $\DD(g_i) \circ f_j$ live in $\Hom(\obstuff{B}, \obstuff{B}')$ and descend to a basis for $\catstuff{I}_{\le y}(\obstuff{B},\obstuff{B}')/\catstuff{I}_{< y}(\obstuff{B},\obstuff{B}')$. Taking the union of these compositions over all $y$, we obtain a right $R$ basis for $\Hom(\obstuff{B}, \obstuff{B}')$. \end{Theorem}

\begin{Remark}
When $B$ and $B'$ are Bott--Samelson bimodules, the basis outlined in \autoref{thm:SHFplus} is very close to the double leaves basis from \cite{Li-light-leaves,ElWi-soergel-calculus}. 
(Technically the double leaves basis factors through $B_{\un{y}}$ instead of $B_y$, where $\un{y}$ is a reduced expression for $y$, but modifying the double leaves construction one can factor through $B_y$ instead.) This is unsurprising as many ``natural'' bases in similar settings are of this type of flavor, see, for example, \cite{El-ladders-clasps,AnStTu-cellular-tilting,ElLa-trace-hecke,An-tensor-q-tilting-modules} and many more references (too many to cite here).\end{Remark}

\begin{proof}(Sketch) It is not the purpose of this article to smooth over this gap in the literature, but let us give a rough outline of how existing results can be patched together to prove this theorem.

The ideals $\catstuff{I}$ above give a filtration on morphism spaces between Soergel bimodules, defined in terms of factorization. Meanwhile, Soergel
\cite{So-kl-bimodules} originally studied filtrations on the bimodules themselves, defined in terms of support. The relationship between these filtrations is studied in \cite[Section
3]{KELP2}, see in particular \cite[Lemma 3.29]{KELP2}; a morphism lies in $\catstuff{I}_{< w}$ if its image has support $< w$. From here it is not too difficult to unravel the original
proof of the hom formula and prove this stronger result. For example, the bimodule-theoretic statement that $B_y$ is isomorphic to (a shift of) the standard module $R_y$ modulo lower
terms in the support filtration can be transferred to prove that $\Hom(B,B_y)/\catstuff{I}_{<y}(B,B_y) \cong \Hom(B,R_y)$ up to shift, see \cite[Lemma 3.32]{KELP2}. Maps between
Soergel bimodules and standard modules are also controlled by Soergel's hom formula \cite[Theorem 5.15]{So-kl-bimodules}, so this morphism space is free over $R$.

That there is a spanning set built from compositions $\DD(g_i) \circ f_j$ is straightforward using induction, and they must be linearly independent or $\Hom(B,B')$ would be
too small. Alternatively, one could prove linear independence by localization to the fraction field $\setstuff{Q}$ of $\setstuff{R}$ (which makes the support filtration
split), see e.g. \cite{Elias2020LocalizedCF}. It is clear that $\catstuff{I}_{< w}$, after localization, agrees with the ideal generated by the identity maps of standard
bimodules $R_y$ indexed by $y<w$. The semisimplicity of Soergel bimodules after localization implies \eqref{strongerSHF} as a statement about ranks of free
$\setstuff{Q}$-modules, from which we can deduce the original rank before localizing.
\end{proof}

We elaborate on the expanded hom formula in the situation of \autoref{lem:deg012initial}.

\begin{Lemma} \label{lem:deg012} Fix $w \in W$. For each $x < w$, choose a basis $\{f^x_1, \ldots, f^x_{\mu(w,x)}\}$ for $\Hom^1(B_w,B_x)$, noting that $\mu(w,x) = 0$ is possible. One obtains $\mu(w,x)^2$ elements of $\End^2(B_w)$ as compositions $\DD(f^x_j) \circ f^x_k$. Let $\{h_1, \ldots, h_{\dim V}\}$ be a basis for $V$. One obtains $\dim V$ elements of $\End^2(B_w)$ as $\id_w \cdot h_j$ (resp. $h_j \cdot \id_w$). Together, the elements just described form a basis for $\End^2(B_w)$ over $\R$. 
\end{Lemma}

\begin{proof} We reiterate the moral from \autoref{rmk:nolowerindeg1}: there are no interesting factorizations of degree $+1$ maps. Thus $\catstuff{I}_{<x}(B_w,B_x)$ is $0$ in degree $+1$, and the basis $\{f^x_1, \ldots, f^x_{\mu(w,x)}\}$ above does descend to a basis for $\Hom^1(B_w,B_x)$ modulo lower terms. Now the result follows from \autoref{thm:SHFplus}. \end{proof}

\subsection{The diagrammatic Hecke category}\label{subsection:hecke-diag}

Extending earlier work of the first author with Khovanov, the reference \cite[Section 5]{ElWi-soergel-calculus} defines a $\Z$-graded 
$\R$-linear monoidal category $\tilde{\catstuff{D}}=\tilde{\catstuff{D}}(W,S)$ by generators and relations. Let $\catstuff{D}$ denote the (additive and) Karoubi envelope
\[ \catstuff{D} :=\mathrm{cl}^{\oplus,\dsum}(\tilde{\catstuff{D}}). \]
Then the main result of \cite{ElWi-soergel-calculus} is an equivalence of categories
\begin{gather}
\catstuff{S} \cong \catstuff{D}.
\end{gather}
More precisely, the category $\tilde{\catstuff{D}}$ is equivalent to the full-subcategory in $\catstuff{S}$ consisting of Bott--Samelson bimodules.

The category of bimodules $\catstuff{S}$ is more easily connected to geometry and representation theory, and came first historically, but the primary advantage of the diagrammatic
category $\tilde{\catstuff{D}}$ is that it is far easier to compute with. We will use $\tilde{\catstuff{D}}$ in our calculations, and refer the reader to \cite{ElWi-soergel-calculus} or
\cite{ElMaThWi-soergel} for more details.

Let us recall that the objects of $\tilde{\catstuff{D}}$ are $\otimes$-generated by $\{\obstuff{B}_{i}\mid i\in S\}$. Thus objects can be identified with expressions $\un{w}$. The endomorphism ring of the monoidal identity is $R$. In addition to polynomials, the generating morphisms are
\begin{gather} \label{generatorsofD}
\left\{
\begin{tikzpicture}[anchorbase,scale=1]
\draw[usual,white] (0,0.5) to (0,1)node[above]{$i$};
\draw[usual,marked=1] (0,0)node[below]{$i$} to (0,0.5);
\end{tikzpicture}
,\;
\begin{tikzpicture}[anchorbase,scale=1]
\draw[usual,white] (0,0) to (0,-0.5)node[below]{$i$};
\draw[usual,marked=0] (0,0) to (0,0.5)node[above]{$i$};
\end{tikzpicture}
,\;
\begin{tikzpicture}[anchorbase,scale=1]
\draw[usual] (0,0)node[below]{$i$} to (0.25,0.5);
\draw[usual] (0.5,0)node[below]{$i$} to (0.25,0.5);
\draw[usual] (0.25,0.5) to (0.25,1)node[above]{$i$};
\end{tikzpicture}
,\;
\begin{tikzpicture}[anchorbase,scale=1]
\draw[usual] (0,0)node[above]{$i$} to (0.25,-0.5);
\draw[usual] (0.5,0)node[above]{$i$} to (0.25,-0.5);
\draw[usual] (0.25,-0.5) to (0.25,-1)node[below]{$i$};
\end{tikzpicture}
,\;
\begin{tikzpicture}[anchorbase,scale=1]
\draw[usual] (0,0)node[below]{$i\phantom{j}\!$} to (0,0.5);
\draw[usuald] (0.25,0)node[below]{$j\phantom{j}\!$} to (0,0.5);
\draw[usual] (0.5,0)node[below]{$i\phantom{j}\!$} to (0,0.5);
\draw[usuald] (-0.25,0)node[below]{$\dots\phantom{j}\!$} to (0,0.5);
\draw[usual] (-0.5,0)node[below]{$i\phantom{j}\!$} to (0,0.5);
\draw[usuald] (0,1)node[above]{$j\phantom{j}\!$} to (0,0.5);
\draw[usual] (0.25,1)node[above]{$i\phantom{j}\!$} to (0,0.5);
\draw[usuald] (0.5,1)node[above]{$j\phantom{j}\!$} to (0,0.5);
\draw[usual] (-0.25,1)node[above]{$\dots\phantom{j}\!$} to (0,0.5);
\draw[usuald] (-0.5,1)node[above]{$j\phantom{j}\!$} to (0,0.5);
\end{tikzpicture}
\text{ or }
\begin{tikzpicture}[anchorbase,scale=1]
\draw[usual,white] (0,0)node[below,black]{$\dots\phantom{j}\!$} to (0,0.5);
\draw[usuald] (0.25,0)node[below]{$j\phantom{j}\!$} to (0,0.5);
\draw[usual] (0.5,0)node[below]{$i\phantom{j}\!$} to (0,0.5);
\draw[usual] (-0.25,0)node[below]{$i\phantom{j}\!$} to (0,0.5);
\draw[usuald] (-0.5,0)node[below]{$j\phantom{j}\!$} to (0,0.5);
\draw[usual,white] (0,1)node[above,black]{$\dots\phantom{j}\!$} to (0,0.5);
\draw[usual] (0.25,1)node[above]{$i\phantom{j}\!$} to (0,0.5);
\draw[usuald] (0.5,1)node[above]{$j\phantom{j}\!$} to (0,0.5);
\draw[usuald] (-0.25,1)node[above]{$j\phantom{j}\!$} to (0,0.5);
\draw[usual] (-0.5,1)node[above]{$i\phantom{j}\!$} to (0,0.5);
\end{tikzpicture}
\Bigg\vert
i,j\in S
\right\}
,
\end{gather}
where the last morphism generator is a $2m_{ij}$-valent vertex, and its coloration depends on the parity of $m_{ij}$ (the pictures above are for $m_{ij}$ odd or even respectively). 
These morphism generators are 
of degrees $1$, $1$, $-1$, $-1$ and $0$, respectively, and
called the \emph{counit or enddot, unit or startdot, multiplication, comultiplication and 
multivalent vertex}. They are
subject to the following relations.

\begin{enumerate}

\item First, $\obstuff{B}_{i}$ for $i\in S$ is a \emph{nilpotent Frobenius algebra object} in $\catstuff{D}$. 
Here nilpotence refers to the fact that multiplication composed with comultiplication 
is zero. For details see \cite[\S 8]{ElMaThWi-soergel}. We recall in particular the \emph{broken trivalent relation}

\begin{equation} \label{eq:brokentrivalent} \begin{tikzpicture}[anchorbase,scale=1]
\draw[usual] (0.6,0) to (0.3,1);
\draw[usual,marked=1] (0,0) to (0.15,0.3);
\end{tikzpicture}
=
\begin{tikzpicture}[anchorbase,scale=1]
\draw[usual] (0,0) to (0.3,1);
\draw[usual,marked=1] (0.6,0) to (0.45,0.3);
\end{tikzpicture}
+
\begin{tikzpicture}[anchorbase,scale=1]
\draw[usual] (0.6,0) to[out=90,in=0] (0.3,0.3) to[out=180,in=90] (0,0);
\draw[usual,marked=1] (0.3,1) to (0.3,0.7);
\end{tikzpicture}
-
\begin{tikzpicture}[anchorbase,scale=1]
\draw[usual] (0,0) to (0.3,0.3) to (0.3,1);
\node[scale=0.8] at (0.8,0.5) {$\alpha_i$};
\draw[usual] (0.6,0) to (0.3,0.3);
\end{tikzpicture} \end{equation}
and the \emph{needle relation} \begin{equation}\label{eq:needle}
\begin{tikzpicture}[anchorbase,scale=1]
\draw[usual] (0,0) to[out=270, in=180] (0.3,-0.3) to[out=0, in=270] (0.6,0) to[out=90,in=0] (0.3,0.3) to[out=180,in=90] (0,0);
\draw[usual] (0.3,-0.7) to (0.3,-0.3);
\end{tikzpicture}
=
0.
\end{equation}

\item Second, relations which manipulate polynomials. The second relation holds for all $f \in R$.
\begin{gather}\label{eq:circle-evaluation}
\begin{tikzpicture}[anchorbase,scale=1]
\draw[usual,marked=0,marked=1] (0,-0.3) to (0,0)node[right,xshift=-2pt]{$i$} to (0,0.3);
\end{tikzpicture}
=
\alpha_{i}
,\qquad
\begin{tikzpicture}[anchorbase,scale=1]
\draw[usual] (0,-0.5)node[below]{$i$} to (0,0.5)node[above]{$i$};
\node at (-.6,0) {$f$};
\end{tikzpicture}
+
\begin{tikzpicture}[anchorbase,scale=1]
\draw[usual] (0,-0.5)node[below]{$i$} to (0,0.5)node[above]{$i$};
\node at (.6,0) {$s_i(f)$};
\end{tikzpicture}
=
\begin{tikzpicture}[anchorbase,scale=1]
\draw[usual,marked=1] (0,-0.5)node[below]{$i$} to (0,-0.25);
\draw[usual,marked=0] (0,0.25) to (0,0.5)node[above]{$i$};
\node at (-.6,0) {$\partial_{i}(f)$};
\end{tikzpicture}
.
\end{gather}
Here we have used the \emph{Demazure operator}, a degree $-2$ operator defined for each $i \in S$ by
\begin{gather}
\partial_{i}\colon
\setstuff{R}\to\setstuff{R}^{i},
\quad
f\mapsto
\tfrac{f-s_i(f)}{\alpha_{i}}.
\end{gather}
Equivalently, one can also define $\partial_i$ by the following two properties. First, as a map from degree two to degree zero in $R$, i.e. from $V$ to $\R$, we have $\partial_i(v) = \alpha_i^\vee(v)$. Second, $\partial_i$ satisfies the \emph{twisted Leibniz rule}
\begin{gather}\label{eq:leibniz}
\partial_{i}(fg)=
\partial_{i}(f)g+s_i(f)\partial_{i}(g).
\end{gather}

A consequence of these relations and the needle relation is that
\begin{equation}\label{eq:needle_polym}
\begin{tikzpicture}[anchorbase,scale=1]
\draw[usual] (0,0) to[out=270, in=180] (0.3,-0.3) to[out=0, in=270] (0.6,0) to[out=90,in=0] (0.3,0.3) to[out=180,in=90] (0,0);
\node at (0.3,0) {$f$};
\draw[usual] (0.3,-0.7) to (0.3,-0.3);
\end{tikzpicture}
=
\partial_i(f)
\cdot
\begin{tikzpicture}[anchorbase,scale=1]
\draw[usual, marked=1] (0.3,-0.7) to (0.3,0);
\end{tikzpicture}
\end{equation}

\item Finally, certain relations involving multivalent vertices, cf. \cite[(5.6)--(5.12)]{ElWi-soergel-calculus}. The most important consequence of these relations is the following. If $\un{w}$ and $\un{w}'$ are two reduced expressions for the same element $w$ which are related by a braid relation, then there is a multivalent vertex giving a morphism $\un{w} \to \un{w}'$, and this morphism is an isomorphism modulo $\catstuff{I}_{< w}$.

\end{enumerate}

The category $\tilde{\catstuff{D}}$ has a duality functor which flips diagrams upside-down, and an adjunction functor which rotates them 180 degrees. The isomorphism $\catstuff{D} \cong \catstuff{S}$ respects the structures of adjunction and duality.



%
%
%
%

We use the following color code for our simple reflections.
\begin{gather}\label{eq:color-code}
\text{arbitrary}
\leftrightsquigarrow
\begin{tikzpicture}[anchorbase,scale=1]
\draw[usual] (0,0) to (0,1);
\end{tikzpicture}
\; ,\quad
1\leftrightsquigarrow
\begin{tikzpicture}[anchorbase,scale=1]
\draw[soergelone] (0,0) to (0,1);
\end{tikzpicture}
\; ,\quad
2\leftrightsquigarrow
\begin{tikzpicture}[anchorbase,scale=1]
\draw[soergeltwo] (0,0) to (0,1);
\end{tikzpicture}
\; ,\quad
3\leftrightsquigarrow
\begin{tikzpicture}[anchorbase,scale=1]
\draw[soergelthree] (0,0) to (0,1);
\end{tikzpicture}
\; ,\quad
4\leftrightsquigarrow
\begin{tikzpicture}[anchorbase,scale=1]
\draw[soergelfour] (0,0) to (0,1);
\end{tikzpicture}
\; ,\quad
5\leftrightsquigarrow
\begin{tikzpicture}[anchorbase,scale=1]
\draw[soergelfive] (0,0) to (0,1);
\end{tikzpicture}
\;.
\end{gather}

\section{On idempotents}\label{section:idempotents}

To the best of our knowledge, the literature does not have a good ``user's guide'' to understanding idempotents in the Hecke category, so the goal of this section is to fill this gap in the literature.

A \emph{mixed category} (in the sense of \cite[Definition 1.11]{We-canonical-bases-higher-rep}) is a graded additive category, linear over a field and possessing a contravariant involutive duality functor, with the following property: there is a set of self-dual indecomposable objects $\{B_w\}_{w \in W}$ indexed by some set $W$, such that: \begin{itemize}
\item every indecomposable object is isomorphic to a shift of some $B_w$,
\item the space $\End^0(B_w)$ is one-dimensional, consisting only of scalar multiples of the identity,
\item the space $\Hom^i(B_w, B_y)$ is zero for $i<0$ and all $w, y \in W$,
\item the space $\Hom^0(B_w,B_y)$ is zero for $w \ne y$. \end{itemize}
Mixedness is a graded analog of semisimplicity. The Hecke category is mixed by the Soergel hom formula.

Most of the techniques we discuss in this section apply more broadly to mixed categories with interesting ideal filtrations (corresponding to a poset structure on $W$), though for simplicity we state them only for the Hecke category. (In this broader context, one should replace rex moves with ``neutral moves'' in that category, e.g. the neutral ladders of \cite{El-ladders-clasps}.)

In the following section we focus on the implications of Kazhdan--Lusztig combinatorics to the computation of idempotents, working with techniques which are specific to the Hecke
category.

\subsection{Decompositions}\label{subsection:decomp}

Below we will discuss objects which are decomposed as direct sums of other objects. By convention, the data of a decomposition includes the choice of inclusion and projection maps.

\begin{Definition} Let $X$ be an object in an additive category. A \emph{decomposition} $X \cong M \oplus N$ is the data of morphisms
\begin{equation} \iota_M \colon M \to X, \quad p_M \colon X \to M, \quad \iota_N \colon N \to X, \quad p_N \colon X \to N, \end{equation}
which satisfy
\begin{equation} \id_X = \iota_M p_M + \iota_N p_N, \quad p_M \iota_N = 0, \quad p_N \iota_M = 0, \quad \id_N = p_N \iota_N, \quad \id_M = p_M \iota_M. \end{equation}
A decomposition of $X$ into more than two summands is defined similarly (or equivalently, by induction from the case of two summands).
\end{Definition}

Again, we reiterate that a decomposition is slightly more data than the choice of an idempotent decomposition $\id_X = e_M + e_N$ into orthogonal idempotents, whose images are
isomorphic to $M$ and $N$ respectively. It also includes the inclusion and projection maps, which factor these idempotents as $e_M = \iota_M p_M$ and $e_N = \iota_N p_N$.

Given a decomposition and an invertible scalar $\lambda_M$, one can replace $p_M$ with $\lambda_M p_M$ and $\iota_M$ with $\lambda_M^{-1} \iota_M$, and get a different decomposition
with the same idempotents. This symmetry is called \emph{rescaling the decomposition}.

\subsection{Top idempotents}\label{subsection:topdef}

Let $w \in W$ and $\un{w}\in \Rex(w)$. Recall that $\obstuff{B}_{w}$ is a direct summand of $\obstuff{B}_{\un{w}}$, the unique summand (up to isomorphism) which is not a shift of $\obstuff{B}_y$ for $y < w$.
Thus there exists a primitive idempotent 
$\morstuff{e}(\un{w})\in\End(\obstuff{B}_{\un{w}})$
whose image is isomorphic to $\obstuff{B}_{w}$. Any such idempotent will be called a \emph{top idempotent}. Top idempotents are not always unique, see \autoref{nastyexample}.

Indeed, the reader should think that $\obstuff{B}_w$, as previously defined and discussed, is better viewed as an isomorphism class of bimodule than as a biomdule per se. The purpose of a top idempotent is to pick out a bimodule in this isomorphism class explicitly, as the image of said idempotent.

More generally we have the following definition.

\begin{Definition} Let $X$ be any self-dual object in $\catstuff{S}$ (or equivalently, in $\catstuff{D}$). We say that $B_w$ is a \emph{top summand} of $X$ if we have chosen a decomposition of $X$ into indecomposable direct summands, such that (a summand isomorphic to) $B_w$ appears once, and all other summands are (isomorphic to) shifts of $B_y$ for $y < w$. That is,
\begin{equation} X = B_w \oplus \bigoplus_{y < w} B_y^{\oplus n_y} \end{equation}
for some graded multiplicities $n_y \in \N[v,v^{-1}]$, possibly zero. When $B_w$ is a top summand of $X$, we call $X$ a \emph{$w$-object}, and implicit in this notation is the choice of an explicit decomposition. In this case, the primitive idempotent $e \in \End(X)$ whose image is isomorphic to $B_w$ will be called the \emph{top idempotent} (associated to this decomposition of $X$). Note that a different decomposition of the same object $X$ might produce a different top idempotent. \end{Definition}

We often denote top idempotents by rectangles, e.g.:
\begin{gather*}
\begin{tikzpicture}[anchorbase,scale=1]
\draw[mor] (0,0) rectangle (1,0.5)node[black,pos=0.5]{$\morstuff{e}$};
\end{tikzpicture}. 
\end{gather*}
The following example is a crucial one to keep in mind, and illustrates the non-uniqueness of top idempotents.

\begin{Example} \label{nastyexample} Let $w = S_4$ with the usual simple reflections $\{1,2,3\}$. Let $\un{w} = (3,1,2,1,3)$, a reduced expression for $w = 31213$, and $y =
31$. Note that $B_y \cong B_1 B_3$. A straightforward computation in the Hecke algebra shows that $b_w = b_3 b_{121} b_3$, and consequently we can describe $B_w$ inside
$X = B_{\un{w}}$ as the image of the following top idempotent $e$. \begin{gather} \label{e31213}
e=
\begin{tikzpicture}[anchorbase,scale=1]
\diagrammaticmorphism{0}{0}{5}{3,1,2,1,3}
\end{tikzpicture}
\; + \;
\begin{tikzpicture}[anchorbase,scale=1]
\draw[soergelthree] (0,0) to (0,0.5);
\draw[soergelone] (0.2,0) to (0.2,0.5);
\draw[soergeltwo,markedtwo=1] (0.4,0) to (0.4,0.15);
\draw[soergeltwo,markedtwo=1] (0.4,0.5) to (0.4,0.35);
\draw[soergelone] (0.6,0) to (0.6,0.5);
\draw[soergelone] (0.2,0.25) to (0.6,0.25);
\draw[soergelthree] (0.8,0) to (0.8,0.5);
\end{tikzpicture} \; .
\end{gather} We have taken the (top) idempotent inside $B_{(1,2,1)}$ projecting to $B_{121}$, and tensored on the left and right with the identity maps of $B_3$.

The Kazhdan--Lusztig polynomial $p_{w,y}$ is equal to $v + v^3$, so there is a degree $+1$ map
$g \colon B_w \to B_y$. Identifying $g \in \Hom(B_w,B_y)$ with $\tilde{g} \in \Hom(B_{\un{w}},B_y)e$, this degree $+1$ map becomes \begin{gather}
\tilde{g} = \;
\begin{tikzpicture}[anchorbase,scale=1]
\draw[mor] (-0.1,-0.5) rectangle (0.9,0)node[black,pos=0.5]{$e$};
\draw[soergelthree] (0,0) to (0,0.6);
\draw[soergelone, markedone=1] (0.2,0) to (0.2,0.15);
\draw[soergeltwo,markedtwo=1] (0.4,0) to (0.4,0.15);
\draw[soergelone] (0.6,0) to (0.6,0.1) to (0.2,0.5) to (0.2,0.6);
\draw[soergelthree] (0.8,0) to (0.8,0.3) to (0,0.3);	
\end{tikzpicture}
\; = \;
\begin{tikzpicture}[anchorbase,scale=1]
\diagrammaticmorphism{-0.1}{-0.4}{4}{3,1,2,1,3}
\draw[soergelthree] (0,0) to (0,0.6);
\draw[soergelone, markedone=1] (0.2,0) to (0.2,0.15);
\draw[soergeltwo,markedtwo=1] (0.4,0) to (0.4,0.15);
\draw[soergelone] (0.6,0) to (0.6,0.1) to (0.2,0.5) to (0.2,0.6);
\draw[soergelthree] (0.8,0) to (0.8,0.3) to (0,0.3);	
\end{tikzpicture}
\; + \;
\begin{tikzpicture}[anchorbase,scale=1]
\draw[soergelthree] (0,-0.4) to (0,0.6);
\draw[soergelone] (0.2,-0.4) to (0.2,-0.2) to (0.6,-0.2);
\draw[soergeltwo,markedtwo=1] (0.4,-0.4) to (0.4,-0.3);
\draw[soergelone] (0.6,-0.4) to (0.6,0.1) to (0.2,0.5) to (0.2,0.6);
\draw[soergelthree] (0.8,-0.4) to (0.8,0.3) to (0,0.3);
\node at (0.3, 0) {\(\alpha_2\)};
\end{tikzpicture} \;.
\end{gather} We also have $b_1 b_3 b_2 b_1 b_3 = b_w +
(v+v^{-1}) b_y$, so there is a degree $-1$ inclusion map $B_y \to B_{\un{w}}$, drawn as\begin{gather}
\iota=
\begin{tikzpicture}[anchorbase,scale=1]
\draw[soergelthree] (0,0) to (0,0.4);
\draw[soergelone] (0.2,0) to (0.2,0.4);
\draw[soergeltwo,markedtwo=1] (0.4,0.4) to (0.4,0.3);
\draw[soergelone] (0.2,0.2) to (0.6,0.2) to (0.6,0.4);
\draw[soergelthree] (0,0.1) to (0.8,0.1) to (0.8,0.4);	
\end{tikzpicture}\; .
\end{gather} Note that $\tilde{g} \circ \iota = 0$, and indeed $e \circ \iota = 0$, which is to be expected since $\Hom^{-1}(B_y,B_w)=0$.

Let $f = \iota \circ \tilde{g}$, a nonzero endomorphism
of $B_{\un{w}}$ living in the degree zero part of $\catstuff{I}_{< w}$ (whose existence was foretold by the Soergel hom formula). One can compute that $fe = f$ since $\tilde{g}e = \tilde{g}$, whereas $ef = 0$ since $e \iota = 0$. Thus $f^2 =
fef = 0$.

As a consequence, for any scalar $\lambda$, we have \begin{equation} (e + \lambda f)^2 = e^2 + \lambda ef + \lambda fe + \lambda^2 f^2 = e + \lambda f. \end{equation} Letting $e' = e +
\lambda f$, we also have $e' f= 0$ and $f e' = f$. We have $e e' e = e$ and $e' e e' = e'$. Thus $e$ and $e'$ have isomorphic images. We have a family $\{e + \lambda f\}$ of
idempotents, all of which are top idempotents.

We call upon this example often below, where we write $e = e_w$ for clarity. \end{Example}

We discuss the general failure of top idempotents to be unique in \autoref{subsection:topunique}, but the above example is prototypical. The main issue is the existence of a nonzero morphism in $\catstuff{I}_{< w}(B_w,X)$ of degree zero, which was possible in the example above because $B_y$ appeared in $X$ with a shift.

\begin{Definition}\label{def:perverse} A $w$-object $X$ is called \emph{perverse} if 
\begin{equation} X \cong B_w \oplus \bigoplus_{y < w} B_y^{\oplus n_y} \end{equation}
where $n_y \in \N$ (i.e. there are no shifts). \end{Definition}

Bott--Samelson objects for reduced expressions are commonly not perverse, but perverse $w$-objects abound.

\begin{Example} \label{Ex:wi} Let $x \in W$ and $i \in S$ such that $w = xi > x$. Then $B_x B_i$ is a perverse $w$-object. This follows from \eqref{inductivedecomp}. \end{Example}

\begin{Example} It is tautological, but $B_w$ itself is a perverse $w$-object. Concretely, consider the setup of \autoref{nastyexample}, and let $X'$ denote the image of the idempotent $e$, a summand of $X=B_{\un{w}}$ isomorphic to $B_w$. There is a unique (top) idempotent in $\End(X') \cong e \End(X) e$, namely $\id_{X'} \leftrightarrow e$, as $\End^0(B_w)$ is one-dimensional. At the same time, $efe = 0$, and there are no nonzero morphisms in $\catstuff{I}_{<w}(X',X')$. \end{Example}

We will show later that a perverse $w$-object has a unique top idempotent, though there are also many non-perverse $w$-objects with a unique top idempotent.

\begin{Example} \label{extendingnastypart1} We continue the previous examples, so that $X = B_{(3,1,2,1,3)}$ is a $w$-object for $w = 31213$, as is $X'$ which is the image of $e_w$. Then both $X B_2$ and $X' B_2$ are $w_0$-objects, for $w_0 = 312132$ the longest element of $S_4$. The top
idempotent $e_{0} \in \End(X' B_2)$ is
\begin{equation} \label{idempotente0}
e_0
=
\begin{tikzpicture}[anchorbase,scale=1]
\diagrammaticmorphism{0}{0.5}{3}{3,1,2,1,3}
\draw[mor] (0,0) rectangle (1,0.5)node[black,pos=0.5]{$e$};
\diagrammaticmorphism{0}{-0.3}{3}{3,1,2,1,3}
\draw[soergeltwo] (1.1,-0.3) to (1.1,0.8);
\end{tikzpicture}
+
\begin{tikzpicture}[anchorbase,scale=1]
\draw[mor] (0,-0.2) rectangle (1,0.3)node[black,pos=0.5]{$e$};
\diagrammaticmorphism{0}{0.3}{3}{3,1,2}
\drawtwomvalentvertex{2}{soergelthree}{soergelone}{0}{0.6}
\drawtwomvalentvertex{3}{soergelthree}{soergeltwo}{0.2}{0.7}
\draw[soergeltwo] (0.5,0.6) to (0.5,0.7);
\draw[soergelone,markedone=1] (0.7,0.3) to (0.7,0.5);
\draw[soergelthree] (0.9,0.3) to (0.9,0.5) to (0.7,0.7);
\draw[soergelone] (0.1,0.7) to (0.1,1);
\diagrammaticmorphism{0}{0.9}{2}{1,2,3,2}
\draw[soergeltwo] (1.1,-0.2) to (1.1,1) to (0.7,1);
\draw[mor] (0,1.1) rectangle (0.8,1.6)node[black,pos=0.5]{$1232$};
\begin{scope}[yscale=-1, shift={(0,-2.7)}]
\draw[mor] (0,-0.2) rectangle (1,0.3)node[black,pos=0.5]{$e$};
\diagrammaticmorphism{0}{0.3}{3}{3,1,2}
\drawtwomvalentvertex{2}{soergelthree}{soergelone}{0}{0.6}
\drawtwomvalentvertex{3}{soergelthree}{soergeltwo}{0.2}{0.7}
\draw[soergeltwo] (0.5,0.6) to (0.5,0.7);
\draw[soergelone,markedone=1] (0.7,0.3) to (0.7,0.5);
\draw[soergelthree] (0.9,0.3) to (0.9,0.5) to (0.7,0.7);
\draw[soergelone] (0.1,0.7) to (0.1,1);
\diagrammaticmorphism{0}{0.9}{2}{1,2,3,2}
\draw[soergeltwo] (1.1,-0.2) to (1.1,1) to (0.7,1);
\end{scope}
\end{tikzpicture}	
+
\begin{tikzpicture}[anchorbase,scale=1]
\draw[mor] (0,0) rectangle (1,0.5)node[black,pos=0.5]{$e$};
\drawtwomvalentvertex{3}{soergelone}{soergeltwo}{0.2}{0.5}
\diagrammaticmorphism{0.2}{0.7}{2}{2,1,2}
\draw[soergelthree,markedthree=1] (0.9,0.5) to (0.9,0.6);
\draw[soergeltwo] (1.1,0) to (1.1,0.7) to (0.7,0.7);
\draw[soergelthree] (0.1,0.5) to (0.1,0.9);
\draw[mor] (0,0.9) rectangle (0.8,1.4)node[black,pos=0.5]{$3212$};
\begin{scope}[yscale=-1, shift={(0,-2.3)}]
\draw[mor] (0,0) rectangle (1,0.5)node[black,pos=0.5]{$e$};
\drawtwomvalentvertex{3}{soergelone}{soergeltwo}{0.2}{0.5}
\diagrammaticmorphism{0.2}{0.7}{2}{2,1,2}
\draw[soergelthree,markedthree=1] (0.9,0.5) to (0.9,0.6);
\draw[soergeltwo] (1.1,0) to (1.1,0.7) to (0.7,0.7);
\draw[soergelthree] (0.1,0.5) to (0.1,0.9);
\end{scope}
\end{tikzpicture}	\;.
\end{equation}
The idempotents for $3212$ and $1232$ which appear above are easy to compute (similar to the idempotent $e$ from \eqref{e31213}) but are actually unnecessary; one can replace them with the identity map and get the same morphism; we include them to indicate to the reader that these morphisms factor through $B_{3212}$ and $B_{1232}$ respectively. The savvy reader can try to verify that \eqref{idempotente0} yields a top idempotent now, though it will be justified in \autoref{ex:justifycoeff1}. It turns out that $e_0$ is the unique top idempotent in $\End(X' B_2)$, which is not surprising since $X' B_2$ is perverse.

Note that $e_{0}$ can be identified with a top idempotent $e_{\un{w}_0} \in
\End(B_{(3,1,2,1,3,2)}) = \End(X B_2)$, via the identification $\End(X' B_2) = (e_w \ot \id_2) \circ \End(X B_2) \circ (e_w \ot \id_2)$. Note that $X B_2$ is not perverse. While these two top idempotents $e_{0}$ and
$e_{\un{w}_0}$ seem the same, they live in different contexts; the unicity of $e_{0}$ in $\End(X' B_2)$ does not imply the unicity of $e_{\un{w}_0}$ in $\End(X B_2)$. As it turns out, $e_{\un{w}_0}$ is also unique, see \autoref{Thm:criterionforclasp}. \end{Example}

\begin{Example} \label{Ex:nonperverseclasp} Let $w = 31213$ as in \autoref{nastyexample}, and let $x = s_1$. Let $Y \cong B_w \oplus B_x(1) \oplus B_x(-1)$. Then $Y$ admits a unique top idempotent, see \autoref{Thm:criterionforclasp}. There is no degree $+1$ map $B_w \to B_x$, contrasting from \autoref{nastyexample} where there was a degree $+1$ map $B_w \to B_y$. \end{Example}

In \autoref{subsection:clasps} we define the concept of a clasp idempotent in $\End(X)$, analogous to the concept of a Jones--Wenzl projector, or a clasp in the theory of webs. Clasp
idempotents are defined by an orthogonality property: they annihilate the degree zero part in the ideal $\catstuff{I}_{< w}(X,X)$. It is easy to show that clasp idempotents are unique
when they exist. Though their definition is not a priori related to that of top idempotents, we prove that a clasp idempotent exists precisely when the top idempotent is unique, in
which case it equals that unique top idempotent. Then we will demonstrate an algorithm to construct clasp idempotents successively for perverse $w$-objects of the form in
\autoref{Ex:wi}.

\subsection{Top idempotents and lower terms}\label{subsection:toplower}

Fix $w \in W$, and fix a $w$-object $X$. Let $\underline{y}$ be an expression of length less than $\len(w)$. We illustrate morphisms $\morstuff{g} \colon X \to \obstuff{B}_{\underline{y}}$ or $\morstuff{f} \colon \obstuff{B}_{\underline{y}} \to X$ using \emph{trapezes}:
\begin{gather*}
\begin{tikzpicture}[anchorbase,scale=1]
\draw[mor] (0,0) to (0.25,0.5) to (0.75,0.5) to (1,0) to (0,0);
\node at (0.5,0.25){$\morstuff{g}_{\underline{y}}$};
\end{tikzpicture}
,\quad
\begin{tikzpicture}[anchorbase,scale=1]
\draw[mor] (0,1) to (0.25,0.5) to (0.75,0.5) to (1,1) to (0,1);
\node at (0.5,0.75){$\morstuff{f}_{\underline{y}}$};
\end{tikzpicture}.
\end{gather*}
These are morphisms in $\catstuff{I}_{< \len(w)}$, and a general morphism in $\catstuff{I}_{< \len(w)}(X,X)$ is a linear combination of morphisms of the form
\begin{gather} \label{eq:downup} \begin{tikzpicture}[anchorbase,scale=1]
\draw[mor] (0,0) to (0.25,0.5) to (0.75,0.5) to (1,0) to (0,0);
\node at (0.5,0.25){$\morstuff{g}_{\underline{y}}$};
\draw[mor] (0,1) to (0.25,0.5) to (0.75,0.5) to (1,1) to (0,1);
\node at (0.5,0.75){$\morstuff{f}_{\underline{y}}$};
\end{tikzpicture}.
\end{gather}
The width of the diagram loosely represents the length of the Bott--Samelson bimodule being factored through.

Similarly, when $y < w$, we might omit the underline to indicate a morphism $X \to B_y$ or $B_y \to X$. A general morphism in $\catstuff{I}_{< w}(X,X)$ is a linear combination of morphisms of the form
\begin{gather} \label{eq:downupredux} \begin{tikzpicture}[anchorbase,scale=1]
\draw[mor] (0,0) to (0.25,0.5) to (0.75,0.5) to (1,0) to (0,0);
\node at (0.5,0.25){$\morstuff{g}_{y}$};
\draw[mor] (0,1) to (0.25,0.5) to (0.75,0.5) to (1,1) to (0,1);
\node at (0.5,0.75){$\morstuff{f}_{y}$};
\end{tikzpicture}.
\end{gather}

\begin{Lemma} \label{Lem:basiswithidentity} Let $X$ be a $w$-object. The hom space $\End(X)$ has a basis, as a right (or left) $\setstuff{R}$-module, where
one element of the basis is $\id_{X}$, and every other element of the basis is in $\catstuff{I}_{< w}$. Moreover, we have $\catstuff{I}_{< w}(X,X) = \catstuff{I}_{< \len(w)}(X,X)$. \end{Lemma}

\begin{proof} This follows immediately from the expanded Soergel hom formula, \autoref{thm:SHFplus}. \end{proof}

\begin{Notation} Let $X$ be a $w$-object.
Any basis of $\End(X)$ as a right $R$-module satisfying the criteria of \autoref{Lem:basiswithidentity} will be called a \emph{(right) basis with identity}. \end{Notation}

\begin{Notation} Let $\morstuff{h} \in \End(X)$, where $X$ is a $w$-object. The \emph{(right) identity coefficient} in $\morstuff{h}$, denoted $\idcoeff{h}$, is the coefficient of $\id_X$ when writing $\morstuff{h}$ with respect to any (right) basis with identity. This is an element of $\setstuff{R}$ of the same degree as $\morstuff{h}$. In particular, when $\morstuff{h}$ has degree zero, then $\idcoeff{h} \in \R = \R \cdot \id_R$. \end{Notation}

\begin{Lemma} The identity coefficient of $\morstuff{h}$ is well-defined, that is, it is independent of the choice of basis with identity. \end{Lemma}

\begin{proof} The quotient of $\End(X)$ by $\catstuff{I}_{< w}(X,X)$ is a free $R$-module of rank 1 spanned by the image of the identity map, and the identity coefficient of $\morstuff{h}$ can be viewed as the image of $\morstuff{h}$ under this quotient. \end{proof}

\begin{Lemma} The identity coefficient of a morphism is preserved by the duality functor (i.e. being flipped upside-down). \end{Lemma}

\begin{proof} The duality functor fixes the identity map and the $R$-module it generates, and sends a basis with identity to another basis with identity, from which the result is clear. \end{proof}

\begin{Remark} There is also the analogous notion of the left identity coefficient of $\morstuff{h}$. The left and right identity coefficients differ by an application of $w$,
since $X$ behaves like the standard bimodule $R_w$ modulo $\catstuff{I}_{< w}$. This statement follows from facts about standard filtrations; for calculations when $X$ is a
Bott--Samelson bimodule for a reduced expression, see \cite[Lemma 12.31]{ElMaThWi-soergel}. In particular, for morphisms of degree zero, the left and right identity
coefficients agree in $\R$. In this paper, we will view morphism spaces primarily as right $R$-modules. \end{Remark}

\begin{Lemma} \label{lem:zerosandwich} Let $X$ be a $w$-object, and let $\morstuff{e} \in \End(X)$ be a top idempotent. For any $\morstuff{h} \in \End^0(X)$ we have $\morstuff{e} \morstuff{h} \morstuff{e} = \idcoeff{h} \morstuff{e}$, that is
\begin{gather}
\begin{tikzpicture}[anchorbase,scale=1]
\draw[mor] (0,0) rectangle (1,0.5)node[black,pos=0.5]{$\morstuff{e}$};
\draw[mor] (0,1) rectangle (1,1.5)node[black,pos=0.5]{$\morstuff{e}$};
\draw[mor] (0,0.5) rectangle (1,1)node[black,pos=0.5]{$h$};
\end{tikzpicture}
=
\kappa(h) \;
\begin{tikzpicture}[anchorbase,scale=1]
\draw[mor] (0,0) rectangle (1,0.5)node[black,pos=0.5]{$\morstuff{e}$};
\end{tikzpicture}\; .
\end{gather}
\end{Lemma}

\begin{proof} Certainly this result is linear in $\morstuff{h}$, so it suffices to prove the result for $\morstuff{h} = \id$ (where it is obvious) and for $\morstuff{h} = \morstuff{f} \circ \morstuff{g}$, where this composition factors through $B_y$ with $y < w$. In this case we wish to show that $\morstuff{e} \morstuff{h} \morstuff{e}=0$. Note that $\deg(\morstuff{f}) + \deg(\morstuff{g}) = 0$, so one of those two morphisms has degree $\le 0$. However, there are no nonzero morphisms of degree $\le 0$ between $B_w$ and $B_y$ in either direction. Since $\Hom(X,B_y) \cdot \morstuff{e} \cong \Hom(B_w,B_y)$ and $\morstuff{e} \cdot \Hom(B_y,X) \cong \Hom(B_y,B_w)$, we deduce that either $\morstuff{e} \morstuff{f} = 0$ or $\morstuff{g} \morstuff{e} = 0$. \end{proof}

\begin{Corollary} \label{cor:topidempotentidcoeff} The identity coefficient of any top idempotent is $1$. \end{Corollary}

In the proof of \autoref{lem:zerosandwich} we argued that either $\morstuff{e} \morstuff{f} = 0$ or $\morstuff{g} \morstuff{e} = 0$, when $f \circ g$ has degree zero, but each one individually can be nonzero. One can not deduce from this lemma that $\morstuff{h} \morstuff{e} = 0$ for $\morstuff{h} \in \catstuff{I}^0_{< w}$, and indeed the result may fail, see \autoref{nastyexample}. Here and below we use $\catstuff{I}^0_{< w}$ to indicate degree zero morphisms in this ideal. For example, $\catstuff{I}^0_{< w}(\obstuff{B}_w,\obstuff{B}_w)= 0$. 

\subsection{Clasp idempotents}\label{subsection:clasps}

Now we are prepared to define our favorite kind of idempotent. We also denote these idempotents as rectangles.

\begin{Definition}\label{definition:clasp-idempotent}
Let $w \in W$ and $X$ be a $w$-object. An element $\morstuff{e} \in \End^0(X)$ is called a \emph{clasp idempotent} if it satisfies the following two properties.

\begin{enumerate}
\item The identity coefficient of $\morstuff{e}$ is $1$. That is, 
\begin{gather}\label{eq:jw3}
\begin{tikzpicture}[anchorbase,scale=1]
\draw[mor] (0,0) rectangle (1,0.5)node[black,pos=0.5]{$\morstuff{e}$};
\end{tikzpicture}
= \id_X
+
\sum \begin{tikzpicture}[anchorbase,scale=1]
\draw[mor] (0,0) to (0.25,0.5) to (0.75,0.5) to (1,0) to (0,0);
\draw[mor] (0,1) to (0.25,0.5) to (0.75,0.5) to (1,1) to (0,1);
\end{tikzpicture}.
\end{gather}

\item Precomposition with $\morstuff{e}$ annihilates $\catstuff{I}^{0}_{< w}(X,X)$. That is,
\begin{gather}\label{eq:annihilate_bot}
\begin{tikzpicture}[anchorbase,scale=1]
\draw[mor] (0,0) to (0.25,0.5) to (0.75,0.5) to (1,0) to (0,0);
\node at (0.5,0.25){$\morstuff{g}$};
\draw[mor] (0,1) to (0.25,0.5) to (0.75,0.5) to (1,1) to (0,1);
\node at (0.5,0.75){$\morstuff{f}$};
\draw[mor] (0,-0.5) rectangle (1,0)node[black,pos=0.5]{$\morstuff{e}$};
\end{tikzpicture}
= 0,
\end{gather}
whenever $\deg(\morstuff{f}) + \deg(\morstuff{g}) = 0$.
\end{enumerate}
Sometimes these are just called clasps for short.
\end{Definition}

\begin{Theorem} Let $w \in W$ and $X$ be a $w$-object. If a clasp idempotent $\morstuff{e} \in \End^0(X)$ exists, then it is unique, and it is an idempotent in the usual sense of the word:\begin{gather}
\begin{tikzpicture}[anchorbase,scale=1]
\draw[mor] (0,0) rectangle (1,0.5)node[black,pos=0.5]{$\morstuff{e}$};
\draw[mor] (0,0.5) rectangle (1,1)node[black,pos=0.5]{$\morstuff{e}$};
\end{tikzpicture}
=
\begin{tikzpicture}[anchorbase,scale=1]
\draw[mor] (0,0) rectangle (1,0.5)node[black,pos=0.5]{$\morstuff{e}$};
\end{tikzpicture}.
\end{gather} Moreover, postcomposition with $\morstuff{e}$ annihilates $\catstuff{I}^{0}_{< w}(X,X)$, that is \begin{gather}\label{eq:annihilate_top}
\begin{tikzpicture}[anchorbase,scale=1]
\draw[mor] (0,1) rectangle (1,1.5)node[black,pos=0.5]{$\morstuff{e}$};
\draw[mor] (0,0) to (0.25,0.5) to (0.75,0.5) to (1,0) to (0,0);
\node at (0.5,0.25){$\morstuff{g}$};
\draw[mor] (0,1) to (0.25,0.5) to (0.75,0.5) to (1,1) to (0,1);
\node at (0.5,0.75){$\morstuff{f}$};
\end{tikzpicture}
=
0.
\end{gather} whenever $\deg(\morstuff{f}) + \deg(\morstuff{g}) = 0$. Finally, $\morstuff{e}$ is preserved by $\DD$, i.e. being flipped upside-down. \begin{gather}
\begin{tikzpicture}[anchorbase,scale=1]
\draw[mor] (0,0) rectangle (1,0.5)node[black,pos=0.5]{$\morstuff{e}$};
\end{tikzpicture}
=
\begin{tikzpicture}[anchorbase,scale=1]
\draw[mor] (0,0) rectangle (1,0.5)node[black,pos=0.5,yscale=-1]{$\morstuff{e}$};
\end{tikzpicture}.
\end{gather} \end{Theorem}

\begin{Remark} Indeed, the proof indicates that we could have replaced \eqref{eq:annihilate_bot} with \eqref{eq:annihilate_top} and it would have given an equivalent definition for clasp idempotents. \end{Remark}

\begin{proof} Let $T \subset \End^0(X)$ be the subset of morphisms satisfying \eqref{eq:annihilate_bot}, and $U\subset \End^0(X)$ be the subset satisfying \eqref{eq:annihilate_top}. By assumption, there exists an element $\morstuff{e} \in T$ with identity coefficient $1$. Flipping $\morstuff{e}$ upside-down, we obtain an element $\morstuff{e}'$ of $U$ with identity coefficient $1$.

Let $\morstuff{g}$ be any element of $U$, having identity coefficient $\kappa \in \R$. Then 
\begin{equation} \morstuff{g} \circ \morstuff{e} = (\kappa \id) \circ \morstuff{e} + \catstuff{I}^{0}_{< w} \circ \morstuff{e} = \kappa \morstuff{e}. \end{equation}
But using the same argument to expand $\morstuff{e}$ rather than $\morstuff{g}$, we see that
\begin{equation} \morstuff{g} \circ \morstuff{e} =  \morstuff{g} \circ \id +  \morstuff{g} \circ \catstuff{I}^{0}_{< w} = \morstuff{g}. \end{equation}
Thus $\morstuff{g} = \kappa \morstuff{e}$, proving that $U$ is contained in the one-dimensional span of $\morstuff{e}$, and an element of $U$ is uniquely determined by its identity coefficient.

By the same argument, $T$ is contained in the span of $\morstuff{e}'$, and an element of $T$ is uniquely determined by its identity coefficient. Consequently, $\morstuff{e} =
\morstuff{e}'$, and $T = U$, and any element of $T$ or $U$ is preserved by being flipped upside-down. Moreover, the computation above implies that $\morstuff{e} \circ \morstuff{e}
= \morstuff{e}$, so it is idempotent. \end{proof}

Because clasp idempotents in $\End^0(X)$ are unique when they exist, we often abuse notation and label our rectangles with $X$ rather than $\morstuff{e}$, as in
\begin{equation*} \begin{tikzpicture}[anchorbase,scale=1]
\draw[mor] (0,0) rectangle (1,0.5)node[black,pos=0.5]{$\morstuff{e}$};
\end{tikzpicture}
\; =: \; \begin{tikzpicture}[anchorbase,scale=1]
\draw[mor] (0,0) rectangle (1,0.5)node[black,pos=0.5]{$X$};
\end{tikzpicture}.
\end{equation*}
When $X = \obstuff{B}_{\un{w}}$ for $\un{w} \in \Rex(w)$, or when $X = B_x B_i$ for $w = xi>x$, we shorten this to
\begin{equation*} \begin{tikzpicture}[anchorbase,scale=1]
\draw[mor] (0,0) rectangle (1,0.5)node[black,pos=0.5]{$\morstuff{e}$};
\end{tikzpicture}
\; =: \; \begin{tikzpicture}[anchorbase,scale=1]
\draw[mor] (0,0) rectangle (1,0.5)node[black,pos=0.5]{$\un{w}$};
\end{tikzpicture}, \qquad
\begin{tikzpicture}[anchorbase,scale=1]
\draw[mor] (0,0) rectangle (1,0.5)node[black,pos=0.5]{$\morstuff{e}$};
\end{tikzpicture}
\; =: \; \begin{tikzpicture}[anchorbase,scale=1]
\draw[mor] (0,0) rectangle (1,0.5)node[black,pos=0.5]{$x \otimes i$};
\end{tikzpicture}.
\end{equation*}
Note however that $X = \obstuff{B}_{\un{w}}$ often does not admit a clasp idempotent.

Let us pause to note that \eqref{eq:jw3} is equivalent to $\morstuff{e}$ being idempotent, given the other condition.

\begin{Lemma}\label{lem:charac_clasp_idemp} Let $X$ be a $w$-object. Let $e \in \End^0(X)$ satisfy \eqref{eq:annihilate_bot}. Then $e$ is a clasp idempotent if and only if $e$ is a nonzero idempotent. \end{Lemma}

\begin{proof} One direction follows from the previous theorem. Now suppose $e$ is a nonzero idempotent. The image of $e$ modulo $\catstuff{I}_{< w}$ is still an idempotent, and consequently $\idcoeff{e}$ is an idempotent in $\R$. However, if $\idcoeff{e} = 0$ then $e \in \catstuff{I}_{< w}$, so $e = e^2 = 0$ by \eqref{eq:annihilate_top}. This is a contradiction, so we must have $\idcoeff{e} = 1$. \end{proof}

\begin{Theorem} \label{thm:alltopareclasp} Let $w \in W$ and $X$ be a $w$-object. Suppose that a clasp idempotent $\morstuff{e}$ exists. Then any top idempotent is equal to $\morstuff{e}$.  \end{Theorem}

\begin{proof} Let $\morstuff{e}'$ be a top idempotent.
Extend $\morstuff{e}'$ to an idempotent decomposition of the identity, so that
\begin{equation} \label{decompwithq} \id_{X} = \morstuff{e}' + \sum \morstuff{q}_j, \end{equation}
where $\morstuff{q}_j$ are (mutually orthogonal) idempotents whose images are isomorphic to $v^k \obstuff{B}_y$ for various $y < w$ and $k \in \Z$.  Then
\begin{equation} \morstuff{e} = \morstuff{e} \circ \id = \morstuff{e} \circ \morstuff{e}' + \sum \morstuff{e} \circ \morstuff{q}_j = \morstuff{e} \circ \morstuff{e}'. \end{equation}
The last equality holds since $\morstuff{q}_j \in \catstuff{I}^0_{< w}$. In particular, this implies that the identity coefficient of $\morstuff{e}'$ is $1$.

By a similar argument (precomposing with $\id$) we have $\morstuff{e} = \morstuff{e}' \morstuff{e}$. Thus
\begin{equation} \morstuff{e} = \morstuff{e} \morstuff{e}' = \morstuff{e}' \morstuff{e} = \morstuff{e}' \morstuff{e} \morstuff{e}'. \end{equation}

We can identify $\morstuff{e}' \End(X) \morstuff{e}'$ with $\End(\obstuff{B}_w)$ in the usual way. By Soergel's hom formula (see \autoref{lem:deg012}),
$\End^0(\obstuff{B}_w)$ is spanned by the identity map, so the only nonzero idempotent in $\End(\obstuff{B}_w)$ is the identity, which is identified with the idempotent $\morstuff{e}'$ in $\End^0(X)$. Hence
$\morstuff{e}'\morstuff{e}\morstuff{e}' = \morstuff{e}'$, from which we deduce $\morstuff{e} = \morstuff{e}'$. \end{proof}

\subsection{Existence of clasp idempotents versus uniqueness of top idempotents}\label{subsection:topunique}

Beyond this point, the theory of clasp idempotents is more complicated than other similar theories in the literature (e.g. clasp idempotents for webs) because of the existence of
morphisms of nonzero degree between distinct indecomposable objects.

Now we discuss the question of uniqueness of top idempotents and the possible non-existence of clasp idempotents, as introduced in \autoref{nastyexample}. Note that none of the idempotents $e' = e + \lambda f$ from \autoref{nastyexample}  are clasp idempotents, since $f e' = f \ne 0$.

Let us explain the issue for a general $w$-object $X$, following the prototype of \autoref{nastyexample}. There are no morphisms of degree $\le 0$ between $B_w$ and $B_y$ for $w \ne
y$. However, there are morphisms of positive degree. Suppose that $B_y$ appears within $X$ with a sufficiently large shift. One might then hypothesize the existence of a degree zero
morphism $B_w \to X$ which factors as a positive degree morphism $g \colon B_w \to B_y$ followed by a negative degree inclusion $\iota \colon B_y \to X$. Precompose $g$ with the
projection $X \to B_w$ induced by a given top idempotent to obtain $\tilde{g}$, and let $f = \iota \circ \tilde{g}$. Because projections (resp. inclusions) are epic (resp. monic), $f$
is nonzero. Then $f$ lives in $\catstuff{I}^0_{< w}$ but is not annihilated by the top idempotent, whence the top idempotent is not a clasp. Using the exact same arguments as in
\autoref{nastyexample}, $\{e + \lambda f\}$ will be a family of top idempotents with isomorphic images.

The following theorems state that this is exactly the obstruction to the existence of clasp idempotents.

\begin{Theorem} \label{thm:claspvsunique} Let $w \in W$ and $X$ be a $w$-object. The following are equivalent. \begin{enumerate}
\item The clasp idempotent exists.
\item $\catstuff{I}^0_{< w}(B_w, X) = 0$.
\item The top idempotent is unique.
\end{enumerate} \end{Theorem}

\begin{proof} \autoref{thm:alltopareclasp} shows that (1) implies (3) . We now prove that (3) implies (2), and (2) implies (1).

If $\catstuff{I}^0_{< w}(B_w, X) \ne 0$ then we can choose a nonzero morphism  $\iota \circ g$ therein, factoring through $B_y$ for some $y < w$. The map $g \colon B_w \to B_y$ must have positive degree by Soergel's hom formula, and $\iota$ therefore must have negative degree. Thus $e \circ \iota = 0$, since $\Hom(B_y,B_w)$ is zero in negative degrees. As above, precompose $\iota \circ g$ with the projection $X \to B_w$ to obtain $f \in \catstuff{I}^0_{< w}(X, X)$ such that $fe = f$ and $ef= 0$. Then $\{e + \lambda f\}$ is another top idempotent for any $\lambda \in \R$. By the contrapositive, (3) implies (2).

Now suppose that $\catstuff{I}^0_{< w}(B_w,X) = 0$. Let $e$ be any top idempotent. Let $g \in \catstuff{I}^0_{< w}(X,X)$. Then $ge \in \End(X)e$ can be identified with some morphism in $\Hom(B_w,X)$ which still lies in the ideal $\catstuff{I}_{< w}$, and hence $ge = 0$. Thus $e$ satisfies the second property of a clasp idempotent. By \autoref{lem:charac_clasp_idemp} we see that $e$ is a clasp idempotent. \end{proof}

Now we point out a crucial feature of this story: whether or not $\catstuff{I}^0_{< w}(B_w, X) = 0$ is determined in the Grothendieck group! Thus we can give a decategorified criterion for the existence of a clasp idempotent.

\begin{Theorem} \label{Thm:criterionforclasp} Let $X$ be a $w$-object. For all $x \le y$ in $W$ let $t_{y,x}$ be the lowest exponent of $v$ such that $v^{t_{y,x}}$ appears with nonzero coefficient in $p_{y,x}$. In particular, $t_{y,y} = 0$, and $t_{y,x} > 0$ for $x < y$. For a given $y < w$ such that some shift of $B_y$ is a summand of $X$, let $a_y \ge 0$ denote the highest exponent of $v$ such that $v^{a_y} B_y$ is a summand of $X$, and such that $a_y$ has the same parity as $t_{w,x} + t_{y,x}$. (Since $X$ is self-dual, $v^{-a_y} B_y$ is also a summand.) Then a clasp idempotent exists if and only if $t_{w,x} + t_{y,x} - a_y > 0$ for all $x \le y < w$ where $y$ appears in the decomposition. \end{Theorem}

\begin{proof} We use \autoref{thm:SHFplus} throughout this proof.

By the previous theorem, a clasp idempotent exists if and only if $\catstuff{I}^0_{< w}(B_w, X) = 0$. We will find a spanning set of $\catstuff{I}_{< w}(B_w,X)$ as a right $R$-module, and analyze the degrees in that spanning set. If some degree is even and non-positive, we can multiply by an element of $R$ to get a degree zero morphism in $\catstuff{I}_{< w}(B_w,X)$, so this space is nonzero. If all degrees are odd or positive, then $\catstuff{I}^0_{< w}(B_w, X) = 0$.

Any morphism $B_w \to X$ is spanned by morphisms which factor as a composition $B_w \to B_y \to X$ for some $y \le w$, where $B_y(a)$ (up to shift) is one of
the summands of $X$, and $B_y \to X$ is the inclusion map of that summand. In particular, the map $B_y \to X$ has degree $-a$. Meanwhile, by the expanded Soergel hom formula, the map
$B_w \to B_y$ is spanned over $R$ by maps factoring as $B_w \to B_x \to B_y$ for some $x \le y$, where the map $B_w \to B_x$ has degree at least $t_{w,x}$ (and has the same parity as
$t_{w,x}$) and the map $B_x \to B_y$ has degree at least $t_{y,x}$ (and has the same parity as $t_{y,x}$). Some nonzero map $B_w \to B_y$ of minimal degree $t_{w,x} + t_{y,x}$ does exist, and since inclusion maps are monic, the composition with $B_y \to X$ is still nonzero.
Thus the overall degree of the composition $B_w \to X$ is at least $t_{w,x} + t_{y,x} - a$, and the minimal degree is attained. Now the
requirement that $t_{w,x} + t_{y,x} - a_y > 0$ is equivalent to the nonexistence of a spanning morphism of even non-positive degree. This morphism is in $\catstuff{I}_{< w}(B_w,X)$ if and only if $x < w$. Note finally that the case $y=w$ can be safely excluded, as $t_{w,x} > 0$ and $a_w = 0$. \end{proof}

\begin{Corollary} \label{Cor:perversehasclasp} Let $X$ be a $w$-object which is perverse (see \autoref{def:perverse}). Then $X$ admits a clasp idempotent. \end{Corollary}

\begin{proof} This follows from the previous theorems, since $a_y = 0$ for all $y$ which appear in the decomposition, and $t_{w,x} > 0$ and $t_{y,x} \ge 0$ for all $x \le y < w$. \end{proof}

\begin{Corollary} \label{Cor:xihasclasp} Let $x \in W$ and $i \in S$ such that $w = xi > x$. Then there is a clasp idempotent in $\End(B_x B_i)$ with image isomorphic to $B_{w}$.\end{Corollary}

\begin{proof} By \eqref{inductivedecomp}, $B_x B_i$ is a perverse $w$-object. \end{proof}	

\subsection{Subtleties with Karoubi envelopes} \label{subsection:trickytop}

An unintuitive consequence of \autoref{Cor:xihasclasp} is that non-clasps can often be used in inductive constructions of clasps.

\begin{Example} We illustrate this phenomenon in a continuation of \autoref{extendingnastypart1}. Recall that $\un{w} = (3,1,2,1,3)$, and the idempotent $e \in \End(B_{\un{w}})$ was not a clasp, as $ge \ne 0$ for a particular endomorphism $g$ in the ideal of lower terms. Concatenating $\un{w}$ with the simple reflection $2$, we get a reduced expression $\un{w}_0$ for the longest element $w_0 = 312132$. Define $e_0$ as in \eqref{idempotente0}. 
We can view $e_0$ as either an element of $\End(B_{\un{w}_0})$ or an element of 
\[ \End(B_w \ot B_2) \cong (e \ot \id_2) \End(B_{\un{w} 2}) (e \ot \id_2). \]
Indeed, $e_0$ is a clasp in either endomorphism ring. One can compute in the Grothendieck group that
\[ B_{(3,1,2,1,3,2)} \cong B_{w_0} \oplus B_{3121} \oplus B_{1323} \oplus B_{132}(1) \oplus B_{132}(-1). \]
Unlike \autoref{nastyexample}, there are no morphisms of degree $1$ from $B_{w_0}$ to $B_{132}$, as $p_{w_0,132} = v^3$. Thus there are no degree zero morphisms from $B_{w_0}$ to this Bott--Samelson which factor through lower terms, and $e_0$ is a clasp in $\End(B_{\un{w}_0})$ by \autoref{Thm:criterionforclasp}.
The reader can verify that the orthogonality issue also disappears, as $(g \ot \id_2) e_0 = 0$. 
\end{Example}

We need to emphasize a very confusing subtlety. Let $e$ be a top idempotent for $B_x$ in some $x$-object $X$, not necessarily a clasp idempotent. Suppose that $w = xi > x$ for some $i \in S$. Then $X B_i$ is a $w$-object, which
need not admit a clasp. However, $B_x B_i$ is a summand, so that $\End(B_x B_i)$ is a summand of $\End(X B_i)$; more precisely, \[ \End(B_x \ot B_i) \cong (e \ot \id_i) \End(X B_i) (e
\ot \id_i) \subset \End(X B_i). \] Thus one can find an idempotent $e_w \in \End(X B_i)$ living inside the summand $\End(B_w B_i)$, such that $e_w$ is a clasp idempotent inside
$\End(B_w B_i)$, though the ``same'' morphism $e_w$ might not be a clasp idempotent within $\End(X B_i)$. In the previous example, it happened to be a clasp idempotent within $X B_i$
as well, but this is not guaranteed.

\begin{Example} The original idempotent $e$ from \autoref{nastyexample} can be constructed as $e = e_{3121} \ot \id_3$ for a clasp idempotent $e_{3121} \in \End(B_{(3,1,2,1)})$. While $e$ is a not a clasp in $\End(B_{(3,1,2,1,3)})$, it is a clasp in another space, namely in $\End(B_{3121} B_3)$. \end{Example}

Whether a top idempotent is a clasp idempotent is determined by the object it is an endomorphism of, purely using the Grothendieck group. One drawback of diagrammatics when
working with Karoubi envelopes is that the source and target of a morphism is not always diagrammatically evident. If $e,f \in \End(X)$ with $e^2 = e$ and $f = efe$, then $f$ can be
interpreted either as an endomorphism of $X$ or an endomorphism of the image of $e$, and this choice is not clear just from the picture, even when idempotents $e$ are explicitly
drawn in the picture. When $f$ is a top idempotent, whether or not $f$ is a clasp may depend on the invisible choice of source and target.

\subsection{Elaborations on orthogonality} \label{subsection:trickytop2}

We wish to provide another criterion for a morphism to be a clasp, and also to ward off a major potential confusion, both arising from the subtly different approach typically taken in
the literature. We have defined clasps in terms of orthogonality to lower terms within $\End(X)$, or more precisely within $\End^0(X)$, see \eqref{eq:annihilate_bot}. Meanwhile, the
literature on clasps (in non-graded contexts) does not focus on endomorphisms, but would define a clasp in terms of orthogonality to lower terms in general. One might ask if
there is a replacement for \eqref{eq:annihilate_bot} where one only insists that $ge = 0$ for maps $g$ (non necessarily endomorphisms) in the ideal of lower terms. Let us demonstrate a few false starts along these lines before stating our theorem.

We say that a morphism $e \in \End^0(X)$ is \emph{truly orthogonal to lower terms} if $ge = 0$ for all $g \in \catstuff{I}_{<w}(X,Y)$ for any object $Y$. It is simply false
that clasp idempotents are truly orthogonal to lower terms! For example, $\id_{B_w}$ is a clasp inside $X = B_w$, but there are nonzero maps (\emph{of positive degree}) between
$B_w$ and $B_y$ for all $y < w$.

Inspired by this, we might say that a morphism $e \in \End^0(X)$ is \emph{orthogonal to non-positive lower terms} if $ge = 0$ for all $g \in \Hom^{\le
0}(X,B_y)$, for all $y < w$. However, every top idempotent satisfies this property, since $\Hom^{\le 0}(B_w, B_y) = 0$, so this is hardly a suitable criterion for a morphism to be a clasp idempotent; it is not a replacement for \eqref{eq:annihilate_bot}. Nonetheless, when a clasp idempotent exists, a top idempotent is a clasp idempotent, so perhaps this is a useful criterion after all.

One might imagine that orthogonality to non-positive lower terms might replace \eqref{eq:annihilate_bot} when the numerical criteria of \autoref{Thm:criterionforclasp} are satisfied. This too is false, as other non-top idempotents satisfy this criterion as well. We give an example of this phenomenon in the setting of \autoref{nastyexample} because it has already been introduced, acknowledging that clasp idempotents do not exist here; a completely analogous construction gives an example in the setting of \autoref{Ex:nonperverseclasp}, where clasp idempotents do exist.

\begin{Example} We continue \autoref{nastyexample}. Let $\id_{\un{w}} = e + q_{-1} + q_1$ be a decomposition of $B_{\un{w}}$ into orthogonal idempotents. Here
$q_{-1} = \rho \circ \iota$, where $\iota$ is the map of degree $-1$ from \autoref{nastyexample}, and $\rho$ has degree $+1$, so that $q_{-1}$ is an idempotent projecting to
$B_y(-1)$. Meanwhile, $q_1 = \bar{\iota} \circ \rho'$ projects to $B_y(1)$, where $\bar{\iota}$ is the vertical flip of $\iota$. We have that $\iota \circ q_1 = 0$. In particular,
$gq_1 = 0$ for all $g \colon B_{\un{w}} \to B_y$ of degree $\le 0$. The idempotent $e + q_1$ has identity coefficient $1$ and is orthogonal to non-positive lower terms, but is
not a top idempotent. \end{Example}

Instead, a genuine criterion for clasphood along these lines can be achieved, but one must be orthogonal to some positive degree lower terms as well.

\begin{Theorem} \label{thm:orthogonalitycriterion} Let $X$ be a $w$-object admitting a clasp idempotent, and define $a_y$ as in \autoref{Thm:criterionforclasp}. Then an endomorphism $\morstuff{e} \in \End^0(X)$ is a clasp idempotent if and only if it has identity coefficient $1$, and $ge = 0$ for all $g \colon X \to B_y$ of degree $\le a_y$, for all $y$ appearing in the decomposition of $X$. \end{Theorem}

\begin{proof} Suppose that $ge = 0$ for all maps $g \colon X \to B_y$ of degree $\le a_y$, as in the statement. Fix a top idempotent $e'$ and a corresponding decomposition of $X$. Then $q_j e = 0$ for all the idempotents $q_j$ in a decomposition \eqref{decompwithq}, since the projection maps to $B_y$ all have degree $\le a_y$. Thus $e = e' e$. Since $e'$ is a clasp idempotent by \autoref{thm:alltopareclasp}, $e' e = \kappa e'$ where $\kappa$ is the identity coefficient of $e$, which is assumed to be $1$. Hence $e = e' e = e'$.

Conversely, suppose that $e$ is a clasp idempotent. Then it is a top idempotent, with image isomorphic to $B_w$. If $g \colon X \to B_y$ then $ge \in \Hom(X,B_y) e \cong \Hom(B_w,
B_y)$ corresponds to a map $B_w \to B_y$ of the same degree. We claim that $\Hom(B_w,B_y)$ is zero in degrees $\le a_y$, which proves that $ge=0$. Define $t_{y,x}$ as in
\autoref{Thm:criterionforclasp}. The expanded Soergel hom formula gives a basis for maps $B_w \to B_y$ where each basis element factors through $B_x$ for some $x \le y$, and has
degree at least $t_{w,x} + t_{y,x}$. Since $t_{w,x} + t_{y,x} > a_y$ by \autoref{Thm:criterionforclasp}, we are done. \end{proof}

\begin{Corollary} \label{Cor:perverseorthocriterion} If $X$ is a perverse $w$-object and $e \in \End^0(X)$ has identity coefficient $1$ and is orthogonal to non-positive lower terms, then $e$ is a clasp idempotent. \end{Corollary}

\begin{proof} A clasp idempotent exists by \autoref{Cor:perversehasclasp}. Now the result follows from the previous theorem, as $a_y = 0$ for all $y$ appearing in the
decomposition. \end{proof}

\subsection{Top summands and $1$-tensors} \label{subsection:top1tensor}

The bimodule $\obstuff{B}_i$ is supported in degrees $\ge -1$. It has an element $1 \otimes 1$ living in degree $-1$, and spanning this negative-most (nonzero) homogeneous component.
We call it the \emph{$1$-tensor}, and denote it $\onetensor$. More generally, a tensor product $\obstuff{B}_{i_1}\obstuff{B}_{i_2}\cdots \obstuff{B}_{i_d}$ has an element in degree $-d$, which is the tensor
product of the various $1$-tensors from each factor. We call this element the \emph{$1$-tensor}, and denote it $\onetensor$. It also spans the negative-most homogeneous component. The theory of 1-tensors was developed and their relationship to the ideals $\catstuff{I}_{< w}$ explored in \cite[Section 3.8]{KELP2}.

\begin{Lemma} \label{Lem:multivalentpreserves1tensor} The multivalent vertices send $\onetensor$ to $\onetensor$. \end{Lemma}

\begin{proof} This is by definition of the functor from diagrammatics to bimodules, see \cite[Definition 4.6]{ElWi-soergel-calculus}. Indeed, the multivalent vertices live in a one-dimensional homogeneous component within their hom space, and are normalized precisely for this lemma to hold. \end{proof}

\begin{Lemma} \label{Lem:onetensorandstuff} Let $\un{w}$ be a reduced expression for $w$. In any decomposition of $\obstuff{B}_{\un{w}}$ into indecomposable summands, the top summand $\obstuff{B}_w$ will contain the 1-tensor, which spans the negative-most homogeneous component of $\obstuff{B}_w$. If $\morstuff{h} \in \Hom^0(\obstuff{B}_{\un{w}},\obstuff{B}_{\un{w}})$, then $\morstuff{h}(\onetensor) = \idcoeff{h} \cdot \onetensor$. Top idempotents preserve the 1-tensor. \end{Lemma}

\begin{proof} The first statement is well-known, and proving it amounts to proving that $\obstuff{B}_w$ is nonzero (and hence one-dimensional) in degree $-\len(w)$. This can be shown in several ways, e.g. using the standard filtration on $\obstuff{B}_w$, see \cite[Satz 6.16]{So-kl-bimodules}.

For the second statement, by additivity we reduce to the case when $\morstuff{h}$ is the identity (obvious), or when $\morstuff{h} \in \catstuff{I}^0_{< w}$. But $\catstuff{I}^0_{< w}(\obstuff{B}_w,\obstuff{B}_w)= 0$, so $\morstuff{h}$ must induce the zero map on the summand $\obstuff{B}_w$, which contains $\onetensor$. The third statement follows from the second statement and \autoref{cor:topidempotentidcoeff}. \end{proof}

One is tempted to define the 1-tensor as a vector in the summand $\obstuff{B}_w$. This is not really a well-defined vector, as it depends on a choice of projection map from
$\obstuff{B}_{\un{w}}$ to $\obstuff{B}_w$. Rescaling the projection and inclusion maps, one will rescale the choice of 1-tensor (this process does not change the top idempotent, and is an issue even when clasp idempotents exist). This seems like an annoyance, but actually, it is the solution to a different annoyance, because this paragraph began with a misleading premise. We now explain.

Recall that $\obstuff{B}_w$ was defined up to isomorphism as the unique summand in a decomposition of $\obstuff{B}_{\un{w}}$ (for a reduced expression $\un{w}$) not appearing
in a shorter Bott--Samelson bimodule, and it was proven by Soergel that any two such summands (for any decomposition, or any reduced expression) are isomorphic. This definition produces an isomorphism class of bimodules, not a specific bimodule. One could fix a representative for this isomorphism class by choosing a particular $\un{w}$ and a particular
decomposition. A better and more flexible fix is to determine a (compatible) family of isomorphisms between the summand produced by any two choices, which we use to identify them with
each other. See e.g. \cite[Definition 2.13 and surrounding discussion]{El-thick-soergel-typea}. Note that $\End^0(\obstuff{B}_w)$ is one-dimensional, spanned by the identity map, so the degree zero morphism space between any two such summands will be one-dimensional, and the isomorphism between them is determined up to invertible scalar.

\begin{Notation} If we are given two reduced expressions $\un{w}$ and
$\un{w}'$ and decompositions of their respective Bott--Samelson bimodules, we always identify their top summands via the unique isomorphism which sends the 1-tensor to the
1-tensor. We call this isomorphism the \emph{canonical isomorphism}. \end{Notation}

In so doing, the $1$-tensor in $\obstuff{B}_w$ becomes a well-defined vector. Said another way, the bimodule $\obstuff{B}_w$ is defined up to isomorphism, while the pair consisting of $\obstuff{B}_w$ and its $1$-tensor is defined up to unique isomorphism.

If two reduced expressions are related by a braid relation, then the multivalent vertex induces a degree zero morphism between their respective top summands. This morphism is
non-zero, because multivalent vertices are isomorphisms modulo $\catstuff{I}_{< w}$. Since it preserves the 1-tensor, the multivalent vertex induces the canonical
isomorphism between these top summands.

We want to extend these concepts to other $w$-objects beyond Bott--Samelson bimodules.

\begin{Definition} A \emph{based $w$-object} will be a $w$-object $X$ (with its tacitly fixed decomposition), together with a homogeneous vector in $X$ called the 1-tensor, such that
the 1-tensor is contained within the top summand of $X$. \end{Definition}

If one has a vector in an $R$-bimodule $M$ which one calls the $1$-tensor, and similarly for $M'$, then the tensor product of these 1-tensors will be called the 1-tensor in $M
\otimes M'$. For example, we have a well-defined $1$-tensor in $\obstuff{B}_w \otimes \obstuff{B}_x$ for any $w, x \in W$. Of particular importance is the 1-tensor in $B_x B_i$, which is a $w$-object when $w = xi > x$.

We now repeat \autoref{Lem:onetensorandstuff} in this context.

\begin{Lemma}  Let $X$ be a based $w$-object. The top summand $\obstuff{B}_w$ will contain the 1-tensor, which spans its negative-most homogeneous component. If $\morstuff{h} \in \Hom^0(X,X)$, then $\morstuff{h}(\onetensor) = \idcoeff{h} \cdot \onetensor$. The top idempotent preserves the 1-tensor. \end{Lemma}

\begin{proof} This is largely tautological or the same as \autoref{Lem:onetensorandstuff}, and left to the reader. \end{proof}

\begin{Definition} Given two based $w$-objects $X$ and $X'$, the \emph{canonical isomorphism} between their top summands is the one which sends the 1-tensor to the 1-tensor. The
\emph{canonical morphism} is the degree zero map $X \to X'$ which projects to the top summand, applies the canonical isomorphism, and then includes from the top summand. \end{Definition}

Recall that a $w$-object has an implicit decomposition, which are used in the definition of the canonical (iso)morphism. Canonical morphisms clearly preserve the 1-tensor. Canonical morphisms are unique, as they live in a one-dimensional space (isomorphic to $\End^0(B_w)$), and are determined within this space by the fact that they preserve the 1-tensor.

\begin{Lemma} \label{Lem:canonmapcrit} A degree zero morphism $N \colon X \to X'$ between based $w$-objects is the canonical morphism if and only if $N = e_{X'} N e_X$ and $N$ preserves the 1-tensor. \end{Lemma}

\begin{proof} A morphism $N \colon X \to X'$ factors through the top summand if and only if it satisfies $N = e_{X'} N e_X$. Since $\Hom^0(B_w,B_w)$ is one-dimensional, any such $N$ is a scalar multiple of the canonical morphism, and that scalar multiple is determined by its action on the 1-tensor. \end{proof}

\begin{Corollary} If a degree zero morphism $N \colon X \to X'$ between based $w$-objects preserves the $1$-tensor, then $e_{X'} N e_X$ is the canonical morphism. \end{Corollary}

\begin{proof} Since top idempotents preserve the 1-tensor, clearly $e_{X'} N e_X$  satisfies the criteria of \autoref{Lem:canonmapcrit}. \end{proof}

In the case of the previous corollary, we might say that $N$ \emph{induces} the canonical morphism.

\begin{Lemma} A composition of canonical (iso)morphisms is a canonical (iso)morphism. \end{Lemma}

\begin{proof} The proof is easy, and we mostly write it down to elaborate on the statement. Suppose for example that $N_1 \colon X \to X'$ and $N_2 \colon X' \to X''$ are canonical morphisms between based $w$-objects. By \autoref{Lem:canonmapcrit}, $e_{X''} N_2 = N_2$ and $N_1 e_X = N_1$, and $N_2 \circ N_1$ preserves the 1-tensor. Again using \autoref{Lem:canonmapcrit}, we deduce that $N_2 \circ N_1$ is a canonical morphism. \end{proof}

In fact, a composition of morphisms which preserve the 1-tensor is another morphism which preserves the 1-tensor, so the lemma also holds for morphisms which induce the canonical morphism.

In this paper, we will typically draw morphisms which induce the canonical morphism (i.e. morphisms which preserve the 1-tensor) as purple strips. These purple strips will typically be composed with idempotents, so that this composition is a genuine canonical morphism. For example, let $N \colon X \to X'$ induce the canonical morphism. Then the following diagram represents the canonical morphism $X \to X'$.
\begin{equation} \begin{tikzpicture}[anchorbase,scale=1]
\draw[mor] (0,-0.5) rectangle (1,0) node[black,pos=0.5]{$X$};
\draw[rex] (0,0) rectangle (1,0.1);
\draw[mor] (0,0.1) rectangle (1,0.6) node[black,pos=0.5]{$X'$};
\end{tikzpicture} \quad := \quad \begin{tikzpicture}[anchorbase,scale=1]
\draw[mor] (0,0) rectangle (1,0.5) node[black,pos=0.5]{$X$};
\draw[mor] (0,0.5) rectangle (1,1) node[color=darktomato,pos=0.5]{$N$};
\draw[mor] (0,1) rectangle (1,1.5) node[black,pos=0.5]{$X^{\prime}$};
\end{tikzpicture}. \end{equation}
Recall that the rectangles labeled $X$ and $X'$ represent their corresponding top idempotents $e_X$ and $e_{X'}$. In the larger context of this diagram the actual morphism represented by the purple strip is irrelevant, since the canonical morphism is unique.

When one of the top idempotents is a clasp, then one can simplify the diagram above.

\begin{Lemma} \label{Lem:ifclaspremovee} When $e_{X}$ is a clasp, and $N$ induces the canonical morphism, then $e_{X'} N e_{X} = N e_{X}$, or diagrammatically:
\begin{equation}
\begin{tikzpicture}[anchorbase,scale=1]
\draw[mor] (0,-0.5) rectangle (1,0) node[black,pos=0.5]{$X$};
\draw[rex] (0,0) rectangle (1,0.1);
\draw[mor] (0,0.1) rectangle (1,0.6) node[black,pos=0.5]{$X^{\prime}$};
\end{tikzpicture} \quad = \quad 
\begin{tikzpicture}[anchorbase,scale=1]
\draw[mor] (0,-0.5) rectangle (1,0) node[black,pos=0.5]{$X$};
\draw[rex] (0,0) rectangle (1,0.1);
\end{tikzpicture}\; .
\end{equation} \end{Lemma}

\begin{proof} Note that $e_{X'} = \id_{X'} - \sum q_j$, where $q_j$ are idempotents projecting to lower summands. Each $q_j \circ N$ lives in $\catstuff{I}^0_{< w}$, so it is orthogonal to $e_{X}$. Hence $e_{X'}$ can be safely replaced with the identity of $X'$. \end{proof}

\subsection{Idempotents and rex moves}\label{subsection:idempotentsrex}

It is essential to note the following issue, for which we return to the setting of Bott--Samelson bimodules. A diagram built only by tensoring and composing multivalent vertices (with sources and targets as in \eqref{generatorsofD}) and identity maps is called a \emph{rex move}. All rex moves preserve the 1-tensor.

\begin{Lemma} \label{Lem:rexgivesnewtop} Suppose $\un{w}, \un{w}' \in \Rex(w)$, and $N_1 \colon \obstuff{B}_{\un{w}}\to \obstuff{B}_{\un{w}'}$ and $N_2 \colon \obstuff{B}_{\un{w}'}\to \obstuff{B}_{\un{w}}$ are rex moves. Let $\morstuff{e}$ be a top idempotent in $\obstuff{B}_{\un{w}}$. Then $e' := N_1 e N_2 \in \End(B_{\un{w}'})$ is a top idempotent. \end{Lemma}

\begin{proof} First note that $e'$ is an idempotent, as $e N_2 N_1 e = e$ by \autoref{lem:zerosandwich} and \autoref{Lem:onetensorandstuff}. It is a top idempotent, because $e' N_1 e$ and $e N_2 e'$ give inverse isomorphisms between the images of $e$ and $e'$. \end{proof}

In this way, given a single top idempotent for one reduced expression, one can produce a large family of top idempotents for all reduced expressions, obtained from the first by applying rex moves.

For example, given any rex move $N$ which is an endomorphism of $B_{\un{w}}$, $Ne$ is another top summand of $B_{\un{w}}$. We have $\morstuff{e} \neq N \morstuff{e}$ in
general! There is no requirement that top idempotents be preserved by rex moves. On the other hand, it was proven in \cite[Lemma 7.4]{ElWi-soergel-calculus} that all rex moves with the same source and target agree modulo lower terms (the cited lemma uses different terminology; their morphisms with strictly negative-positive decompositions can be rephrased as our $\catstuff{I}_{< \len(w)}$), so in particular $N$ agrees with the identity modulo lower
terms (i.e. the identity coefficient of $N$ is $1$). If a clasp idempotent $\morstuff{e}$ exists, then $N \morstuff{e} = \morstuff{e}$ (which we also knew from the uniqueness of top
idempotents)! Thus rex move endomorphisms preserve clasp idempotents.

More generally, if $N_1'$ and $N_2'$ are different rex moves with the same sources and targets as $N_1$ and $N_2$ above, then it is possible that $N_1 e N_2 \neq N_1' e N_2'$ for a top idempotent $e$. If $e$ is a clasp idempotent, then $N_1 e N_2 = N_1' e N_2'$.

In the setting of \autoref{Lem:rexgivesnewtop}, when $e$ is a clasp idempotent one might expect that $N_1 e N_2$ is also a clasp idempotent. This is not always the case, for the
reason that clasp idempotents might not exist! Since $N_1 e N_2$ is a top idempotent, it is a clasp when the clasp exists.

\begin{Example} Continue the setting of \autoref{nastyexample}.  Let $\un{w} = (3,1,2,1,3)$ and $\un{w}' = (3,2,1,2,3)$. A calculation shows that $B_{\un{w}'}$ is a perverse $w$-object, and thus admits a clasp, whereas we have seen that $B_{\un{w}}$ does not admit a clasp. \end{Example}

Finally, we note the relationship between rex moves and canonical morphisms.

\begin{Lemma} Suppose $\un{w}, \un{w}' \in \Rex(w)$, and $N \colon \obstuff{B}_{\un{w}}\to \obstuff{B}_{\un{w}'}$. Let $e_{\un{w}}$ and $e_{\un{w}'}$ be the top idempotents. Then $e_{\un{w}'} N e_{\un{w}} \colon \obstuff{B}_{\un{w}}\to \obstuff{B}_{\un{w}'}$ is the canonical morphism, and is thus independent of the choice of rex move $N$. \end{Lemma}

\begin{proof} This follows immediately from \autoref{Lem:canonmapcrit}. \end{proof}

When rex moves function as (morphisms which induce) canonical morphisms, we will denote them as purple strips as above. For example, when $\un{w}$ and $\un{w}'$ are reduced expressions for $w$ and we have chosen top summands within their respective Bott--Samelson bimodules, then for any rex move $N$ between these Bott--Samelsons we can denote it as a purple strip, as in:
\begin{gather}\label{eq:rex}
\begin{tikzpicture}[anchorbase,scale=1]
\draw[mor] (0,-0.5) rectangle (1,0) node[black,pos=0.5]{$\un{w}$};
\draw[rex] (0,0) rectangle (1,0.1);
\draw[mor] (0,0.1) rectangle (1,0.6) node[black,pos=0.5]{$\un{w}^{\prime}$};
\end{tikzpicture}
\quad := \quad
\begin{tikzpicture}[anchorbase,scale=1]
\draw[mor] (0,0) rectangle (1,0.5) node[black,pos=0.5]{$\un{w}$};
\draw[mor] (0,0.5) rectangle (1,1) node[color=darktomato,pos=0.5]{$N$};
\draw[mor] (0,1) rectangle (1,1.5) node[black,pos=0.5]{$\un{w}^{\prime}$};
\end{tikzpicture}
.
\end{gather}
We can compose any existing purple strip with a rex move to get another valid purple strip, which gives us the flexibility of changing the reduced expression, for example:
\begin{equation} \begin{tikzpicture}[anchorbase,scale=1]
\draw[mor] (0,-0.5) rectangle (1,0) node[black,pos=0.5]{$X$};
\draw[rex] (0,0) rectangle (1,0.1);
\draw[mor] (0,0.1) rectangle (1,0.6) node[black,pos=0.5]{$\un{w}$};
\draw[rex] (0,0.6) rectangle (1,0.7);
\draw[mor] (0,0.7) rectangle (1,1.2) node[black,pos=0.5]{$\un{w}^{\prime}$};
\end{tikzpicture}
=
\begin{tikzpicture}[anchorbase,scale=1]
\draw[mor] (0,-0.5) rectangle (1,0) node[black,pos=0.5]{$X$};
\draw[rex] (0,0) rectangle (1,0.1);
\draw[mor] (0,0.1) rectangle (1,0.6) node[black,pos=0.5]{$\un{w}^{\prime}$};
\end{tikzpicture} \; .
\end{equation} 
Equality holds again because the morphism between the extremal idempotents is irrelevant so long as it preserves the 1-tensor.
Note that \autoref{Lem:ifclaspremovee} holds mutatis mutandis when $e_X$ is a clasp.

\begin{Example}\label{nastyexamplenew} Let us redraw the idempotent $e_0$ from \autoref{extendingnastypart1} using purple strip notation.
\begin{equation}
e_0
=
\begin{tikzpicture}[anchorbase,scale=1]
\diagrammaticmorphism{0}{0.5}{3}{3,1,2,1,3}
\draw[mor] (0,0) rectangle (1,0.5)node[black,pos=0.5]{$e$};
\diagrammaticmorphism{0}{-0.3}{3}{3,1,2,1,3}
\draw[soergeltwo] (1.1,-0.3) to (1.1,0.8);
\end{tikzpicture}
+
\begin{tikzpicture}[anchorbase,scale=1]
\draw[mor] (0,1.2) rectangle (1,1.7)node[black,pos=0.5]{$e$};
\draw[rex] (0,1.1) rectangle (1,1.2);
\diagrammaticmorphism{0}{0.8}{3}{1,2,3,2}
\draw[mor] (0,0.3) rectangle (0.8,0.8)node[black,pos=0.5]{$1232$};
\diagrammaticmorphism{0}{0}{3}{1,2,3,2}
\draw[rex] (0,-0.1) rectangle (1,0);
\draw[mor] (0,-0.6) rectangle (1,-0.1)node[black,pos=0.5]{$e$};
\draw[soergelone,markedone=1] (0.9,0) to (0.9,0.1);
\draw[soergelone,markedone=1] (0.9,1.1) to (0.9,1);
\draw[soergeltwo] (0.7,0.9) to (1.2,0.9) to (1.2,1.7);
\draw[soergeltwo] (0.7,0.2) to (1.2,0.2) to (1.2,-0.6);
\end{tikzpicture}	
+
\begin{tikzpicture}[anchorbase,scale=1]
\draw[mor] (0,1.2) rectangle (1,1.7)node[black,pos=0.5]{$e$};
\draw[rex] (0,1.1) rectangle (1,1.2);
\diagrammaticmorphism{0}{0.8}{3}{3,2,1,2}
\draw[mor] (0,0.3) rectangle (0.8,0.8)node[black,pos=0.5]{$3212$};
\diagrammaticmorphism{0}{0}{3}{3,2,1,2}
\draw[rex] (0,-0.1) rectangle (1,0);
\draw[mor] (0,-0.6) rectangle (1,-0.1)node[black,pos=0.5]{$e$};
\draw[soergelthree,markedthree=1] (0.9,0) to (0.9,0.1);
\draw[soergelthree,markedthree=1] (0.9,1.1) to (0.9,1);
\draw[soergeltwo] (0.7,0.9) to (1.2,0.9) to (1.2,1.7);
\draw[soergeltwo] (0.7,0.2) to (1.2,0.2) to (1.2,-0.6);
\end{tikzpicture}.
\end{equation}
In the next section we will introduce additional diagrammatic language, and we will be able to redraw this idempotent as
\begin{equation}
e_0
=
\begin{tikzpicture}[anchorbase,scale=1]
\diagrammaticmorphism{0}{0.5}{3}{3,1,2,1,3}
\draw[mor] (0,0) rectangle (1,0.5)node[black,pos=0.5]{$e$};
\diagrammaticmorphism{0}{-0.3}{3}{3,1,2,1,3}
\draw[soergeltwo] (1.1,-0.3) to (1.1,0.8);
\end{tikzpicture}
+
\begin{tikzpicture}[anchorbase,scale=1]
\draw[mor] (0,0.6) rectangle (1,1.1)node[black,pos=0.5]{$e$};
\draw[rex] (0,0.5) rectangle (1,0.6);
\draw[mor] (0,0) rectangle (0.8,0.5)node[black,pos=0.5]{$1232$};
\draw[rex] (0,-0.1) rectangle (1,0);
\draw[mor] (0,-0.6) rectangle (1,-0.1)node[black,pos=0.5]{$e$};
\draw[soergelone,markedone=1] (0.9,0) to (0.9,0.15);
\draw[soergelone,markedone=1] (0.9,0.5) to (0.9,0.35);
\draw[soergeltwo] (0.8,0.25) to (1.2,0.25) to (1.2,1.1);
\draw[soergeltwo] (1.2,0.25) to (1.2,-0.6);
\end{tikzpicture}	
+
\begin{tikzpicture}[anchorbase,scale=1]
\draw[mor] (0,0.6) rectangle (1,1.1)node[black,pos=0.5]{$e$};
\draw[rex] (0,0.5) rectangle (1,0.6);
\draw[mor] (0,0) rectangle (0.8,0.5)node[black,pos=0.5]{$3212$};
\draw[rex] (0,-0.1) rectangle (1,0);
\draw[mor] (0,-0.6) rectangle (1,-0.1)node[black,pos=0.5]{$e$};
\draw[soergelthree,markedthree=1] (0.9,0) to (0.9,0.15);
\draw[soergelthree,markedthree=1] (0.9,0.5) to (0.9,0.35);
\draw[soergeltwo] (0.8,0.25) to (1.2,0.25) to (1.2,1.1);
\draw[soergeltwo] (1.2,0.25) to (1.2,-0.6);
\end{tikzpicture} \; .
\end{equation}
More details are below.
\end{Example}

\subsection{Stacking top summands} \label{subsection:topstack}

\begin{Lemma} \label{Lem:topstack}
Let $y, w \in W$ with $y<w$. Let $Y$ be a $y$-object, and $X$ be a $w$-object of the form $M \ot Y \ot N$ for some objects $M, N$. Let $e_w$ and $e_y$ denote top idempotents in $X$ and $Y$ respectively. If $e_w$ is a clasp then
\begin{equation}
\begin{tikzpicture}[anchorbase,scale=1]
\draw[mor] (0,0) rectangle (1,0.5) node[black,pos=0.5]{$y$};
\draw[mor] (-0.25,0.5) rectangle (1.25,1) node[black,pos=0.5]{$w$};
\end{tikzpicture}
=
\begin{tikzpicture}[anchorbase,scale=1]
\draw[mor] (-0.25,0.5) rectangle (1.25,1) node[black,pos=0.5]{$w$};
\end{tikzpicture}
=
\begin{tikzpicture}[anchorbase,scale=1]
\draw[mor] (0,0) rectangle (1,0.5) node[black,pos=0.5]{$y$};
\draw[mor] (-0.25,-0.5) rectangle (1.25,0) node[black,pos=0.5]{$w$};
\end{tikzpicture},
\end{equation} where in the empty regions we have $\id_M$ and $\id_N$ respectively.
\end{Lemma}

\begin{proof} The proof is similar to that of \autoref{Lem:ifclaspremovee}, and starts by writing $e_y = \id_{\un{y}} - \sum q_j$, where $q_j$ are idempotents projecting to lower summands. To argue that $\id_M \ot q_j \ot \id_N$ is still in $\catstuff{I}^0_{< w}$, we use that $\catstuff{I}^0_{< w} = \catstuff{I}^0_{< \len(w)}$. \end{proof}

Whenever $y < w$, we can find a reduced expression $\un{y}$ for $y$ as a contiguous subword of some reduced expression $\un{w}$ for $w$, so we might set $Y = B_{\un{y}}$ and $X = B_{\un{w}}$. This gives an example of the above phenomenon, when $X$ happens to admit a clasp. Since Bott--Samelson bimodules often do not admit clasps, it helps to have a slight generalization of the above result.

\begin{Lemma} \label{Lem:topstack2}
Let $y, w \in W$ with $y<w$. Let $X$ be a $w$-object, and choose a reduced expression $\un{w}$ for $w$ which has a reduced expression $\un{y}$ for $y$ as a contiguous subword. Let the purple strip induce the canonical morphism $X \to B_{\un{w}}$, or vice versa. If $e_X$ is a clasp, then
\begin{equation}
\begin{tikzpicture}[anchorbase,scale=1]
\draw[mor] (0,0) rectangle (1,0.5) node[black,pos=0.5]{$y$};
\draw[rex] (-0.25,-0.1) rectangle (1.25,0);
\draw[mor] (-0.25,-0.6) rectangle (1.25,-0.1) node[black,pos=0.5]{$X$};
\end{tikzpicture} 
= 
\begin{tikzpicture}[anchorbase,scale=1]
\draw[mor] (-0.25,0.5) rectangle (1.25,1) node[black,pos=0.5]{$X$};
\draw[rex] (-0.25,1) rectangle (1.25,1.1);
\end{tikzpicture}
\; , \qquad 
\begin{tikzpicture}[anchorbase,scale=1]
\draw[mor] (0,0) rectangle (1,0.5) node[black,pos=0.5]{$y$};
\draw[rex] (-0.25,0.5) rectangle (1.25,0.6);
\draw[mor] (-0.25,0.6) rectangle (1.25,1.1) node[black,pos=0.5]{$X$};
\end{tikzpicture}
=
\begin{tikzpicture}[anchorbase,scale=1]
\draw[mor] (-0.25,0.5) rectangle (1.25,1) node[black,pos=0.5]{$X$};
\draw[rex] (-0.25,0.4) rectangle (1.25,0.5);
\end{tikzpicture} \; .
\end{equation}

\end{Lemma}

\begin{proof}
Easy and omitted.
\end{proof}

\section{Computing idempotents and traces}\label{section:idempotentsinhecke}

Now we focus on the implications of Kazhdan--Lusztig combinatorics.

\subsection{Generalized trivalent vertices}\label{subsection:idempotentstrivalent}

Suppose that $w \in W$ and $i\in S$ are such that $wi< w$. Then $B_w B_i \cong B_w(1) \oplus B_w(-1)$ by \eqref{downdecomp}. As a consequence, Soergel's hom formula implies that $\Hom^{-1}(B_w, B_w B_i)$ is one-dimensional. As we will see below, there is an element in this one-dimensional space which sends $\onetensor$ to $\onetensor$, that we call the \emph{$(w,i)$-trivalent vertex}. Consequently, an element of $\Hom^{-1}(B_w, B_w B_i)$ is determined by its action on $\onetensor$. Let us construct the $(w,i)$-trivalent vertex explicitly. We do this using the $(i,i)$-trivalent vertex $T_i \colon B_i \to B_i B_i$, which is one of the generators in the Soergel calculus.

\begin{Definition} \label{Def:witri} Let $x, w \in W$ and $i \in S$, such that $w = xi > x$. Let $B_w$ be the image of a top idempotent $e_w$ inside a based $w$-object $X$. Let $B_x$ be the image of a top idempotent $e_x$ inside a based $x$-object $Y$. Then $Y B_i$ is another based $w$-object.
Let $\phi \colon X \to Y B_i$ denote the canonical morphism, and let $\psi$ denote the canonical morphism $Y B_i \to X$, both of which we denote with a purple strip below. Then we define the  \emph{$(w,i)$-trivalent vertex} inside 
\[ \Hom^{-1}(B_w, B_w B_i) = (e_w \ot \id_i) \circ \Hom^{-1}(X,X B_i) \circ e_w \]
as the composition
\begin{equation} (\psi \ot \id_{B_i}) \circ (e_x \otimes T_i) \circ \phi, \end{equation}
or in pictures
\begin{equation} \label{eq:witridef}
\begin{tikzpicture}[anchorbase,scale=1]
\draw[mor] (0,0) rectangle (1,0.5) node[black,pos=0.5]{$e_w$};
\draw[usual] (1,0.25) to (1.3,0.25) to (1.3,0.5);
\node at (1.3,0.7) {$i$};
\end{tikzpicture}
\coloneqq 
\begin{tikzpicture}[anchorbase,scale=1]
\draw[mor] (0,-0.5) rectangle (1.2,0) node[black,pos=0.5]{$e_w$};
\draw[rex] (0,0) rectangle (1.2,0.1);
\draw[mor] (0,0.1) rectangle (.8,0.6) node[black,pos=0.5]{$e_x$};
\draw[usual] (1.1,0.1) to (1.1,0.6);
\draw[usual] (1.1,0.35) to (1.4,0.35) to (1.4,1.2);
\draw[rex] (0,0.6) rectangle (1.2,0.7);
\node at (1.5,1.2) {$i$};
\draw[mor] (0,0.7) rectangle (1.2,1.2) node[black,pos=0.5]{$e_w$};
\end{tikzpicture}  \; .
\end{equation}
We use this shorthand notation throughout.
\end{Definition}


As a special case of a based $w$-object we might have $X = B_w$, where $B_w$ is viewed as an abstract bimodule with fixed element $\onetensor$. In this case $e_w$ is the identity map. Traditional diagrammatics does not have a way of drawing the identity map of $B_w$ except as the image of an idempotent inside a Bott--Samelson bimodule, but it helps to draw more abstract bimodules schematically, so we represent $\id_{B_w}$ as a rectangle labeled by $w$ instead of $e_w$. Similarly we might have $Y = B_x$ and $e_x = \id_{B_x}$. In this abstract setting, the trivalent vertex would be drawn as
\begin{equation} \label{eq:witridefabstract}
\begin{tikzpicture}[anchorbase,scale=1]
\draw[mor] (0,0) rectangle (1,0.5) node[black,pos=0.5]{$w$};
\draw[usual] (1,0.25) to (1.3,0.25) to (1.3,0.5);
\node at (1.3,0.7) {$i$};
\end{tikzpicture}
\coloneqq 
\begin{tikzpicture}[anchorbase,scale=1]
\draw[mor] (0,-0.5) rectangle (1.2,0) node[black,pos=0.5]{$w$};
\draw[rex] (0,0) rectangle (1.2,0.1);
\draw[mor] (0,0.1) rectangle (.8,0.6) node[black,pos=0.5]{$x$};
\draw[usual] (1.1,0.1) to (1.1,0.6);
\draw[usual] (1.1,0.35) to (1.4,0.35) to (1.4,1.2);
\draw[rex] (0,0.6) rectangle (1.2,0.7);
\node at (1.5,1.2) {$i$};
\draw[mor] (0,0.7) rectangle (1.2,1.2) node[black,pos=0.5]{$w$};
\end{tikzpicture}  \; .
\end{equation}

We give examples shortly, after proving that the idempotent $e_x$ above is redundant.

\begin{Lemma} The $(w,i)$-trivalent vertex is the unique element of $\Hom^{-1}(B_w, B_w B_i)$ which preserves $\onetensor$. It is independent of the choice of $x$-object $Y$. \end{Lemma}

\begin{proof} See the discussion before the definition. It remains to note that canonical morphisms and ordinary $(i,i)$-trivalent vertices preserve the one-tensor, and thus so does the $(w,i)$-trivalent vertex. \end{proof}

\begin{Lemma} \label{Lem:redundantidempotent1} Continue the setup of \autoref{Def:witri}. Then
\begin{equation}\label{eq:redundantidempotent}
\begin{tikzpicture}[anchorbase,scale=1]
\draw[mor] (0,0) rectangle (1,0.5) node[black,pos=0.5]{$e_w$};
\draw[usual] (1,0.25) to (1.3,0.25) to (1.3,0.5);
\node at (1.3,0.7) {$i$};
\end{tikzpicture}
\; = \;
\begin{tikzpicture}[anchorbase,scale=1]
\draw[mor] (0,-0.5) rectangle (1.2,0) node[black,pos=0.5]{$e_w$};
\draw[rex] (0,0) rectangle (1.2,0.1);
\draw[mor] (0,0.1) rectangle (.8,0.6) node[black,pos=0.5]{$\id_Y$};
\draw[usual] (1.1,0.1) to (1.1,0.6);
\draw[usual] (1.1,0.35) to (1.4,0.35) to (1.4,1.2);
\draw[rex] (0,0.6) rectangle (1.2,0.7);
\node at (1.5,1.2) {$i$};
\draw[mor] (0,0.7) rectangle (1.2,1.2) node[black,pos=0.5]{$e_w$};
\end{tikzpicture}  \; ,
\end{equation}
where we have replaced the idempotent $e_x \in \End(Y)$ with the identity map.
\end{Lemma}

\begin{proof} This follows from the previous lemma, as the right-hand side preserves $\onetensor$. \end{proof}

We now give some examples of $(w,i)$-trivalent vertices, when $X$ is a Bott--Samelson bimodule associated to a reduced expression $\un{w}$ of $w$. In practice, the canonical morphism drawn in purple in \eqref{eq:witridef} could be a rex move between $X$ and a reduced expression $\un{w}'$ ending in $i$, i.e. a reduced expression $\un{w}' = \un{y} i$; the object $Y$ would be $B_{\un{y}}$. Notably, the $(w,i)$-trivalent vertex is independent of the choice of rex move. 

\begin{Example} \label{Ex:witrivalents12321} Continuing \autoref{nastyexample}, the element $w = 13231=31213 = 32123$ has several reduced expressions, corresponding to Bott--Samelson bimodules which are $w$-objects. We draw several examples of the $(w,1)$-trivalent vertex, for different possible reduced expressions, using \eqref{eq:redundantidempotent}.
\begin{equation}
\begin{tikzpicture}[anchorbase,scale=1]
\draw[mor] (0,0) rectangle (1,0.5) node[black,pos=0.5]{$13231$};
\draw[soergelone] (1,0.25) to (1.1,0.25) to (1.1,0.5);
\end{tikzpicture}
=
\begin{tikzpicture}[anchorbase,scale=1]
\draw[mor] (0,0) rectangle (1,0.5) node[black,pos=0.5]{$13231$};
\diagrammaticmorphism{0}{0.5}{5}{1,3,2,3,1}
\draw[mor] (0,1) rectangle (1,1.5) node[black,pos=0.5]{$13231$};
\draw[soergelone] (0.9,0.75) to (1.1,0.75) to (1.1,1.5);
\end{tikzpicture}
, \quad
\begin{tikzpicture}[anchorbase,scale=1]
\draw[mor] (0,0) rectangle (1,0.5) node[black,pos=0.5]{$31213$};
\draw[soergelone] (1,0.25) to (1.1,0.25) to (1.1,0.5);
\end{tikzpicture}
=
\begin{tikzpicture}[anchorbase,scale=1]
\draw[mor] (0,0) rectangle (1,0.5) node[black,pos=0.5]{$31213$};
\draw[rex] (0,0.5) rectangle (1,0.6);
\draw[rex] (0,0.9) rectangle (1,1);
\draw[mor] (0,1) rectangle (1,1.5) node[black,pos=0.5]{$31213$};
\diagrammaticmorphism{0}{0.6}{3}{3,1,2,3,1}
\draw[soergelone] (0.9,0.75) to (1.1,0.75) to (1.1,1.5);
\end{tikzpicture}
=
\begin{tikzpicture}[anchorbase,scale=1.5, every node/.style={transform shape}]
\draw[mor] (0,0) rectangle (1,0.5) node[black,pos=0.5]{$31213$};
\draw[mor] (0,0.9) rectangle (1,1.4) node[black,pos=0.5]{$31213$};
\diagrammaticmorphism{0}{0.5}{4}{3,1,2}
\draw[soergelone] (0.7,0.5) to (0.9,0.7) to (0.7,0.9);
\draw[soergelone] (0.9,0.7) to (1.1,0.7) to (1.1,1.4);
\draw[soergelthree] (0.9,0.5) to (0.7,0.7) to (0.9,0.9);
\end{tikzpicture}
,\quad \begin{tikzpicture}[anchorbase,scale=1]
\draw[mor] (0,0) rectangle (1,0.5) node[black,pos=0.5]{$32123$};
\draw[soergelone] (1,0.25) to (1.1,0.25) to (1.1,0.5);
\end{tikzpicture}
=
\begin{tikzpicture}[anchorbase,scale=1]
\draw[mor] (0,0) rectangle (1,0.5) node[black,pos=0.5]{$32123$};
\draw[rex] (0,0.5) rectangle (1,0.6);
\draw[rex] (0,0.9) rectangle (1,1);
\draw[mor] (0,1) rectangle (1,1.5) node[black,pos=0.5]{$32123$};
\diagrammaticmorphism{0}{0.6}{3}{3,1,2,3,1}
\draw[soergelone] (0.9,0.75) to (1.1,0.75) to (1.1,1.5);
\end{tikzpicture}
=
\begin{tikzpicture}[anchorbase,scale=1.5, every node/.style={transform shape}]
\draw[mor] (0,0) rectangle (1,0.5) node[black,pos=0.5]{$32123$};
\draw[mor] (0,1.1) rectangle (1,1.6) node[black,pos=0.5]{$32123$};
\drawtwomvalentvertex{3}{soergeltwo}{soergelone}{0.2}{0.5}
\drawtwomvalentvertex{2}{soergelone}{soergelthree}{0.6}{0.7}
\drawtwomvalentvertex{2}{soergelthree}{soergelone}{0.6}{0.8}
\drawtwomvalentvertex{3}{soergelone}{soergeltwo}{0.2}{0.9}
\draw[soergelone] (0.9,0.8) to (1.1,0.8) to (1.1,1.6);
\draw[soergelthree] (0.1,0.5) to (0.1,1.1);
\diagrammaticmorphism{0.2}{0.7}{2}{1,2}
\draw[soergelthree] (0.9,0.5) to (0.9,0.7);
\draw[soergelthree] (0.9,0.9) to (0.9,1.1);
\end{tikzpicture}		\;.
\end{equation}
We often simply use the left-hand pictures.
\end{Example}

Let us give an alternate construction of the $(w,i)$-trivalent vertex. If $wi < w$, then the deletion condition implies that, for any $\un{w} \in \Rex(w)$, right-multiplication by $i$ agrees with removing some simple reflection $j$ from $\un{w}$. A choice of $\un{w}$ yields $x_1, x_2 \in W$ and $j \in S$ with 
\begin{equation} \label{deletioncomposition} w = x_1 j x_2, \quad \len(w) = \len(x_1) + 1 + \len(x_2), \quad wi = x_1 x_2, \quad x_2 i = j x_2. \end{equation}
In particular, $\len(j x_2) = \len(x_2 i) = \len(x_2) + 1$ or else the composition for $w$ would not be reduced. We call the triple $(x_1, j, x_2)$ a \emph{deletion composition} for the pair $(w,i)$.

\begin{Lemma} \label{Lem:alternatetri} Suppose $wi<w$, and pick a deletion composition $(x_1, j, x_2)$ for $w$. Then we have
\begin{equation}
\includegraphics[valign=c]{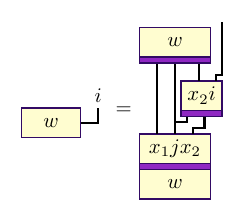}
.
\end{equation}
The purple strips represent rex moves from $jx_2$ to $x_2 i$, or the canonical map from $B_w$ to $B_{\un{w}}$ where $\un{w}$ a concatenation of reduced expressions for $x_1$, $j$, and $x_2$.
\end{Lemma}

\begin{proof} The right-hand side preserves the 1-tensor. \end{proof}

\begin{Example} We now draw the $(w,i)$-trivalent vertices from \autoref{Ex:witrivalents12321} anew.
\begin{equation}
\begin{tikzpicture}[anchorbase,scale=1.5, every node/.style={transform shape}]
\draw[mor] (0,0) rectangle (1,0.5) node[black,pos=0.5]{$31213$};
\draw[mor] (0,0.9) rectangle (1,1.4) node[black,pos=0.5]{$31213$};
\diagrammaticmorphism{0}{0.5}{4}{3,1,2}
\draw[soergelone] (0.7,0.5) to (0.9,0.7) to (0.7,0.9);
\draw[soergelone] (0.9,0.7) to (1.1,0.7) to (1.1,1.4);
\draw[soergelthree] (0.9,0.5) to (0.7,0.7) to (0.9,0.9);
\end{tikzpicture}
=
\begin{tikzpicture}[anchorbase,scale=1.5, every node/.style={transform shape}]
\draw[mor] (0,0) rectangle (1,0.5) node[black,pos=0.5]{$31213$};
\draw[mor] (0,0.7) rectangle (1,1.2) node[black,pos=0.5]{$31213$};
\diagrammaticmorphism{0}{0.5}{2}{3,1,2,1,3}
\draw[soergelone] (0.7,0.6) to (1.1,0.6) to (1.1,1.2);
\end{tikzpicture}
,\quad\includegraphics[valign=c]{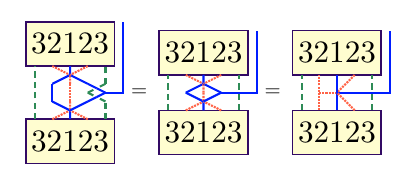}.
\end{equation} 
For each equality above, the left-hand side is the original description from \eqref{eq:redundantidempotent}, the right-hand side is the the alternate description of \autoref{Lem:alternatetri}, and each equality follows from a two-color associativity relation. 
\end{Example}

\begin{Example} When $w = w_0$ is the longest element of $S_n$, the $(w,i)$-trivalent vertex was described explicitly for all $i \in S$ in \cite[Definition 3.34]{El-thick-soergel-typea}. The basic properties stated below are analogs of \cite[Claim 3.37, Proposition 3.41]{El-thick-soergel-typea}. \end{Example}

Here are two basic properties of the $(w,i)$-trivalent vertex, analogous to the usual associativity and unit relations (which are special cases when $w=s_i$).

\begin{Lemma} We have
\begin{equation} \label{triassoc}	\begin{tikzpicture}[anchorbase,scale=1]
\draw[mor] (0,0) rectangle (1,0.5) node[black,pos=0.5]{$w$};
\draw[mor] (0,0.5) rectangle (1,1) node[black,pos=0.5]{$w$};
\draw[usual] (1,0.25) to (1.5,0.25) to (1.5,1);
\draw[usual] (1,0.75) to (1.3,0.75) to (1.3,1);
\end{tikzpicture}
=
\begin{tikzpicture}[anchorbase,scale=1]
\draw[mor] (0,0) rectangle (1,0.5) node[black,pos=0.5]{$w$};
\draw[usual] (1,0.25) to (1.5,0.25) to (1.5,0.5);
\draw[usual] (1.3,0.25) to (1.3,0.5);
\end{tikzpicture}.
\end{equation}
\begin{equation}\label{triunit}	\begin{tikzpicture}[anchorbase,scale=1]
\draw[mor] (0,0) rectangle (1,0.5) node[black,pos=0.5]{$w$};
\draw[mor] (0,0.5) rectangle (1,1) node[black,pos=0.5]{$w$};
\draw[usual] (1.3,0.75) to (1.3,0.25) to (1,0.25);
\fill[black] (1.3,0.75) circle (.055cm);
\end{tikzpicture}
=
\begin{tikzpicture}[anchorbase,scale=1]
\draw[mor] (0,0) rectangle (1,0.5) node[black,pos=0.5]{$w$};
\end{tikzpicture}.
\end{equation}
\end{Lemma}

\begin{proof} Both $\Hom^{-2}(B_w,B_wB_iB_i)$ and $\Hom^0(B_w,B_w)$ are one-dimensional. All morphisms drawn preserve the one-tensor. \end{proof}

Here is a third crucial property of trivalent vertices.

\begin{Lemma} The $(w,i)$-trivalent vertex slides past idempotents. That is, using the canonical morphism to and from $Y B_i$, we have
\begin{equation}
\label{trislide} \begin{tikzpicture}[anchorbase,scale=1]
\draw[mor] (0,0) rectangle (1,0.5) node[black,pos=0.5]{$w$};
\draw[rex]  (0,0.5) rectangle (1,0.6);
\draw[rex]  (0,-0.1) rectangle (1,0);
\draw[usual] (0.9,-0.5) to (0.9,-0.1);
\draw[usual] (0.9,0.6) to (0.9,1);
\draw[usual] (0.9,-0.25) to (1.3,-0.25) to (1.3,1);
\end{tikzpicture}
= 	\begin{tikzpicture}[anchorbase,scale=1]
\draw[mor] (0,0) rectangle (1,0.5) node[black,pos=0.5]{$w$};
\draw[rex]  (0,0.5) rectangle (1,0.6);
\draw[rex]  (0,-0.1) rectangle (1,0);
\draw[usual] (1,0.25) to (1.3,0.25) to (1.3,0.5);
\end{tikzpicture} =
\begin{tikzpicture}[anchorbase,scale=1]
\draw[mor] (0,0) rectangle (1,0.5) node[black,pos=0.5]{$w$};
\draw[rex]  (0,0.5) rectangle (1,0.6);
\draw[rex]  (0,-0.1) rectangle (1,0);
\draw[usual] (0.9,-0.5) to (0.9,-0.1);
\draw[usual] (0.9,0.6) to (0.9,1);
\draw[usual] (0.9,0.75) to (1.3,0.75) to (1.3,1);
\end{tikzpicture}.
\end{equation} We have omitted $\id_Y$ several times in this diagram.
\end{Lemma}

\begin{proof} Again, all three diagrams preserve $\onetensor$. \end{proof}

The three relations \eqref{triassoc} and \eqref{triunit} and \eqref{trislide} we call the \emph{trivalent package}. We chose above to prove the trivalent package by applying the
functor from diagrams to Soergel bimodules, and examining what both sides do to the one-tensor. One can ask for a purely diagrammatic proof of these statements, without recourse to
bimodules and elements therein. The next two remarks discuss alternative proofs of the trivalent package.

\begin{Remark} \label{justwithdiagramsremark} Let $w = xi > x$. In \autoref{Thm:inductiveclasp} we provide the formula \eqref{exidiag} for the clasp idempotent $e_w \in \End(B_x
B_i)$. From this explicit description of $e_w$ it is completely straightforward to prove the trivalent package for the $(w,i)$-trivalent vertex, as a consequence of the one-color
relations. Meanwhile, the proof of \autoref{Thm:inductiveclasp} relies on the trivalent package for $(y,j)$-trivalent vertices, for various $y < w$ and $j \in S$. Thus one can prove both \autoref{Thm:inductiveclasp} and the trivalent package by a simultaneous induction. \end{Remark}

\begin{Remark} \label{ssbimtriremark} Here is another perspective on trivalent vertices, coming from singular Soergel bimodules. Let $wi < w$, and let $[w]$ denote the coset containing $w$ in $W/W_{\{i\}}$. There is an indecomposable singular Soergel $(R,R^i)$-bimodule $B_{[w]}$. As a special case of \cite[Proposition 5.4.3(1)]{Wi-sing-soergel},
the indecomposable $(R,R)$-bimodule $B_w$ is the induction of $B_{[w]}$, i.e. $B_w \cong B_{[w]} \ot_{R^i} R$.

The $(w,i)$-trivalent vertex really does not depend on $w$ or $[w]$ in any measurable way! It really comes from the existence of degree $-1$ natural transformation
\begin{equation} \label{witrinatltrans} (-) \otimes_{R^i} R  \to (-) \otimes_{R^i} R \otimes_R B_i, \end{equation}
(this is a natural transformation between two functors from $(R,R^i)$-bimodules to $(R,R)$-bimodules). In turn this comes from an $(R^i,R)$-bimodule morphism
\begin{equation} \label{whencewitri} R^i \otimes_{R^i} R  \to R^i \otimes_{R^i} R \otimes_R B_i \end{equation}
which comes directly from adjunction. That is, recall that $B_i \cong R \otimes_{R^i} R(1)$. Before inducing there is already a degree $-1$ $(R^i,R^i)$-bimodule morphism
\begin{equation} R^i \to R^i \otimes_{R^i} R \otimes_{R^i} R^i(1) \end{equation}
which is just the ring inclusion of $R^i$ into $R$, and part of the Frobenius extension structure.

Using the diagrammatic calculus for singular Soergel bimodules from \cite[\S 24]{ElMaThWi-soergel}, the morphism \eqref{whencewitri} would be drawn as
\begin{gather}
\begin{tikzpicture}[anchorbase]
\draw[dashed, gray] (-1, 1) -- (1.5, 1);  
\draw[dashed, gray] (-1, -1) -- (1.5, -1); 
\draw[tomato, postaction={decoration={markings, mark=at position 0.5 with {\arrow{<}}}, decorate}, thick] (0.5,-1) .. controls (0.6,-0.3) and (1.2,0.3) .. (1.3,1);		
\draw[tomato, postaction={decoration={markings, mark=at position 0.5 with {\arrow{<}}}, decorate}, thick] (1, 1) arc (0:-180:0.6);		
\fill[tomato, opacity=0.5] (-1, -1) --  (0.5,-1) .. controls (0.6,-0.3) and (1.2,0.3) .. (1.3,1) -- (1,1) arc (0:-180:.6) -- (-1, 1) -- cycle; 
\end{tikzpicture} \; .
\end{gather} From this perspective, the trivalent package follows immediately from isotopy relations in singular diagrammatics. \end{Remark}

\subsection{Maps of degree one}\label{subsection:idempotentsdots}

\begin{Notation} For $y, w \in W$ we write $y \muless w$ if $\mu(w,y) \ne 0$ and $y < w$. Equivalently, the coefficient of $v^1$ in the Kazhdan--Lusztig polynomial $p_{w,y}$ is nonzero. \end{Notation}

When $y < w$ note that $\mu(w,y)$ is the dimension of both $\Hom^1(B_w,B_y)$ and $\Hom^1(B_y, B_w)$, and when $y$ is incomparable with $w$ then $\Hom^1(B_w, B_y) = 0$, see \autoref{lem:deg012initial}. We typically draw a nonzero degree $+1$ map $B_w \to B_y$ resp. $B_y \to B_w$ as a trapezoid labeled with ``$+1$'', as below:
\[ \begin{tikzpicture}[anchorbase,scale=1]
\draw[mor] (0,0) rectangle (1.2,0.5) node[black,pos=0.5]{$w$};
\draw[mor] (0,0.5) to (1.2,0.5) to (1.1,1) to (0.1,1) to (0,0.5) node[black] at (0.6,0.75){${+1}$};
\draw[mor] (0.1,1) rectangle (1.1,1.5) node[black,pos=0.5]{$y$};
\end{tikzpicture} \qquad \text{resp.} \qquad \begin{tikzpicture}[anchorbase,scale=1]
\draw[mor] (0.1,1) rectangle (1.1,1.5) node[black,pos=0.5]{$y$};
\draw[mor] (0.1,1.5) to (1.1,1.5) to (1.2,2) to (0,2) to (0.1,1.5) node[black] at (0.6,1.75) {${+1}$};
\draw[mor] (0,2) rectangle (1.2,2.5) node[black,pos=0.5]{$w$};
\end{tikzpicture} \; .\]
Later we will normalize some degree $+1$ maps, to give a more specific meaning to this diagram.

Note that $\ell(y) = \ell(w) - 1$ and $y < w$ implies that $y \muless w$, and indeed $\mu(w,y) = 1$. There is a dichotomy in the kinds of pairs $(y,w)$ with $y \muless w$: the case where $\ell(y) = \ell(w)-1$ and the case where $\ell(y) < \ell(w)-1$ (though $\ell(y)$ and $\ell(w)-1$ will always have the same parity). In the first case, morphisms of degree $+1$ are fairly straightforward. In the latter case morphisms of degree $+1$ can be fairly complicated.

We focus first on the case when $y<w$ and $\ell(y) = \ell(w) - 1$. In this case, a degree $+1$ morphism is the correct degree to send $\onetensor \in B_w$ to $\onetensor \in B_y$. A degree $+1$ morphism $B_w \to B_y$ which preserves $\onetensor$ will be called the $(w,y)$-dot. The $(w,y)$-dot is unique if it exists, since $\Hom^1(B_w,B_y)$ is one-dimensional. We construct it next, thus proving existence. Recall that for any $j \in S$, the dot map $B_j \to \munit$ is one of the generators of the Soergel category, and it preserves the one-tensor.

\begin{Definition} Let $y,w \in W$ be such that $y<w$ and $\ell(y) = \ell(w)-1$. Let $\un{w}$ be any reduced expression for $w$. There is some simple reflection $j$ in $\un{w}$ such that removing $j$ yields a reduced expression $\un{y}$ for $y$. By tensoring the $j$-colored dot with identity maps, we get a morphism $B_{\un{w}} \to B_{\un{y}}$ of degree $+1$ which preserves the 1-tensor. Using the canonical morphism from $B_w$ to $B_{\un{w}}$, and the canonical morphism from $B_{\un{y}}$ to $B_y$, we obtain a morphism we call the \emph{$(w,y)$-dot} or the \emph{$(w,y)$-enddot}, and it preserves the 1-tensor.
\begin{equation} \label{wydotdef}	\begin{tikzpicture}[anchorbase,scale=1]
\draw[mor] (0,0) rectangle (1.2,0.5) node[black,pos=0.5]{$w$};
\draw[mor] (0,0.5) to (1.2,0.5) to (1.1,1) to (0.1,1) to (0,0.5) node[black] at (0.6,0.75){${+1}$};
\draw[mor] (0.1,1) rectangle (1.1,1.5) node[black,pos=0.5]{$y$};
\end{tikzpicture}	
\coloneqq 
\begin{tikzpicture}[anchorbase,scale=1]
\draw[mor] (-0.2,-0.1) rectangle (1.2,0.4) node[black,pos=0.5]{$w$};
\draw[rex] (-0.2,0.4) rectangle (1.2,0.5);
\draw[usual] (-0.1,0.5) to (0,1.4);
\draw[usual] (0.3,0.5) to (0.4,1.4);
\draw[usual] (0.7,0.5) to (0.6,1.4);
\draw[usual] (1.1,0.5) to (1,1.4);
\draw[usual] (0.5,0.5) to (0.5,0.7);
\fill[black] (0.5,0.7) circle (.055cm);
\node at (0.5,1) {$j$};
\draw[rex] (-0.1,1.4) rectangle (1.1,1.5);
\draw[mor] (-0.1,1.5) rectangle (1.1,2) node[black,pos=0.5]{$y$};
\end{tikzpicture}.
\end{equation}

We define the \emph{(upside-down) $(w,y)$-dot} or the \emph{$(w,y)$-startdot} to be the dual morphism, namely
\begin{equation} \label{eq:upside-wy-dot} 
\begin{tikzpicture}[anchorbase,scale=1]
\begin{scope}[yscale=-1]
\draw[mor] (0,0) rectangle (1.2,0.5) node[black,pos=0.5]{$w$};
\draw[mor] (0,0.5) to (1.2,0.5) to (1.1,1) to (0.1,1) to (0,0.5) node[black] at (0.6,0.75){${+1}$};
\draw[mor] (0.1,1) rectangle (1.1,1.5) node[black,pos=0.5]{$y$};
\end{scope}
\end{tikzpicture}	
\coloneqq 
\begin{tikzpicture}[anchorbase,scale=1]
\begin{scope}[yscale=-1]
\draw[mor] (-0.2,-0.1) rectangle (1.2,0.4) node[black,pos=0.5]{$w$};
\draw[rex] (-0.2,0.4) rectangle (1.2,0.5);
\draw[usual] (-0.1,0.5) to (0,1.4);
\draw[usual] (0.3,0.5) to (0.4,1.4);
\draw[usual] (0.7,0.5) to (0.6,1.4);
\draw[usual] (1.1,0.5) to (1,1.4);
\draw[usual] (0.5,0.5) to (0.5,0.7);
\fill[black] (0.5,0.7) circle (.055cm);
\node at (0.5,1) {$j$};
\draw[rex] (-0.1,1.4) rectangle (1.1,1.5);
\draw[mor] (-0.1,1.5) rectangle (1.1,2) node[black,pos=0.5]{$y$};
\end{scope}
\end{tikzpicture} \; .
\end{equation}

In this context, trapezoids labeled $+1$ will always be normalized as in \eqref{wydotdef} and \eqref{eq:upside-wy-dot}. The $(w,y)$-enddot is normalized so that $\onetensor$ is
preserved. The $(w,y)$-startdot does not preserve $\onetensor$ for degree reasons, but it is normalized by being dual to the $(w,y)$-enddot. \end{Definition}


\begin{Lemma} Let $y,w \in W$ be such that $y<w$ and $\ell(y) = \ell(w)-1$. The $(w,y)$-dot, as a morphism $B_w \to B_y$, is independent of the choice of $\un{w}$, and is the unique morphism of degree $+1$ which preserves $\onetensor$. \end{Lemma}

\begin{proof} The independence of $\un{w}$ follows from unicity, which was proven before the definition. \end{proof}

\begin{Example} For each $j$ in the right descent set of $w$, let $y = wj < w$. Then there is a reduced expression $\un{w}$ for $w$ which ends in $j$. Consequently, we can set
\begin{equation} \label{wydotdefspecialcase} \begin{tikzpicture}[anchorbase,scale=1]
\draw[mor] (0,0) rectangle (1.2,0.5) node[black,pos=0.5]{$w$};
\draw[mor] (0,0.5) to (1.2,0.5) to (1.1,1) to (0.1,1) to (0,0.5) node[black] at (0.6,0.75){${+1}$};
\draw[mor] (0.1,1) rectangle (1.1,1.5) node[black,pos=0.5]{$y$};
\end{tikzpicture}	
\coloneqq 
\begin{tikzpicture}[anchorbase,scale=1]
\draw[mor] (0,0) rectangle (1.2,0.5) node[black,pos=0.5]{$w$};
\draw[rex] (0,0.5) rectangle (1.2,0.6);
\diagrammaticmorphism{0}{0.6}{3}{0,0,0,0,0}
\draw[rex] (0,0.9) rectangle (1,1);
\draw[mor] (0,1) rectangle (1,1.5) node[black,pos=0.5]{$y$};
\draw[usual] (1.1,0.6) to (1.1,0.75);
\fill[black] (1.1,0.75) circle (.055cm);
\node at (1.3,0.9) {$j$};
\end{tikzpicture}. \end{equation}
It is often convenient in our calculations to choose this description of the $(w,y)$-dot when possible. \end{Example}

Now we discuss the interaction between $(w,y)$-dots and $(w,i)$-trivalent vertices. There are two cases, when $yi > y$ and when $yi<y$. The first case is very restrictive; in the following lemma we do not assume that $\len(y) = \len(w) - 1$.

\begin{Lemma} \label{yiisw} If $y \muless w$ and $wi < w$ and $yi>y$ then $yi = w$. \end{Lemma}

\begin{proof} This is a standard result in Kazhdan--Lusztig combinatorics. Suppose that $y<w$ and $wi<w$ and $y<yi$. By \cite[Theorem 6.6(c)]{Lu-hecke-book} we have $yi \le w$ and $p_{y,w} = v p_{yi,w}$. Thus the coefficient $\mu(w,y)$ of $v$ in $p_{y,w}$ equals the constant coefficient of $p_{yi,w}$, which is zero unless $yi = w$. \end{proof}

\begin{Remark} Another way to think about the previous lemma is the following. Suppose that $y \muless w$ and $wi< w$ and $yi>y$ and $yi \ne w$. Then the morphism
\[ \begin{tikzpicture}[anchorbase,scale=1]
\draw[mor] (0,0) rectangle (1.2,0.5) node[black,pos=0.5]{$w$};
\draw[mor] (0,0.5) to (1.2,0.5) to (1.1,1) to (0.1,1) to (0,0.5) node[black] at (0.6,0.75){${+1}$};
\draw[mor] (0.1,1) rectangle (1.1,1.5) node[black,pos=0.5]{$y$};
\draw[usual] (1.2,0.25) to (1.4,0.25) to (1.4,2);
\end{tikzpicture} \]
is a degree zero morphism $B_w \to B_y B_i$. But $B_y B_i$ is perverse, and all summands are $\le yi$ in the Bruhat order. Thus the diagram is zero for degree reasons. This is a contradiction (adding a dot to the $i$-colored output, and using \eqref{triunit}, this would imply that the original degree $+1$ map is also zero). \end{Remark}

\begin{Example} Suppose that $w = yi > y$. Note that $y \muless w$. Then we have the following computation.
\begin{equation} \begin{tikzpicture}[anchorbase,scale=1]
\draw[mor] (0,0) rectangle (1.2,0.5) node[black,pos=0.5]{$w$};
\draw[mor] (0,0.5) to (1.2,0.5) to (1.1,1) to (0.1,1) to (0,0.5) node[black] at (0.6,0.75){${+1}$};
\draw[mor] (0.1,1) rectangle (1.1,1.5) node[black,pos=0.5]{$y$};
\draw[usual] (1.2,0.25) to (1.4,0.25) to (1.4,2);
\end{tikzpicture} \; = \; \begin{tikzpicture}[anchorbase,scale=1]
\draw[mor] (0,0) rectangle (1.2,0.5) node[black,pos=0.5]{$w$};
\draw[rex] (0,0.5) rectangle (1.2,0.6);
\diagrammaticmorphism{0}{0.6}{3}{0,0,0,0,0}
\draw[rex] (0,0.9) rectangle (1,1);
\draw[mor] (0,1) rectangle (1,1.5) node[black,pos=0.5]{$y$};
\draw[usual] (1.1,0.6) to (1.1,0.75);
\fill[black] (1.1,0.75) circle (.055cm);
\draw[usual] (1.2,0.25) to (1.6,0.25) to (1.6,2);
\node at (1.25,0.9) {$i$};
\end{tikzpicture} \; = \; \begin{tikzpicture}[anchorbase,scale=1]
\draw[mor] (0,0) rectangle (1.2,0.5) node[black,pos=0.5]{$w$};
\draw[rex] (0,0.5) rectangle (1.2,0.6);
\diagrammaticmorphism{0}{0.6}{3}{0,0,0,0,0}
\draw[rex] (0,0.9) rectangle (1,1);
\draw[mor] (0,1) rectangle (1,1.5) node[black,pos=0.5]{$y$};
\draw[usual] (1.1,0.6) to (1.1,1);
\fill[black] (1.1,1) circle (.055cm);
\draw[usual] (1.1,0.8) to (1.6,0.8) to (1.6,2);
\node at (1.3,1.2) {$i$};
\end{tikzpicture} \; = \; 
\begin{tikzpicture}[anchorbase,scale=1]
\draw[mor] (0,0) rectangle (1.2,0.5) node[black,pos=0.5]{$w$};
\draw[rex] (0,0.5) rectangle (1.2,0.6);
\diagrammaticmorphism{0}{0.6}{3}{0,0,0,0,0}
\draw[rex] (0,0.9) rectangle (1,1);
\draw[mor] (0,1) rectangle (1,1.5) node[black,pos=0.5]{$y$};
\draw[usual] (1.1,0.6) to (1.1,1.5);
\end{tikzpicture} \; .
\end{equation} We used \eqref{trislide} and the unit axiom in this computation.
\end{Example}

Now we consider the other case, when $yi<y$.

\begin{Example} Suppose that $y \muless w$ and $\ell(y) = \ell(w) - 1$, and $wi < w$ and $yi < y$. Choose a reduced expression $\un{w}$ ending in $i$. There is some simple reflection which can be removed from $\un{w}$ to obtain a reduced expression for $y$, and it is not the final $i$ (because $w = yi$ forces $yi>y$). Thus we have the following picture.
\begin{equation} \label{trislidepastdot} \begin{tikzpicture}[anchorbase,scale=1]
\draw[mor] (0,0) rectangle (1.2,0.5) node[black,pos=0.5]{$w$};
\draw[mor] (0,0.5) to (1.2,0.5) to (1.1,1) to (0.1,1) to (0,0.5) node[black] at (0.6,0.75){${+1}$};
\draw[mor] (0.1,1) rectangle (1.1,1.5) node[black,pos=0.5]{$y$};
\draw[usual2] (1.2,0.25) to (1.4,0.25) to (1.4,2);
\node at (1.5,2.1) {$i$};
\end{tikzpicture} \; = \; \begin{tikzpicture}[anchorbase,scale=1]
\draw[mor] (-0.2,-0.1) rectangle (1.2,0.4) node[black,pos=0.5]{$w$};
\draw[rex] (-0.2,0.4) rectangle (1.2,0.5);
\draw[usual] (-0.1,0.5) to (0,1.4);
\draw[usual] (0.3,0.5) to (0.4,1.4);
\draw[usual] (0.7,0.5) to (0.6,1.4);
\draw[usual2] (1.1,0.5) to (1,1.4);
\draw[usual] (0.5,0.5) to (0.5,0.7);
\fill[black] (0.5,0.7) circle (.055cm);
\node at (0.5,1) {$j$};
\draw[rex] (-0.1,1.4) rectangle (1.1,1.5);
\draw[mor] (-0.1,1.5) rectangle (1.1,2) node[black,pos=0.5]{$y$};
\draw[usual2] (1.2,0.25) to (1.4,0.25) to (1.4,2);
\node at (1.5,2.1) {$i$};
\end{tikzpicture} \; = \; \begin{tikzpicture}[anchorbase,scale=1]
\draw[mor] (-0.2,-0.1) rectangle (1.2,0.4) node[black,pos=0.5]{$w$};
\draw[rex] (-0.2,0.4) rectangle (1.2,0.5);
\draw[usual] (-0.1,0.5) to (0,1.4);
\draw[usual] (0.3,0.5) to (0.4,1.4);
\draw[usual] (0.7,0.5) to (0.6,1.4);
\draw[usual2] (1.1,0.5) to (1,1.4);
\draw[usual] (0.5,0.5) to (0.5,0.7);
\fill[black] (0.5,0.7) circle (.055cm);
\node at (0.5,1) {$j$};
\draw[rex] (-0.1,1.4) rectangle (1.1,1.5);
\draw[mor] (-0.1,1.5) rectangle (1.1,2) node[black,pos=0.5]{$y$};
\draw[usual2] (1.05,0.95) to (1.4,0.95) to (1.4,2);
\node at (1.5,2.1) {$i$};
\end{tikzpicture} \; = \; \begin{tikzpicture}[anchorbase,scale=1]
\draw[mor] (-0.2,-0.1) rectangle (1.2,0.4) node[black,pos=0.5]{$w$};
\draw[rex] (-0.2,0.4) rectangle (1.2,0.5);
\draw[usual] (-0.1,0.5) to (0,1.4);
\draw[usual] (0.3,0.5) to (0.4,1.4);
\draw[usual] (0.7,0.5) to (0.6,1.4);
\draw[usual2] (1.1,0.5) to (1,1.4);
\draw[usual] (0.5,0.5) to (0.5,0.7);
\fill[black] (0.5,0.7) circle (.055cm);
\node at (0.5,1) {$j$};
\draw[rex] (-0.1,1.4) rectangle (1.1,1.5);
\draw[mor] (-0.1,1.5) rectangle (1.1,2) node[black,pos=0.5]{$y$};
\draw[usual2] (1.1,1.75) to (1.4,1.75) to (1.4,2);
\node at (1.5,2.1) {$i$};
\end{tikzpicture} \; = \; \begin{tikzpicture}[anchorbase,scale=1]
\draw[mor] (0,0) rectangle (1.2,0.5) node[black,pos=0.5]{$w$};
\draw[mor] (0,0.5) to (1.2,0.5) to (1.1,1) to (0.1,1) to (0,0.5) node[black] at (0.6,0.75){${+1}$};
\draw[mor] (0.1,1) rectangle (1.1,1.5) node[black,pos=0.5]{$y$};
\draw[usual2] (1.1,1.25) to (1.3,1.25) to (1.3,2);
\node at (1.4,2.1) {$i$};
\end{tikzpicture} \; .
\end{equation} 
So the strand slides.
\end{Example}

The conclusion of the previous example will also hold when $\ell(y) < \ell(w) - 1$, see below.

Now we discuss the case when $y \muless w$ but $\ell(y) < \ell(w) - 1$. In this case, any morphism $B_w \to B_y$ of degree $+1$ will kill the $1$-tensor for degree reasons, and
there is no particularly good normalization (that we are aware of) for any particular degree $+1$ map. When $\mu(w,y) = 1$, we henceforth fix some nonzero degree $+1$ map (which therefore forms a basis for $\Hom^1(B_w,B_y)$) and call it the \emph{$(w,y)$-dot}. We draw it as a trapezoid labeled $+1$ as before, and draw its dual upside-down. When $\mu(w,y) > 1$ and thus $\ell(y) < \ell(w) - 1$, one can pick a basis for $\Hom^1(B_w,B_y)$ and call them the $(w,y)$-dot family. What we do with these morphisms below will then require more complicated linear algebra.

\begin{Example} \label{otherdeg1ex} Here is a morphism of degree $1$ from $B_{2132}$ to $B_2$, where $\{1,2,3\}$ are the usual simple reflections in $S_4$.
\begin{equation}\begin{tikzpicture}[anchorbase,scale=1]
\draw[mor] (0,0) rectangle (1.2,0.5) node[black,pos=0.5]{$2132$};
\draw[mor] (0,0.5) to (1.2,0.5) to (1.1,1) to (0.1,1) to (0,0.5) node[black] at (0.6,0.75){${+1}$};
\draw[mor] (0.1,1) rectangle (1.1,1.5) node[black,pos=0.5]{$2$};
\end{tikzpicture}	 \; = \; 
\begin{tikzpicture}[anchorbase,scale=1]
\draw[mor] (0,0) rectangle (0.8,0.5) node[black,pos=0.5]{$2132$};
\draw[mor] (.3,1.2) rectangle (0.5,1.7) node[black,pos=0.5]{$2$};
\draw[soergeltwo] (0.1,0.5) to (0.1,.9);
\draw[soergelone,markedone=1] (0.3,0.5) to (0.3,0.75);
\draw[soergelthree,markedthree=1] (0.5,0.5) to (0.5,0.75);
\draw[soergeltwo] (0.7,0.5) to (0.7,0.9) to (0.1,0.9);
\draw[soergeltwo] (0.4,0.9) to (0.4,1.2);
\end{tikzpicture}.
\end{equation} This follows from a direct computation.\end{Example}

\begin{Example} \label{yetanotherdeg1ex} Here is a morphism of degree $1$ from $B_{2413524}$ to $B_{24}$ (both indecomposables are Bott--Samelsons), where $\{1,2,3,4,5\}$ are the usual simple reflections in $S_6$.
\begin{equation}
\includegraphics[valign=c]{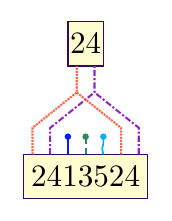}.
\end{equation}This also follows from a direct computation.
\end{Example}

\begin{Theorem} \label{Thm:trislideoverdot} Let $y \muless w$ and let the trapezoid below denote an arbitrary map $B_w \to B_y$ of degree $+1$. If $wi < w$ and $yi < y$ for some $i \in S$, then
\begin{equation} \label{trislideoverplusone} 	\begin{tikzpicture}[anchorbase,scale=1]
\draw[mor] (0,0) rectangle (1.2,0.5) node[black,pos=0.5]{$w$};
\draw[mor] (0,0.5) to (1.2,0.5) to (1.1,1) to (0.1,1) to (0,0.5) node[black] at (0.6,0.75){${+1}$};
\draw[mor] (0.1,1) rectangle (1.1,1.5) node[black,pos=0.5]{$y$};
\draw[usual] (1.2,0.25) to (1.4,0.25) to (1.4,2);
\end{tikzpicture}
\; = \;
\begin{tikzpicture}[anchorbase,scale=1]
\draw[mor] (0,0) rectangle (1.2,0.5) node[black,pos=0.5]{$w$};
\draw[mor] (0,0.5) to (1.2,0.5) to (1.1,1) to (0.1,1) to (0,0.5) node[black] at (0.6,0.75){${+1}$};
\draw[mor] (0.1,1) rectangle (1.1,1.5) node[black,pos=0.5]{$y$};
\draw[usual] (1.1,1.25) to (1.3,1.25) to (1.3,2);
\end{tikzpicture} \; . \end{equation}
\end{Theorem}

The best proof of this result uses singular Soergel bimodules, see \autoref{ssbimtriremark}, and we sketch this proof now. We also prove it as a quick corollary of \autoref{howdotfactorswhenkindescent} below.

\begin{proof} (Sketch)  Both $B_w$ and $B_y$ are induced from $(R,R^i)$-bimodules $B_{[w]}$ and $B_{[y]}$ respectively by \cite[Prop. 5.4.3(1)]{Wi-sing-soergel}, and one can compute that the Kazhdan--Lusztig polynomial $p_{w,y}$ agrees with the singular polynomial $p_{[w],[y]}$ in degree $1$ using \cite[Prop. 2.2.4(1)]{Wi-sing-soergel}. Consequently, the functor $(-) \ot_{R^i} R$ induces an isomorphism $\Hom^1(B_{[w]},B_{[y]}) \to \Hom^1(B_w,B_y)$. There is a degree $+1$ morphism $B_{[w]} \to B_{[y]}$ which we might call the $([w],[y])$-dot, from which the $(w,y)$-dot is obtained by applying the functor $(-) \ot_{R^i} R$. The equality \eqref{trislideoverplusone} above is merely an elaboration of the naturality the natural transformation \eqref{witrinatltrans}, as it applies to the $([w],[y])$-dot. \end{proof}

\subsection{Interlude: more properties of degree one maps}\label{subsection:moreonplusone}

We do not use the results of this section, but this seems like a natural place to record what will surely be the helpful observation \autoref{classifydeg1mapswhenkindescent}.
We continue to study the $(w,y)$-dot when $wi < w$ and $yi<y$ for some $i \in S$. We set $x = wi$ and $z = yi$, and assume that $\ell(x) - \ell(z)$ is odd (or else there could be no $(w,y)$-dot).

\begin{Lemma} Fix $x, z \in W$ and $i \in S$ with $x<xi$ and $z<zi$. For any morphism $f \colon B_{xi} \to B_{zi}$ of degree $1$, there is a morphism $g \colon B_{x} \to B_{zi}$ of degree $+2$ such that
\begin{equation} \label{howdotfactorswhenkindescent} \begin{tikzpicture}[anchorbase,scale=1]
\draw[mor] (0,0) rectangle (1.2,0.5) node[black,pos=0.5]{$xi$};
\draw[mor] (0,0.5) to (1.2,0.5) to (1.1,1) to (0.1,1) to (0,0.5) node[black] at (0.6,0.75){${+1}$};
\draw[mor] (0.1,1) rectangle (1.1,1.5) node[black,pos=0.5]{$zi$};
\end{tikzpicture} \; = \; \begin{tikzpicture}[anchorbase,scale=1]
\draw[mor] (0,-0.5) rectangle (1.3,0) node[black,pos=0.5]{$xi$};
\draw[mor] (0,0) rectangle (1,0.5) node[black,pos=0.5]{$x$};
\draw[mor,rounded corners=0.1cm] (-0.1,0.5) rectangle (1.1,1) node[black, pos = 0.5] {${+2}$};
\draw[mor] (0,1) rectangle (1,1.5) node[black,pos=0.5]{$zi$};
\draw[usual] (1, 1.25) to (1.2,1.25) to (1.2,0);
\end{tikzpicture} \; .
\end{equation}
Here the $+1$ trapezoid represents $f$, while the $+2$ rounded rectangle represents $g$.
\end{Lemma}

\begin{proof} By including $B_{xi}$ into $B_x \ot B_i$ and $B_{zi}$ into $B_z \ot B_i$, the left-hand side is induced from some morphism in $\Hom^1(B_x \ot B_i, B_z \ot B_i)$. By adjunction this hom space is isomorphic to $\Hom^1(B_x B_i B_i, B_z)$. Since $B_i B_i \cong v B_i \oplus v^{-1} B_i$, the hom space in question is isomorphic to $\Hom^0(B_x B_i, B_z) \oplus \Hom^2(B_x B_i, B_z)$, and the isomorphism involves composition with the two different projection maps from $B_i B_i$ to $B_i$: up to scalar, they are the trivalent vertex, and the trivalent vertex with $\alpha_i$ below it. Consequently we deduce that
\begin{equation} \begin{tikzpicture}[anchorbase,scale=1]
\draw[mor] (0,0) rectangle (1.2,0.5) node[black,pos=0.5]{$xi$};
\draw[mor] (0,0.5) to (1.2,0.5) to (1.1,1) to (0.1,1) to (0,0.5) node[black] at (0.6,0.75){${+1}$};
\draw[mor] (0.1,1) rectangle (1.1,1.5) node[black,pos=0.5]{$zi$};
\end{tikzpicture}
=
\begin{tikzpicture}[anchorbase,scale=1]
\draw[mor] (0,-0.5) rectangle (1.3,0) node[black,pos=0.5]{$xi$};
\draw[mor] (0,0) rectangle (1,0.5) node[black,pos=0.5]{$x$};
\draw[mor,rounded corners=0.1cm] (-0.1,0.5) rectangle (1.4,1) node[black, pos = 0.5] {${+1}$};
\draw[mor] (0,1) rectangle (1,1.5) node[black,pos=0.5]{$z$};
\draw[mor] (0,1.5) rectangle (1.3,2) node[black,pos=0.5]{$zi$};
\draw[usual] (1.2,0) to (1.2,0.5);
\draw[usual] (1.2,1) to (1.2,1.5);
\end{tikzpicture} 
=
\begin{tikzpicture}[anchorbase,scale=1]
\draw[mor] (0,-0.5) rectangle (1.5,0) node[black,pos=0.5]{$xi$};
\draw[mor] (0,0) rectangle (1.2,0.5) node[black,pos=0.5]{$x$};
\draw[mor,rounded corners=0.1cm] (-0.1,0.5) rectangle (1.3,1) node[black, pos = 0.5] {${+2}$};
\draw[mor] (0,1) rectangle (1,1.5) node[black,pos=0.5]{$z$};
\draw[mor] (0,1.5) rectangle (1.3,2) node[black,pos=0.5]{$zi$};
\draw[usual] (1.4,0) to (1.4,1) to (1.2,1.25) to (1.2,1.5);
\draw[usual] (1.2,1.25) to (1.2,1);
\end{tikzpicture} 
+
\begin{tikzpicture}[anchorbase,scale=1]
\draw[mor] (0,-0.5) rectangle (1.5,0) node[black,pos=0.5]{$xi$};
\draw[mor] (0,0) rectangle (1.2,0.5) node[black,pos=0.5]{$x$};
\draw[mor,rounded corners=0.1cm] (-0.1,0.5) rectangle (1.3,1) node[black, pos = 0.5] {${+0}$};
\draw[mor] (0,1) rectangle (1,1.5) node[black,pos=0.5]{$z$};
\draw[mor] (0,1.5) rectangle (1.3,2) node[black,pos=0.5]{$zi$};
\draw[usual] (1.4,0) to (1.4,1) to (1.2,1.25) to (1.2,1.5);
\draw[usual] (1.2,1.25) to (1.2,1);
\end{tikzpicture}\alpha_i
=
\begin{tikzpicture}[anchorbase,scale=1]
\draw[mor] (0,-0.5) rectangle (1.3,0) node[black,pos=0.5]{$xi$};
\draw[mor] (0,0) rectangle (1,0.5) node[black,pos=0.5]{$x$};
\draw[mor,rounded corners=0.1cm] (-0.1,0.5) rectangle (1.1,1) node[black, pos = 0.5] {${+2}$};
\draw[mor] (0,1) rectangle (1,1.5) node[black,pos=0.5]{$zi$};
\draw[usual] (1, 1.25) to (1.2,1.25) to (1.2,0);
\end{tikzpicture} 
+
\begin{tikzpicture}[anchorbase,scale=1]
\draw[mor] (0,-0.5) rectangle (1.3,0) node[black,pos=0.5]{$xi$};
\draw[mor] (0,0) rectangle (1,0.5) node[black,pos=0.5]{$x$};
\draw[mor,rounded corners=0.1cm] (-0.1,0.5) rectangle (1.1,1) node[black, pos = 0.5] {${+0}$};
\draw[mor] (0,1) rectangle (1,1.5) node[black,pos=0.5]{$zi$};
\draw[usual] (1, 1.25) to (1.2,1.25) to (1.2,0);
\end{tikzpicture} 
\alpha_i,
\end{equation}
for some morphisms of degree $2$ and degree $0$ respectively.
Note that $zi \ne x$ since $x < xi$, and thus $\Hom^0(B_x,B_{zi}) = 0$. Thus the second diagram on the right-hand side is zero, and the result is proven.
\end{proof}

As promised, this gives a second proof of \autoref{Thm:trislideoverdot}.

\begin{proof}[Proof of \autoref{Thm:trislideoverdot}]\label{proof:second_proof} Given \eqref{howdotfactorswhenkindescent}, the result is now a straightforward consequence of \eqref{triassoc} and \eqref{trislide}. \end{proof}

With the same assumptions as the previous lemma, suppose further that $zi < xi$, which holds if and only if $z<x$. Both $zi < x$ and $zi \not< x$ are possible, but clearly $zi \ne x$. Let us better understand the degree $+2$ morphism appearing in
\eqref{howdotfactorswhenkindescent}. Let $p_{x,i,zi}$ denote the coefficient of $h_{zi}$ in $b_x b_i$. The rule $h_{xi}b_i = h_{x}+v^{-1}h_{xi}$ (see \cite[Equation 3.24]{ElMaThWi-soergel}) implies that \begin{equation} p_{x,i,zi} =
v^{-1} p_{x,zi} + p_{x,z} \in v\N[v] \end{equation} is an odd polynomial with non-negative coefficients. The coefficient of $v$ in $p_{x,i,zi}$ is thus the sum of the coefficient of
$v$ in $p_{x,z}$ and the coefficient of $v^2$ in $p_{x,zi}$.
As $b_{xi}$ is a summand of $b_x b_i$, its coefficient of $h_{zi}$ (namely the KL polynomial $p_{xi,zi}$) is obtained from
$p_{x,i,zi}$ by subtracting some of these positive terms. Note that $p_{x,zi} = 0$ if $zi \not< x$. Note also that $zi \muless xi$ does not imply that $z \muless x$, so it is also possible that the term $p_{x,z}$ does not contribute in degree $+1$.

Let us categorify this observation. Given any morphism in $\Hom^1(B_x,B_z)$ (i.e. a $(x,z)$-dot), we get a morphism in $\Hom^1(B_{xi},B_{zi})$ as follows.
\begin{equation}\begin{tikzpicture}[anchorbase,scale=1]
\draw[mor] (0,0) rectangle (1.2,0.5) node[black,pos=0.5]{$x$};
\draw[mor] (0,0.5) to (1.2,0.5) to (1.1,1) to (0.1,1) to (0,0.5) node[black] at (0.6,0.75){${+1}$};
\draw[mor] (0.1,1) rectangle (1.1,1.5) node[black,pos=0.5]{$z$};
\end{tikzpicture}
\mapsto
\begin{tikzpicture}[anchorbase,scale=1]
\draw[mor] (0,-0.5) rectangle (1.5,0) node[black,pos=0.5]{$xi$};
\draw[mor] (0,0) rectangle (1.2,0.5) node[black,pos=0.5]{$x$};
\draw[mor] (0,0.5) to (1.2,0.5) to (1.1,1) to (0.1,1) to (0,0.5) node[black] at (0.6,0.75){${+1}$};
\draw[mor] (0.1,1) rectangle (1.1,1.5) node[black,pos=0.5]{$z$};
\draw[mor] (0.1,1.5) rectangle (1.4,2) node[black,pos=0.5]{$zi$};
\draw[usual] (1.4,0) to (1.3,1.5);
\end{tikzpicture}.
\end{equation}
Given any morphism in $\Hom^2(B_x, B_{zi})$, we get a morphism in $\Hom^1(B_{xi},B_{zi})$ via \eqref{howdotfactorswhenkindescent}:
\begin{equation}
\begin{tikzpicture}[anchorbase,scale=1]
\draw[mor] (0,0) rectangle (1,0.5) node[black,pos=0.5]{$x$};
\draw[mor,rounded corners=0.1cm] (-0.1,0.5) rectangle (1.1,1) node[black, pos = 0.5] {${+2}$};
\draw[mor] (0,1) rectangle (1,1.5) node[black,pos=0.5]{$zi$};
\end{tikzpicture} 
\mapsto 
\begin{tikzpicture}[anchorbase,scale=1]
\draw[mor] (0,-0.5) rectangle (1.3,0) node[black,pos=0.5]{$xi$};
\draw[mor] (0,0) rectangle (1,0.5) node[black,pos=0.5]{$x$};
\draw[mor,rounded corners=0.1cm] (-0.1,0.5) rectangle (1.1,1) node[black, pos = 0.5] {${+2}$};
\draw[mor] (0,1) rectangle (1,1.5) node[black,pos=0.5]{$zi$};
\draw[usual] (1, 1.25) to (1.2,1.25) to (1.2,0);
\end{tikzpicture} .
\end{equation} However, we are not interested in all morphisms in $\Hom^2(B_x,B_{zi})$, but only those which correspond to the coefficient of $v^2$ in $p_{x,zi}$, namely morphisms which descend to a basis of $\Hom^2(B_x,B_{zi})$ modulo $\catstuff{I}_{< zi}$. We indicate elements of such a basis with trapezoids labeled by $+2$. This distinction was not important for morphisms of degree $+1$ by \autoref{rmk:nolowerindeg1}. In conclusion, we expect a spanning set for $\Hom^1(B_{xi},B_{zi})$ consisting of morphisms of the following form:
\begin{equation} \label{spanningsetfordeg1} \begin{tikzpicture}[anchorbase,scale=1]
\draw[mor] (0,-0.5) rectangle (1.5,0) node[black,pos=0.5]{$xi$};
\draw[mor] (0,0) rectangle (1.2,0.5) node[black,pos=0.5]{$x$};
\draw[mor] (0,0.5) to (1.2,0.5) to (1.1,1) to (0.1,1) to (0,0.5) node[black] at (0.6,0.75){${+1}$};
\draw[mor] (0.1,1) rectangle (1.1,1.5) node[black,pos=0.5]{$z$};
\draw[mor] (0.1,1.5) rectangle (1.4,2) node[black,pos=0.5]{$zi$};
\draw[usual] (1.4,0) to (1.3,1.5);
\end{tikzpicture} \quad \text{and} \quad \begin{tikzpicture}[anchorbase,scale=1]
\draw[mor] (0,-0.5) rectangle (1.5,0) node[black,pos=0.5]{$xi$};
\draw[mor] (0,0) rectangle (1.2,0.5) node[black,pos=0.5]{$x$};
\draw[mor] (0,0.5) to (1.2,0.5) to (1.1,1) to (0.1,1) to (0,0.5) node[black] at (0.6,0.75){${+2}$};
\draw[mor] (0.1,1) rectangle (1.1,1.5) node[black,pos=0.5]{$zi$};
\draw[usual] (1.1, 1.3) to (1.4,1.3) to (1.4,0);
\end{tikzpicture}. \end{equation}
Diagrams of the first form do not arise unless $z \muless x$, and diagrams of the second form do not arise unless $zi < x$. When $B_{xi} = B_x B_i$, one expects such morphisms to form a basis (as the trapezoids run over bases). When $B_{xi}$ is a proper summand of $B_x B_i$, there may be some linear dependencies, arising from those linear combinations of diagrams which become zero when precomposed with the top idempotent $e_{xi}$.

Let us prove this expectation, focusing only on the question of spanning.

\begin{Lemma} \label{classifydeg1mapswhenkindescent} Fix $x, z \in W$ and $i \in S$ with $x<xi$ and $z<zi$ and $zi < xi$. Then $\Hom^1(B_{xi},B_{zi})$ is spanned by morphisms of the form in \eqref{spanningsetfordeg1}. \end{Lemma}

\begin{proof} The hom space in question is spanned by morphisms of the form of the right-hand side of \eqref{howdotfactorswhenkindescent}. Meanwhile, by the expanded Soergel hom formula, $\Hom^2(B_x,B_{zi})$ is spanned by two kinds of morphisms: degree $+2$ trapezoids straight to $zi$, and compositions of maps of degree $+1$. Loosely, we have
\begin{equation}
\begin{tikzpicture}[anchorbase,scale=1]
\draw[mor] (0,0) rectangle (1,0.5) node[black,pos=0.5]{$x$};
\draw[mor,rounded corners=0.1cm] (-0.1,0.5) rectangle (1.1,1) node[black, pos = 0.5] {${+2}$};
\draw[mor] (0,1) rectangle (1,1.5) node[black,pos=0.5]{$zi$};
\end{tikzpicture} 
=
\begin{tikzpicture}[anchorbase,scale=1]
\draw[mor] (0,0) rectangle (1.2,0.5) node[black,pos=0.5]{$x$};
\draw[mor] (0,0.5) to (1.2,0.5) to (1.1,1) to (0.1,1) to (0,0.5) node[black] at (0.6,0.75){${+2}$};
\draw[mor] (0.1,1) rectangle (1.1,1.5) node[black,pos=0.5]{$zi$};
\end{tikzpicture}
+
\sum_{u}
\begin{tikzpicture}[anchorbase,scale=1]
\draw[mor] (0,0) rectangle (1.2,0.5) node[black,pos=0.5]{$x$};
\draw[mor] (0,0.5) to (1.2,0.5) to (1.1,1) to (0.1,1) to (0,0.5) node[black] at (0.6,0.75){${+1}$};
\draw[mor] (0.1,1) rectangle (1.1,1.5) node[black,pos=0.5]{$u$};
\draw[mor] (0.1,1.5) to (1.1,1.5) to (1.2,2) to (0,2) to (0.1,1.5) node[black] at (0.6,1.75) {${+1}$};
\draw[mor] (0,2) rectangle (1.2,2.5) node[black,pos=0.5]{$zi$};
\end{tikzpicture} \; ,
\end{equation} where $u$ in the ``sum'' satisfy $u \muless x$ and $u \muless zk$.
Thus the hom space in question is spanned by morphisms of the second kind from \eqref{spanningsetfordeg1}, and morphisms of the form
\begin{equation}\label{goesthroughuthistime}
\begin{tikzpicture}[anchorbase,scale=1]
\draw[mor] (0,-0.5) rectangle (1.4,0) node[black,pos=0.5]{$xi$};
\draw[mor] (0,0) rectangle (1.2,0.5) node[black,pos=0.5]{$x$};
\draw[mor] (0,0.5) to (1.2,0.5) to (1.1,1) to (0.1,1) to (0,0.5) node[black] at (0.6,0.75){${+1}$};
\draw[mor] (0.1,1) rectangle (1.1,1.5) node[black,pos=0.5]{$u$};
\draw[mor] (0.1,1.5) to (1.1,1.5) to (1.2,2) to (0,2) to (0.1,1.5) node[black] at (0.6,1.75) {${+1}$};
\draw[mor] (0,2) rectangle (1.2,2.5) node[black,pos=0.5]{$zi$};
\draw[usual] (1.3,0) to (1.3,2.25) to (1.2,2.25);
\end{tikzpicture} \; .
\end{equation}
A morphism as in \eqref{goesthroughuthistime} factors through a morphism in $\Hom^0(B_u B_i, B_{zi})$.

If $ui < u$ then by \autoref{Thm:trislideoverdot}, the diagram above is equal to
\begin{equation}		\begin{tikzpicture}[anchorbase,scale=1]
\draw[mor] (0,-0.5) rectangle (1.4,0) node[black,pos=0.5]{$xi$};
\draw[mor] (0,0) rectangle (1.2,0.5) node[black,pos=0.5]{$x$};
\draw[mor] (0,0.5) to (1.2,0.5) to (1.1,1) to (0.1,1) to (0,0.5) node[black] at (0.6,0.75){${+1}$};
\draw[mor] (0.1,1) rectangle (1.1,1.5) node[black,pos=0.5]{$u$};
\draw[mor] (0.1,1.5) to (1.1,1.5) to (1.2,2) to (0,2) to (0.1,1.5) node[black] at (0.6,1.75) {${+1}$};
\draw[mor] (0,2) rectangle (1.2,2.5) node[black,pos=0.5]{$zi$};
\draw[usual] (1.3,0) to (1.3,2.25) to (1.2,2.25);
\end{tikzpicture}
=
\begin{tikzpicture}[anchorbase,scale=1]
\draw[mor] (0,-0.5) rectangle (1.4,0) node[black,pos=0.5]{$xi$};
\draw[mor] (0,0) rectangle (1.2,0.5) node[black,pos=0.5]{$x$};
\draw[mor] (0,0.5) to (1.2,0.5) to (1.1,1) to (0.1,1) to (0,0.5) node[black] at (0.6,0.75){${+1}$};
\draw[mor] (0.1,1) rectangle (1.1,1.5) node[black,pos=0.5]{$u$};
\draw[mor] (0.1,1.5) to (1.1,1.5) to (1.2,2) to (0,2) to (0.1,1.5) node[black] at (0.6,1.75) {${+1}$};
\draw[mor] (0,2) rectangle (1.2,2.5) node[black,pos=0.5]{$zi$};
\draw[usual] (1.3,0) to (1.3,1.25) to (1.1,1.25);
\end{tikzpicture}, \end{equation}
which factors through a degree zero morphism $B_{xi} \to B_u$. Since $u \ne xi$ (because $u < x$), such a morphism is zero.

If $ui > u$ then $B_u B_i$ is perverse, and admits no degree zero maps to $B_{zi}$ unless $B_{zi}$ is one of the summands of $B_u B_i$. But since $u \muless zi$, this is only possible if $u = z$, in which case the $(zi,u)$-dot is given as in \eqref{wydotdefspecialcase}. In this case we compute that
\begin{equation}
\begin{tikzpicture}[anchorbase,scale=1]
\draw[mor] (0,-0.5) rectangle (1.4,0) node[black,pos=0.5]{$xi$};
\draw[mor] (0,0) rectangle (1.2,0.5) node[black,pos=0.5]{$x$};
\draw[mor] (0,0.5) to (1.2,0.5) to (1.1,1) to (0.1,1) to (0,0.5) node[black] at (0.6,0.75){${+1}$};
\draw[mor] (0.1,1) rectangle (1.1,1.5) node[black,pos=0.5]{$u$};
\draw[mor] (0.1,1.5) to (1.1,1.5) to (1.2,2) to (0,2) to (0.1,1.5) node[black] at (0.6,1.75) {${+1}$};
\draw[mor] (0,2) rectangle (1.2,2.5) node[black,pos=0.5]{$zi$};
\draw[usual] (1.3,0) to (1.3,2.25) to (1.2,2.25);
\end{tikzpicture}
=
\begin{tikzpicture}[anchorbase,scale=1]
\draw[mor] (0,-0.5) rectangle (1.6,0) node[black,pos=0.5]{$xi$};
\draw[mor] (0,0) rectangle (1.2,0.5) node[black,pos=0.5]{$x$};
\draw[mor] (0,0.5) to (1.2,0.5) to (1.1,1) to (0.1,1) to (0,0.5) node[black] at (0.6,0.75){${+1}$};
\draw[mor] (0.1,1) rectangle (1.1,1.5) node[black,pos=0.5]{$z$};
\draw[usual, marked=1] (1.3,1.5) to (1.3,1.2);
\draw[mor] (0.1,1.5) rectangle (1.4,2) node[black,pos=0.5]{$zi$};
\draw[usual] (1.5,0) to (1.5,1.75) to (1.4,1.75);
\end{tikzpicture}
=
\begin{tikzpicture}[anchorbase,scale=1]
\draw[mor] (0,-0.5) rectangle (1.5,0) node[black,pos=0.5]{$xi$};
\draw[mor] (0,0) rectangle (1.2,0.5) node[black,pos=0.5]{$x$};
\draw[mor] (0,0.5) to (1.2,0.5) to (1.1,1) to (0.1,1) to (0,0.5) node[black] at (0.6,0.75){${+1}$};
\draw[mor] (0.1,1) rectangle (1.1,1.5) node[black,pos=0.5]{$z$};
\draw[mor] (0.1,1.5) rectangle (1.4,2) node[black,pos=0.5]{$zi$};
\draw[usual] (1.4,0) to (1.3,1.5);
\end{tikzpicture}, \end{equation}
as desired. \end{proof}

\subsection{Inductive construction of clasp idempotents}\label{subsection:induction}

Suppose $w = xi > x$. We remind the reader of \eqref{inductivedecomp} which stated that
\begin{equation} B_x B_i = B_w \oplus \bigoplus_{y \muless x, yi<y} B_y^{\oplus \mu(x,y)}. \end{equation}
In particular, $B_y$ is a summand of $B_x B_i$ if and only if $y \muless x$ and $yi < y$.

To construct the clasp idempotent projecting to $B_w$, we can take the identity of $B_x B_i$ and subtract the idempotents projecting to $B_y$. In order to construct these lower idempotents, we need to construct degree zero morphisms $B_x B_i \to B_y$ to serve as projection maps, and degree zero morphisms backward to serve as inclusions. Since $yi < y$, we have already constructed a $(y,i)$-trivalent vertex $B_y \to B_y B_i$, which under adjunction corresponds to a degree $-1$ morphism $B_y B_i \to B_y$. Composing this with the $(x,y)$-dot, we obtain our desired morphism.

We start by exploring the case when $\mu(x,y) = 1$.

\begin{Definition} \label{Defn:xiyproject} Let $x,w,y \in W$ and $i \in S$ such that $w = xi > x$ and $yi < y < x$ and $\mu(x,y) = 1$. Let $p_{(x,i,y)}$, the \emph{$(x,i,y)$-projection}, be the composition
\begin{equation} \label{eq:pxiy} p_{(x,i,y)} := \begin{tikzpicture}[anchorbase,scale=1]
\draw[mor] (0,0) rectangle (1.2,0.5) node[black,pos=0.5]{$x$};
\draw[mor] (0,0.5) to (1.2,0.5) to (1.1,1) to (0.1,1) to (0,0.5) node[black] at (0.6,0.75){${+1}$};
\draw[mor] (0.1,1) rectangle (1.1,1.5) node[black,pos=0.5]{$y$};
\draw[usual] (1.1, 1.3) to (1.4,1.3) to (1.4,0);
\end{tikzpicture} \;, \end{equation}
living in $\Hom^0(B_x B_i,B_y)$. Similarly, let $\iota_{(x,i,y)} \in \Hom^0(B_y,B_x B_i)$ denote the dual map
\begin{equation} \iota_{(x,i,y)} := \begin{tikzpicture}[anchorbase,scale=1]
\begin{scope}[yscale=-1]
\draw[mor] (0,0) rectangle (1.2,0.5) node[black,pos=0.5]{$x$};
\draw[mor] (0,0.5) to (1.2,0.5) to (1.1,1) to (0.1,1) to (0,0.5) node[black] at (0.6,0.75){${+1}$};
\draw[mor] (0.1,1) rectangle (1.1,1.5) node[black,pos=0.5]{$y$};
\draw[usual] (1.1, 1.3) to (1.4,1.3) to (1.4,0);
\end{scope}
\end{tikzpicture} \;, \end{equation}
called \emph{$(x,i,y)$-inclusion}.
\end{Definition}

\begin{Lemma} The morphisms $p_{(x,i,y)}$ and $\iota_{(x,i,y)}$ are nonzero, and span their respective morphism spaces in degree zero. \end{Lemma}

\begin{proof} The morphism spaces are one-dimensional since $\mu(x,y) = 1$ by assumption. Composing $p_{(x,i,y)}$ with an $i$-colored dot yields the $(x,y)$-dot, which is nonzero.
\begin{equation} \begin{tikzpicture}[anchorbase,scale=1]
\draw[mor] (0,0) rectangle (1.2,0.5) node[black,pos=0.5]{$x$};
\draw[mor] (0,0.5) to (1.2,0.5) to (1.1,1) to (0.1,1) to (0,0.5) node[black] at (0.6,0.75){${+1}$};
\draw[mor] (0.1,1) rectangle (1.1,1.5) node[black,pos=0.5]{$y$};
\draw[usual,marked=1] (1.1, 1.3) to (1.4,1.3) to (1.4,0);
\end{tikzpicture}
=
\begin{tikzpicture}[anchorbase,scale=1]
\draw[mor] (0,0) rectangle (1.2,0.5) node[black,pos=0.5]{$x$};
\draw[mor] (0,0.5) to (1.2,0.5) to (1.1,1) to (0.1,1) to (0,0.5) node[black] at (0.6,0.75){${+1}$};
\draw[mor] (0.1,1) rectangle (1.1,1.5) node[black,pos=0.5]{$y$};
\end{tikzpicture}.
\end{equation}The proof is complete. \end{proof}

Note that $p_{(x,i,y)}$ kills the one-tensor of $B_x B_i$ (which lives in the summand $B_w$ and not in $B_y$).

In any direct sum decomposition of $B_x B_i$, the projection map to $B_y$ must be a scalar multiple of $p_{(x,i,y)}$, and the inclusion map a multiple of $\iota_{(x,i,y)}$, in such a way that their composition is the identity of $B_y$. To understand the scalars which appear, we call upon the local composition pairing or local intersection pairing, see \cite[Definition 11.70]{ElMaThWi-soergel}. This is the bilinear pairing 
\[ \Hom^0(B_x B_i,B_y) \times \Hom^0(B_y,B_x B_i) \to \R, \qquad (p,i) \mapsto \idcoeff{p \circ i}. \]
That is, one sends two morphisms to the coefficient of the identity map in their degree zero composition, which is an element of $\End^0(B_y) = \R \cdot \id$. Since we have fixed bases for both hom spaces in question, the local intersection pairing becomes a $1 \times 1$ matrix, which we interpret simply as a scalar.

\begin{Definition} Continue the notation of \autoref{Defn:xiyproject} The \emph{local intersection form} of $y$ at $(x,i)$ is the scalar $\LIF(x,i,y) := \idcoeff{p_{(x,i,y)} \circ \iota_{(x,i,y)}} \in \C$. Pictorially we have
\begin{equation} \label{LIFpicture}
\begin{tikzpicture}[anchorbase,scale=1]
\draw[mor] (0,0) rectangle (1.2,0.5) node[black,pos=0.5]{$x$};
\draw[mor] (0,0.5) to (1.2,0.5) to (1.1,1) to (0.1,1) to (0,0.5) node[black] at (0.6,0.75){${+1}$};
\draw[mor] (0.1,1) rectangle (1.1,1.5) node[black,pos=0.5]{$y$};
\draw[usual] (1.1, 1.3) to (1.4,1.3) to (1.4,0);
\begin{scope}[yscale=-1, shift={(0,0)}]
\draw[mor] (0,0) rectangle (1.2,0.5) node[black,pos=0.5]{$x$};
\draw[mor] (0,0.5) to (1.2,0.5) to (1.1,1) to (0.1,1) to (0,0.5) node[black] at (0.6,0.75){${+1}$};
\draw[mor] (0.1,1) rectangle (1.1,1.5) node[black,pos=0.5]{$y$};
\draw[usual] (1.1, 1.3) to (1.4,1.3) to (1.4,0);
\end{scope}
\end{tikzpicture} 
=
\LIF(x,i,y) \cdot \id_{B_y}.
\end{equation}
By convention, we set $\LIF(x,i,xi) = 1$.
\end{Definition}

\begin{Lemma} Continue the notation of \autoref{Defn:xiyproject}. Then $\LIF(x,i,y) \ne 0$. \end{Lemma}

\begin{proof} If $\LIF(x,i,y) = 0$ then the composition of $\Hom^0(B_y,B_x B_i)$ with $\Hom^0(B_x B_i,B_y)$ is identically zero, and it impossible for $B_y$ to be a summand of $B_x B_i$. \end{proof}

Since $p_{(x,i,y)} \circ \iota_{(x,i,y)} = \LIF(x,i,y) \cdot \id_{B_y}$, we deduce that
\begin{equation} \label{defqy} q_y := \frac{1}{\LIF(x,i,y)} \iota_{(x,i,y)} \circ p_{(x,i,y)}, \end{equation}
is an idempotent projecting to $B_y$ inside $B_x B_i$. In pictures we have
\begin{equation} \LIF(x,i,y) \cdot q_y := \;
\begin{tikzpicture}[anchorbase,scale=1]
\draw[mor] (0,0) rectangle (1.2,0.5) node[black,pos=0.5]{$w$};
\draw[mor] (0,0.5) to (1.2,0.5) to (1.1,1) to (0.1,1) to (0,0.5) node[black] at (0.6,0.75){${+1}$};
\draw[mor] (0.1,1) rectangle (1.1,1.5) node[black,pos=0.5]{$y$};
\draw[mor] (0.1,1.5) rectangle (1.1,2) node[black,pos=0.5]{$y$};
\draw[mor] (0.1,2) to (1.1,2) to (1.2,2.5) to (0,2.5) to (0.1,2) node[black] at (0.6,2.25) {${+1}$};
\draw[mor] (0,2.5) rectangle (1.2,3) node[black,pos=0.5]{$w$};
\draw[usual] (1.1, 1.3) to (1.4,1.3) to (1.4,0);
\draw[usual] (1.1, 1.7) to (1.4,1.7) to (1.4,3);
\end{tikzpicture}
\; = \; \begin{tikzpicture}[anchorbase,scale=1]
\draw[mor] (0,0) rectangle (1.2,0.5) node[black,pos=0.5]{$x$};
\draw[mor] (0,0.5) to (1.2,0.5) to (1.1,1) to (0.1,1) to (0,0.5) node[black] at (0.6,0.75){${+1}$};
\draw[mor] (0.1,1) rectangle (1.1,1.5) node[black,pos=0.5]{$y$};
\draw[mor] (0.1,1.5) to (1.1,1.5) to (1.2,2) to (0,2) to (0.1,1.5) node[black] at (0.6,1.75) {${+1}$};
\draw[mor] (0,2) rectangle (1.2,2.5) node[black,pos=0.5]{$x$};
\draw[usual] (1.1,1.25) to (1.4,1.25);
\draw[usual] (1.4,0) to (1.4,2.5);
\node at (1.5,1.25) {$i$};
\end{tikzpicture}, \end{equation}
where for the second equality we have used \eqref{triassoc}.

\begin{Theorem} \label{Thm:inductiveclasp} Suppose that $w, x \in W$ and $i \in S$ are such that $w = xi > x$, and $\mu(x,y) = 1$ for all $y$ such that $yi < y < x$. For each such $y$ let $q_y$ be defined as in \eqref{defqy}. Then 
\begin{equation} \label{defnexi} e_{x,i} := \id_{B_x B_i} - \sum_{yi < y < x} q_y, \end{equation}
agrees with the clasp in $\End(B_x B_i)$. In diagrams we have
\begin{equation}\label{exidiag} \begin{tikzpicture}[anchorbase,scale=1]
\draw[mor] (0,0) rectangle (1,0.5)node[black,pos=0.5]{$x \otimes i$};
\end{tikzpicture} \; = \; \raisebox{0.05cm}{\begin{tikzpicture}[anchorbase,scale=1]
\draw[mor] (0,0) rectangle (1,0.5)node[black,pos=0.5]{$x$};
\draw[usual] (1.2,0) to (1.2,0.5);
\node at (1.3,0.4) {$i$};
\end{tikzpicture}} \;
- \sum_{y \muless x, yi < y} \LIF(x,i,y)^{-1} \begin{tikzpicture}[anchorbase,scale=1]
\draw[mor] (0,0) rectangle (1.2,0.5) node[black,pos=0.5]{$x$};
\draw[mor] (0,0.5) to (1.2,0.5) to (1.1,1) to (0.1,1) to (0,0.5) node[black] at (0.6,0.75){${+1}$};
\draw[mor] (0.1,1) rectangle (1.1,1.5) node[black,pos=0.5]{$y$};
\draw[mor] (0.1,1.5) to (1.1,1.5) to (1.2,2) to (0,2) to (0.1,1.5) node[black] at (0.6,1.75) {${+1}$};
\draw[mor] (0,2) rectangle (1.2,2.5) node[black,pos=0.5]{$x$};
\draw[usual] (1.1,1.25) to (1.4,1.25);
\draw[usual] (1.4,0) to (1.4,2.5);
\node at (1.5,1.25) {$i$};
\end{tikzpicture}. \end{equation}
\end{Theorem}

\begin{proof} We have argued above that $q_y$ is an idempotent projecting to $B_y$. It is clear that $q_y q_{y'} = 0$ for $y \ne y'$, since $\Hom^0(B_y,B_{y'}) = 0$. Thus $\sum q_y$
is an idempotent projecting to the sum of all lower summands, and $\id - \sum q_y$ is a top idempotent. Since a clasp idempotent exists, it is equal to this top idempotent. \end{proof}

\begin{Remark} It may also be pedagogically helpful to see an alternative proof which uses \autoref{Cor:perverseorthocriterion}. We need only prove that $e_{x,i}$ is orthogonal to
lower terms of non-positive degree. As $B_x B_i$ is perverse, the only nonzero morphisms $B_x B_i \to B_y$ of non-positive degree occur when $B_y$ is a summand of $B_x B_i$, and such
maps are spanned by $p_{(x,i,y)}$. So we need to show that $p_{(x,i,y)} e_{x,i} = 0$ for all $y$ with $yi < y < x$. We have $p_{(x,i,y)} q_{y'} = 0$ since $\Hom^0(B_y,B_{y'}) = 0$. It
is also straightforward to show that $p_{(x,i,y)} q_y = p_{(x,i,y)}$, and from there to deduce the result. \end{Remark}

We pause for a quick corollary.

\begin{Corollary} \label{Cor:putsomedotsonit} With the same assumptions as \autoref{Thm:inductiveclasp}, we have
\begin{equation}
\begin{tikzpicture}[anchorbase,scale=1]
\draw[mor] (0,0) rectangle (1,0.5)node[black,pos=0.5]{$x \otimes i$};
\draw[usual] (.9,.5) to (.9,.7);
\fill[black] (0.9,.7) circle (.055cm);
\draw[usual] (.9,0) to (.9,-.2);
\fill[black] (.9,-.2) circle (.055cm);
\end{tikzpicture} \; = \;
\begin{tikzpicture}[anchorbase,scale=1]
\draw[mor] (0,0) rectangle (1,0.5)node[black,pos=0.5]{$x$};
\end{tikzpicture} \cdot \alpha_i \;
- \sum_{y \muless x, yi < y} \LIF(x,i,y)^{-1} \begin{tikzpicture}[anchorbase,scale=1]
\draw[mor] (0,0) rectangle (1.2,0.5) node[black,pos=0.5]{$x$};
\draw[mor] (0,0.5) to (1.2,0.5) to (1.1,1) to (0.1,1) to (0,0.5) node[black] at (0.6,0.75){${+1}$};
\draw[mor] (0.1,1) rectangle (1.1,1.5) node[black,pos=0.5]{$y$};
\draw[mor] (0.1,1.5) to (1.1,1.5) to (1.2,2) to (0,2) to (0.1,1.5) node[black] at (0.6,1.75) {${+1}$};
\draw[mor] (0,2) rectangle (1.2,2.5) node[black,pos=0.5]{$x$};
\end{tikzpicture}\; . \end{equation}
\end{Corollary}

\begin{proof} This follows from \eqref{exidiag} together with the fact that
\begin{equation} \begin{tikzpicture}[anchorbase,scale=1]
\draw[mor] (0,0) rectangle (1.2,0.5) node[black,pos=0.5]{$w$};
\draw[mor] (0,0.5) to (1.2,0.5) to (1.1,1) to (0.1,1) to (0,0.5) node[black] at (0.6,0.75) {${+1}$};
\draw[mor] (0.1,1) rectangle (1.1,1.5) node[black,pos=0.5]{$y$};
\draw[mor] (0.1,1.5) to (1.1,1.5) to (1.2,2) to (0,2) to (0.1,1.5) node[black] at (0.6,1.75) {${+1}$};
\draw[mor] (0,2) rectangle (1.2,2.5) node[black,pos=0.5]{$w$};
\draw[usual] (1.1,1.25) to (1.3,1.25);
\draw[usual] (1.3,1.25) to (1.3,2);
\draw[usual] (1.3,1.25) to (1.3,0.5);
\node at (1.4,1.25) {$i$};
\fill[black] (1.3,2) circle (.055cm);
\fill[black] (1.3,0.5) circle (.055cm);
\end{tikzpicture} \; = \; \begin{tikzpicture}[anchorbase,scale=1]
\draw[mor] (0,0) rectangle (1.2,0.5) node[black,pos=0.5]{$w$};
\draw[mor] (0,0.5) to (1.2,0.5) to (1.1,1) to (0.1,1) to (0,0.5) node[black] at (0.6,0.75) {${+1}$};
\draw[mor] (0.1,1) rectangle (1.1,1.5) node[black,pos=0.5]{$y$};
\draw[mor] (0.1,1.5) to (1.1,1.5) to (1.2,2) to (0,2) to (0.1,1.5) node[black] at (0.6,1.75) {${+1}$};
\draw[mor] (0,2) rectangle (1.2,2.5) node[black,pos=0.5]{$w$};
\draw[usual] (1.1,1.25) to (1.3,1.25);
\fill[black] (1.3,1.25) circle (.055cm);
\end{tikzpicture} \; = \;
\begin{tikzpicture}[anchorbase,scale=1]
\draw[mor] (0,0) rectangle (1.2,0.5) node[black,pos=0.5]{$w$};
\draw[mor] (0,0.5) to (1.2,0.5) to (1.1,1) to (0.1,1) to (0,0.5) node[black] at (0.6,0.75) {${+1}$};
\draw[mor] (0.1,1) rectangle (1.1,1.5) node[black,pos=0.5]{$y$};
\draw[mor] (0.1,1.5) to (1.1,1.5) to (1.2,2) to (0,2) to (0.1,1.5) node[black] at (0.6,1.75) {${+1}$};
\draw[mor] (0,2) rectangle (1.2,2.5) node[black,pos=0.5]{$w$};
\end{tikzpicture} \;.
\end{equation}
The first equality is the usual unit axiom, and the second equality is \eqref{triunit}. \end{proof}

To use \autoref{Thm:inductiveclasp} effectively we need a method to compute local intersection forms. To compute \eqref{LIFpicture}, we would want to expand the idempotent $e_x$ describing $B_x$
(within some reasonable $x$-object like $B_{\un{x}}$, or $B_z B_j$ with $x = zj > z$) as a linear combination of simpler diagrams. If this idempotent $e_x$ is constructed inductively
in the same fashion, there is a chance for this process to be done algorithmically, in special cases. This is the algorithm we pursue in the rest of this section, starting with the next subsection.

It is not uncommon that $B_x B_i$ is itself indecomposable, which happens when $yi>y$ for all $y \muless x$. In this case $B_x B_i \cong B_w$, and the complicated part of
\eqref{exidiag} is an empty sum. This is in some sense the ``base case'' of our inductive computation. In particular, when $x$ is the identity element of $W$, then $B_x B_i$ is
indecomposable.

Now let us discuss the situation when $\mu = \mu(x,y) > 1$ for some $y < x$ with $xi > x$ and $yi < y$. Choose a basis $\phi_1, \ldots, \phi_{\mu}$ for $\Hom^1(B_x,B_y)$, and call it the $(w,y)$-dot family. Then $\Hom^0(B_x B_i, B_y)$ also has dimension $\mu$, and one can find a basis by composing the $(w,y)$ dot
family with the $(y,i)$ trivalent vertex.

\begin{Definition} Suppose that $xi > x$ and $yi < y < x$ and $\mu(x,y) \ge 1$. For any morphism $\phi \colon B_x \to B_y$ of degree $+1$, we define the corresponding projection $p_{\phi,i}$ as in \eqref{eq:pxiy}, but with $\phi$ replacing the $+1$-labeled trapezoid. One defines $\iota_{\phi,i}$ similarly.
\end{Definition}

Now the local intersection pairing is a $\mu \times \mu$ matrix, whose $(k,\ell)$-entry is the coefficient of $\id_y$ in the composition $p_{\phi_k} \circ \iota_{\phi_{\ell}}$. One
can construct orthogonal idempotents by finding dual bases for the local intersection pairing; indeed, the pairing must be non-degenerate or else $B_y$ would not appear as a summand
of $B_x B_i$ with multiplicity $\mu$. See \cite[Corollary 11.71]{ElMaThWi-soergel} for more details. Computing dual bases is certainly possible in examples, but doing it
algorithmically is a headache. Once these dual bases are found, one can construct the idempotent $e_{x,i}$ as a linear combination of diagrams just like in \eqref{exidiag}, but with
the trapezoids ranging over dual bases for each $y$. The form of $e_{x,i}$ will be important for further arguments, so let us record this conclusion in a theorem.

\begin{Theorem} \label{Thm:inductiveclaspalt} Suppose that $w, x \in W$ and $i \in S$ are such that $w = xi > x$. Then 
\begin{equation} \label{defnexialt} e_{x,i} := \id_{B_x B_i} - \sum_{yi < y < x} \sum_{k,\ell=1}^{\mu(x,y)} L^{-1}_{k \ell} \iota_{\phi_k,i} \circ p_{\phi_{\ell},i}, \end{equation}
where $L^{-1}_{k \ell}$ represents the $(k,\ell)$-entry of the inverse matrix of the local intersection pairing. We have the schematic picture
\begin{equation}\label{exidiagalt} \begin{tikzpicture}[anchorbase,scale=1]
\draw[mor] (0,0) rectangle (1,0.5)node[black,pos=0.5]{$x \otimes i$};
\end{tikzpicture} \; = \; \raisebox{0.05cm}{\begin{tikzpicture}[anchorbase,scale=1]
\draw[mor] (0,0) rectangle (1,0.5)node[black,pos=0.5]{$x$};
\draw[usual] (1.2,0) to (1.2,0.5);
\node at (1.3,0.4) {$i$};
\end{tikzpicture}} \;
- \sum \text{scalar} \; \cdot \;
\begin{tikzpicture}[anchorbase,scale=1]
\draw[mor] (0,0) rectangle (1.2,0.5) node[black,pos=0.5]{$x$};
\draw[mor] (0,0.5) to (1.2,0.5) to (1.1,1) to (0.1,1) to (0,0.5) node[black] at (0.6,0.75){${+1}$};
\draw[mor] (0.1,1) rectangle (1.1,1.5) node[black,pos=0.5]{$y$};
\draw[mor] (0.1,1.5) to (1.1,1.5) to (1.2,2) to (0,2) to (0.1,1.5) node[black] at (0.6,1.75) {${+1}$};
\draw[mor] (0,2) rectangle (1.2,2.5) node[black,pos=0.5]{$x$};
\draw[usual] (1.1,1.25) to (1.4,1.25);
\draw[usual] (1.4,0) to (1.4,2.5);
\node at (1.5,1.25) {$i$};
\end{tikzpicture}. \end{equation}
where $y$ ranges over elements $y \muless x$ with $yi<y$, and the trapezoids represent various morphisms of degree $1$.
\end{Theorem}

\subsection{A closed formula for some local intersection forms}\label{subsection:LIF}

In this section we compute $\LIF(x,i,y)$ whenever there exists $j \in S$ such that $yj = x$. We postpone the statement of the theorem to the end of the section, but one can view this
entire section as the proof. First we will produce a recursive formula for this local intersection form, and then we will solve that recursion.

\begin{Definition} Let $x, y \in W$ and $i \in S$. We call a triple $(x,i,y)$ \emph{recursible} if $y<x<xi$ and $B_y$ is a summand of $B_x B_i$, and there exists $j \in S$ such that $yj = x$. Though $j$ is determined by $x$ and $y$, we may include $j$ in the notation, and say that the quadruple $(x,i,y,j)$ is recursible. \end{Definition}

\begin{Remark} We will also use recursible triples later when studying recursive formulas for partial traces. \end{Remark}

Note that $\mu(x,y) = 1$ whenever $(x,i,y)$ is recursible. We state the following lemma mostly to establish notation for this section.

\begin{Lemma} \label{Lemma:recursiblenotation} Suppose that $(x,i,y,j)$ is recursible. Then $i \ne j$ and $m_{ij} > 2$, and there exist  unique $z,w \in W$ such that $z < zi = y < yj = x < xi = w$. \end{Lemma}

\begin{proof} If $B_y$ is a summand of $B_x B_i$ then $yi<y$. Since $y \muless x$ we have $y < x = yj$, so $i \ne j$. Let $z = yi < y$. Since $z < y < x$ and $x = zij$ we must have $\ell(z) = \ell(x)-2$, so that $zij$ is a reduced composition. If $m_{ij}=2$ then $x = zji$ has $i$ in its right descent set, a contradiction. \end{proof}

For the rest of the section, fix $(x,i,y,j)$ recursible and let $z, w \in W$ be as in the lemma.

We have described the $(x,y)$-dot in \eqref{wydotdefspecialcase}. Thus
\begin{equation} \label{LIFcomputationspecial1}
\begin{tikzpicture}[anchorbase,scale=1]
\draw[mor] (0,0) rectangle (1.2,0.5) node[black,pos=0.5]{$x$};
\draw[mor] (0,0.5) to (1.2,0.5) to (1.1,1) to (0.1,1) to (0,0.5) node[black] at (0.6,0.75){${+1}$};
\draw[mor] (0.1,1) rectangle (1.1,1.5) node[black,pos=0.5]{$y$};
\draw[usual] (1.1, 1.3) to (1.4,1.3) to (1.4,0);
\begin{scope}[yscale=-1, shift={(0,0)}]
\draw[mor] (0,0) rectangle (1.2,0.5) node[black,pos=0.5]{$x$};
\draw[mor] (0,0.5) to (1.2,0.5) to (1.1,1) to (0.1,1) to (0,0.5) node[black] at (0.6,0.75){${+1}$};
\draw[mor] (0.1,1) rectangle (1.1,1.5) node[black,pos=0.5]{$y$};
\draw[usual] (1.1, 1.3) to (1.4,1.3) to (1.4,0);
\end{scope}
\end{tikzpicture}
= 
\begin{tikzpicture}[anchorbase,scale=1]
\draw[rex] (0,0.5) rectangle (1.2,0.6);
\diagrammaticmorphism{0}{0.6}{3}{0,0,0,0,0}
\draw[rex] (0,0.9) rectangle (1,1);
\draw[mor] (0,1) rectangle (1,1.5) node[black,pos=0.5]{$y$};
\draw[usual] (1.1,0.6) to (1.1,0.75);
\fill[black] (1.1,0.75) circle (.055cm);
\node at (1.25,0.9) {$j$};
\draw[usual] (1.0, 1.25) to (1.4,1.25) to (1.4,0);
\node at (1.5,.6) {$i$};
\begin{scope}[yscale=-1,shift={(0,-0.5)}]
\draw[mor] (0,0) rectangle (1.2,0.5) node[black,pos=0.5]{$x$};
\draw[rex] (0,0.5) rectangle (1.2,0.6);
\diagrammaticmorphism{0}{0.6}{3}{0,0,0,0,0}
\draw[rex] (0,0.9) rectangle (1,1);
\draw[mor] (0,1) rectangle (1,1.5) node[black,pos=0.5]{$y$};
\draw[usual] (1.1,0.6) to (1.1,0.75);
\fill[black] (1.1,0.75) circle (.055cm);
\node at (1.25,0.9) {$j$};
\draw[usual] (1.0, 1.3) to (1.4,1.3) to (1.4,0);
\end{scope}
\end{tikzpicture}
.
\end{equation} 
We can then describe $B_x$ as the image of the clasp idempotent in $B_y B_j$, using  \eqref{exidiagalt}. We let $u$ (rather than $y$) be our variable indicating a lower summand inside $B_y B_j$. Thus
\begin{equation} \label{LIFcomputationspecial2}	\begin{tikzpicture}[anchorbase,scale=1]
\draw[rex] (0,0.5) rectangle (1.2,0.6);
\diagrammaticmorphism{0}{0.6}{3}{0,0,0,0,0}
\draw[rex] (0,0.9) rectangle (1,1);
\draw[mor] (0,1) rectangle (1,1.5) node[black,pos=0.5]{$y$};
\draw[usual] (1.1,0.6) to (1.1,0.75);
\fill[black] (1.1,0.75) circle (.055cm);
\node at (1.25,0.9) {$j$};
\draw[usual] (1.0, 1.25) to (1.4,1.25) to (1.4,0);
\begin{scope}[yscale=-1,shift={(0,-0.5)}]
\draw[mor] (0,0) rectangle (1.2,0.5) node[black,pos=0.5]{$x$};
\draw[rex] (0,0.5) rectangle (1.2,0.6);
\diagrammaticmorphism{0}{0.6}{3}{0,0,0,0,0}
\draw[rex] (0,0.9) rectangle (1,1);
\draw[mor] (0,1) rectangle (1,1.5) node[black,pos=0.5]{$y$};
\draw[usual] (1.1,0.6) to (1.1,0.75);
\fill[black] (1.1,0.75) circle (.055cm);
\node at (1.25,0.9) {$j$};
\draw[usual] (1.0, 1.3) to (1.4,1.3) to (1.4,0);
\end{scope}
\end{tikzpicture}
=
\begin{tikzpicture}[anchorbase,scale=1]
\diagrammaticmorphism{0}{0.5}{5}{0,0,0,0,0}
\draw[mor] (0,1) rectangle (1,1.5) node[black,pos=0.5]{$y$};
\draw[usual] (1.1,0) to (1.1,0.75);
\fill[black] (1.1,0.75) circle (.055cm);
\node at (1.25,0.9) {$j$};
\draw[usual] (1, 1.25) to (1.4,1.25) to (1.4,0);
\begin{scope}[yscale=-1,shift={(0,-0.5)}]
\draw[mor] (0,0) rectangle (1,0.5) node[black,pos=0.5]{$y$};
\diagrammaticmorphism{0}{0.5}{5}{0,0,0,0,0}
\draw[mor] (0,1) rectangle (1,1.5) node[black,pos=0.5]{$y$};
\draw[usual] (1.1,0.5) to (1.1,0.75);
\fill[black] (1.1,0.75) circle (.055cm);
\node at (1.25,0.9) {$j$};
\draw[usual] (1, 1.3) to (1.4,1.3) to (1.4,0);
\end{scope}
\end{tikzpicture}
-  \raisebox{-1ex}{\scalebox{2}{$\Sigma$}},
\end{equation}
where $\Sigma$ is a linear combination of diagrams of the form
\begin{equation} \label{LIFcomputationspecial3} \begin{tikzpicture}[anchorbase,scale=1]
\draw[mor] (0,0) rectangle (1.2,0.5) node[black,pos=0.5]{$y$};
\draw[mor] (0,0.5) to (1.2,0.5) to (1.1,1) to (0.1,1) to (0,0.5);
\node[black] at (0.6,0.75) {${+1}$};
\draw[mor] (0.1,1) rectangle (1.1,1.5) node[black,pos=0.5]{$u$};
\draw[mor] (0.1,1.5) to (1.1,1.5) to (1.2,2) to (0,2) to (0.1,1.5); 
\node[black] at (0.6,1.75) {${+1}$};
\draw[mor] (0,2) rectangle (1.2,2.5) node[black,pos=0.5]{$y$};
\draw[usual] (1.2,0.25) to (1.4,0.25) to (1.4,2.25) to (1.2,2.25);
\node at (1.5, 1.25) {$i$};
\end{tikzpicture}, \; \end{equation}
for various $u \muless y$ with $uj<u$.


The first diagram on the right of \eqref{LIFcomputationspecial2} simplifies using \eqref{triassoc} and \eqref{triunit} and the needle relation as
\begin{equation}	\begin{tikzpicture}[anchorbase,scale=1]
\diagrammaticmorphism{0}{0.5}{5}{0,0,0,0,0}
\draw[mor] (0,1) rectangle (1,1.5) node[black,pos=0.5]{$y$};
\draw[usual] (1.1,0) to (1.1,0.75);
\fill[black] (1.1,0.75) circle (.055cm);
\node at (1.25,0.9) {$j$};
\draw[usual] (1, 1.25) to (1.4,1.25) to (1.4,0);
\begin{scope}[yscale=-1,shift={(0,-0.5)}]
\draw[mor] (0,0) rectangle (1,0.5) node[black,pos=0.5]{$y$};
\diagrammaticmorphism{0}{0.5}{5}{0,0,0,0,0}
\draw[mor] (0,1) rectangle (1,1.5) node[black,pos=0.5]{$y$};
\draw[usual] (1.1,0.5) to (1.1,0.75);
\fill[black] (1.1,0.75) circle (.055cm);
\node at (1.25,0.9) {$j$};
\draw[usual] (1, 1.3) to (1.4,1.3) to (1.4,0);
\end{scope}
\end{tikzpicture}
=
\begin{tikzpicture}[anchorbase,scale=1]
\diagrammaticmorphism{0}{0.5}{5}{0,0,0,0,0}
\draw[mor] (0,1) rectangle (1,1.5) node[black,pos=0.5]{$y$};
\draw[usual] (1, 1.25) to (1.4,1.25) to (1.4,0);
\node at (1.5,0.6) {$i$};
\begin{scope}[yscale=-1,shift={(0,-0.5)}]
\diagrammaticmorphism{0}{0}{10}{0,0,0,0,0}
\draw[mor] (0,1) rectangle (1,1.5) node[black,pos=0.5]{$y$};
\node at (1.2,0.25) {$\alpha_j$};
\draw[usual] (1, 1.3) to (1.4,1.3) to (1.4,0);
\end{scope}
\end{tikzpicture}
=
\begin{tikzpicture}[anchorbase,scale=1]
\draw[mor] (0,1) rectangle (1,1.5) node[black,pos=0.5]{$y$};
\draw[usual] (1, 1.25) to (1.25,1.25);
\draw[usual] (1.6,1.25) circle (.35cm);
\node at (1.6,1.25) {$\alpha_j$};
\end{tikzpicture}
=
\begin{tikzpicture}[anchorbase,scale=1]
\draw[mor] (0,0) rectangle (1,0.5) node[black,pos=0.5]{$y$};
\end{tikzpicture} \; \cdot \partial_i(\alpha_j). \end{equation}
In particular, if $B_y B_j$ is indecomposable (the sum over $u$ is empty), then we have just computed $\LIF(x,i,y)$, and it equals $\partial_i(\alpha_j)$.

Now focus on the diagram in \eqref{LIFcomputationspecial3} for a given $u \muless y$. The following lemma computes this diagram.

\begin{Lemma} Let $y, u \in W$ and $i \in S$. Suppose that $yi<y$ and $u \muless y$. We have
\begin{equation} \label{uz} \begin{tikzpicture}[anchorbase,scale=1]
\draw[mor] (0,0) rectangle (1.2,0.5) node[black,pos=0.5]{$y$};
\draw[mor] (0,0.5) to (1.2,0.5) to (1.1,1) to (0.1,1) to (0,0.5);
\node[black] at (0.6,0.75) {${+1}$};
\draw[mor] (0.1,1) rectangle (1.1,1.5) node[black,pos=0.5]{$u$};
\draw[mor] (0.1,1.5) to (1.1,1.5) to (1.2,2) to (0,2) to (0.1,1.5); 
\node[black] at (0.6,1.75) {${+1}$};
\draw[mor] (0,2) rectangle (1.2,2.5) node[black,pos=0.5]{$y$};
\draw[usual] (1.2,0.25) to (1.4,0.25) to (1.4,2.25) to (1.2,2.25);
\node at (1.5, 1.25) {$i$};
\end{tikzpicture} \; = \delta_{u,z} \id_y. \end{equation}
To elaborate, the result is zero (for any trapezoidal morphisms of degree $+1$) unless $u = yi < y$. When $u = yi$, we use our convention that the trapezoid represents the $(y,u)$ dot as in \eqref{wydotdefspecialcase}, in which case the result is $\id_y$.
\end{Lemma}

\begin{proof} There are two cases to consider, the case $ui > u$ and the case $ui<u$. First we treat the case $ui>u$.

The diagram in question is the composition of a degree zero morphism $B_y \to B_u B_i$ and a degree zero morphism $B_u B_i \to B_y$. If $ui > u$ then $B_u B_i$ is perverse, so it admits a nonzero degree zero morphism to $B_y$ if and only if $B_y$ is a direct summand of $B_u B_i$. Every summand of $B_u B_i$ has the form $B_v$ for $\ell(v) \le \ell(u) + 1$, with equality holding only for $v = ui$. Since $u < y$, the only way that $\ell(y) \le \ell(u) + 1$ is if equality holds, and thus $y = ui$ and $u=z$.

If $u=z$ then the $(y,u)$-dot is described by \eqref{wydotdefspecialcase}. Using \eqref{trislide} and the unit axiom we have
\begin{equation} \begin{tikzpicture}[anchorbase,scale=1]
\draw[mor] (0,0) rectangle (1.2,0.5) node[black,pos=0.5]{$y$};
\draw[mor] (0,0.5) to (1.2,0.5) to (1.1,1) to (0.1,1) to (0,0.5);
\node[black] at (0.6,0.75) {${+1}$};
\draw[mor] (0.1,1) rectangle (1.1,1.5) node[black,pos=0.5]{$z$};
\draw[mor] (0.1,1.5) to (1.1,1.5) to (1.2,2) to (0,2) to (0.1,1.5); 
\node[black] at (0.6,1.75) {${+1}$};
\draw[mor] (0,2) rectangle (1.2,2.5) node[black,pos=0.5]{$y$};
\draw[usual] (1.2,0.25) to (1.4,0.25) to (1.4,2.25) to (1.2,2.25);
\node at (1.5, 1.25) {$i$};
\end{tikzpicture}
=
\begin{tikzpicture}[anchorbase,scale=1]
\draw[mor] (0,0) rectangle (1.2,0.5) node[black,pos=0.5]{$y$};
\draw[usual] (1.2,0.25) to (1.5,0.25) to (1.5,2.25) to (1.2,2.25);
\draw[rex] (0,0.5) rectangle (1.2,0.6);
\diagrammaticmorphism{0}{0.6}{3}{0,0,0,0,0}
\draw[rex] (0,0.9) rectangle (1,1);
\draw[rex] (0,1.5) rectangle (1,1.6);
\draw[usual] (1.1,0.6) to (1.1,0.75);
\fill[black] (1.1,0.75) circle (.055cm);
\node at (1.25,0.9) {$i$};
\begin{scope}[yscale=-1,shift={(0,-2.5)}]
\draw[mor] (0,0) rectangle (1.2,0.5) node[black,pos=0.5]{$y$};
\draw[rex] (0,0.5) rectangle (1.2,0.6);
\diagrammaticmorphism{0}{0.6}{3}{0,0,0,0,0}
\draw[rex] (0,0.9) rectangle (1,1);
\draw[mor] (0,1) rectangle (1,1.5) node[black,pos=0.5]{$z$};
\draw[rex] (0,1.5) rectangle (1,1.6);
\draw[usual] (1.1,0.6) to (1.1,0.75);
\fill[black] (1.1,0.75) circle (.055cm);
\node at (1.25,0.9) {$i$};
\end{scope}
\end{tikzpicture}
=
\begin{tikzpicture}[anchorbase,scale=1]
\draw[mor] (0,0) rectangle (1.2,0.5) node[black,pos=0.5]{$y$};
\draw[usual] (1.1,0.7) to (1.5,0.7) to (1.5,1.8) to (1.1,1.8);
\draw[rex] (0,0.5) rectangle (1.2,0.6);
\diagrammaticmorphism{0}{0.6}{3}{0,0,0,0,0}
\draw[rex] (0,0.9) rectangle (1,1);
\draw[rex] (0,1.5) rectangle (1,1.6);
\draw[usual] (1.1,0.6) to (1.1,0.75);
\fill[black] (1.1,0.75) circle (.055cm);
\node at (1.25,0.9) {$i$};
\begin{scope}[yscale=-1,shift={(0,-2.5)}]
\draw[mor] (0,0) rectangle (1.2,0.5) node[black,pos=0.5]{$y$};
\draw[rex] (0,0.5) rectangle (1.2,0.6);
\diagrammaticmorphism{0}{0.6}{3}{0,0,0,0,0}
\draw[rex] (0,0.9) rectangle (1,1);
\draw[mor] (0,1) rectangle (1,1.5) node[black,pos=0.5]{$z$};
\draw[rex] (0,1.5) rectangle (1,1.6);
\draw[usual] (1.1,0.6) to (1.1,0.75);
\fill[black] (1.1,0.75) circle (.055cm);
\node at (1.25,0.9) {$i$};
\end{scope}
\end{tikzpicture}
=
\begin{tikzpicture}[anchorbase,scale=1]
\draw[mor] (0,0) rectangle (1.2,0.5) node[black,pos=0.5]{$y$};
\draw[usual] (1.1,0.6) to (1.1,1.9);
\draw[rex] (0,0.5) rectangle (1.2,0.6);
\diagrammaticmorphism{0}{0.6}{3}{0,0,0,0,0}
\draw[rex] (0,0.9) rectangle (1,1);
\draw[rex] (0,1.5) rectangle (1,1.6);
\node at (1.25,0.9) {$i$};
\begin{scope}[yscale=-1,shift={(0,-2.5)}]
\draw[mor] (0,0) rectangle (1.2,0.5) node[black,pos=0.5]{$y$};
\draw[rex] (0,0.5) rectangle (1.2,0.6);
\diagrammaticmorphism{0}{0.6}{3}{0,0,0,0,0}
\draw[rex] (0,0.9) rectangle (1,1);
\draw[mor] (0,1) rectangle (1,1.5) node[black,pos=0.5]{$z$};
\draw[rex] (0,1.5) rectangle (1,1.6);
\end{scope}
\end{tikzpicture}
=
\begin{tikzpicture}[anchorbase,scale=1]
\draw[mor] (0,0) rectangle (1.2,0.5) node[black,pos=0.5]{$y$};
\end{tikzpicture}
\end{equation}

If $ui<z$ then $B_u B_i$ is not perverse and we must use a different argument. Using \eqref{trislideoverplusone} and \eqref{triassoc} and the needle relation we have
\begin{equation} \begin{tikzpicture}[anchorbase,scale=1]
\draw[mor] (0,0) rectangle (1.2,0.5) node[black,pos=0.5]{$y$};
\draw[mor] (0,0.5) to (1.2,0.5) to (1.1,1) to (0.1,1) to (0,0.5);
\node[black] at (0.6,0.75) {${+1}$};
\draw[mor] (0.1,1) rectangle (1.1,1.5) node[black,pos=0.5]{$u$};
\draw[mor] (0.1,1.5) to (1.1,1.5) to (1.2,2) to (0,2) to (0.1,1.5); 
\node[black] at (0.6,1.75) {${+1}$};
\draw[mor] (0,2) rectangle (1.2,2.5) node[black,pos=0.5]{$y$};
\draw[usual] (1.2,0.25) to (1.4,0.25) to (1.4,2.25) to (1.2,2.25);
\node at (1.5, 1.25) {$i$};
\end{tikzpicture} \; = \;
\begin{tikzpicture}[anchorbase,scale=1]
\draw[mor] (0,0) rectangle (1.2,0.5) node[black,pos=0.5]{$y$};
\draw[mor] (0,0.5) to (1.2,0.5) to (1.1,1) to (0.1,1) to (0,0.5) node[black] at (0.6,0.75){${+1}$};
\draw[mor] (0.1,1) rectangle (1.1,1.5) node[black,pos=0.5]{$u$};
\draw[mor] (0.1,1.5) rectangle (1.1,2) node[black,pos=0.5]{$u$};
\draw[mor] (0.1,2) to (1.1,2) to (1.2,2.5) to (0,2.5) to (0.1,2) node[black] at (0.6,2.25) {${+1}$};
\draw[mor] (0,2.5) rectangle (1.2,3) node[black,pos=0.5]{$y$};
\draw[usual] (1.1,1.25) to (1.4,1.25) to (1.4,1.75) to (1.1,1.75);
\end{tikzpicture}
\; = \;
\begin{tikzpicture}[anchorbase,scale=1]
\draw[mor] (0,0) rectangle (1.2,0.5) node[black,pos=0.5]{$y$};
\draw[mor] (0,0.5) to (1.2,0.5) to (1.1,1) to (0.1,1) to (0,0.5);
\node[black] at (0.6,0.75) {${+1}$};
\draw[mor] (0.1,1) rectangle (1.1,1.5) node[black,pos=0.5]{$u$};
\draw[mor] (0.1,1.5) to (1.1,1.5) to (1.2,2) to (0,2) to (0.1,1.5); 
\node[black] at (0.6,1.75) {${+1}$};
\draw[mor] (0,2) rectangle (1.2,2.5) node[black,pos=0.5]{$y$};
\draw[usual] (1.1,1.25) to (1.4,1.25) to (1.4,1.5) to (1.6,1.5) to (1.6,1) to (1.4,1) to (1.4,1.25);
\end{tikzpicture}
\; = 0. \end{equation}The proof is complete.
\end{proof}

In conclusion, while $B_y B_j$ may have many summands $B_u$, each giving a term in \eqref{LIFcomputationspecial2}, at most one of these terms is nonzero: a term indexed by $u=z$. We know that $z \muless y$. Hence $B_z$ will be a summand of $B_y B_j$ if and only if $zj<z$.

\begin{Theorem} \label{Thm:LIFrecursion} Let $(x,i,y,j)$ be recursible and let $z = yi$. Then
\begin{equation} \label{eq:LIFrecursion} \LIF(x,i,y) = \begin{cases} \partial_i(\alpha_j) & \text{ if } zj>z, \\ \partial_i(\alpha_j) - \LIF(y,j,z)^{-1} & \text{ if } zj<z. \end{cases} \end{equation}
\end{Theorem}

\begin{proof} Given the previous lemma, it remains to note that when $u = z$ the diagram \eqref{LIFcomputationspecial3} appears in the sum $\Sigma$ from \eqref{LIFcomputationspecial2} with coefficient $\LIF(y,j,z)^{-1}$ (see e.g. \autoref{Thm:inductiveclasp}). \end{proof}

Before solving this recursive formula, we give some illustrative examples.

\begin{Example} \label{ex:justifycoeff1}
We continue \autoref{extendingnastypart1} and \autoref{nastyexamplenew}. We claim that $\LIF(31213,2,1232) = \partial_{2}(\alpha_1)=-1$. Here we have $j=1$ and $z = 123$, and note that $1232 = zi$ and $31213 = 1232j$ and $zj>z$. Similarly, $\LIF(31213,2,3212) = \partial_2(\alpha_3)=-1$. This time we have $j=3$ and $z = 321$, and again $zj > z$. Plugging these numbers $-1$ into the formula from \autoref{Thm:inductiveclasp}, we verify the computation of $e_0$ from \autoref{extendingnastypart1}.
\end{Example}

\begin{Example}
Consider the dihedral group of type $G_2$, and the realization where $\partial_1(\alpha_2) = -3$ and $\partial_2(\alpha_1) = -1$. We verify the known computations of the top idempotents, which in this case come from colored Jones--Wenzl projectors, see \cite{El-two-color-soergel}. The Bott--Samelson bimodules $B_1$ and $B_1 B_2$ are indecomposable.

We have $B_{12} \ot B_1 \cong B_{121} \oplus B_1$, so we must compute $\LIF(12,1,1)$. Then $j=2$ and $z = \id_W$ so $zj>z$. Hence
\begin{equation}
\LIF(12,1,1) = \partial_{1}(\alpha_2) = -3.
\end{equation}

Next we have $B_{121} \ot B_2 \cong B_{1212} \oplus B_{12}$, so we must compute $\LIF(121,2,12)$. Then $j=1$ and $z = 1$ and $zj < z$. Therefore
\begin{equation}
\LIF(121,2,12) = \partial_{2}(\alpha_1) - \LIF(12,1,1)^{-1} = -1 + \frac{1}{3} = \frac{-2}{3}.
\end{equation}

Next we have $B_{1212} \ot B_1 \cong B_{12121} \oplus B_{121}$, so we must compute $\LIF(1212,1,121)$. Then $j=2$ and $z = 12$ and $zj < z$. Therefore
\begin{equation}
\LIF(1212,1,121) = \partial_{1}(\alpha_2) - \LIF(121,2,12)^{-1} = -3 + \frac{3}{2} = \frac{-3}{2}.
\end{equation}

Finally we have $B_{12121} \ot B_2 \cong B_{121212} \oplus B_{1212}$, so we must compute $\LIF(12121,2,1212)$. Then $j=1$ and $z = 121$ and $zj < z$. Therefore
\begin{equation}
\LIF(12121,2,1212) = \partial_{2}(\alpha_1) - \LIF(1212,1,121)^{-1} = -1 + \frac{2}{3} = \frac{-1}{3}.
\end{equation}So we have $-1/3$.
\end{Example}

The closed form solution to $\LIF(x,i,y)$ involves two-colored quantum numbers, originally defined in \cite[Appendix A.1]{El-two-color-soergel}. For the $G_2$ example above, the two-colored quantum numbers in question form the sequence $1, 3, 2, 3, 1$, and the local intersection forms are signed ratios of successive two-colored quantum numbers.


\begin{Definition} Let $i \ne j \in S$. Define $[2]_{i,j} := -\partial_i(\alpha_j)$ and $[2]_{j,i} = -\partial_j(\alpha_i)$. Let $[1]_{i,j} = [1]_{j,i} = 1$ and $[0]_{i,j} = [0]_{j,i} = 0$. For $k > 2$ set
\begin{equation} [k]_{i,j} = [2]_{i,j} [k-1]_{j,i} - [k-2]_{i,j}. \end{equation}
\end{Definition}

Recall that $W_{i,j}$ is the parabolic subgroup generated by $i$ and $j$, which may be infinite.

\begin{Theorem} \label{Thm:LIFrecursionWORKS} Let $(x,i,y,j)$ be recursible. Let $v$ be the minimal representative of the right $W_{\{i,j\}}$-coset containing $x$, and let $x', y' \in W_{\{i,j\}}$ be the unique elements such that $x = vx'$ and $y = vy'$. Then $(x',i,y',j)$ is recursible, and $\LIF(x,i,y) = \LIF(x',i,y')$. Moreover, we have
\begin{equation} \label{eq:LIFviaquantumnumbers} \LIF(x',i,y') = (-1) \cdot \frac{[\len(x')]_{i,j}}{[\len(x')-1]_{j,i}}, \end{equation}
a ratio of successive two-colored quantum numbers.
\end{Theorem}

\begin{proof} Suppose $(x,i,y,j)$ is recursible and let $z = yi$. Then \eqref{eq:LIFrecursion} computes $\LIF(x,i,y)$, but in the case
when $zj < z$ we need to be able to compute $\LIF(y,j,z)$ to use this formula. However, we have $z < y < yj$ with $B_z$ being a summand of $B_y B_j$, and $i \in S$ satisfies $zi = y$.
Thus $(y,j,z,i)$ is recursible, and is computable using \eqref{eq:LIFrecursion} as well. This recursive process continues, staying within the right $W_{\{i,j\}}$-coset containing $x$. If it terminates when computing $\LIF(x,i,y)$, then it terminates because $zj > z$. But if $zi>z$ and $zj > z$ then $z = v$ is the minimal element in its right $W_{\{i,j\}}$-coset.

Let us be explicit about this recursive process. Let $n = \len(x')$, and $u_n = x'$. Note that $x' i > x'$ and $x' j < x'$, so $1 \le \len(x') < m_{ij}$. For $0 \le k \le n$ let $u_k$ represent the unique element in $W_{\{i,j\}}$ with length $k$ which is less than $x'$ in the weak right Bruhat order; in other words, extending $u_k$ by an alternating word $(\ldots,i,j)$ of length $n-k$ gives a reduced expression for $x'$. For example $u_0$ is the identity element, and $u_1$ is either $i$ or $j$, whichever is in the left descent set of $x'$. Let $s_k \in \{i,j\}$ be such that $u_k s_k = u_{k+1}$. The formula \eqref{eq:LIFrecursion} applies to both $(x',i,y')$ and $(x,i,y)$. It implies that for all $k \ge 3$ we have 
\begin{equation} \LIF(u_k,s_k,u_{k-1}) = \partial_{s_k}(\alpha_{s_{k-1}}) - \LIF(u_{k-1},s_{k-1},u_{k-2})^{-1}, \end{equation}
\begin{equation} \LIF(vu_k,s_k,vu_{k-1}) = \partial_{s_k}(\alpha_{s_{k-1}}) - \LIF(vu_{k-1},s_{k-1},vu_{k-2})^{-1}, \end{equation}
and for $k=2$ we have
\begin{equation} \LIF(u_2,s_2,u_1) = \partial_{s_2}(\alpha_{s_1}), \quad \LIF(vu_2,s_2,vu_1) = \partial_{s_2}(\alpha_{s_1}). \end{equation}
(There is no need to consider $\LIF(u_1,s_1,u_0)$ since $u_0 < u_0 s_1$ so $B_{u_0}$ cannot be a summand of $B_{u_1} B_{s_1}$.)
In particular, the same recursion applies to compute $\LIF(x,i,y)$ as to compute $\LIF(x',i,y')$, with the same base for the recursion, so these values are equal.
Note that $\partial_{s_k}(\alpha_{s_{k-1}})$ is either $-[2]_{i,j}$ or $-[2]_{j,i}$ depending on parity. That \eqref{eq:LIFviaquantumnumbers} solves the recursion follows directly
from the definition of the two-colored quantum numbers. \end{proof}

\begin{Remark} Note that we do not assert that $B_x \cong B_v B_{x'}$, which is false in general, nor do we assert that the decomposition of $B_x B_i$ or $B_y B_j$ into direct summands mimics the decomposition of $B_{x'} B_i$ or $B_{y'} B_j$. There may be many more summands in $B_y B_j$ than in $B_{y'} B_j$, but as shown in \eqref{uz}, these other summands do not contribute to the particular local intersection form being studied. \end{Remark}

\begin{Remark} In \cite{ElLi-universal-soergel}, the first author and Libedinsky study universal Coxeter groups (where $m_{ij} = \infty$ for all $i \ne j \in S$). In the language of the present paper, they proved that $B_y$ is a summand of $B_x B_i$ (for $y \ne xi > x$) if and only if $(x,i,y)$ is recursible, and they gave an independent proof that the local intersection forms are equal to ratios of successive two-colored quantum numbers. \end{Remark}

\subsection{Other local intersection forms}\label{subsection:the rest}

To compute the clasp idempotent for $w$ in $B_x B_i$ we need to compute all $\LIF(x,i,y)$ for $yi<y$, $y \muless x$.  Of course, many such triples $(x,i,y)$ are not recursible, including all triples where $\ell(y) < \ell(x) - 1$, or where $\ell(y) = \ell(x)-1$ and $y<x$ in the Bruhat order but not in the weak right Bruhat order.

\begin{Example} Let $w = S_4$ with the usual simple reflections $S = \{1,2,3\}$. Let $x = 2132$ and $y = 212$. Then $B_y$ is a summand of $B_x B_1$, but there is no $j \in S$ such
that $yj = x$. The local intersection form $\LIF(x,1,y)$ is not computable using \eqref{eq:LIFrecursion}, but it is easy enough to compute by hand, especially since $B_x \cong B_2 B_1
B_3 B_2$ is a Bott--Samelson object. We leave this exercise to the interested reader. In this case we have that \begin{equation} \LIF(x,1,y) = \partial_1(s_2(\alpha_3)) = -1.
\end{equation}as one easily computes. \end{Example}

When some $(x,i,y)$ is not recursible, we can still compute the local intersection form, but it takes creativity and/or additional work. To do this, the steps are to:
\begin{enumerate}
\item Compute the idempotent describing $B_x$ inside some reasonable $x$-object, like a Bott--Samelson bimodule.
\item Find a suitable morphism to serve as the $(x,y)$-dot.
\item Compute the composition \eqref{LIFpicture}, or at least determine its coefficient of the identity.
\end{enumerate}
The first step is done inductively by this same process, so let us presume that idempotents for all elements less than $xi$ are already computed. More precisely, the idempotent formula \eqref{exidiag} gives a formula for $B_x$ inside $B_y B_j$ for some $y$, rather than an idempotent inside a Bott--Samelson bimodule $\un{x}$. However, we also have a formula for the idempotent describing $B_y$ inside $B_z B_i$, etc. Iterating this process, for $\un{x} \in \Rex(x)$ we could construct an idempotent $e_{\un{x}} \in \End(B_{\un{x}})$ projecting to $B_x$ (not necessarily a clasp). Practically speaking, the choice of $\un{x}$ can make this iterative process easier or more difficult.

The remaining steps can be done by hand, but one can also make use of computer software. Existing software \cite{Gi-ASLoc-code} can transform Soergel diagrams (morphisms between Bott--Samelson bimodules) into sparse matrices of Laurent polynomials, and multiply those matrices to compose morphisms (or do something else funky and explicit to tensor morphisms). Thus if the problem is rephrased in terms of a computation involving Soergel diagrams, computers can perform this computation for us, in theory. We now focus on this rephrasing, having already described all smaller idempotents as linear combinations of diagrams.

When $\ell(y) = \ell(x)-1$ we always have a straightforward choice for the $(x,y)$-dot, as in \eqref{wydotdef}. When $\ell(y) < \ell(x)-1$ we must find one. Given $\un{x} \in
\Rex(x)$ and $\un{y} \in \Rex(y)$, the light leaves construction \cite[\S 6]{ElWi-soergel-calculus} gives a basis for all morphisms from $B_{\un{x}}$ to $B_{\un{y}}$ modulo
$\catstuff{I}_{< y}$. It is a finite amount of work to look for all light leaves $\phi$ of degree $+1$, and precompose with the idempotent $e_{\un{x}}$ projecting to $B_x$. So long
as the result is nonzero, then $e_y \circ \phi \circ e_{\un{x}}$ is a valid $(x,y)$-dot. If $\mu(x,y) > 1$ then find all light leaves of degree $+1$, and do linear algebra to find a
set of size $\mu(x,y)$ which remain linearly independent after composing with $e_{\un{x}}$.

Once one has found an $(x,y)$-dot (or a family thereof) then the composition \eqref{LIFpicture} is a composition of linear combinations of diagrams, and can be done by computer. Conveniently, the
identity coefficient is one of the matrix coefficients in the matrix representation of a morphism.

Thus computing local intersection forms is possible by computer, and this has roughly been the process used by Williamson, Jensen, Gibson, and others to compute $p$-KL polynomials
\cite{GiJeWi}. Whether it is practical, given the limitations of time and memory, is quite another story. For a Bott--Samelson associated to an expression of length $d$, the matrices
involved have size $2^d$, and though they are sparse, this becomes prohibitive very quickly. We were happy to be able to bypass many of these computations with
\autoref{Thm:LIFrecursionWORKS}.

\begin{Remark} In the sequel \cite{ElRoTu-ah-2}, those local intersections we were not able to compute with \autoref{Thm:LIFrecursionWORKS} we found either by hand or by computer as above, with numerous practical simplifications to make the problem tractable. \end{Remark}

\subsection{Branching graphs}

We now define a graph-theoretical bookkeeping device that we will heavily use.

\begin{Definition} Fix a Coxeter system $(W,S)$. The \emph{branching graph} $\Gamma(W,S)$ is an oriented graph whose edges are labeled/colored by the simple reflections $S$. The vertices are in bijection with $W$, such that there is an $i$-colored edge from $x$ to $y$ if $xi>x$ and $B_y$ is a direct summand of $B_x B_i$.

When we discuss a path in the branching graph, the default assumption is that the path begins at the identity element of $W$, a vertex we typically denote as $\emptyset$. Given a path in the branching graph, the labels on the edges give a sequence of simple reflections, the \emph{expression} underlying the path. Conversely, given an expression $\un{w}$, there are potentially many paths (starting at the identity) whose underlying expression is $\un{w}$, and we say those paths are \emph{subordinate} to $\un{w}$. The union of all these we might denote by $\Gamma(\un{w})$.

When $\un{w} = (s_1, \ldots, s_d)$ is a reduced expression for $w$, there is a unique subordinate path which ends in $w$, and we call this path the \emph{spine} of $\Gamma(\un{w})$. To spell this out more precisely, if $\un{k}$ is the element expressed by $(s_1, \ldots, s_k)$ (with $\un{0}$ being $\emptyset$), then the spine has $d$ total edges, with an $s_k$-colored edge from $\un{k-1}$ to $\un{k}$ for all $1 \le k \le d$. Finally, the \emph{superspine} is the subgraph of $\Gamma(\un{w})$ which includes the spine, as well as all $s_k$-colored edges from $\un{k-1}$ which leave the spine (and their targets). Loosely, the superspine includes the main reduced path along $\un{w}$ as well as the first step along any false starts. We denote the superspine by $\something(\un{w})$. \end{Definition}

\begin{Example}\label{example:p-graph}
If $m_{12}\geq 5$ and $\un{w} = (1,2,1,2,1)$, then 
\begin{gather} \label{12121graph}
\Gamma(\un{w}) = \something(\un{w}) =
\begin{tikzcd}[ampersand replacement=\&,row sep=scriptsize,column sep=scriptsize,arrows={shorten >=-0.5ex,shorten <=-0.5ex},labels={inner sep=0.05ex},arrow style=tikz]
\emptyset \ar[r,yshift=0.1cm,soergelone]
\&
1\ar[r,yshift=0.1cm,soergeltwo]
\&
12\ar[r,yshift=0.1cm,soergelone]
\ar[l,yshift=-0.1cm,soergelone]
\&
121\ar[r,yshift=0.1cm,soergeltwo]
\ar[l,yshift=-0.1cm,soergeltwo]
\&
1212\ar[r,yshift=0.1cm,soergelone]
\ar[l,yshift=-0.1cm,soergelone]
\&
12121
\end{tikzcd}
,
\end{gather}
and the graph is a line.
\end{Example}

\begin{Example} \label{ex:H3nastygraph} In type $H_3$, with $m_{12} = 5$, $m_{23} = 3$, $m_{13} = 2$ and $\un{w} = (1,2,1,2,1,3,2,1,2,1,3,2,1,2)$, the graph $\Gamma(\un{w})$ begins as in the previous example. We now illustrate how the branching graph continues, starting at $\un{5} = 12121$.
\begin{gather}\label{eq:h3c5-pgraphBen}
\begin{tikzcd}[ampersand replacement=\&,row sep=scriptsize,column sep=scriptsize,arrows={shorten >=-0.5ex,shorten <=-0.5ex},labels={inner sep=0.05ex},arrow style=tikz]
\underline{5}\ar[r,yshift=0.1cm,soergelthree]
\& 
\underline{6}\ar[r,yshift=0.1cm,soergeltwo]\ar[l,yshift=-0.1cm,soergeltwo]
\& 
\underline{7}\ar[r,yshift=0.1cm,soergelone]\ar[l,yshift=-0.1cm,soergelone]
\& 
\underline{8}\ar[r,yshift=0.1cm,soergeltwo]\ar[l,yshift=-0.1cm,soergeltwo]
\& 
\underline{9}\ar[r,yshift=0.1cm,soergelone]\ar[l,yshift=-0.1cm,soergelone]
\& 
\underline{10}\ar[r,yshift=0.1cm,soergelthree]
\& 
\underline{11}\ar[r,yshift=0.1cm,soergeltwo]\ar[l,yshift=-0.1cm,soergeltwo]\ar[ld,soergeltwo, yshift=-0.1cm]
\& 
\underline{12}\ar[r,yshift=0.1cm,soergelone]\ar[l,yshift=-0.1cm,soergelone]
\& 
\underline{13}\ar[r,yshift=0.1cm,soergeltwo]\ar[l,yshift=-0.1cm,soergeltwo]
\& 
\underline{14} \\
\& \& \& \& \& x=\underline{9}3\ar[ru,yshift=0.1cm,soergelone]
\& \& \& \& \\
\end{tikzcd}.
\end{gather}
Any path subordinate to $\un{w}$ whose first five steps go from $\emptyset$ straight to $\un{5}$ will continue along a path in \eqref{eq:h3c5-pgraphBen}. The subgraph $\something(\un{w})$ is almost identical to the union of \eqref{eq:h3c5-pgraphBen} and \eqref{12121graph}; one need only leave out the blue edge from $x$ to $\un{11}$.
The graph $\Gamma(\un{w})$ is larger than the union of \eqref{eq:h3c5-pgraphBen} and \eqref{12121graph}, though.  If the first five steps of a path subordinate to $\un{w}$ end up at $1$ or $121$, then the sixth step will go to $13$ or $1213$. The superspine does not include the edge from $1$ to $13$ or from $121$ to $1213$, even though these edges leave the spine, since they leave the spine after taking several steps along a false start.

Let $\un{w}' = (1,2,1,2,1,3,2,1,2,3,1,2,1,2)$, differing from $\un{w}$ in locations $10$ and $11$ where the commuting reflections $1$ and $3$ have been swapped. The spine of $\Gamma(\un{w}')$ factors through $x$ instead of $\un{10}$. The union of $\something(\un{w})$ and $\something(\un{w}')$ is the union of \eqref{12121graph} and the more symmetric graph
\begin{gather}\label{eq:h3c5-pgraph0early}
\begin{tikzcd}[ampersand replacement=\&,row sep=scriptsize,column sep=scriptsize,arrows={shorten >=-0.5ex,shorten <=-0.5ex},labels={inner sep=0.05ex},arrow style=tikz]
\underline{5}\ar[r,yshift=0.1cm,soergelthree]
\& 
\underline{6}\ar[r,yshift=0.1cm,soergeltwo]\ar[l,yshift=-0.1cm,soergeltwo]
\& 
\underline{7}\ar[r,yshift=0.1cm,soergelone]\ar[l,yshift=-0.1cm,soergelone]
\& 
\underline{8}\ar[r,yshift=0.1cm,soergeltwo]\ar[l,yshift=-0.1cm,soergeltwo]
\& 
\underline{9}\ar[r,yshift=0.1cm,soergelone]\ar[l,yshift=-0.1cm,soergelone]\ar[rd,yshift=0.1cm,soergelthree]\ar[l,yshift=-0.2cm,soergelthree]
\& 
\underline{10}\ar[r,yshift=0.1cm,soergelthree]
\& 
\underline{11}\ar[r,yshift=0.1cm,soergeltwo]\ar[l,yshift=-0.1cm,soergeltwo]\ar[ld,tomato, yshift=-0.1cm]
\& 
\underline{12}\ar[r,yshift=0.1cm,soergelone]\ar[l,yshift=-0.1cm,soergelone]
\& 
\underline{13}\ar[r,yshift=0.1cm,soergeltwo]\ar[l,yshift=-0.1cm,soergeltwo]
\& 
\underline{14} \\
\& \& \& \& \& x=\underline{9}3\ar[ru,yshift=0.1cm,soergelone]
\& \& \& \& \\
\end{tikzcd}.
\end{gather}
Both graphs are already more involved than the line graph in 
\autoref{example:p-graph}.
\end{Example}

\begin{Lemma} In $\Gamma(W)$, a vertex $y$ is the source of an $i$-colored edge if and only if $yi > y$. A vertex $z$ is the target of an $i$-colored edge if and only if $zi < z$. Consequently, no vertex can be both the source and target of an $i$-colored edge, and every vertex is either the source or target of at least one $i$-colored edge. \end{Lemma}

\begin{proof} By definition of the graph $\Gamma(W)$, there cannot be an $i$-colored edge leaving $y$ unless $yi>y$. If $yi>y$ then there is an $i$-colored edge $y \to yi$, and if $zi < z$ then there is an $i$-colored edge $zi \to z$. Moreover, if there is an $i$-colored edge $y \to z$, then $B_z$ is a summand of $B_y B_i$, implying that $zi < z$. \end{proof}

\begin{Remark} Most vertices in \eqref{eq:h3c5-pgraph0early} abut arrows of all three colors, so one can easily determine their right descent sets. Vertices $\un{7}$, $\un{12}$,
and $\un{13}$ do not abut a $3$-colored arrow, and they all have $3$ in their right descent set; this means they are a target of a $3$-colored arrow in the larger graph
$\Gamma(W)$. \end{Remark}

\begin{Example} \label{Ex:D4branching} In type $D_4$, where $2$ is the center vertex of the Dynkin diagram, let $\un{d} = 143$ and $\un{w} = \un{d}2341$. The superspine of $\Gamma(\un{w})$ is
\begin{gather}
\begin{tikzcd}[ampersand replacement=\&,row sep=scriptsize,column sep=scriptsize,arrows={shorten >=-0.5ex,shorten <=-0.5ex},labels={inner sep=0.05ex},arrow style=tikz]
\emptyset \ar[r,yshift=0.1cm,soergelone]
\&
1 \ar[r,yshift=0.1cm,soergelfour]
\&
14 \ar[r,yshift=0.1cm,soergelthree]
\&
\underline{d}\ar[r,yshift=0.2cm,soergeltwo]
\&
\underline{d}2 \ar[r,yshift=0.1cm,soergelthree] \ar[l,yshift=0cm,soergelthree]
\&
\underline{d}23 \ar[r,yshift=0.1cm,soergelfour]
\&
\underline{d}234 \ar[r,yshift=0.1cm,soergelone]\ar[lll,bend left=15,yshift=-0.2cm,soergelone,"{2}"]
\&
\un{w} 
\end{tikzcd}
\end{gather} 
The ``long edge'' going backwards is labeled with the corresponding local intersection form.\end{Example}

\begin{Example} \label{Ex:D5branching} Here is an example in type $D_5$, with vertices $4$ and $5$ being the two tines of the Dynkin diagram ($m_{45} = 2$). Let $\underline{d}=1535435$ and $\un{w} = \underline{d}234531$. The superspine of $\Gamma(\un{w})$ starting with $\un{d}$ is
\begin{gather}
\begin{tikzcd}[ampersand replacement=\&,row sep=scriptsize,column sep=scriptsize,arrows={shorten >=-0.5ex,shorten <=-0.5ex},labels={inner sep=0.05ex},arrow style=tikz]
\underline{d} \ar[r,yshift=0.1cm,soergeltwo]
\&
\underline{d}2 \ar[r,yshift=0.2cm,soergelthree]\ar[l,yshift=-0.1cm,soergelthree]
\&
\underline{d}23 \ar[r,yshift=0.1cm,soergelfour]\ar[l,yshift=0cm,soergelfour]\ar[rd,yshift=0.1cm,soergelfive]\ar[l,yshift=-0.2cm,soergelfive]
\&
\underline{d}234 \ar[r,yshift=0.1cm,soergelfive]
\&
\underline{d}2345 \ar[r,yshift=0.1cm,soergelthree]\ar[l,yshift=-0.1cm,soergelthree]\ar[ld,yshift=-0.2cm,soergelthree]
\&
\underline{d}23453 \ar[r,yshift=0.1cm,soergelone]\ar[lllll,bend right=15,yshift=0.2cm,soergelone,"{2}"]
\&
\un{w}
\\
\&
\&
\&
\underline{d}235\ar[ru,yshift=0cm,soergelfour]
\&
\&
\&
\end{tikzcd}
\end{gather}
Again, the ``long edge'' going backwards has been labeled with the corresponding local intersection form. All other local intersection forms are all equal to $-1 = \partial_i(\alpha_{i \pm 1})$, as guaranteed by \autoref{Thm:LIFrecursionWORKS}. The value $2$ for the long edge (and the long edge in \autoref{Ex:D4branching}) was calculated by computer.
\end{Example}

\begin{Remark}\label{remark:p-graph} Calculating branching graphs only involves calculations on the Grothendieck level and can easily been done via computer, if necessary. See
\autoref{subsection:hecke-computer}. \end{Remark}

\begin{Definition} \label{def:recursible} We say that a reduced expression $\un{w}$ is \emph{recursible} if every edge which is in the superspine but not the spine corresponds to a triple $(x,i,y)$ which is recursible. Here $x$ is the source of the edge, $y$ the target, and $i$ the label. We say that $w \in W$ is \emph{recursible} if it has a recursible reduced expression.

We say that a reduced expression $\un{w}$ is \emph{linear} if $\something(\un{w})$ is a \emph{line}. More precisely, we mean that there is at most one edge leaving the $k$-th vertex of the spine which is not in the spine, and it goes to the $(k-1)$-st vertex of the spine. We say that $w \in W$ is \emph{linear} if it has a linear reduced expression. \end{Definition}

A linear reduced expression is a particularly easy example of a recursible reduced expression. The element $12121$ from \autoref{example:p-graph} is linear. The expression $\un{w}$ from \autoref{ex:H3nastygraph} is not linear, but it is recursible. Obviously if $\un{w}$ is linear/recursible, and $\un{k}$ is the expression given by first $k$ letters of $\un{w}$, then $\un{k}$ is linear/recursible.  The expressions from \autoref{Ex:D4branching} and \autoref{Ex:D5branching} are not recursible, because of their ``long edges.''

If $\un{w}$ is recursible, then every local intersection form along $\something(\un{w})$ is computable by \autoref{Thm:LIFrecursionWORKS}. By iterating \autoref{Thm:inductiveclasp}, we obtain a formula for the idempotent inside $B_{\un{w}}$ projecting to $B_w$.

\subsection{Partial traces}\label{subsection:traces}

It will be important when computing categorical dimensions to consider partial traces of clasp idempotents, which we now define.

\begin{Definition} Let $f \in R$, and let $x, w \in W$ and $i \in S$ be such that $w = xi > x$. Let $e_w$ denote the clasp idempotent within $B_x B_i$. The diagram below is an endomorphism of $B_x$, called the \emph{partial trace} $\ptr_{x,i}(f)$.
\begin{equation} \label{eq:ptrdef} \begin{tikzpicture}[anchorbase,scale=1]
\draw[mor] (0,-0.5) rectangle (0.6,0) node[black,pos=0.5]{$x$};
\draw[mor] (0,0) rectangle (1,0.5) node[black,pos=0.5]{$w$};
\draw[mor] (0,0.5) rectangle (0.6,1) node[black,pos=0.5]{$x$};
\draw[soergelone] (0.75,0) to[out=270,in=180] (1,-0.25) to[out=0,in=270] (1.25,0) to (1.25,0.5) to[out=90,in=0] (1,0.75) to[out=180,in=90] (0.75,0.5);
\node[scale=0.7] at (1.1,0.25) {$f$};
\end{tikzpicture}\;. \end{equation}

More generally, if $\phi$ is any endomorphism of $B_w$, then its \emph{partial trace} $\ptr_{x,i}(\phi) \in \End(B_x)$ is the following morphism:
\begin{equation} \label{eq:ptrphidef} \begin{tikzpicture}[anchorbase,scale=1]
\draw[mor] (0,-0.5) rectangle (0.6,0) node[black,pos=0.5]{$x$};
\draw[mor] (0,0) rectangle (1,0.5) node[black,pos=0.5]{$w$};
\draw[mor] (0,0.5) rectangle (1,1) node[black,pos=0.5]{$\phi$};
\draw[mor] (0,1) rectangle (1,1.5) node[black,pos=0.5]{$w$};
\draw[mor] (0,1.5) rectangle (0.6,2) node[black,pos=0.5]{$x$};
\draw[soergelone] (0.75,0) to[out=270,in=180] (1,-0.25) to[out=0,in=270] (1.25,0) to (1.25,1.5) to[out=90,in=0] (1,1.75) to[out=180,in=90] (0.75,1.5);
\end{tikzpicture}\;. \end{equation}
Thus, $\ptr_{x,i}(f) = \ptr_{x,i}(\id_w \cdot f)$.
\end{Definition}

By \autoref{Lem:topstack} the idempotents $e_x$ in \eqref{eq:ptrdef} are redundant, and we often omit them. We displayed them above to emphasize that the diagram is an endomorphism of $B_x$.

\begin{Lemma} \label{lemma:AtoB} For any $x \in W$ and $i \in S$ with $w = xi > x$ we have $\ptr_{x,i}(\id_w) = 0$. More generally we have
\begin{equation} 
\begin{tikzpicture}[anchorbase,scale=1]
\draw[mor] (0.6,0) rectangle (1.6,0.5) node[black,pos=0.5]{\raisebox{-0.06cm}{$x \otimes i$}};
\draw[usual] (1.5,0) to[out=270,in=180] (1.7,-0.2) to[out=0,in=270] (1.9,0) to (1.9,0.5)to[out=90,in=0] (1.7,0.7) to[out=180,in=90] (1.5,0.5);
\node[scale=0.8] at (1.7,0.25) {$f$};	
\end{tikzpicture} \; = \; 	\begin{tikzpicture}[anchorbase,scale=1]
\draw[mor] (0.6,0) rectangle (1.6,0.5) node[black,pos=0.5]{\raisebox{-0.06cm}{$x \otimes i$}};
\draw[usual,marked=1] (1.5,0.5) to (1.5,0.9);
\draw[usual,marked=1] (1.5,0) to (1.5,-0.4);
\node[scale=0.8] at (2,0.25) {$\partial_i(f)$};	
\end{tikzpicture} \; . \end{equation}
\end{Lemma}

\begin{proof} We compute:
\begin{gather}
\begin{tikzpicture}[anchorbase,scale=1]
\draw[mor] (0.6,0) rectangle (1.6,0.5) node[black,pos=0.5]{\raisebox{-0.06cm}{$x \otimes i$}};
\draw[usual] (1.5,0) to[out=270,in=180] (1.7,-0.2) to[out=0,in=270] (1.9,0) to (1.9,0.5)to[out=90,in=0] (1.7,0.7) to[out=180,in=90] (1.5,0.5);
\node[scale=0.8] at (1.7,0.25) {$f$};	
\end{tikzpicture}	
=
\begin{tikzpicture}[anchorbase,scale=1]
\draw[mor] (0.6,0) rectangle (1.6,0.5) node[black,pos=0.5]{\raisebox{-0.06cm}{$x \otimes i$}};
\draw[usual,marked=1] (1.5,0.5) to (1.5,0.9);
\draw[usual,marked=1] (1.5,0) to (1.5,-0.4);
\draw[usual] (1.5,-0.2) to[out=0,in=270] (1.9,0) to (1.9,0.5)to[out=90,in=0] (1.5,0.7);
\node[scale=0.8] at (1.7,0.25) {$f$};	
\end{tikzpicture}	
=
\begin{tikzpicture}[anchorbase,scale=1]
\draw[mor] (0.6,0) rectangle (1.6,0.5) node[black,pos=0.5]{\raisebox{-0.06cm}{$x \otimes i$}};
\draw[mor] (0.6,0.5) rectangle (1.6,1) node[black,pos=0.5]{\raisebox{-0.06cm}{$x \otimes i$}};
\draw[usual,marked=1] (1.5,1) to (1.5,1.4);
\draw[usual,marked=1] (1.5,0) to (1.5,-0.4);
\draw[usual] (1.6,0.25) to[out=0,in=270] (1.9,0.5) to[out=90,in=0] (1.6,0.75);
\node[scale=0.8] at (1.7,0.5) {$f$};	
\end{tikzpicture}
=
\begin{tikzpicture}[anchorbase,scale=1]
\draw[mor] (0.6,0) rectangle (1.6,0.5) node[black,pos=0.5]{\raisebox{-0.06cm}{$x \otimes i$}};
\draw[usual,marked=1] (1.5,0.5) to (1.5,0.9);
\draw[usual,marked=1] (1.5,0) to (1.5,-0.4);
\node[scale=0.8] at (2.2,0.25) {$f$};	
\draw[usual] (2.2,0.25) circle (0.25);
\draw[usual]  (1.6,0.25) to (1.95,0.25);
\end{tikzpicture}
=
\begin{tikzpicture}[anchorbase,scale=1]
\draw[mor] (0.6,0) rectangle (1.6,0.5) node[black,pos=0.5]{\raisebox{-0.06cm}{$x \otimes i$}};
\draw[usual,marked=1] (1.5,0.5) to (1.5,0.9);
\draw[usual,marked=1] (1.5,0) to (1.5,-0.4);
\node[scale=0.8] at (2.2,0.25) {$\partial_i(f)$};	
\draw[usual,marked=1]  (1.6,0.25) to (1.8,0.25);
\end{tikzpicture}
=
\begin{tikzpicture}[anchorbase,scale=1]
\draw[mor] (0.6,0) rectangle (1.6,0.5) node[black,pos=0.5]{\raisebox{-0.06cm}{$x \otimes i$}};
\draw[usual,marked=1] (1.5,0.5) to (1.5,0.9);
\draw[usual,marked=1] (1.5,0) to (1.5,-0.4);
\node[scale=0.8] at (2,0.25) {$\partial_i(f)$};	
\end{tikzpicture} \; ,
\end{gather}
and we are done.
\end{proof}

The previous lemma combines with \autoref{Cor:putsomedotsonit} when the assumptions of \autoref{Thm:inductiveclasp} hold to give a precise formula for $\ptr_{x,i}(f)$.

We are also interested in iterated partial traces. For example, in the setting of \autoref{Thm:LIFrecursion}, we might want to consider $\ptr_{z,i}(\ptr_{y,j}(\ptr_{x,i}(f)))$, which we would abbreviate to $\ptr_{z,(i,j,i)}(f)$.

\begin{Definition} Let $\un{d}$ and $\un{t} = (t_1, \ldots, t_m)$ be reduced expressions (for elements $d$ and $t$) such that the concatenation $\un{dt}$ is reduced. Then for $\phi \in \End(B_{dt})$ set 
\begin{equation} \ptr_{d,\un{t}}(\phi) = \ptr_{d,t_1}(\ptr_{dt_1,t_2}(\cdots(\ptr_{dt_1 \cdots t_{m-1},t_m}(\phi))\cdots )), \end{equation}
and call it an \emph{iterated partial trace}. For example, when $\un{t} = (1,2,3)$ we have
\begin{equation}
\ptr_{d,\un{t}}(f) = \;
\begin{tikzpicture}[anchorbase,scale=1]
\draw[mor] (0,0.5) rectangle (0.5,1) node[black,pos=0.5]{\raisebox{-0.06cm}{$d$}};
\draw[mor] (0,0) rectangle (1.1,0.5) node[black,pos=0.5]{\raisebox{-0.06cm}{$\underline{d}123$}};
\draw[mor] (0,-0.5) rectangle (0.5,0) node[black,pos=0.5]{\raisebox{-0.06cm}{$d$}};
\draw[soergelthree] (1,0) to[out=270,in=180] (1.4,-0.4) to[out=0,in=270] (1.8,0) to (1.8,0.5)to[out=90,in=0] (1.4,0.9) to[out=180,in=90] (1,0.5);
\node[scale=0.8] at (1.4,0.25) {$f$};
\draw[soergeltwo] (0.8,0) to[out=270,in=180] (1.4,-0.6) to[out=0,in=270] (2,0) to (2,0.5)to[out=90,in=0] (1.4,1.1) to[out=180,in=90] (0.8,0.5);
\draw[soergelone] (0.6,0) to[out=270,in=180] (1.4,-0.8) to[out=0,in=270] (2.2,0) to (2.2,0.5)to[out=90,in=0] (1.4,1.3) to[out=180,in=90] (0.6,0.5);
\end{tikzpicture} \; .
\end{equation}
(Here we iterate three times.)
\end{Definition}

\subsection{Linear branching graphs and iterated partial traces} \label{subsec:linearpartialtraces}

Let $\un{d}$ and $\un{t} = (t_1, \ldots, t_m)$ be reduced expressions (for elements $d$ and $t$) such that the concatenation $\un{dt}$ is reduced. In this section we provide an algorithm to compute partial traces of polynomials when the branching graph of $d\un{t}$ (i.e. starting at $d$) is linear.

For this purpose, let us provide notation for a number of commonly-used diagrams. For $1 \le k \le m$ let
\begin{equation} A_k(f) = 
\begin{tikzpicture}[anchorbase,scale=1]
\draw[mor] (0,0) rectangle (1.6,0.5) node[black,pos=0.5]{\raisebox{-0.06cm}{$d~ t_1\ldots t_k$}};
\draw[mor] (0,0.5) rectangle (0.4,1) node[black,pos=0.5]{\raisebox{-0.06cm}{$d$}};
\draw[mor] (0,-0.5) rectangle (0.4,0) node[black,pos=0.5]{\raisebox{-0.06cm}{$d$}};
\draw[usual] (1.5,0) to[out=270,in=180] (1.7,-0.2) to[out=0,in=270] (1.9,0) to (1.9,0.5)to[out=90,in=0] (1.7,0.7) to[out=180,in=90] (1.5,0.5);
\node[scale=0.8] at (1.75,0.25) {$f$};	
\draw[usual] (0.5,0) to[out=270,in=180] (1.7,-0.5) to[out=0,in=270] (2.3,0) to (2.3,0.5)to[out=90,in=0] (1.7,1) to[out=180,in=90] (0.5,0.5);
\node[scale=0.8] at (1,-0.25) {$...$};
\node[scale=0.8] at (1,0.65) {$...$};
\node[scale=0.8] at (2.1,0.25) {$...$};
\end{tikzpicture}
\quad \in \End(B_d).	 \end{equation}
Let
\begin{equation} B_{k,1}(f) = 		\begin{tikzpicture}[anchorbase,scale=1]
\draw[mor] (0,0) rectangle (1.6,0.5) node[black,pos=0.5]{\raisebox{-0.06cm}{$d~ t_1\ldots t_k$}};
\draw[mor] (0,0.5) rectangle (0.4,1) node[black,pos=0.5]{\raisebox{-0.06cm}{$d$}};
\draw[mor] (0,-0.5) rectangle (0.4,0) node[black,pos=0.5]{\raisebox{-0.06cm}{$d$}};
\draw[usual, marked=1] (1.5,0.5) to (1.5,0.7);
\draw[usual, marked=1] (1.5,0) to (1.5,-0.2);
\draw[usual] (0.5,0) to[out=270,in=180] (1.7,-0.5) to[out=0,in=270] (2.3,0) to (2.3,0.5)to[out=90,in=0] (1.7,1) to[out=180,in=90] (0.5,0.5);
\node[scale=0.8] at (1,-0.25) {$...$};
\node[scale=0.8] at (1,0.65) {$...$};
\node[scale=0.8] at (2,0.25) {$f$};
\end{tikzpicture} \;, \end{equation}
which differs from $A_k(f)$ only in that the innermost loop is ``broken.'' More generally, for $a \ge 1$ let
\begin{equation} B_{k,a}(f) = 	\begin{tikzpicture}[anchorbase,scale=1]
\draw[mor] (0,0.5) rectangle (0.4,1) node[black,pos=0.5]{\raisebox{-0.06cm}{$d$}};
\draw[mor] (0,0) rectangle (1.6,0.5) node[black,pos=0.5]{\raisebox{-0.06cm}{$d~ t_1\ldots ~t_k$}};
\draw[mor] (0,-1) rectangle (1.6,-0.5) node[black,pos=0.5]{\raisebox{-0.06cm}{$d~ t_1\ldots ~t_k$}};
\draw[mor] (0,-2) rectangle (1.6,-1.5) node[black,pos=0.5]{\raisebox{-0.06cm}{$d~ t_1\ldots ~t_k$}};
\draw[mor] (0,-2.5) rectangle (0.4,-2) node[black,pos=0.5]{\raisebox{-0.06cm}{$d$}};	
\draw[usual] (1.3,-2) to[out=270,in=180] (1.7,-2.4) to[out=0,in=270] (2.1,-2) to (2.1,0.5)to[out=90,in=0] (1.7,0.9) to[out=180,in=90] (1.3,0.5);
\node[scale=0.8] at (1.8,-0.75) {$f$};	
\draw[usual] (0.5,-2) to[out=270,in=180] (1.7,-2.8) to[out=0,in=270] (2.7,-2) to (2.7,0.5)to[out=90,in=0] (1.7,1.3) to[out=180,in=90] (0.5,0.5);
\draw[usual,marked=1] (1.5,0.5) to (1.5,0.6);
\draw[usual,marked=1] (1.5,-2) to (1.5,-2.1);
\node[scale=0.8] at (0.8,-1.25) {$\vdots$};
\node[scale=0.8] at (1,0.65) {$...$};
\node[scale=0.8] at (2.4,-0.75) {$...$};
\draw[usual] (0.2,0) to (0.2,-0.5);	
\draw[usual] (1.3,0) to (1.3,-0.5);
\draw[usual] (0.5,0) to (0.5,-0.5);
\draw[usual,marked=1] (1.5,0) to (1.5,-0.1);
\draw[usual,marked=1] (1.5,-0.5) to (1.5,-0.4);
\end{tikzpicture} \; . \end{equation}
There are $2a$ dots in this picture. Meanwhile, let
\begin{equation} C_{k,a} = 	\begin{tikzpicture}[anchorbase,scale=1]
\draw[mor] (0,0.5) rectangle (0.4,1) node[black,pos=0.5]{\raisebox{-0.06cm}{$d$}};
\draw[mor] (0,0) rectangle (1.6,0.5) node[black,pos=0.5]{\raisebox{-0.06cm}{$d~ t_1\ldots ~t_k$}};
\draw[mor] (0,-1) rectangle (1.6,-0.5) node[black,pos=0.5]{\raisebox{-0.06cm}{$d~ t_1\ldots ~t_k$}};
\draw[mor] (0,-2) rectangle (1.6,-1.5) node[black,pos=0.5]{\raisebox{-0.06cm}{$d~ t_1\ldots ~t_k$}};
\draw[mor] (0,-2.5) rectangle (0.4,-2) node[black,pos=0.5]{\raisebox{-0.06cm}{$d$}};	
\draw[usual] (1.5,-2) to[out=270,in=180] (1.8,-2.4) to[out=0,in=270] (2.1,-2) to (2.1,0.5)to[out=90,in=0] (1.8,0.9) to[out=180,in=90] (1.5,0.5);
\node[scale=0.8] at (1.8,-0.75) {$f$};	
\draw[usual] (0.5,-2) to[out=270,in=180] (1.8,-2.8) to[out=0,in=270] (2.7,-2) to (2.7,0.5)to[out=90,in=0] (1.8,1.3) to[out=180,in=90] (0.5,0.5);
\node[scale=0.8] at (0.8,-1.25) {$\vdots$};
\node[scale=0.8] at (1,0.65) {$...$};
\node[scale=0.8] at (2.4,-0.75) {$...$};
\draw[usual] (0.2,0) to (0.2,-0.5);	
\draw[usual] (0.5,0) to (0.5,-0.5);
\draw[usual,marked=1] (1.5,0) to (1.5,-0.1);
\draw[usual,marked=1] (1.5,-0.5) to (1.5,-0.4);
\end{tikzpicture} \; .	 \end{equation}
There are $2a$ dots in this picture. By convention, $C_{k,0}(f) = A_k(f)$. By convention $A_0(f) = \id_{d} \cdot f$. To reiterate, all of these diagrams live in $\End(B_d)$.

The iterated partial trace we wish to compute is $A_m(f)$ for some polynomial $f$. The result is an endomorphism of $B_d$ of the same degree as $f$. We are most interested in the coefficient of the identity of this endomorphism, or equivalently, the value of this morphism modulo the ideal $\catstuff{I}_{< d}$. First we recall the construction of idempotents in this case, mostly to set our notation.

\begin{Theorem} \label{thm:linearidempotents} Let $\un{d}$ and $\un{t} = (t_1, \ldots, t_m)$ be reduced expressions (for elements $d$ and $t$) such that the concatenation $\un{dt}$ is reduced. Let $x_k$ be the element expressed by $\un{d} t_1 t_2 \cdots t_k$, so that $x_0 = d$ and $x_m = dt$. Assume that its branching graph between $d$ and $dt$ is linear, which is to say that
\begin{equation} B_{x_k} B_{t_{k+1}} \cong \begin{cases} B_{x_{k+1}} & \text{ if } x_{k-1} t_{k+1} > x_{k-1}, \\  B_{x_{k+1}} \oplus B_{x_{k-1}} & \text{ if } x_{k-1} t_{k+1} < x_{k-1}, \end{cases} \end{equation}
for all $1 \le k \le m-1$. We set $\delta_{k} = 1$ if $x_{k-1} t_{k+1} < x_{k-1}$, and $\delta_k = 0$ otherwise. Then 
\begin{equation} \label{exk1}
e_{x_{k+1}}
\; = \; \begin{tikzpicture}[anchorbase,scale=1]
\draw[mor] (0,0) rectangle (1,0.5)node[black,pos=0.5]{$x_k$};
\draw[usual] (1.2,0) to (1.2,0.5);
\node at (1.6,0.25) {$\scriptstyle t_{k+1}$};
\end{tikzpicture} \;
- \delta_k \LIF(x_k,t_{k+1},x_{k-1})^{-1} \; \cdot \;
\begin{tikzpicture}[anchorbase,scale=1]
\draw[mor] (0,0) rectangle (1.2,0.5) node[black,pos=0.5]{$x_k$};
\draw[usual,marked=1] (1.1,0.5) to (1.1,0.7);
\draw[usual] (0.8,0.5) to (0.85,1);
\node at (0.5,0.75) {$\scriptstyle \cdots$};
\draw[usual] (0.1,0.5) to (0.15,1);
\draw[mor] (0.1,1) rectangle (1.1,1.5) node[black,pos=0.5]{$x_{k-1}$};
\draw[usual,marked=1] (1.1,2) to (1.1,1.8);
\draw[usual] (0.8,2) to (0.85,1.5);
\draw[usual] (0.1,2) to (0.15,1.5);
\node at (0.5,1.75) {$\scriptstyle \cdots$};
\draw[mor] (0,2) rectangle (1.2,2.5) node[black,pos=0.5]{$x_k$};
\draw[usual] (1.1,1.25) to (1.4,1.25);
\draw[usual] (1.4,0) to (1.4,2.5);
\node at (1.8,1.25) {$\scriptstyle t_{k+1}$};
\end{tikzpicture}. \end{equation}
and consequently
\begin{equation} \label{exk1dots}
\begin{tikzpicture}[anchorbase,scale=1]
\draw[mor] (0,0) rectangle (1,0.5)node[black,pos=0.5]{$x_{k+1}$};
\draw[usual] (.9,.5) to (.9,.7);
\fill[black] (0.9,.7) circle (.055cm);
\draw[usual] (.9,0) to (.9,-.2);
\fill[black] (.9,-.2) circle (.055cm);
\end{tikzpicture} \; = \;
\begin{tikzpicture}[anchorbase,scale=1]
\draw[mor] (0,0) rectangle (1,0.5)node[black,pos=0.5]{$x_k$};
\end{tikzpicture} \cdot \alpha_i \;
- \delta_k \LIF(x_k,t_{k+1},x_{k-1})^{-1} \begin{tikzpicture}[anchorbase,scale=1]
\draw[mor] (0,0) rectangle (1.1,0.5) node[black,pos=0.5]{$x_k$};
\draw[mor] (0,1) rectangle (1.1,1.5) node[black,pos=0.5]{$x_k$};
\draw[marked=1] (1,0.5) to (1,0.65);
\draw[usual] (0.25,0.5) to (0.25,1);
\draw[usual] (0.5,0.5) to (0.5,1);
\draw[marked=1] (1,1) to (1,0.85);
\end{tikzpicture}\; . \end{equation}

\end{Theorem}

We have used a trick of notation above: the second diagram on the RHS of  \eqref{exk1} does not make sense when $x_{k-1} t_{k+1} > x_{k-1}$, but we accept this since we have set its coefficient to be zero in this case.

\begin{proof} This is an immediate application of \autoref{Thm:inductiveclasp} and \autoref{Cor:putsomedotsonit}. \end{proof}

\begin{Theorem} \label{thm:linearpartialtracerecursion} Continue the hypotheses and notation of \autoref{thm:linearidempotents}. For all $k \ge 1$ we have
\begin{equation} \label{Arecursion} A_{k}(f) = B_{k,1}(\partial_{t_k}(f)). \end{equation}
For all $k \ge 0$ and $a \ge 1$ we have
\begin{equation} \label{Crecursion} C_{k+1,a}(f) = B_{k+1,a+1}(\partial_{t_{k+1}}(f)) + \sum_{b=0}^{a-1} (-1)^b \binom{a}{b+1} B_{k+1,a-b}(\alpha_{t_{k+1}}^{b} \cdot s_{t_{k+1}}(f)). \end{equation}
For all $k \ge 1$ and $a \ge 1$ we have
\begin{equation} \label{Brecursion} B_{k+1,a}(f) = A_{k}(\alpha_{t_{k+1}}^a f) + \delta_{k} \sum_{b=1}^a \binom{a}{b} (-1)^b \LIF(x_{k},t_{k+1},x_{k-1})^{-b} C_{k,b}(\alpha_{t_{k+1}}^{a-b} f). \end{equation}
Note that $A_{k}(\alpha_{t_{k+1}}^a f) = C_{k,0}(\alpha_{t_{k+1}}^a f)$ agrees with the $b=0$ term in the sum, but we separate it because the $b>0$ terms only appear when $B_{x_{k-1}}$ is a summand of $B_{x_{k}} B_{t_{k+1}}$. Meanwhile for $a \ge 1$ we have
\begin{equation} \label{Bbasecase} B_{1,a}(f) \equiv A_0(\alpha_{t_1}^{a} f) + \catstuff{I}_{< d}, \end{equation}
and no lower terms appear when $B_d B_{t_1} = B_{dt_1}$.
\end{Theorem}

These recursions replace $A_k(-)$ with a linear combination of $B_{k,-}(-)$, then with $C_{k-1,-}(-)$ with $B_{k-1,-}(-)$ etcetera, on down to $B_{1,-}(-)$ being replaced with $A_0(-)$ modulo lower terms (these lower terms can be thought of as $C_{0,a}(-)$ for $a > 0$). Thus, one can recursively write $A_k(f)$ (modulo lower terms) as $A_0(g)$ for some polynomial $g$ determined from $f$.

\begin{proof} The first equality \eqref{Arecursion} is precisely \autoref{lemma:AtoB}.

Now fix $k \ge 1$. The essential computations take place in $\End(B_{x_k})$. Write $\ptr_k$ for the iterated partial trace operation which takes an endomorphism of $B_{x_k}$ to an endomorphism of $B_d = B_{x_0}$. Instead of working with $A, B, C$ we work with replacements in $\End(B_{x_k})$ whose partial trace gives the original quantity. For example, we set $A'_k = \id_{x_k}$ since $\ptr_k(\id_{x_k} \cdot f) = A_k(f)$.

The equation \eqref{exk1dots} also takes place in $\End(B_{x_k})$. Let us give names to various elements of this endomorphism ring. The left-hand side of \eqref{exk1dots} will be called the \emph{dotted clasp}, and will be denoted $B'_1 = B'_{k+1,1}$. Let $B'_{k+1,a} = (B'_1)^a \in \End(B_{x_{k}})$ denote the $a$-th power of the dotted clasp, so that $\ptr_{k}(B'_{k+1,a} \cdot f) = B_{k+1,a}(f)$. The right-hand side of  \eqref{exk1dots} is a linear combination of one or two diagrams, the \emph{barbell diagram} and (if $\delta_k \ne 0$) the \emph{broken diagram}, in that order. We denote the broken diagram by $C'_1 = C'_{k,1}$, and let $C'_{k,a} = (C'_1)^a \in \End(B_{x_k})$ denote the $a$-th power of the broken diagram, so that $\ptr_k(C'_{k,a} \cdot f) = C_{k,a}(f)$.

Now we will prove \eqref{Brecursion}, and the method is simple: apply \eqref{exk1dots} to resolve each of the $a$ idempotents in $B'_{k+1,a}$. If $\delta_k = 0$ then only barbell diagrams appear, so
\begin{equation} B'_{k+1,a} = \id_{x_k} \cdot \alpha_{t_{k+1}}^a = A'_{k} \cdot \alpha_{t_{k+1}}^a. \end{equation}
Applying $\ptr_k$ we get the desired result. This also proves \eqref{Bbasecase} when $B_d B_{t_1} = B_{dt_1}$.

Now suppose $\delta_k = 1$. The barbell diagram is in the center of $\End(B_{x_{k}})$, as right multiplication by any polynomial must be, and Thus, it commutes with the broken diagram. This commutativity permits us to use the binomial theorem when computing $(B'_1)^a$ via \eqref{exk1dots}. Let $L$ denote $\LIF(x_k,t_{k+1},x_{k-1})$. Thus
\begin{equation} B'_{k+1,a} = (\alpha_{t_{k+1}} - L^{-1} C'_1)^a = \sum_{b=0}^a \binom{a}{b} (-1)^b \alpha_{t_{k+1}}^{a-b} L^{-b} (C'_1)^b = \sum_{b=0}^a \binom{a}{b} (-1)^b \alpha_{t_{k+1}}^{a-b} L^{-b} C'_{k,b}. \end{equation}
Applying $\ptr_k$ yields \eqref{Brecursion}. Note that the case $b=0$ accounts for the term $A_{k}(\alpha_{t_{k+1}}^{a} f)$.

We illustrate this when $a=2$, after applying the partial trace.

\begin{gather}
B_{k+1,2}(f) = \; \begin{tikzpicture}[anchorbase,scale=1]
\draw[mor] (-0.4,0.5) rectangle (0,1) node[black,pos=0.5]{\raisebox{-0.06cm}{${d}$}};
\draw[mor] (-0.4,0) rectangle (1.6,0.5) node[black,pos=0.5]{\raisebox{-0.06cm}{${d}~ t_1\ldots ~t_{k+1}$}};
\draw[mor] (-0.4,-1) rectangle (1.6,-0.5) node[black,pos=0.5]{\raisebox{-0.06cm}{${d}~ t_1\ldots ~t_{k+1}$}};
\draw[mor] (-0.4,-1.5) rectangle (0,-1) node[black,pos=0.5]{\raisebox{-0.06cm}{${d}$}};	
\draw[usual] (1.3,-1) to[out=270,in=180] (1.7,-1.4) to[out=0,in=270] (2.1,-1) to (2.1,0.5)to[out=90,in=0] (1.7,0.9) to[out=180,in=90] (1.3,0.5);
\node[scale=0.8] at (1.8,-0.25) {$f$};	
\draw[usual] (0.5,-1) to[out=270,in=180] (1.7,-1.8) to[out=0,in=270] (2.7,-1) to (2.7,0.5)to[out=90,in=0] (1.7,1.3) to[out=180,in=90] (0.5,0.5);
\draw[usual,marked=1] (1.5,0.5) to (1.5,0.6);
\draw[usual,marked=1] (1.5,-1) to (1.5,-1.1);
\node[scale=0.8] at (1,-1.25) {$...$};
\node[scale=0.8] at (1,0.65) {$...$};
\node[scale=0.8] at (2.4,-0.25) {$...$};
\draw[usual] (0.4,0) to (0.4,-0.5);	
\draw[usual] (1.3,0) to (1.3,-0.5);
\draw[usual] (-0.1,0) to (-0.1,-0.5);
\draw[usual,marked=1] (1.5,0) to (1.5,-0.1);
\draw[usual,marked=1] (1.5,-0.5) to (1.5,-0.4);
\end{tikzpicture}	\quad = \quad
\\ \nonumber
\begin{tikzpicture}[anchorbase,scale=1]
\draw[mor] (0,0.5) rectangle (0.4,1) node[black,pos=0.5]{\raisebox{-0.06cm}{${d}$}};
\draw[mor] (0,0) rectangle (1.4,0.5) node[black,pos=0.5]{\raisebox{-0.06cm}{${d}\ldots t_k$}};
\draw[mor] (0,-1) rectangle (1.4,-0.5) node[black,pos=0.5]{\raisebox{-0.06cm}{${d}\ldots t_k$}};
\draw[mor] (0,-1.5) rectangle (0.4,-1) node[black,pos=0.5]{\raisebox{-0.06cm}{${d}$}};	
\draw[usual] (1.3,-1) to[out=270,in=180] (1.7,-1.4) to[out=0,in=270] (2.1,-1) to (2.1,0.5)to[out=90,in=0] (1.7,0.9) to[out=180,in=90] (1.3,0.5);
\node[scale=0.8] at (1.8,-0.25) {$f$};	
\draw[usual] (0.5,-1) to[out=270,in=180] (1.7,-1.8) to[out=0,in=270] (2.7,-1) to (2.7,0.5)to[out=90,in=0] (1.7,1.3) to[out=180,in=90] (0.5,0.5);
\draw[usual,marked=1] (1.5,0) to (1.5,0.6);
\draw[usual,marked=1] (1.5,-0.5) to (1.5,-1.1);
\node[scale=0.8] at (1,-1.25) {$...$};
\node[scale=0.8] at (1,0.65) {$...$};
\node[scale=0.8] at (2.4,-0.25) {$...$};
\draw[usual] (0.2,0) to (0.2,-0.5);	
\draw[usual] (1.3,0) to (1.3,-0.5);
\draw[usual] (0.5,0) to (0.5,-0.5);
\draw[usual,marked=1] (1.5,0) to (1.5,-0.1);
\draw[usual,marked=1] (1.5,-0.5) to (1.5,-0.4);
\end{tikzpicture}
-L^{-1}
\begin{tikzpicture}[anchorbase,scale=1]
\draw[mor] (0,0.5) rectangle (0.4,1) node[black,pos=0.5]{\raisebox{-0.06cm}{${d}$}};
\draw[mor] (0,0) rectangle (1.4,0.5) node[black,pos=0.5]{\raisebox{-0.06cm}{${d}\ldots t_k$}};
\draw[mor] (0,-1) rectangle (1.4,-0.5) node[black,pos=0.5]{\raisebox{-0.06cm}{${d}\ldots t_k$}};
\draw[mor] (0,-2) rectangle (1.4,-1.5) node[black,pos=0.5]{\raisebox{-0.06cm}{${d}\ldots t_k$}};
\draw[mor] (0,-2.5) rectangle (0.4,-2) node[black,pos=0.5]{\raisebox{-0.06cm}{${d}$}};	
\draw[usual] (1.3,-2) to[out=270,in=180] (1.7,-2.4) to[out=0,in=270] (2.1,-2) to (2.1,0.5)to[out=90,in=0] (1.7,0.9) to[out=180,in=90] (1.3,0.5);
\node[scale=0.8] at (1.8,-0.75) {$f$};	
\draw[usual] (0.5,-2) to[out=270,in=180] (1.7,-2.8) to[out=0,in=270] (2.7,-2) to (2.7,0.5)to[out=90,in=0] (1.7,1.3) to[out=180,in=90] (0.5,0.5);
\draw[usual,marked=1] (1.5,-2) to (1.5,-2.1);
\draw[usual,marked=1] (1.5,-2) to (1.5,-1.4);
\draw[usual] (0.2,-1) to (0.2,-1.5);
\draw[usual] (0.5,-1) to (0.5,-1.5);
\draw[usual] (1.3,-1) to (1.3,-1.5);
\node[scale=0.8] at (1,0.65) {$...$};
\node[scale=0.8] at (2.4,-0.75) {$...$};
\draw[usual] (0.2,0) to (0.2,-0.5);	
\draw[usual] (0.5,0) to (0.5,-0.5);
\draw[usual,marked=1] (1.3,0) to (1.3,-0.1);
\draw[usual,marked=1] (1.3,-0.5) to (1.3,-0.4);
\end{tikzpicture}
- L^{-1}
\begin{tikzpicture}[anchorbase,scale=1]
\begin{scope}[yscale=-1]
\draw[mor] (0,0.5) rectangle (0.4,1) node[black,pos=0.5]{\raisebox{-0.06cm}{${d}$}};
\draw[mor] (0,0) rectangle (1.4,0.5) node[black,pos=0.5]{\raisebox{-0.06cm}{${d}\ldots t_k$}};
\draw[mor] (0,-1) rectangle (1.4,-0.5) node[black,pos=0.5]{\raisebox{-0.06cm}{${d}\ldots t_k$}};
\draw[mor] (0,-2) rectangle (1.4,-1.5) node[black,pos=0.5]{\raisebox{-0.06cm}{${d}\ldots t_k$}};
\draw[mor] (0,-2.5) rectangle (0.4,-2) node[black,pos=0.5]{\raisebox{-0.06cm}{${d}$}};	
\draw[usual] (1.3,-2) to[out=270,in=180] (1.7,-2.4) to[out=0,in=270] (2.1,-2) to (2.1,0.5)to[out=90,in=0] (1.7,0.9) to[out=180,in=90] (1.3,0.5);
\node[scale=0.8] at (1.8,-0.75) {$f$};	
\draw[usual] (0.5,-2) to[out=270,in=180] (1.7,-2.8) to[out=0,in=270] (2.7,-2) to (2.7,0.5)to[out=90,in=0] (1.7,1.3) to[out=180,in=90] (0.5,0.5);
\draw[usual,marked=1] (1.5,-2) to (1.5,-2.1);
\draw[usual,marked=1] (1.5,-2) to (1.5,-1.4);
\draw[usual] (0.2,-1) to (0.2,-1.5);
\draw[usual] (0.5,-1) to (0.5,-1.5);
\draw[usual] (1.3,-1) to (1.3,-1.5);
\node[scale=0.8] at (1,0.65) {$...$};
\node[scale=0.8] at (2.4,-0.75) {$...$};
\draw[usual] (0.2,0) to (0.2,-0.5);	
\draw[usual] (0.5,0) to (0.5,-0.5);
\draw[usual,marked=1] (1.3,0) to (1.3,-0.1);
\draw[usual,marked=1] (1.3,-0.5) to (1.3,-0.4);
\end{scope}
\end{tikzpicture}
+
L^{-2}
\begin{tikzpicture}[anchorbase,scale=1]
\draw[mor] (0,0.5) rectangle (0.4,1) node[black,pos=0.5]{\raisebox{-0.06cm}{${d}$}};
\draw[mor] (0,0) rectangle (1.4,0.5) node[black,pos=0.5]{\raisebox{-0.06cm}{${d}\ldots t_k$}};
\draw[mor] (0,-1) rectangle (1.4,-0.5) node[black,pos=0.5]{\raisebox{-0.06cm}{${d}\ldots t_k$}};
\draw[mor] (0,-2) rectangle (1.4,-1.5) node[black,pos=0.5]{\raisebox{-0.06cm}{${d}\ldots t_k$}};
\draw[mor] (0,-3) rectangle (1.4,-2.5) node[black,pos=0.5]{\raisebox{-0.06cm}{${d}\ldots t_k$}};
\draw[mor] (0,-3.5) rectangle (0.4,-3) node[black,pos=0.5]{\raisebox{-0.06cm}{${d}$}};	
\draw[usual] (1.3,-3) to[out=270,in=180] (1.7,-3.4) to[out=0,in=270] (2.1,-3) to (2.1,0.5)to[out=90,in=0] (1.7,0.9) to[out=180,in=90] (1.3,0.5);
\node[scale=0.8] at (1.8,-1.25) {$f$};	
\draw[usual] (0.5,-3) to[out=270,in=180] (1.7,-3.8) to[out=0,in=270] (2.7,-3) to (2.7,0.5)to[out=90,in=0] (1.7,1.3) to[out=180,in=90] (0.5,0.5);
\draw[usual] (0.2,-1) to (0.2,-1.5);
\draw[usual] (0.5,-1) to (0.5,-1.5);
\draw[usual] (1.3,-1) to (1.3,-1.5);
\node[scale=0.8] at (1,0.65) {$...$};
\node[scale=0.8] at (2.4,-1.25) {$...$};
\draw[usual] (0.2,0) to (0.2,-0.5);	
\draw[usual] (0.5,0) to (0.5,-0.5);
\draw[usual,marked=1] (1.3,0) to (1.3,-0.1);
\draw[usual,marked=1] (1.3,-0.5) to (1.3,-0.4);
\draw[usual,marked=1] (1.3,-2) to (1.3,-2.1);
\draw[usual,marked=1] (1.3,-2.5) to (1.3,-2.4);
\end{tikzpicture}
\\ \nonumber =
C_{k,0}(f\alpha_i^2)
- 2 L^{-1} C_{k-1,1}(f\alpha_i)
+
L^{-2} C_{k_1,2}(f).
\end{gather}
Now we prove \eqref{Crecursion}, focusing on 
\[ \ptr_{x_k,t_{k+1}} (C'_{k+1,a} \cdot f) \in \End(B_{x_k}). \]
The first step is to use polynomial forcing: 

\begin{equation}
\begin{tikzpicture}[anchorbase,scale=1]
\draw[mor] (-0.2,0) rectangle (1.6,0.5) node[black,pos=0.5]{\raisebox{-0.06cm}{${d} \ldots t_{k} t_{k+1}$}};
\draw[mor] (-0.2,-1) rectangle (1.6,-0.5) node[black,pos=0.5]{\raisebox{-0.06cm}{${d} \ldots t_{k}t_{k+1}$}};
\draw[mor] (-0.2,-2) rectangle (1.6,-1.5) node[black,pos=0.5]{\raisebox{-0.06cm}{${d} \ldots t_{k} t_{k+1}$}};
\draw[usual] (1.3,-2) to[out=270,in=180] (1.8,-2.5) to[out=0,in=270] (2.3,-2) to (2.3,0.5)to[out=90,in=0] (1.8,1) to[out=180,in=90] (1.3,0.5);
\node at (1.8,-0.75) {$f$};	
\node[scale=0.8] at (0.7,-1.25) {$\vdots$};
\draw[usual] (0.2,0) to (0.2,-0.5);	
\draw[usual] (0.5,0) to (0.5,-0.5);
\draw[usual,marked=1] (1.3,0) to (1.3,-0.1);
\draw[usual,marked=1] (1.3,-0.5) to (1.3,-0.4);
\end{tikzpicture}	
=	
\begin{tikzpicture}[anchorbase,scale=1]
\draw[mor] (-0.2,0) rectangle (1.6,0.5) node[black,pos=0.5]{\raisebox{-0.06cm}{${d} \ldots t_{k} t_{k+1}$}};
\draw[mor] (-0.2,-1) rectangle (1.6,-0.5) node[black,pos=0.5]{\raisebox{-0.06cm}{${d} \ldots t_{k}t_{k+1}$}};
\draw[mor] (-0.2,-2) rectangle (1.6,-1.5) node[black,pos=0.5]{\raisebox{-0.06cm}{${d} \ldots t_{k} t_{k+1}$}};
\node at (2.4,-0.75) {${t_{k+1}}f$};	
\draw[usual] (1.5,-2) to[out=270,in=180] (1.7,-2.3) to[out=0,in=270] (1.9,-2) to (1.9,0.5)to[out=90,in=0] (1.7,0.8) to[out=180,in=90] (1.5,0.5);
\node[scale=0.8] at (0.8,-1.25) {$\vdots$};
\draw[usual] (0.2,0) to (0.2,-0.5);	
\draw[usual] (0.5,0) to (0.5,-0.5);
\draw[usual,marked=1] (1.5,0) to (1.5,-0.1);
\draw[usual,marked=1] (1.5,-0.5) to (1.5,-0.4);
\end{tikzpicture}	
+
\begin{tikzpicture}[anchorbase,scale=1]
\draw[mor] (-0.2,0) rectangle (1.6,0.5) node[black,pos=0.5]{\raisebox{-0.06cm}{${d} \ldots t_{k} t_{k+1}$}};
\draw[mor] (-0.2,-1) rectangle (1.6,-0.5) node[black,pos=0.5]{\raisebox{-0.06cm}{${d} \ldots t_{k}t_{k+1}$}};
\draw[mor] (-0.2,-2) rectangle (1.6,-1.5) node[black,pos=0.5]{\raisebox{-0.06cm}{${d} \ldots t_{k} t_{k+1}$}};
\node at (2.3,-0.75) {$\partial_{t_{k+1}}(f)$};	
\draw[usual,marked=1] (1.5,0.5) to (1.5,0.6);
\draw[usual,marked=1] (1.5,-2) to (1.5,-2.1);
\node[scale=0.8] at (0.8,-1.25) {$\vdots$};
\draw[usual] (0.2,0) to (0.2,-0.5);	
\draw[usual] (1.3,0) to (1.3,-0.5);
\draw[usual] (0.5,0) to (0.5,-0.5);
\draw[usual,marked=1] (1.5,0) to (1.5,-0.1);
\draw[usual,marked=1] (1.5,-0.5) to (1.5,-0.4);
\end{tikzpicture} \; .
\end{equation}
The second term on the right-hand side is $(B'_1)^{a+1} \cdot \partial_{t_{k+1}}(f)$. For the first term, we use the following manipulation, like in \autoref{lemma:AtoB}.
\begin{equation} \label{lalalalala}
\begin{tikzpicture}[anchorbase,scale=1]
\draw[mor] (-0.2,0) rectangle (1.6,0.5) node[black,pos=0.5]{\raisebox{-0.06cm}{${d} \ldots t_{k} t_{k+1}$}};
\draw[mor] (-0.2,-1) rectangle (1.6,-0.5) node[black,pos=0.5]{\raisebox{-0.06cm}{${d} \ldots t_{k}t_{k+1}$}};
\draw[mor] (-0.2,-2) rectangle (1.6,-1.5) node[black,pos=0.5]{\raisebox{-0.06cm}{${d} \ldots t_{k} t_{k+1}$}};
\draw[usual] (1.5,-2) to[out=270,in=180] (1.8,-2.4) to[out=0,in=270] (2.1,-2) to (2.1,0.5)to[out=90,in=0] (1.8,0.9) to[out=180,in=90] (1.5,0.5);
\node[scale=0.8] at (0.8,-1.25) {$\vdots$};
\draw[usual] (0.2,0) to (0.2,-0.5);	
\draw[usual] (0.5,0) to (0.5,-0.5);
\draw[usual,marked=1] (1.5,0) to (1.5,-0.1);
\draw[usual,marked=1] (1.5,-0.5) to (1.5,-0.4);
\draw[usual,marked=1] (1.5,-1) to (1.5,-1.1);
\draw[usual,marked=1] (1.5,-1.5) to (1.5,-1.4);
\end{tikzpicture}	
=
\begin{tikzpicture}[anchorbase,scale=1]
\draw[mor] (-0.2,0) rectangle (1.6,0.5) node[black,pos=0.5]{\raisebox{-0.06cm}{${d} \ldots t_{k} t_{k+1}$}};
\draw[mor] (-0.2,-1) rectangle (1.6,-0.5) node[black,pos=0.5]{\raisebox{-0.06cm}{${d} \ldots t_{k}t_{k+1}$}};
\draw[mor] (-0.2,-2) rectangle (1.6,-1.5) node[black,pos=0.5]{\raisebox{-0.06cm}{${d} \ldots t_{k} t_{k+1}$}};
\draw[usual] (1.6,-1.75) to[out=0,in=270] (2,-1.5) to (2,0)to[out=90,in=0] (1.6,0.25);
\node[scale=0.8] at (0.8,-1.25) {$\vdots$};
\draw[usual] (0.2,0) to (0.2,-0.5);	
\draw[usual] (0.5,0) to (0.5,-0.5);
\draw[usual,marked=1] (1.5,0.5) to (1.5,0.9);
\draw[usual,marked=1] (1.5,0) to (1.5,-0.1);
\draw[usual,marked=1] (1.5,-0.5) to (1.5,-0.4);
\draw[usual,marked=1] (1.5,-1) to (1.5,-1.1);
\draw[usual,marked=1] (1.5,-1.5) to (1.5,-1.4);
\draw[usual,marked=1] (1.5,-2) to (1.5,-2.4);
\end{tikzpicture}\;.
\end{equation}
We define a family of diagrams $D'_a = D'_{k+1,a} \in \End(B_{x_k})$ as
\begin{equation} D'_a := \;	\begin{tikzpicture}[anchorbase,scale=1]
\draw[mor] (-0.2,0) rectangle (1.6,0.5) node[black,pos=0.5]{\raisebox{-0.06cm}{${d} \ldots t_{k} t_{k+1}$}};
\draw[mor] (-0.2,-1) rectangle (1.6,-0.5) node[black,pos=0.5]{\raisebox{-0.06cm}{${d} \ldots t_{k}t_{k+1}$}};
\draw[mor] (-0.2,-2) rectangle (1.6,-1.5) node[black,pos=0.5]{\raisebox{-0.06cm}{${d} \ldots t_{k} t_{k+1}$}};
\draw[usual] (1.6,-1.75) to[out=0,in=270] (2,-1.5) to (2,0)to[out=90,in=0] (1.6,0.25);
\node[scale=0.8] at (0.8,-1.25) {$\vdots$};
\draw[usual] (0.2,0) to (0.2,-0.5);	
\draw[usual] (0.5,0) to (0.5,-0.5);
\draw[usual,marked=1] (1.5,0.5) to (1.5,0.9);
\draw[usual,marked=1] (1.5,0) to (1.5,-0.1);
\draw[usual,marked=1] (1.5,-0.5) to (1.5,-0.4);
\draw[usual,marked=1] (1.5,-1) to (1.5,-1.1);
\draw[usual,marked=1] (1.5,-1.5) to (1.5,-1.4);
\draw[usual,marked=1] (1.5,-2) to (1.5,-2.4);
\end{tikzpicture}\;.
\end{equation}
So we have
\begin{equation} \label{CtoD} \ptr_{x_k,t_{k+1}} (C'_{k+1,a} \cdot f) = (B'_1)^{a+1} \cdot \partial_{t_{k+1}}(f) + D'_{a} \cdot s_{t_{k+1}}(f). \end{equation}

To evaluate $D'_{a}$, we use the broken trivalent relation \eqref{eq:brokentrivalent}, which we recall here.
\begin{equation} \begin{tikzpicture}[anchorbase,scale=1]
\draw[usual] (0.6,0) to (0.3,1);
\draw[usual,marked=1] (0,0) to (0.15,0.3);
\end{tikzpicture}
=
\begin{tikzpicture}[anchorbase,scale=1]
\draw[usual] (0,0) to (0.3,1);
\draw[usual,marked=1] (0.6,0) to (0.45,0.3);
\end{tikzpicture}
+
\begin{tikzpicture}[anchorbase,scale=1]
\draw[usual] (0.6,0) to[out=90,in=0] (0.3,0.3) to[out=180,in=90] (0,0);
\draw[usual,marked=1] (0.3,1) to (0.3,0.7);
\end{tikzpicture}
-
\begin{tikzpicture}[anchorbase,scale=1]
\draw[usual] (0,0) to (0.3,0.3) to (0.3,1);
\node[scale=0.8] at (0.8,0.5) {$\alpha_i$};
\draw[usual] (0.6,0) to (0.3,0.3);
\end{tikzpicture}. \end{equation}
We apply this relation to $D'_a$ as follows.
\begin{equation} \label{Dprimereduction}	\begin{tikzpicture}[anchorbase,scale=1]
\draw[mor] (-0.2,0) rectangle (1.6,0.5) node[black,pos=0.5]{\raisebox{-0.06cm}{${d} \ldots t_{k} t_{k+1}$}};
\draw[mor] (-0.2,-1) rectangle (1.6,-0.5) node[black,pos=0.5]{\raisebox{-0.06cm}{${d} \ldots t_{k}t_{k+1}$}};
\draw[mor] (-0.2,-2) rectangle (1.6,-1.5) node[black,pos=0.5]{\raisebox{-0.06cm}{${d} \ldots t_{k} t_{k+1}$}};
\draw[usual] (1.6,-1.75) to[out=0,in=270] (2,-1.5) to (2,0)to[out=90,in=0] (1.6,0.25);
\node[scale=0.8] at (0.8,-1.25) {$\vdots$};
\draw[usual] (0.2,0) to (0.2,-0.5);	
\draw[usual] (0.5,0) to (0.5,-0.5);
\draw[usual,marked=1] (1.5,0.5) to (1.5,0.9);
\draw[usual,marked=1] (1.5,0) to (1.5,-0.1);
\draw[usual,marked=1] (1.5,-0.5) to (1.5,-0.4);
\draw[usual,marked=1] (1.5,-1) to (1.5,-1.1);
\draw[usual,marked=1] (1.5,-1.5) to (1.5,-1.4);
\draw[usual,marked=1] (1.5,-2) to (1.5,-2.4);
\end{tikzpicture}	
=
\begin{tikzpicture}[anchorbase,scale=1]
\draw[mor] (-0.2,0) rectangle (1.6,0.5) node[black,pos=0.5]{\raisebox{-0.06cm}{${d} \ldots t_{k} t_{k+1}$}};
\draw[mor] (-0.2,-1) rectangle (1.6,-0.5) node[black,pos=0.5]{\raisebox{-0.06cm}{${d} \ldots t_{k}t_{k+1}$}};
\draw[mor] (-0.2,-2) rectangle (1.6,-1.5) node[black,pos=0.5]{\raisebox{-0.06cm}{${d} \ldots t_{k} t_{k+1}$}};
\draw[usual,marked=1] (1.6,-1.75) to[out=0,in=270] (2,-1.5) to (2,-0.5);
\node[scale=0.8] at (0.8,-1.25) {$\vdots$};
\draw[usual] (0.2,0) to (0.2,-0.5);	
\draw[usual] (0.5,0) to (0.5,-0.5);
\draw[usual,marked=1] (1.5,0.5) to (1.5,0.9);
\draw[usual,marked=1] (1.5,0) to (1.5,-0.1);
\draw[usual] (1.5,-0.5) to (1.5,-0.4) to[out=90, in=270] (1.8,-0.1) to [out=90, in=0] (1.6,0.25);
\draw[usual,marked=1] (1.5,-1) to (1.5,-1.1);
\draw[usual,marked=1] (1.5,-1.5) to (1.5,-1.4);
\draw[usual,marked=1] (1.5,-2) to (1.5,-2.4);
\end{tikzpicture}
+
\begin{tikzpicture}[anchorbase,scale=1]
\draw[mor] (-0.2,0) rectangle (1.6,0.5) node[black,pos=0.5]{\raisebox{-0.06cm}{${d} \ldots t_{k} t_{k+1}$}};
\draw[mor] (-0.2,-1) rectangle (1.6,-0.5) node[black,pos=0.5]{\raisebox{-0.06cm}{${d} \ldots t_{k}t_{k+1}$}};
\draw[mor] (-0.2,-2) rectangle (1.6,-1.5) node[black,pos=0.5]{\raisebox{-0.06cm}{${d} \ldots t_{k} t_{k+1}$}};
\draw[usual] (1.6,-1.75) to[out=0,in=270] (2,-1.5) to (2,-0.5) to[out=90, in=0] (1.8,-0.2) to[out=180, in=90] (1.5,-0.5);
\node[scale=0.8] at (0.8,-1.25) {$\vdots$};
\draw[usual] (0.2,0) to (0.2,-0.5);	
\draw[usual] (0.5,0) to (0.5,-0.5);
\draw[usual,marked=1] (1.5,0.5) to (1.5,0.9);
\draw[usual,marked=1] (1.5,0) to (1.5,-0.1);
\draw[usual,marked=1] (1.6,.25) to (1.8,0.25);
\draw[usual,marked=1] (1.5,-1) to (1.5,-1.1);
\draw[usual,marked=1] (1.5,-1.5) to (1.5,-1.4);
\draw[usual,marked=1] (1.5,-2) to (1.5,-2.4);
\end{tikzpicture}
-
\begin{tikzpicture}[anchorbase,scale=1]
\draw[mor] (-0.2,0) rectangle (1.6,0.5) node[black,pos=0.5]{\raisebox{-0.06cm}{${d} \ldots t_{k} t_{k+1}$}};
\draw[mor] (-0.2,-1) rectangle (1.6,-0.5) node[black,pos=0.5]{\raisebox{-0.06cm}{${d} \ldots t_{k}t_{k+1}$}};
\draw[mor] (-0.2,-2) rectangle (1.6,-1.5) node[black,pos=0.5]{\raisebox{-0.06cm}{${d} \ldots t_{k} t_{k+1}$}};
\draw[usual] (1.6,-1.75) to[out=0,in=270] (2,-1.5) to (2,-0.5);
\node[scale=0.8] at (0.8,-1.25) {$\vdots$};
\draw[usual] (0.2,0) to (0.2,-0.5);	
\draw[usual] (0.5,0) to (0.5,-0.5);
\draw[usual,marked=1] (1.5,0.5) to (1.5,0.9);
\draw[usual,marked=1] (1.5,0) to (1.5,-0.1);
\draw[usual] (1.6,-1.75) to[out=0,in=270] (2,-1.5) to (2,-0.5) to[out=90, in=0] (1.8,-0.2) to[out=180, in=90] (1.5,-0.5);
\draw[usual] (1.8,-0.2) to (1.8,-0.1) to [out=90, in=0] (1.6,0.25);
\draw[usual,marked=1] (1.5,-1) to (1.5,-1.1);
\draw[usual,marked=1] (1.5,-1.5) to (1.5,-1.4);
\draw[usual,marked=1] (1.5,-2) to (1.5,-2.4);
\node at (2.4,-0.75) {$\alpha_i$};
\end{tikzpicture}\;. \end{equation}
Each of these diagrams can be simplified using \eqref{triassoc} and \eqref{trislide} and \eqref{triunit}, like in the following example:
\begin{equation} 	\begin{tikzpicture}[anchorbase,scale=1]
\draw[mor] (0,0) rectangle (1,0.5) node[black,pos=0.5]{\raisebox{-0.06cm}{}};
\draw[mor] (0,-1) rectangle (1,-0.5) node[black,pos=0.5]{\raisebox{-0.06cm}{}};
\draw[usual](1.4,-1.1) to (1.4,-0.5) to[out=90, in=0] (1.2,-0.2) to[out=180, in=90] (0.9,-0.5);
\draw[usual] (1.2,-0.2) to (1.2,-0.1) to [out=90, in=0] (1,0.25);
\draw[usual, marked=1] (0.9,0) to (0.9,-0.1);
\draw[usual, marked=1] (0.9,-1) to (0.9,-1.1);
\end{tikzpicture}
=
\begin{tikzpicture}[anchorbase,scale=1]
\draw[mor] (0,0) rectangle (1,0.5) node[black,pos=0.5]{\raisebox{-0.06cm}{}};
\draw[mor] (0,-1) rectangle (1,-0.5) node[black,pos=0.5]{\raisebox{-0.06cm}{}};
\draw[usual](1.4,-1.1) to (1.4,-0.5) to[out=90, in=0] (1.2,-0.2) to[out=180, in=90] (0.9,-0.5);
\draw[usual] (1.2,-0.2) to[out=180, in=270] (0.9,0);
\draw[usual, marked=1] (0.9,-1) to (0.9,-1.1);
\end{tikzpicture}
=
\begin{tikzpicture}[anchorbase,scale=1]
\draw[mor] (0,-1) rectangle (1,-0.5) node[black,pos=0.5]{\raisebox{-0.06cm}{}};
\draw[usual](1.4,-1.2) to[out=90, in=0] (1,-0.75);
\draw[usual, marked=1] (0.9,-1) to (0.9,-1.1);
\end{tikzpicture}
=
\begin{tikzpicture}[anchorbase,scale=1]
\draw[mor] (0,-1) rectangle (1,-0.5) node[black,pos=0.5]{\raisebox{-0.06cm}{}};
\draw[usual] (0.9,-1) to (0.9,-1.4);
\end{tikzpicture}\;. \end{equation}
Doing similar reductions, one can verify that the first diagram on the right-hand side of \eqref{Dprimereduction} is $(B'_1)^a$. The second diagram is a $B'_1 \circ D'_{a-1}$. The third diagram is $D'_{a-1} \cdot \alpha_{t_{k+1}}$. Thus
\begin{equation}\label{Dprimerecursion} D'_{a} = (B'_{1})^a + B'_{1} \circ D'_{a-1} - D'_{a-1} \cdot \alpha_{t_{k+1}}. \end{equation}

Now we solve this recursion for $D'$ and prove by induction on $a$ that 
\begin{equation} \label{Dprimesolved} D'_{a} = \sum_{b=0}^{a-1} (-1)^b \binom{a}{b+1} (B'_1)^{a-b} \cdot \alpha_{t_{k+1}}^{b} \end{equation}
Note that the base cases $D'_0 = 0$ and $D'_1 = B'_1$ hold. If the statement holds for $D'_{a-1}$ then \eqref{Dprimerecursion} says
\begin{equation} D'_a = (B'_{1})^a  + \sum_{b=0}^{a-2} (-1)^b \binom{a-1}{b+1} (B'_1)^{a-b} \cdot \alpha_{t_{k+1}}^{b} - \sum_{b=0}^{a-2} (-1)^b \binom{a-1}{b+1} (B'_1)^{a-1-b} \cdot \alpha_{t_{k+1}}^{b+1}. \end{equation}
The coefficient of $(B'_1)^a$ is $1 + \binom{a-1}{1} = \binom{a}{1}$. The coefficient of $(B'_1)^1 \cdot \alpha_{t_{k+1}}^{a-1}$ is $-(-1)^{a-2} \binom{a-1}{a-1} = (-1)^{a-1} \binom{a}{a}$. For all $0 < b < a-1$, a little reindexing shows that the coefficient of $(B'_1)^{a-b} \cdot \alpha_{t_{k+1}}^{b}$ is
\begin{equation} (-1)^b \binom{a-1}{b+1} + (-1)^b \binom{a-1}{b} = (-1)^b \binom{a}{b+1}, \end{equation}
using the Pascal triangle rule for adding binomials. Thus, we have matched all coefficients with \eqref{Dprimesolved}.

Now we plug \eqref{Dprimesolved} back into \eqref{CtoD} and apply $\ptr_{k}$; the result is \eqref{Crecursion}.
\end{proof}

There is actually a simpler recursion which bypasses the morphisms $B_{k,a}$ entirely. The following result supercedes \autoref{thm:linearpartialtracerecursion} and is easier to use.

\begin{Theorem} \label{thm:linearpartialtracerecursionbetter} Continue the hypotheses and notation of \autoref{thm:linearidempotents}. Recall that $C_{k,0} = A_k$. Fix some $k \ge 1$. Let $L = \LIF(x_k,t_{k+1},x_{k-1})$ and $\alpha = \alpha_{t_{k+1}}$ and $\partial = \partial_{t_{k+1}}$. Then for all $a > 0$ we have
\begin{equation} \label{betterCrecursion} C_{k+1,a}(f) = C_{k,0}(\alpha^a f) - \delta_k (-L)^{-a} C_{k,a}(t_{k+1} f) + \delta_k (-L)^{a+1} C_{k,a+1}(\partial f) + \delta_k \sum_{c=1}^a (-L)^{-c} \binom{a+1}{c} C_{k,c}(\alpha^{a-c} f). \end{equation}
The formula \eqref{betterCrecursion} holds for $a=0$ too when $\delta_k = 1$. A formula which always works for $a = 0$ is
\begin{equation} \label{betterCrecursiona0} C_{k+1,0}(f) = C_{k,0}(\alpha \partial f) + \delta_k (-L)^{-1} C_{k,1}(\partial f). \end{equation}
\end{Theorem}

As before, we continue to ignore $C_{0,a}$ for $a > 0$, as it lies within the ideal $\catstuff{I}_{< d}$.

\begin{proof} Combining \eqref{Arecursion} with \eqref{Brecursion} (for $a=1$) we get
\begin{equation} A_{k+1}(f) = B_{k+1,1}(\partial f) = A_k(\alpha \partial f) + \delta_k (-L)^{-1} C_{k,1}(\partial f). \end{equation}
Recalling that $C_{k,0}(f) = A_k(f)$, this is exactly \eqref{betterCrecursiona0}. Recall that $\alpha \partial f = f - t_{k+1} f$ by the definition of Demazure operators. Thus, $A_k(\alpha \partial f) = C_{k,0}(f) - C_{k,0}(t_{k+1} f)$, which agrees with the sum of the first two terms in \eqref{betterCrecursion} when $\delta_k = 1$. The third term of \eqref{betterCrecursion} matches the last term of \eqref{betterCrecursiona0}, and the sum over $c$ in \eqref{betterCrecursion} is empty when $a=0$. We conclude that \eqref{betterCrecursion} and \eqref{betterCrecursiona0} agree when $\delta_k = 1$.

Now fix $a \ge 1$. Combining \eqref{Crecursion} and \eqref{Brecursion} we get
\begin{align} C_{k+1,a}(f) &= B_{k+1,a+1}(\partial f) + \sum_{b=0}^{a-1} (-1)^b \binom{a}{b+1} B_{k+1,a-b}(\alpha^{b} t_{k+1} f) \nonumber \\
&= A_k(\alpha^{a+1} \partial f) + \delta_k \sum_{b=1}^{a+1} \binom{a+1}{b} (-L)^{-b} C_{k,b}(\alpha^{a+1-b} \partial f) \nonumber\\
& \scalebox{0.99}{$+ \sum_{b=0}^{a-1} (-1)^b \binom{a}{b+1} A_k(\alpha^{a-b} \alpha^b t_{k+1} f) + \delta_k \sum_{b=0}^{a-1} \sum_{c=1}^{a-b} (-1)^b \binom{a-b}{c} \binom{a}{b+1} (-L)^{-c} C_{k,c}(\alpha^{a-b-c} \alpha^b t_{k+1} f)$} \\ 
&= A_k(\alpha^{a+1} \partial f) + \sum_{b=0}^{a-1} (-1)^b \binom{a}{b+1} A_k(\alpha^a t_{k+1} f) + \nonumber \\
& + \delta_k \sum_{c=1}^{a+1} \binom{a+1}{c} (-L)^{-c} C_{k,c}(\alpha^{a+1-c} \partial f) \\&+ \delta_k \sum_{c=1}^{a} \sum_{b=0}^{a-c} (-1)^b \binom{a-b}{c} \binom{a}{b+1} (-L)^{-c} C_{k,c}(\alpha^{a-c} t_{k+1} f). \nonumber   \end{align}
In the last step we reindexed both sums: in one sum we merely replaced $b$ with $c$; in the other sum we did not change the meaning of $b$ and $c$, but found a different way to sum over the same choices of the pair $(b,c)$.

Now we observe that for $a \ge 1$ and $0 \le c < a$ we have \begin{equation} \label{funkybinomial} \sum_{b=0}^{a-1} (-1)^b \binom{a}{b+1} = 1, \qquad \sum_{b=0}^{a-c} (-1)^b \binom{a-b}{c} \binom{a}{b+1} =
\binom{a+1}{c}. \end{equation} The first formula is a classic fact about Pascal's triangle, coming from the expansion of $0 = (1-1)^{a}$, and it is the $c=0$ case of the second formula. The second property is not a variant of any of the standard binomial formulas we are familiar with, so we prove it in \autoref{Lemma:binomials}. When $c=a$ then \eqref{funkybinomial} fails, as the left-hand side equals $a$ while the right-hand side equals $a+1$.  Below, we view $a$ as $(a+1)-1$ and view the $-1$ as a separate term.

Using these binomial simplications we have
\begin{align} C_{k+1,a}(f) &= C_{k,0}(\alpha^{a+1} \partial f + \alpha^a t_{k+1} f) + \delta_k \sum_{c=1}^{a+1} \binom{a+1}{c} (-L)^{-c} C_{k,c}(\alpha^{a+1-c} \partial f) \\
&+ \delta_k \sum_{c=1}^{a} \binom{a+1}{c} (-L)^{-c} C_{k,c}(\alpha^{a-c} t_{k+1} f) - \delta_k (-L)^{-a} C_{k,a}(t_{k+1} f). \end{align}
Now we use again that $\alpha \partial f = f - t_{k+1} f$, so that $\alpha^{a+1} \partial f + \alpha^a t_{k+1} f = \alpha^a f$. Treating separately the $c=a$ and $c=a+1$ terms, we get
\begin{equation} C_{k+1,a}(f) = C_{k,0}(\alpha^a f) + \sum_{c=1}^{a} \binom{a+1}{c} (-L)^{-c} C_{k,c}(\alpha^{a-c} f) - \delta_k (-L)^{-a} C_{k,a}(t_{k+1} f) + \delta_k (-L)^{-(a+1)} C_{k,a+1}(\partial f). \end{equation}
This is exactly \eqref{betterCrecursion}.
\end{proof}

The following binomial identity appears to be an alternating-sign variant of the Chu--Vandermonde identity, of which there are many in the literature. However, we were unable to locate this specific form, so we include a short proof for completeness.

\begin{Lemma} \label{Lemma:binomials} For $0 \le c < a$ we have
\begin{equation} \sum_{b=0}^{a-c} (-1)^b \binom{a-b}{c} \binom{a}{b+1} = \binom{a+1}{c}. \end{equation}
\end{Lemma}

\begin{proof} First we note that the sum can be taken over all $0 \le b \le a$, because for $a-c < b \le a$ we have $\binom{a-b}{c} = 0$. We may actually sum over all $0 \le b \le a-1$, because $\binom{a}{b+1} = 0$ when $b=a$.
We have noted above that this formula fails when $c=a$, being off by exactly $1$. It also fails for $c=a+1$, being off by exactly $1$ again. Now consider the following generating function.
\begin{equation} f(u) = u^a + u^{a+1} + \sum_{c=0}^{a+1} \sum_{b=0}^{a-1} (-1)^b \binom{a-b}{c} \binom{a}{b+1} u^c. \end{equation}
The coefficient of $u^c$ in $f(u)$ is the left-hand side above, plus one in the cases $c=a$ and $c=a+1$, so we wish to prove that the coefficient is precisely $\binom{a+1}{c}$. This will be true if and only if $f(u) = (1+u)^{a+1}$, by the binomial theorem.

Swapping the order of summation and using the binomial theorem, we have
\begin{equation} \sum_{c=0}^{a+1} \binom{a-b}{c} u^c = (1+u)^{a-b}, \qquad f(u) = u^a + u^{a+1} + \sum_{b=0}^{a-1} (-1)^b \binom{a}{b+1} (1+u)^{a-b}. \end{equation}
We add and subtract a term corresponding to $b=-1$, and then reindex to replace $b+1$ with $b$. We get
\begin{equation} f(u) = u^a + u^{a+1} - (-1) \cdot (1+u)^{a+1} + \sum_{b=0}^{a} (-1)^{b-1} \binom{a}{b} (1+u)^{a-b+1}. \end{equation}
Factoring out $(-1)\cdot (1+u)^{a+1}$ we get
\begin{equation} \sum_{b=0}^a (-1)^b \binom{a}{b} (1+u)^{-b} = (1 - \frac{1}{1+u})^a = (\frac{u}{1+u})^a, \qquad f(u) = u^a + u^{a+1} + (1+u)^{a+1} - (1+u)^{a+1}(\frac{u}{1+u})^a. \end{equation}
Now we simplify and get the desired formula
\begin{equation} f(u) = u^a(1+u) + (1+u)^{a+1} - u^a (1+u) = (1+u)^{a+1}. \end{equation}The proof is complete.
\end{proof}

\subsection{Examples of linear partial traces}\label{subsection:somepartialtraces}

\begin{Example}\label{example:smoothH3}
We work in type $H_{3}$ with $m_{12} = 5$ and $m_{23} = 3$ and $m_{13} = 2$. We use the color code from \eqref{eq:color-code}. Let $\underline{d}=232$ and $\un{t} = (1,2,3)$. We have a linear branching graph
\begin{gather}\label{eq:h3-pgraph}
\Gamma(\underline{d}123)=
\begin{tikzcd}[ampersand replacement=\&,row sep=scriptsize,column sep=scriptsize,arrows={shorten >=-0.5ex,shorten <=-0.5ex},labels={inner sep=0.05ex},arrow style=tikz]
\emptyset \ar[r,yshift=0.1cm,soergeltwo]
\&
2\ar[r,yshift=0.1cm,soergelthree]
\&
23\ar[r,yshift=0.1cm,soergeltwo]
\ar[l,yshift=-0.1cm,soergeltwo]
\&
\underline{d}\ar[r,yshift=0.1cm,soergelone]
\& 
\underline{d}1\ar[r,yshift=0.1cm,soergeltwo]
\ar[l,yshift=-0.1cm,soergeltwo]
\& 
\underline{d}12\ar[r,yshift=0.1cm,soergelthree]
\ar[l,yshift=-0.1cm,soergelthree]
\& 
\underline{d}123
\end{tikzcd}.\end{gather}
The formula \eqref{eq:LIFviaquantumnumbers} gives $\LIF(d1,2,d) = \partial_2(\alpha_1)$ which we denote as $-[2]_{2,1}$, and $\LIF(d12,3,d1) = \partial_3(\alpha_2)$ which we denote as $-[2]_{3,2}$. (Typically one has $[2]_{3,2} = 1$ and $[2]_{2,1} = \phi$, the golden ratio, though other (non-symmetric) Cartan matrices are possible.) Though not essential for our algorithm, we can compute the idempotents as
\begin{gather}
\begin{tikzpicture}[anchorbase,scale=1]
\draw[mor] (0,0) rectangle (1,0.5) node[black,pos=0.5]{$\underline{d}1$};
\end{tikzpicture}
=
\begin{tikzpicture}[anchorbase,scale=1]
\draw[mor] (0,0) rectangle (1,0.5) node[black,pos=0.5]{$\underline{d}$};
\draw[soergelone] (1.25,0) to (1.25,0.5);
\end{tikzpicture}
,\quad
\begin{tikzpicture}[anchorbase,scale=1]
\draw[mor] (0,0) rectangle (1,0.5) node[black,pos=0.5]{$\underline{d}12$};
\end{tikzpicture}
=
\begin{tikzpicture}[anchorbase,scale=1]
\draw[mor] (0,0) rectangle (1,0.5) node[black,pos=0.5]{$\underline{d}1$};
\draw[soergeltwo] (1.25,0) to (1.25,0.5);
\end{tikzpicture}
+\tfrac{1}{[2]_{2,1}}\cdot
\begin{tikzpicture}[anchorbase,scale=1]
\draw[mor] (0,0) rectangle (1.1,0.5) node[black,pos=0.5]{$\underline{d}1$};
\draw[mor] (0,1) rectangle (1.1,1.5) node[black,pos=0.5]{$\underline{d}$};
\draw[mor] (0,2) rectangle (1.1,2.5) node[black,pos=0.5]{$\underline{d}1$};
\draw[soergeltwo] (1.25,0) to (1.25,0.5) to (1,0.8);
\draw[soergeltwo] (0.75,0.5) to (1,0.8) to (1,1);
\draw[soergelone,markedone=1] (1,0.5) to (1,0.65);
\draw[usual] (0.25,0.5) to (0.25,1);
\draw[usual] (0.5,0.5) to (0.5,1);
\draw[soergeltwo] (1.25,2.5) to (1.25,2) to (1,1.65);
\draw[soergeltwo] (0.75,2) to (1,1.65) to (1,1.5);
\draw[soergelone,markedone=1] (1,2) to (1,1.85);
\draw[usual] (0.25,1.5) to (0.25,2);
\draw[usual] (0.5,1.5) to (0.5,2);
\end{tikzpicture}
,\quad
\begin{tikzpicture}[anchorbase,scale=1]
\draw[mor] (0,0) rectangle (1,0.5) node[black,pos=0.5]{$\underline{d}123$};
\end{tikzpicture}
=
\begin{tikzpicture}[anchorbase,scale=1]
\draw[mor] (0,0) rectangle (1,0.5) node[black,pos=0.5]{$\underline{d}12$};
\draw[soergelthree] (1.25,0) to (1.25,0.5);
\end{tikzpicture}
+\tfrac{1}{[2]_{3,2}}\cdot
\begin{tikzpicture}[anchorbase,scale=1]
\draw[mor] (0,0) rectangle (1.1,0.5) node[black,pos=0.5]{$\underline{d}12$};
\draw[rex] (0,0.5) rectangle (1.1,0.6);
\draw[rex] (0,1.1) rectangle (1.1,1.2);
\draw[mor] (0,1.2) rectangle (1.1,1.7) node[black,pos=0.5]{$\underline{d}1$};
\draw[rex] (0,1.7) rectangle (1.1,1.8);
\draw[rex] (0,2.3) rectangle (1.1,2.4);
\draw[mor] (0,2.4) rectangle (1.1,2.9) node[black,pos=0.5]{$\underline{d}12$};
\draw[soergelthree] (1.25,0) to (1.25,0.6) to (1,0.9);
\draw[soergelthree] (0.75,0.6) to (1,0.9) to (1,1.1);
\draw[soergeltwo,markedtwo=1] (1,0.6) to (1,0.75);
\draw[usual] (0.25,0.6) to (0.25,1.1);
\draw[usual] (0.5,0.6) to (0.5,1.1);
\draw[soergelthree] (1.25,2.9) to (1.25,2.3) to (1,1.95);
\draw[soergelthree] (0.75,2.3) to (1,1.95) to (1,1.8);
\draw[soergeltwo,markedtwo=1] (1,2.3) to (1,2.15);
\draw[usual] (0.25,1.8) to (0.25,2.3);
\draw[usual] (0.5,1.8) to (0.5,2.3);
\end{tikzpicture}
.
\end{gather}
For the final step note that $232123=323123=321323$, which indicates the rex moves needed in the purple strips above.

Now let $f_{3,0} \in R$ be any polynomial. The recursion from \autoref{thm:linearpartialtracerecursionbetter} computes the partial trace $A_3(f_{3,0}) = C_{3,0}(f_{3,0})$ as
\begin{gather*}
\begin{tikzcd}[ampersand replacement=\&]
\begin{tikzpicture}[anchorbase,scale=1]
\draw[mor] (0,0) rectangle (1.1,0.5) node[black,pos=0.5]{\raisebox{-0.06cm}{$\underline{d}123$}};
\draw[soergelthree] (1,0) to[out=270,in=180] (1.4,-0.4) to[out=0,in=270] (1.8,0) to (1.8,0.5)to[out=90,in=0] (1.4,0.9) to[out=180,in=90] (1,0.5);
\node[scale=0.8] at (1.4,0.25) {$f_{3,0}$};
\draw[soergeltwo] (0.8,0) to[out=270,in=180] (1.4,-0.6) to[out=0,in=270] (2,0) to (2,0.5)to[out=90,in=0] (1.4,1.1) to[out=180,in=90] (0.8,0.5);
\draw[soergelone] (0.6,0) to[out=270,in=180] (1.4,-0.8) to[out=0,in=270] (2.2,0) to (2.2,0.5)to[out=90,in=0] (1.4,1.3) to[out=180,in=90] (0.6,0.5);
\end{tikzpicture}	
\ar[d,swap,"{\alpha_{3}\partial_{3}(\placeholder)}"]
\ar[dr,"{-\frac{\partial_{3}(\placeholder)}{[2]_{3,2}}}",orchid]
\&
\&
\\
\begin{tikzpicture}[anchorbase,scale=1]
\draw[mor] (0,0) rectangle (1.1,0.5) node[black,pos=0.5]{\raisebox{-0.06cm}{$\underline{d}12$}};
\draw[soergeltwo] (1,0) to[out=270,in=180] (1.4,-0.4) to[out=0,in=270] (1.8,0) to (1.8,0.5)to[out=90,in=0] (1.4,0.9) to[out=180,in=90] (1,0.5);
\node[scale=0.8] at (1.4,0.25) {$f_{2,0}$};
\draw[soergelone] (0.8,0) to[out=270,in=180] (1.4,-0.6) to[out=0,in=270] (2,0) to (2,0.5)to[out=90,in=0] (1.4,1.1) to[out=180,in=90] (0.8,0.5);
\end{tikzpicture}
\ar[d,swap,"{\alpha_{2}\partial_{2}(\placeholder)}"]
\ar[dr,swap,near start,"{-\frac{\partial_{2}(\placeholder)}{[2]_{2,1}}}",orchid]
\&
\begin{tikzpicture}[anchorbase,scale=1]
\draw[mor] (0,0) rectangle (1.1,0.5) node[black,pos=0.5]{\raisebox{-0.06cm}{$\underline{d}12$}};
\draw[mor] (0,0.9) rectangle (1.1,1.4) node[black,pos=0.5]{\raisebox{-0.06cm}{$\underline{d}12$}};
\draw[soergeltwo, markedtwo=1] (1,0.5) to (1,0.6);
\draw[soergeltwo, markedtwo=1] (1,0.9) to (1,0.8);
\draw[usual] (0.5,0.5) to (0.5,0.9);
\draw[soergeltwo] (1,0) to[out=270,in=180] (1.4,-0.4) to[out=0,in=270] (1.8,0) to (1.8,1.4)to[out=90,in=0] (1.4,2) to[out=180,in=90] (1,1.4);
\node[scale=0.8] at (1.4,0.7) {$f_{2,1}$};
\draw[soergelone] (0.8,0) to[out=270,in=180] (1.4,-0.6) to[out=0,in=270] (2,0) to (2,1.4)to[out=90,in=0] (1.4,2.2) to[out=180,in=90] (0.8,1.4);
\end{tikzpicture}	
\ar[dl,near start,"{\alpha_{2}\cdot}",spinach]
\ar[d,"{-\frac{2}{[2]_{2,1}}}\cdot"]
\ar[dr,"{\frac{1}{[2]_{2,1}^{2}}\partial_{2}(\placeholder)}",orchid]
\&
\\
\begin{tikzpicture}[anchorbase,scale=1]
\draw[mor] (0,0) rectangle (1.1,0.5) node[black,pos=0.5]{\raisebox{-0.06cm}{$\underline{d}1$}};
\draw[soergelone] (1,0) to[out=270,in=180] (1.4,-0.4) to[out=0,in=270] (1.8,0) to (1.8,0.5)to[out=90,in=0] (1.4,0.9) to[out=180,in=90] (1,0.5);
\node[scale=0.8] at (1.4,0.25) {$f_{1,0}$};
\end{tikzpicture}
\ar[d,swap,"{\alpha_{1}\partial_{1}(\placeholder)}"]	
\&
\begin{tikzpicture}[anchorbase,scale=1]
\draw[mor] (0,0) rectangle (1.1,0.5) node[black,pos=0.5]{\raisebox{-0.06cm}{$\underline{d}1$}};
\draw[mor] (0,0.9) rectangle (1.1,1.4) node[black,pos=0.5]{\raisebox{-0.06cm}{$\underline{d}1$}};
\draw[soergelone, markedone=1] (1,0.5) to (1,0.6);
\draw[soergelone, markedone=1] (1,0.9) to (1,0.8);
\draw[usual] (0.5,0.5) to (0.5,0.9);
\draw[soergelone] (1,0) to[out=270,in=180] (1.4,-0.4) to[out=0,in=270] (1.8,0) to (1.8,1.4)to[out=90,in=0] (1.4,1.8) to[out=180,in=90] (1,1.4);
\node[scale=0.8] at (1.4,0.7) {$f_{1,1}$};
\end{tikzpicture}	
\ar[dl,"{\alpha_{1}\cdot}",spinach]
\&
\begin{tikzpicture}[anchorbase,scale=1]
\draw[mor] (0,0) rectangle (1.1,0.5) node[black,pos=0.5]{\raisebox{-0.06cm}{$\underline{d}1$}};
\draw[mor] (0,0.9) rectangle (1.1,1.4) node[black,pos=0.5]{\raisebox{-0.06cm}{$\underline{d}1$}};
\draw[mor] (0,1.8) rectangle (1.1,2.3) node[black,pos=0.5]{\raisebox{-0.06cm}{$\underline{d}1$}};
\draw[soergelone, markedone=1] (1,0.5) to (1,0.6);
\draw[soergelone, markedone=1] (1,0.9) to (1,0.8);
\draw[usual] (0.5,0.5) to (0.5,0.9);
\draw[soergelone, markedone=1] (1,1.4) to (1,1.5);
\draw[soergelone, markedone=1] (1,1.8) to (1,1.7);
\draw[usual] (0.5,1.4) to (0.5,1.8);
\draw[soergelone] (1,0) to[out=270,in=180] (1.4,-0.4) to[out=0,in=270] (1.8,0) to (1.8,2.3)to[out=90,in=0] (1.4,2.7) to[out=180,in=90] (1,2.3);
\node[scale=0.8] at (1.4,1.15) {$f_{1,2}$};
\end{tikzpicture}
\ar[dll,bend left=15,"{\alpha_{1}^{2}\cdot}",spinach]	
\\
\begin{tikzpicture}[anchorbase,scale=1]
\draw[mor] (0,0) rectangle (1.1,0.5) node[black,pos=0.5]{\raisebox{-0.06cm}{$\underline{d}$}};
\node[scale=0.8] at (1.4,0.25) {$f_{0,0}$};
\end{tikzpicture}	
\&
\&
\end{tikzcd}
.
\end{gather*}
This is to be read as follows. We have $f_{2,0}=\alpha_{3}\partial_{3}(f_{3,0})$. 
If more than one arrow points to a cell we take the sum, i.e. $f_{0,0}=\alpha_{1}\partial_{1}(f_{1,0})+\alpha_{1}f_{1,1}+\alpha_{1}^{2}f_{1,2}$. The polynomial $f_{k,a}$ is the input to $C_{k,a}$ at various steps along the recursive computation.

The polynomial $f_{0,0}$ is the partial trace.
\end{Example}


\begin{Example} \label{Ex:F4} We work in type $F_4$, with $m_{23} = 4$. Let $\un{d} = 1434$ and $\un{w} = \un{d}2341$.  The superspine of $\Gamma(\underline{w})$ is:

\begin{gather}
\begin{tikzcd}[ampersand replacement=\&,row sep=scriptsize,column sep=scriptsize,arrows={shorten >=-0.5ex,shorten <=-0.5ex},labels={inner sep=0.05ex},arrow style=tikz]
\emptyset \ar[r,yshift=0.1cm,soergelone]
\&
\underline{1}\ar[r,yshift=0.1cm,soergelfour]
\&
14  \ar[r,yshift=0.1cm,soergelthree]
\&
143 \ar[r,yshift=0.1cm,soergelfour]\ar[l,yshift=-0.1cm,soergelfour]
\&
\underline{d} \ar[r,yshift=0.1cm,soergeltwo]
\&
\underline{d}2 \ar[r,yshift=0.1cm,soergelthree]\ar[l,yshift=-0.1cm,soergelthree]
\&
\underline{d}23  \ar[r,yshift=0.1cm,soergelfour] \ar[l,yshift=-0.1cm,soergelfour]
\&
\underline{d}234 \ar[r,yshift=0.1cm,soergelone] 
\&
\underline{w}
\end{tikzcd}\;.
\end{gather}

Write $[2] = -\partial_3(\alpha_2)$ and $[2]' = -\partial_2(\alpha_3)$, and we know that $[2] [2]' = 2$ since $m_{23} = 4$. We assume that $\partial_3(\alpha_4) = -1 = \partial_4(\alpha_3)$. The relevant local intersection forms are computed to be
\begin{equation} \LIF(d23,4,d2) = -1, \qquad \LIF(d2,3,d) = [2]. \end{equation}
Following the algorithm above, we get that
\begin{equation} \ptr_{d,(2,3,4,1)}(f) = \alpha_2 D(\alpha_1 \partial_1 f), \end{equation}
where $D$ is the operator
\begin{align} \nonumber D(g) = & \partial_2(\alpha_3 \partial_3(\alpha_4 \partial_4(g))) + \partial_2(\alpha_3 \partial_4(g)) + \frac{1}{[2]} \partial_3(\alpha_4 \partial_4(g)) \\
& + \frac{2}{[2]} \partial_4(g) - \frac{1}{[2]} s_3(\partial_4(g)) + \frac{1}{[2]^2} \alpha_2 \partial_3(\partial_4(g)). \end{align}
Using the Leibniz rule one can reduce this with effort (and some cancellation) to
\begin{equation} \label{F4calc1} \ptr_{d,(2,3,4,1)}(f) = \frac{-1}{[2]^2} \alpha_2 Q \partial_{341} f + \alpha_2 s_2(\alpha_1) s_2(\alpha_3) s_2 s_4 (\alpha_3) \partial_{2341}(f), \end{equation}
where $Q$ is the quadratic polynomial
\begin{equation} 2 [2] \alpha_1 \alpha_3 + [2] \alpha_1 \alpha_4 + 3 \alpha_1 \alpha_2 + [2]^2 s_2(\alpha_3) s_2 s_4(\alpha_3). \end{equation}
One can verify easily that both $Q$ and $\alpha_2 s_2(\alpha_1) s_2(\alpha_3) s_2 s_4 (\alpha_3)$ are invariant under $s_4$.

In the next section we will be interested in $\partial_{d}(\ptr_{d,(2,3,4,1)}(f))$, where $\partial_d = \partial_{1434} = \partial_1 \partial_4 \partial_3 \partial_4$. 
As a first step, we apply $\partial_4$ to $\ptr_{d,(2,3,4,1)}(f)$, and since most terms are $s_4$-invariant we get
\begin{equation} \partial_4 \ptr_{d,(2,3,4,1)}(f) = \frac{-1}{[2]^2} \alpha_2 Q \partial_{4341} f + \alpha_2 s_2(\alpha_1) s_2(\alpha_3) s_2 s_4 (\alpha_3) \partial_{42341}(f). \end{equation}
Consequently, if $f$ is a quartic polynomial (so that the end result is a scalar) (technically, $\deg f = 8$), then the second term vanishes, and $\partial_{4341}(f) = \partial_d(f)$ is in the invariant subring $R^{134}$. Applying $\partial_{143}$ we get
\begin{equation} \partial_d \ptr_{d,(2,3,4,1)}(f) = \frac{-1}{[2]^2} \partial_{143}(\alpha_2 Q) \partial_d(f). \end{equation}

The operators $\partial_d \ptr_{d,(2,3,4,1)}(-)$ and $\partial_d(-)$ both live in the high-dimensional space of linear functionals on polynomials of a given degree (quartics). The phenomenon we wish to emphasize here is that these operators $\partial_d \ptr_{d,(2,3,4,1)}(-)$ and $\partial_d(-)$ are secretly colinear, which was not a priori obvious! We leave it as a (tedious) exercise to verify that
\begin{equation} \frac{-1}{[2]^2} \partial_{143}(\alpha_2 Q) = 3, \qquad \partial_d \ptr_{d,(2,3,4,1)}(f) = 3 \partial_d(f). \end{equation}
The number $3$ ends up being the categorical dimension of $w$ in its asymptotic Hecke category, as we explain in the next section.

If we knew somehow that these operators would be colinear, then we could compute their ratio by simply evaluating both operators on a suitably generic polynomial of the appropriate degree. This obviates the need for complicated (yet darkly fascinating) arguments using the Leibniz rule.  \end{Example}

\section{Computing in asymptotic Hecke categories}  \label{section:asymptotichecke}

\subsection{Cell theory basics}\label{subsection:hecke-cell}

Previously in \autoref{subsection:filtration} we were interested in two-sided ideals $\catstuff{I}_{\setstuff{X}}$ generated by identity maps of
various objects. Now we use the monoidal structure, and are interested in three- and four-sided ideals.

\begin{Notation} For a monoidal Krull--Schmidt category $\catstuff{C}$, let $\setstuff{X}$ be a set of objects in $\catstuff{C}$. Let $\catstuff{I}_{\setstuff{X}}^L$ be the left-monoidal ideal inside $\catstuff{C}$
generated by the identity maps of all objects in $\setstuff{X}$. That is, $\catstuff{I}_{\setstuff{X}}^L$ contains the identity maps of these objects, and is closed under composition with arbitrary morphisms, and is closed under the operation $\morstuff{id}_{\obstuff{Z}} \otimes (-)$ for any object $Z$. Similarly, let $\catstuff{I}_{\setstuff{X}}^R$ be the right-monoidal ideal (it is closed under $(-) \otimes \morstuff{id}_{\obstuff{Z}}$), and $\catstuff{I}_{\setstuff{X}}^{LR}$ the two-sided monoidal ideal, generated by the identity maps of all objects in $\setstuff{X}$.
\end{Notation}

Now we focus on $\catstuff{S}$.

\begin{Notation} For $w, x \in W$ we say that $w \le_{L} x$ if $\catstuff{I}_{\obstuff{B}_w}^L \subset \catstuff{I}_{\obstuff{B}_x}^L$. It is clear that $\le_L$ is a preorder, and the equivalence classes under $\le_L$ are called \emph{left cells}. We write $w \sim_L x$ if they are in the same left cell. We define preorders $\le_{R}$ and $\le_{LR}$ and equivalence relations $\sim_R$ and $\sim_{LR}$ similarly, and their equivalence classes are called \emph{right cells} and \emph{two-sided cells} respectively. \end{Notation}

Thus, $w \le_{LR} x$ if $B_w$ is a direct summand of $M \otimes B_x \otimes N$ for some objects $M, N$ in $\catstuff{S}$. Note that the preorder $\le_{LR}$ descends to a partial order $\le_{LR}$ on two-sided cells.

\begin{Remark} The preorders $\le_{?}$ are determined on the Grothendieck group, and match the definitions of Kazhdan--Lusztig from \cite{KaLu-reps-coxeter-groups} since $[\obstuff{B}_w] = b_w$. In particular, the span of $\{b_x\}_{x \le_{L} w}$ is a left ideal in the Hecke algebra, and similarly $\{b_x\}_{x \le_{LR} w}$ is a two-sided ideal. \end{Remark}

\begin{Notation} For $x \in W$ we write $\catstuff{I}_{\le x}^{LR}$ for the monoidal ideal generated by the identity map of $B_x$. This agrees with $\catstuff{I}_{\le w}^{LR}$ for any $w$ with $w \sim_{LR} x$, so that we often index this ideal not by elements of $W$ but by the two-sided cells themselves, i.e. we write $\catstuff{I}^{LR}_{\le \mathcal{J}} := \catstuff{I}_{\le x}^{LR}$ when $\mathcal{J}$ is the two-sided cell containing $x$. The corresponding ideal in the Hecke algebra is denoted $I_{\le \mathcal{J}}^{LR}$. We often omit the superscript $LR$ when the ideal is indexed by a two-sided cell $\mathcal{J}$, so that the two-sided nature is understood. \end{Notation}

\begin{Notation} The quotient category $\catstuff{S}/\catstuff{I}_{< \mathcal{J}}$ will be denoted $\catstuff{S}_{\not< \mathcal{J}}$, and we write morphism spaces between objects $B$ and $B'$ as $\Hom_{\not< \mathcal{J}}(B,B')$. \end{Notation}

The ideals $\catstuff{I}_{\le \mathcal{J}}$ are closed under duality (because the set of generating objects is closed under duality) and under adjunction (any four-sided ideal is
closed under adjunction, because adjunction arises from composition with cups and caps). Thus, $\catstuff{S}_{\not< \mathcal{J}}$ inherits both the duality and rigidity structures from
$\catstuff{S}$.

By the expanded Soergel hom formula, it is easy to compute a \textit{spanning set} for morphisms in $\catstuff{S}_{\not< \mathcal{J}}$ as a right $R$-module. If $\Hom_{\catstuff{S}}(B,B')$ has a basis of size $\sum q_y r_y$ as in \autoref{thm:SHFplus}, where the basis elements factor as $B \to B_y \to B'$,
then one obtains a spanning set for $\Hom_{\not< \mathcal{J}}(B,B')$ by ignoring all basis elements where $y < \mathcal{J}$. However, this is not a basis as a right $R$-module, and
morphism spaces in the quotient are not free!

\begin{Example} The monoidal identity is the unique element in the top cell. The endomorphism ring of the monoidal identity in $\catstuff{S}$ is $R$, with the identity map as its sole basis element. However, multiplication by $\alpha_i \in R$ factors through $B_i$, so it becomes zero in $\catstuff{S}_{\not< \munit}$. Indeed, $\End_{\not< \munit}(\munit) = \R \cdot \id$ whenever $R$ is generated by the simple roots. \end{Example}

\begin{Example} More generally, for any element $x \in \mathcal{J}$ in any two-sided cell, there are nonzero morphisms $\munit \to B_x \to \munit$ (corresponding to the nonzero KL polynomial $p_{x,1}$) which compose to nonzero polynomials in $R$. These polynomials vanish in $\End_{\not\le \mathcal{J}}(\munit)$. The ideals $\catstuff{I}_{\le \mathcal{J}}(\munit, \munit) \subset R$ are scheme-theoretically interesting. For example, in \cite[\S 3]{El-dia-tl-categorification} it is proven that for the cell in $S_n$ associated to the partition $(n-2,1,1)$, the ideal cuts out the union of those lines which are obtained as the intersections of reflection hyperplanes in $V^*$. \end{Example}

\subsection{Cell theory and asymptotics}\label{subsection:hecke-cell-asymp}

For this section we refer to \cite[\S 13-15]{Lu-hecke-book}, and statements given without proof are proven there. A list of properties\footnote{The book \cite{Lu-hecke-book} also
covers Hecke algebras with unequal parameters, where the properties (P1) through (P15) are conjectural.} (P1) through (P15) is given in \cite[Conjecture 14.2]{Lu-hecke-book}, and
they are proven in \cite[\S 15]{Lu-hecke-book}. One must assume a boundedness conjecture \cite[13.4]{Lu-hecke-book}, which Lusztig proves for all finite and affine Weyl groups; it was recently proven in full generality in \cite[Theorem 1.2]{Ch-Hu-boundness-Lusztig-a}.

We recall \emph{Lusztig's $a$-function}, defined as follows. For a given $w \in W$, let $a(w)$ be the maximal integer $k$ such that there exist $y, z \in W$ with $B_w(k)$
being a direct summand of $B_y \otimes B_z$. By duality, $-a(w)$ is also the minimal such integer. Lusztig proved (P4) that $w \sim_{LR} x$ implies $a(w) = a(x)$,
so that $a$ is often viewed as a function on two-sided cells rather than a function on elements. We write $a(\mathcal{J}) := a(w)$ for $w \in \mathcal{J}$.

\begin{Example}\label{example:soergel2}
We continue \autoref{example:soergel1}. Using \eqref{B1B1} and properties of the monoidal identity, one observes that $a(\id) = 0$ and $a(s_1) = 1$. The elements $\id$ and $s_1$ are each in their own two-sided cell. \end{Example}

Let us use the following notation for structure coefficients\footnote{In Lusztig's book \cite[13.1]{Lu-hecke-book}, $m^w_{x,y}$ is denoted as $h_{x,y,z^{-1}}$.} in the Hecke algebra:
\begin{equation} \label{structurecoeff} b_x b_y = \sum_{w \in W} m^w_{x,y} b_w, \qquad B_x \otimes B_y \cong \bigoplus_{w \in W} B_w^{\oplus m^w_{x,y}}. \end{equation}

Fix a two-sided cell $\mathcal{J}$. The definition of the $a$-function guarantees that $m^w_{x,y} \in \N[v,v^{-1}]$ is supported in degrees between $v^{-a(\mathcal{J})}$ and
$v^{a(\mathcal{J})}$, for $x, y, w \in \mathcal{J}$. It is only the structure coefficients within the same cell that are guaranteed to have this degree bound; $m^w_{x,y}$ will often
involve even more positive exponents when $x,y \in \mathcal{J}$ but $w < \mathcal{J}$ is in strictly lower cells.

\begin{Example} \label{Ex:B2cellcomp} Let $W$ have type $B_2$. There is a two-sided cell $\mathcal{J} = \{1,2, 12, 21, 121, 212\}$ consisting of the six elements with a unique reduced expression, and $a(\mathcal{J}) = 1$. We have
\begin{equation} b_{121} b_{121} = (v+v^{-1}) b_1 + (v^{-2} + 2 + v^2)b_{1212}. \end{equation}
The coefficient of $b_1$ has exponents $\le 1$, though the coefficient of $b_{1212}$ does not, which is acceptable since $1212$ lies in a lower cell. For further examples we have
\begin{equation} b_{12} b_{21} = (v+v^{-1}) b_{121} + (v+v^{-1}) b_1, \qquad b_{12} b_1 = b_{121} + b_1. \end{equation}
The last example demonstrates that the exponent $v^{a(\mathcal{J})}$ need not be achieved, i.e. the coefficient of $v^{a(\mathcal{J})}$ might be zero.
\end{Example}

One obtains a quotient algebra from the Hecke algebra by killing the ideal $I^{LR}_{< \mathcal{J}}$, and in this algebra no exponents larger than $v^{a(\mathcal{J})}$ appear in the structure coefficients. The asymptotic Hecke algebra comes from this quotient by taking the limit (within a given cell $\mathcal{J}$) as $v^{-1}$ approaches zero, which
is eventually dominated by the maximal power of $v$ which appears.

\begin{Definition}
Fix a two-sided cell $\mathcal{J}$. For $x, y, w \in \mathcal{J}$ let us write $\gamma^w_{x,y} \in \N$ for the coefficient\footnote{The integer $\gamma^w_{x,y}$ is denoted $\gamma_{x,y,w^{-1}}$ in \cite[13.6]{Lu-hecke-book}.} of $v^{a(\mathcal{J})}$ inside $m^w_{x,y}$ from \eqref{structurecoeff}. Define a $\Z$-algebra $A_{\mathcal{J}}$ with basis $\{a_w\}_{w \in \mathcal{J}}$ and multiplication rule
\begin{equation} a_x a_y = \sum_w \gamma^w_{x,y} a_w. \end{equation}
This is the \emph{asymptotic Hecke algebra} in cell $\mathcal{J}$. \end{Definition}

\begin{Example} Continuing \autoref{Ex:B2cellcomp} we have
\begin{equation} a_{121} a_{121} = a_1, \qquad a_{12} a_{21} = a_{121} + a_1, \qquad a_{12} a_1 = 0, \end{equation}as one easily checks.
\end{Example}

For a proof that $A_{\mathcal{J}}$ is an associative algebra, see \cite[18.3]{Lu-hecke-book}. In fact, Lusztig proves that $A_{\mathcal{J}} \ot_{\Z} \C$ is a semisimple algebra \cite[21.9]{Lu-hecke-book}. Note that $A_{\mathcal{J}}$ might not be unital, but it is locally unital, see below.

\begin{Lemma} ((P7) and \cite[Prop 13.9]{Lu-hecke-book}) \label{Lem:gammasymmetries} We have $\gamma_{x,y}^w = \gamma_{y,w^{-1}}^{x^{-1}} = \gamma_{y^{-1},x^{-1}}^{w^{-1}}$. \end{Lemma}

By the expanded Soergel hom formula, $\gamma_{x,y}^w$ is the dimension of $\Hom^{-a}(B_w,B_x B_y)$ in the quotient category $\catstuff{S}/\catstuff{I}_{< \mathcal{J}}$. The symmetries
in \autoref{Lem:gammasymmetries} come from the symmetries of $\catstuff{S}$, all of which preserve $\catstuff{I}_{< \mathcal{J}}$. That is, adjunction will give an isomorphism
\begin{equation} \Hom^{-a}(B_z,B_x B_y) \cong \Hom^{-a}(\rdual{B_x}, B_y \rdual{B_z}) = \Hom^{-a}(B_{x^{-1}},B_y B_{z^{-1}}). \end{equation} The horizontal flip will give an
isomorphism \begin{equation} \Hom^{-a}(B_z,B_x B_y) \cong \Hom^{-a}(B_{z^{-1}},B_{y^{-1}}B_{x^{-1}}). \end{equation}

An element $d \in \mathcal{J}$ is called a \emph{Duflo involution} if $p_{d,1}$ has a nonzero coefficient of $v^{a(\mathcal{J})}$ (c.f. (P1)).

\begin{Lemma} ((P6), (P13)) \label{lem:duflobasic} There is exactly one Duflo involution in each left cell (resp. in each right cell). Each Duflo involution is indeed an involution. \end{Lemma}

\begin{Lemma} ((P13), (P2), (P5), (P8)) Let $d \in \mathcal{J}$ be a Duflo involution. For each $x, y \in \mathcal{J}$ we have $\gamma_{d,x}^y = \delta_{x,y} \delta_{x \sim_R d}$. That is, $\gamma_{d,x}^y$ is zero unless $x$ is in the same right cell as $d$ and unless $x=y$, in which case it is $1$. We also have $\gamma_{x,y}^d = \delta_{x,y^{-1}} \delta_{x \sim_R d}$ and $\gamma_{x,d}^y = \delta_{x,y} \delta_{x \sim_L d}$. \end{Lemma}

In particular, $\gamma_{d,d}^{d} = 1$. Thus, $a_d \in A_{\mathcal{J}}$ is an idempotent, and $a_d, a_{d'}$ are orthogonal when $d \ne d'$. Moreover, 
\begin{equation} a_d a_x = \delta_{x \sim_R d} a_x.  \end{equation}
As a consequence, the sum over all Duflo involutions $\sum_d a_d$ in $\mathcal{J}$ is an idempotent decomposition of the identity element of $A_{\mathcal{J}}$. If there are infinitely many left cells in $\mathcal{J}$ then this sum is not well-defined and $A_{\mathcal{J}}$ is not unital, but the sum serves instead as a local unit. Lusztig proves there are finitely many left cells in a finite or affine Weyl group in \cite[Theorem 18.2]{Lu-hecke-book}.

Recall that $\obstuff{B}_w$ has adjoint $\rdual{\obstuff{B}_w} = \obstuff{B}_{w^{-1}}$. As one easily checks, $\rdual{\placeholder}$
interchanges left and right cells.  Associated to a left cell $\mathcal{L}$ we have a \emph{diagonal cell} defined as
\begin{gather}
\mathcal{H}=\mathcal{H}(\mathcal{L}):=\mathcal{L}\cap\mathcal{L}^{\star}=\mathcal{L}\cap\mathcal{L}^{-1}.
\end{gather}
Diagonal cells are neither left nor right nor two-sided cells, see \autoref{Ex:typeAintroMorita}.

\begin{Lemma} Each diagonal cell has a unique Duflo involution, and each Duflo involution lives in a unique diagonal cell. \end{Lemma}

\begin{proof} This is immediate from \autoref{lem:duflobasic}. \end{proof}

The $\Z$-span of elements $\{a_w\}_{w \in \mathcal{H}}$ is an algebra we denote\footnote{Lusztig studies these algebras frequently in \cite{Lu-hecke-book} where they are denoted $J^{\Gamma \cap \Gamma^{-1}}$. Here $\Gamma$ is a left cell.} $A_{\mathcal{H}}$, and it agrees with the idempotent truncation $a_d \cdot A_{\mathcal{J}} \cdot a_d$. It has identity element $a_d$ for the unique Duflo
involution of $\mathcal{H}$.

\begin{Example} For $S_n$, the two-sided cells are in bijection with partitions of $n$. Fix such a partition $\lambda$. The Robinson-Schensted correspondence gives a bijection between pairs $(P,Q)$ of standard tableaux of shape $\lambda$, and elements of the cell $\lambda$; we write $w_{(P,Q)}$ for the corresponding permutation, and $b_{(P,Q)}$ or $a_{(P,Q)}$ for the corresponding basis elements. The left cells in $\lambda$ are parametrized by standard tableaux, and have the form $\{w_{(-,Q)}\}$ for a given $Q$. Each diagonal cell consists of a single element $w_{(P,P)}$, which is an involution.  One has
\begin{equation} a_{(P,Q)} a_{(U,V)} = \delta_{Q,U} a_{(P,V)}. \end{equation}
In particular, $A_{\lambda}$ is a matrix algebra with rows and columns indexed by standard tableaux of shape $\lambda$. The diagonal elements $w_{(P,P)}$ are precisely the involutions, and every involution of $S_n$ is a Duflo involution in its diagonal cell. For the diagonal cell $(Q,Q)$ we have $A_{(Q,Q)} \cong \Z$. \end{Example}

\begin{Example} \label{S4middle} As a sub-example, for $S_4$ the middle cell (associated to the partition $(2,2)$) has four elements $\{s_1 s_3, s_1 s_3 s_2, s_2 s_1 s_3, s_2 s_1 s_3 s_2\}$, and the $a$-function is $2$. The involutions are $s_1 s_3$ and $s_2 s_1 s_3 s_2$, and correspond to the diagonal elements of a $2 \times 2$ matrix algebra. The elements $s_1 s_3 s_2$ and $s_2 s_1 s_3$ are the off diagonal elements. For example,
\begin{equation} a_{s_1 s_3} a_{s_1 s_3 s_2} = a_{s_1 s_3 s_2}, \qquad a_{s_1 s_3} a_{s_2 s_1 s_3} = 0, \qquad a_{s_2 s_1 s_3} a_{s_1 s_3 s_2} = a_{s_2 s_1 s_3 s_2}, \end{equation}by direct computation.
\end{Example}

\begin{Example} In type $H_2$, let $\mathcal{H}$ be the diagonal cell containing the Duflo involution $s_1$. Then $\mathcal{H} = \{s_1, s_1 s_2 s_1\}$. We have that $a_1$ is the identity of $A_{\mathcal{H}}$, and $a_{121}^2 = a_{121} + a_1$. Thus, $A_{\mathcal{H}} \cong \Z[\phi]$, where $\phi$ is the golden ratio. \end{Example}

Many examples can be found in the sequel  \cite{ElRoTu-ah-2}.

As a quick aside, when $W$ is finite, recall that multiplication by the longest word $w_{0}\in W$ sends 
left, right and two-sided cells to left, right and two-sided cells, 
and reverses the corresponding cell order \cite[Remark 3.3a]{KaLu-reps-coxeter-groups}. The notation ${\placeholder}^{\prime}$ is used 
for this symmetry, e.g. for all $\mathcal{H}=\mathcal{L}\cap\mathcal{L}^{-1}\subset\mathcal{J}$, there is a \emph{$w_0$-dual}
\begin{gather*}
\mathcal{H}^{\prime}
:=(\mathcal{L}w_{0})\cap(w_{0}\mathcal{L}^{-1}),
\end{gather*}
which satisfies $\#\mathcal{H}=\#\mathcal{H}^{\prime}$.

\subsection{Asymptotic Hecke categories: loose definition}\label{subsection:ahecke-examples}

In this section we recall the definition of the asymptotic Hecke category as typically given in the literature, see e.g. \cite[Section 2]{Lu-cells-tensor-cats} for the original construction (but in the setting of perverse sheaves), and \cite{ElWi-relative-lefschetz} or \cite[\S 18.15]{Lu-hecke-book} for a version with
Soergel bimodules. The essential input to this construction is the idea that one can canonically pick out the summands with the largest shift within a bimodule, though there are
subtleties in this idea which we need to explain.

\begin{Definition} Suppose $a \in \Z$. Suppose that $B$ is a Soergel bimodule, and that no summand of $B$ is isomorphic to $B_w(k)$ for $k>a$. Then there is a decomposition $B \cong X
\oplus Y$, where $X$ is a direct sum of various $B_w(a)$ with multiplicity, and $Y$ is a direct sum of various $B_w(k)$ for $k < a$. We call $X$ the \emph{$a$-degree submodule} of
$B$, noting that it is possible that $X = 0$. When $X$ is nonzero we call $X$ the \emph{max degree submodule} of $B$. \end{Definition}

\begin{Lemma} The $a$-degree submodule is canonical, i.e. the set $X \subset B$ is independent of the choice of decomposition $B \cong X \oplus Y$. \end{Lemma}

\begin{proof} By the Soergel hom formula, there are no morphisms of negative degree between indecomposable Soergel bimodules, and as a consequence $\Hom^0(B_w(a),Y) = 0$ for all $w \in W$. It is straightforward to verify that $X \subset B$ can be abstractly described as the span of the images of all degree zero maps $B_w(a) \to B$. \end{proof}

\begin{Remark} The construction of $X$ as a span of images also implies functoriality of the $a$-degree submodule under degree zero maps. \end{Remark}

\begin{Remark} Applying the same analysis to $Y$, one obtains a canonical submodule $Z \subset Y$ consisting of all bimodules with shift $a-1$. One can repeat this process, giving
rise to a canonical filtration on $B$ called the \emph{perverse filtration}. For details see \cite[\S 2.2]{ElWi-relative-lefschetz}. \end{Remark}

\begin{Remark} By similar arguments, there is a canonical quotient of $B$ corresponding to the summands with minimal (rather than maximal) shift. \end{Remark}

Note that $X$ is a summand of $B$ (whence it is a Soergel bimodule), but we have deliberately called it a submodule. The key subtlety is that the summand $Y \subset B$, complementary
to the max degree submodule, is not canonical. The decomposition of $B$, i.e. the choice of inclusion and projection maps to and from $X$ and $Y$, is not canonical, nor is the
idempotent projecting to $X$. It is better to think that $X$ as a \textit{submodule} is canonical, an element-theoretic statement, but that $X$ as a \textit{summand} is not canonical.

\begin{Example} \label{B1B1maxdegree} For $S_2$ we have $B_1 B_1 \cong B_1(1) \oplus B_1(-1)$. The inclusion map $B_1(1) \to B_1 B_1$ is canonical up to scalar, living in the one-dimensional space
$\Hom^{-1}(B_1, B_1B_1)$, and its image is the canonical subset, the max degree submodule. The projection map $B_1 B_1 \to B_1(-1)$ is also canonical up to scalar, its kernel being
the max degree submodule. However, the projection map $B_1B_1 \to B_1(1)$ is not canonical; the space $\Hom^1(B_1B_1,B_1)$ is multidimensional, and many non-colinear maps (with
different kernels) can pair with the canonical inclusion map to serve as a valid projection map. \end{Example}

Now fix a two-sided cell $\mathcal{J}$ and let $a = a(\mathcal{J})$. In the quotient $\catstuff{S}_{\not< \mathcal{J}}$, the full subcategory whose objects are $\{B_w(k)\}$ for $w \in
\mathcal{J}$ and $k \in \Z$ is closed under tensor product, and forms a (non-unital) monoidal category. Moreover, for any $w, x \in \mathcal{J}$, no summand of $B_w B_x$ in the
quotient has shift greater than $a$. Define a new tensor structure on the unshifted objects $B_w$, where $B_w \bullet B_x$ is defined as the $a$-degree submodule of $B_w \ot B_x$. (This is not necessarily the max degree submodule, since the $a$-degree submodule may be zero.) This new tensor structure is
associative, as $(B_w \bullet B_x) \bullet B_y$ agrees canonically with the $2a$-degree submodule of $(B_w B_x) B_y$, and the ordinary tensor product is associative. (More precisely, Lusztig uses these identifications of submodules to produce an associator.)

\begin{Definition} Let $\catstuff{A} = \catstuff{A}_{\mathcal{J}}$ be the monoidal category with objects $A_w$ associated to $w \in \mathcal{J}$. We think of $A_w$ as being the image
of the self-dual object $B_w$ (without shift) in the cell quotient from the previous paragraph. We set $A_w \star A_x := B_w \bullet B_x(-a)$ to be the $a$-degree submodule of $B_w
\ot B_x$, shifted so that it too is self-dual. Under this monoidal structure, $\bigoplus_d A_d$ is the (local) monoidal identity, as $d$ ranges over Duflo involutions in the cell. We
set morphisms in $\catstuff{A}_{\mathcal{J}}$ to agree with degree zero morphisms in $\catstuff{S}_{\not< \mathcal{J}}$. For example we have \begin{equation}
\Hom_{\catstuff{A}_{\mathcal{J}}}(A_w, A_x \star A_y) = \Hom^0_{\not< \mathcal{J}}(B_w, B_x \bullet B_y(-a)). \end{equation} Note that we can view $\Hom(A_w,A_x \star A_y)$ as living
inside $\Hom^{-a}_{\not< \mathcal{J}}(B_w, B_x \ot B_y)$.

For a diagonal cell $\mathcal{H} \subset \mathcal{J}$ set $\catstuff{A}_{\mathcal{H}} \subset \catstuff{A}_{\mathcal{J}}$ to be the full subcategory with objects $A_w$ for $w \in \mathcal{H}$.\end{Definition}

There are many subtleties hidden in this definition, so let us try to elucidate it, with a focus on being able to compute with it.

\subsection{Asymptotic Hecke categories in the Karoubi envelope}\label{subsection:karoubibullshit}

What kind of thing is the max-degree submodule exactly?

We are faced here with a rare situation where a piece of technical categorical nonsense is actually highly significant to one's
understanding. The essential issue is this: the definition of the $a$-degree submodule makes sense in the category of abstract Soergel bimodules but not in the Karoubi envelope of Bott--Samelson bimodules, and we need to do some work to bridge this gap.

Recall that the category $\catstuff{S}$ is defined by taking the category of Bott--Samelson bimodules (tensor products of $B_i$), taking its closure under direct sums and shifts (call
this $\catstuff{S}'$ temporarily), and then taking the closure under direct summands. The latter operation can be considered from two different perspectives. Within the ambient
category of $(R,R)$-bimodules, closure under direct summands might mean ``include any bimodule isomorphic to a direct summand of a bimodule in $\catstuff{S}'$'', and the resulting
category will be closed under isomorphism. We call these \emph{abstract Soergel bimodules}. Meanwhile, for an abstract additive category which is not embedded in a larger context,
direct summands don't make sense per se, and instead the standard approach is to take the Karoubi envelope $\Kar(\catstuff{S}')$, which formally adds as new objects the images of all
idempotents. Henceforth we make the distinction between these two approaches, and define $\catstuff{S} := \Kar(\catstuff{S}')$. After all, the Karoubi envelope description is the only one which technically applies to the diagrammatic category, where we do computations.

The problem with the definition of the max-degree submodule from the previous section (and with the definition of $\catstuff{A}_{\mathcal{J}}$) is that the max-degree submodule is
canonically an abstract Soergel bimodule, but it is not canonically an element of $\catstuff{S}$. We do not have a canonical idempotent whose image is the max-degree submodule. Our
goal is to be able to compute in the category $\catstuff{A}_{\mathcal{J}}$, and so we need to interpret the definition in terms of $\catstuff{S}$, where we can use diagrammatics and
idempotents. The phenomenon from \autoref{B1B1maxdegree} is completely general: inclusion maps from summands with shift $a(\mathcal{J})$ are canonical (up to a choice of basis)
whereas projection maps to such summands are not. To pin down an idempotent, we need a good way to pick out projection maps.

The \emph{local intersection pairing} or \emph{LIP} (of $B_z(k)$ at $B_x B_y$) is \begin{equation} \Hom^k(B_x B_y, B_z) \times
\Hom^{-k}(B_z,B_xB_y) \to \End^0(B_z) = \C \cdot \id_z \cong \C, \end{equation} a bilinear pairing which sends $(f,g) \mapsto f \circ g$. The rank of the LIP is precisely the
multiplicity of $B_z(k)$ as a summand of $B_x B_y$, and dual sets $\{\proj_i\} \subset \Hom^{k}(B_x B_y,B_z)$ and $\{\incl_i\} \subset \Hom^{-k}(B_z,B_xB_y)$ give (orthogonal) inclusion and projection maps (see \cite[Corollary 11.71]{ElMaThWi-soergel}). (We say dual sets rather than dual bases because the form might be degenerate, so $\{p_i\}$ may not be a basis for the hom space.)

Thus, one way to get a precise handle on $B_x \bullet B_y$, as an object of the Karoubi envelope $\catstuff{S}$, is to fix dual sets of size $\gamma^z_{x,y}$ for the specific local
intersection pairing
\begin{equation} \Hom^{a}(B_x B_y, B_z) \times \Hom^{-a}(B_z,B_xB_y) \to \C, \end{equation} for each $x,y,z$ in cell $\mathcal{J}$. If we fix this data once and for all, for each $x,y,z$, we get a well-defined operation $B_x \bullet B_y$ and a well-defined asymptotic Hecke category. Different choices of decompositions of $B_x B_y$ will give rise to different (but equivalent) constructions of the asymptotic Hecke category.

The asymptotic Hecke category is monoidal, and a monoidal category has certain structure maps. In a monoidal category one has an associator $\alpha_{x,y,z} \colon (A_x \star A_y)
\star A_z \to A_x \star (A_y \star A_z)$ for all objects $x, y, z$, an isomorphism satisfying the pentagon axiom. One also has the left unitor $lu_x \colon A_x \to \munit \star A_x$
and the right unitor $ru_x \colon A_x \to A_x \star \munit$, satisfying the triangle axiom. Lusztig's definition of the associator makes sense for abstract Soergel bimodules, but does
not translate so easily to the Karoubi envelope, as there may not be a single decomposition of $B_x B_y B_z$ which is compatible with the chosen idempotents for both $(B_x B_y) B_z$
and $B_x (B_y B_z)$. We discuss the associators and unitors in detail in the sequel  \cite{ElRoTu-ah-2}.

In the meantime, we discuss the practical matter of finding dual sets for the local intersection pairing. The first observation is that for $x, y, z$ in a two-sided cell
$\mathcal{J}$, when $a = a(\mathcal{J})$, the kernel of the LIP intersects $\Hom^{-a}(B_z,B_xB_y)$ to be precisely $\catstuff{I}^{-a}_{< \mathcal{J}}(B_z,B_xB_y)$, which follows
without much effort from the expanded Soergel hom formula. In particular, a basis of $\Hom_{\not< \mathcal{J}}^{-a}(B_z,B_x B_y)$ gives a basis of inclusion maps, for which we need
to find a dual basis of projection maps.

One solution to this problem was outlined in \cite[\S 5]{ElWi-relative-lefschetz}. One can define the \emph{multiplicity space} of $B_z(-k)$ in $B_x B_y$ to be the quotient of
$\Hom^k(B_z,B_x B_y)$ by the kernel of the LIF. The main result of \cite{ElWi-relative-lefschetz} is that the LIF on multiplicity spaces satisfies the Hodge-Riemann bilinear
relations with respect to a certain Lefschetz operator. We now describe the implications of this result for $a$-degree summands.

\begin{Definition} An element $\rho\in V$ is called \emph{dominant regular} if 
$\partial_{i}(\rho)>0$ for all $i \in S$. The \emph{standard dominant regular element} satisfies $\partial_i(\rho) = 1$ for all $i \in S$. \end{Definition}

On the bimodule $B_x B_y$ we have three actions of $R$: the usual right and left action, and the middle action $M$. For $b \in B_x$ and $b' \in B_y$ and $\rho \in V$ we write 
\begin{equation} M_{\rho}(b \ot b') := b \rho \ot b' = b \ot \rho b'. \end{equation} Then $M_{\rho}$ is a degree $+2$ bimodule endomorphism of $B_xB_y$, and pre-composition with $M_\rho$ gives a map from $\Hom^k(B_x B_y,B_z)$ to $\Hom^{k+2}(B_x B_y,B_z)$.  We write $M_{\rho}^a$ for the $a$-th power of $M_{\rho}$, which is the same as middle multiplication by $\rho^a \in R$. Then $M_{\rho}^a$ is an operator of degree $+2a$. The operator $M_{\rho}$ preserves $\catstuff{I}_{< \mathcal{J}}$, so it descends to the quotient category.

\begin{Definition} Fix $\rho \in V$. The \emph{Lefschetz form} on $\Hom^{-a}(B_z,B_x B_y)$ is the bilinear form $(-,-)_{\rho}$ defined by
\begin{equation} (f,g)_{\rho} := \DD(f) \circ M_{\rho}^a \circ g = \DD(M_{\rho}^a \circ f) \circ g, \end{equation}
where this composition is a local intersection pairing as above. \end{Definition}

The Lefschetz form is computed in $\catstuff{S}$, but also descends to a form on the quotient category $\catstuff{S}_{\not< \mathcal{J}}$.

\begin{Theorem} If $\rho$ is dominant regular, then the Lefschetz form $(-,-)_{\rho}$ on $\Hom_{\not< \mathcal{J}}^{-a}(B_z,B_x B_y)$ is positive definite. \end{Theorem}

\begin{proof}
This is \cite[Theorem 3.3, Proof of Proposition 5.5]{ElWi-relative-lefschetz}.
\end{proof}

In particular, that the Lefschetz form is non-degenerate implies that a dual basis for $\{\incl_i\}$ can be found within the span of $\{\DD(\incl_i) M_{\rho}^a\}$. One can solve for the dual basis explicitly by inverting the Gram matrix of the Lefschetz form. When the multiplicity space is one-dimensional, one can avoid linear algebra entirely.

\begin{Corollary} When $B_z(a)$ appears in $B_x B_y$ with multiplicity one, and $\incl$ is a possible inclusion map $B_z \to B_x B_y$ of degree $-a$, then $\proj = 
\lambda \DD(\incl) \circ M_{\rho}^a$ is the dual projection map, for some scalar $\lambda \in \R_{> 0}$. More precisely, $\lambda = (\incl,\incl)_{\rho}^{-1}$. \end{Corollary}

The multiplicity one setting arises frequently for important maps in the asymptotic Hecke category. Here are the two most important examples for this paper. For both of these, let $d$ be a distinguished involution, and $x$ an element in the same right cell. We work in $\catstuff{S}_{\not< \mathcal{J}}$.
\begin{enumerate}
\item We have that $B_d \bullet B_x \cong B_x(a)$. Inclusion and projection maps for this decomposition of $B_d B_x$ provide the unitor (and its inverse isomorphism, the counitor).
\item We have that $B_d(a)$ is a summand with multiplicity $1$ in $B_x \bullet B_{x^{-1}}$. Inclusion and projection maps for this summand, up to scalar, will provide the cups and caps for a rigid structure on the asymptotic Hecke category, as proven in \cite[Proposition 5.5]{ElWi-relative-lefschetz}.
\end{enumerate}

\subsection{Computing dimensions}\label{subsection:dimensions}

The goal of this paper is to compute categorical dimensions in the fusion category $\catstuff{A}_{\mathcal{H}}$ for various diagonal cells $\mathcal{H}$. In this section we state the computation we will do in the diagrammatic Hecke category. In  \cite{ElRoTu-ah-2} we justify that this computation actually agrees with the desired categorical dimension, which requires a more careful exploration of the associators and unitors.

Throughout this section we work in the quotient category $\catstuff{D}_{\not< \mathcal{J}}$, where $\mathcal{J}$ is the two-sided cell containing $\mathcal{H}$. We argue using the
one-dimensionality of certain hom spaces, which is only guaranteed in the quotient. Thus, the symmetries we discuss below are not guaranteed to hold in $\catstuff{D}$ itself, only in
the quotient.

Fix a diagonal cell $\mathcal{H} \subset \mathcal{J}$ with Duflo involution $d$. Fix some morphism $\incl_d = \incl_d^{d,d} \colon B_d \to B_d \ot B_d$ of degree $-a$, which we might
draw schematically as 
\begin{equation} \trisplitCW{d}{d}{d}{1}\;. \end{equation}

The one-dimensional morphism space $\Hom_{\not< \mathcal{J}}^{-a}(B_d, B_d B_d)$ has a number of symmetries: rotation by 120 degrees (via adjunction) and the horizontal flip. These symmetries must send $\incl_d$ to a scalar multiple of itself. Rotation has order 3 since all morphisms in $\catstuff{D}$ are cyclic (rotation by 360 preserves the morphism), and horizontal flip has order $2$ (it comes from rotation and duality, commuting involutions). Because our morphism is defined (or can be defined) over $\R$, the rotational symmetry acts trivially, meaning that we have
\begin{equation}
\begin{tikzpicture}[anchorbase,scale=1]
\begin{scope}[yscale=-1]
\draw[mor] (0,0) rectangle (1,0.5) node[black,pos=0.5]{$d$};
\draw[mor] (2,0) rectangle (3,0.5) node[black,pos=0.5]{$d$};
\draw[mor] (1,1.5) rectangle (2,2) node[black,pos=0.5]{$d$};
\draw[usual] (0.5,0.5) to (1.5,1) to (1.5,1.5);
\draw[usual] (2.5,0.5) to (1.5,1);
\varrotcircleCCWforflip{1.5}{1}{1}
\draw[usual] (0.5,0)to[out=270,in=0] (0,-0.5) to[out=180,in=270] (-0.5,0) to (-0.5,2);
\draw[usual] (1.5,2)to[out=90,in=180] (2.5,2.5) to[out=0,in=90] (3.5,2) to (3.5,0);
\end{scope}
\end{tikzpicture}
=
\trisplitCW{d}{d}{d}{1}. \end{equation} 
This rotation-invariance was the reason we were permitted to represent $\incl_d$ with a symbol that had rotational symmetry. Meanwhile, there is some sign $\ep_d \in \{\pm 1\}$ such that
\begin{equation} \trisplitCW{d}{d}{d}{1} = \ep_d \trisplitCCW{d}{d}{d}{1} \; . \end{equation}
%
The reason that we include small oriented circles in our symbol for $\incl_d$ is because the morphism might not be flip-symmetric, so we should not use a flip-symmetric symbol. When $\ep_d = 1$ we can remove the circles from our notation, as flip-symmetry does hold.

Using adjunction or duality one gets two potentially different morphisms $B_dB_d \to B_d$,
\begin{equation}
\trimergeCW{d}{d}{d}{1}
\coloneqq
\begin{tikzpicture}[anchorbase,scale=1]
\begin{scope}[yscale=-1]
\draw[mor] (0,0) rectangle (1,0.5) node[black,pos=0.5]{$d$};
\draw[mor] (2,0) rectangle (3,0.5) node[black,pos=0.5]{$d$};
\draw[mor] (1,1.5) rectangle (2,2) node[black,pos=0.5]{$d$};
\draw[usual] (0.5,0.5) to (1.5,1) to (1.5,1.5);
\draw[usual] (2.5,0.5) to (1.5,1);
\varrotcircleCCWforflip{1.5}{1}{1}
\draw[usual] (2.5,0)to[out=270,in=180] (3,-0.5) to[out=360,in=270] (3.5,0) to (3.5,2);
\end{scope}
\end{tikzpicture}\;.
\end{equation}
and
\begin{equation}
\trimergeCCW{d}{d}{d}{1} 
\coloneqq
\DD\left(
\trisplitCW{d}{d}{d}{1}
\right) \;.
\end{equation}
Because the horizontal flip is just the composition of duality and rotation by 180 degrees, we have
\begin{equation} 
\trimergeCW{d}{d}{d}{1} = \ep_d \trimergeCCW{d}{d}{d}{1} \;. \end{equation}
When $\ep_d = 1$ we can safely remove the circles from our diagrammatic notation.

\begin{Remark} When $\len(d) = a(d)$, one can normalize the map $\incl_d$ so that $\onetensor$ is preserved. When $\len(d) \ne a(d)$, the morphism $\incl_d$ is not the right degree to preserve the one-tensor, and there is not an obvious normalization. Similarly, the map $\incl_x$ below typically does not have the correct degree to preserve the one-tensor. \end{Remark}

Let $x$ be another element in $\mathcal{J}$ in the same right cell as $d$. Fix some morphism $\incl_x = \incl_d^{x,x^{-1}} \colon B_d \to B_x \ot B_{x^{-1}}$ of degree $-a$, which we might draw schematically as
\begin{equation}
\incl_x = \trisplitCW{x}{x^{-1}}{d}{1} \; .
\end{equation}
There is no longer any question of rotational symmetry beyond cyclicity. Using adjunction, $\incl_x$ leads to a total of six morphisms of note, the other five being
\begin{equation}\label{6rotates}
lu_x := 
\trisplitCW{d}{x}{x}{1}, 
ru_{x^{-1}} := \trisplitCW{x^{-1}}{d}{x^{-1}}{1}, \end{equation}
\begin{equation} \trimergeCW{x}{x^{-1}}{d}{1}, \trimergeCW{d}{x}{x}{1}, \trimergeCW{x^{-1}}{d}{x^{-1}}{1}. \end{equation}
Using duality one obtains an additional six morphisms, drawn with counter-clockwise circles.	

Note that $\incl_x$ lives in a one-dimensional hom space in the quotient $\catstuff{S}_{\not< \mathcal{J}}$, since $\gamma_{x,x^{-1}}^d = 1$. Once again, the horizontal flip acts on this one-dimensional space by some involution. (This is true regardless of whether or not $x$ is an involution, as the horizontal flip sends $x$ to $x^{-1}$. The typical case when $x$ is also an involution does not produce additional symmetry.) Thus, there is some sign $\ep_x \in \{\pm 1\}$ such that 
\begin{equation} \label{signflippinginclx} \trisplitCW{x}{x^{-1}}{d}{1} = \ep_x \trisplitCCW{x}{x^{-1}}{d}{1}. \end{equation}

\begin{Lemma} \label{Lem:sameep} For all $x$ in the same diagonal cell as $d$, we have $\ep_x = \ep_d$. \end{Lemma}

This Lemma is proven in the sequel \cite{ElRoTu-ah-2}, though see \autoref{rem:epwepdspecialcase} for an important special case. Because of this Lemma, it turns out that the sign $\ep_x$ is not relevant to the computation of categorical dimensions. We wished to discuss it anyway, because it explains why we included oriented circles in our diagrammatic notation for $\incl_x$.

\begin{Theorem} \label{thm:whatisdim} Fix a dominant regular element $\rho$. Choose $x\in W$ within a diagonal cell $\mathcal{H}$, and let $d \in \mathcal{H}$ be the Duflo involution. Fix $\incl_x \colon B_d \to B_x \ot B_{x^{-1}}$ of degree $-a$, and define $lu_x$ as in \eqref{6rotates}. Compute $\lambda_x := \ep_x (lu_x,lu_x)_{\rho} \in \R_{\ge 0}$. Compute $\mu_x := \ep_x (\incl_x,\incl_x)_{\rho} \in \R_{\ge 0}$. Then $\lambda_x$ and $\mu_x$ are nonzero if and only if $\incl_x \notin \catstuff{I}_{< \mathcal{J}}$, in which case the ratio $\mu_x/\lambda_x$ is independent of the choice of $\rho$, and is equal to the categorical dimension of $A_x$ in $\catstuff{A}_{\mathcal{H}}$. \end{Theorem}

This Theorem is proven in the sequel  \cite{ElRoTu-ah-2}.

In diagrams, we define $\lambda_x$ and $\mu_x$ so that
\begin{gather}
\begin{tikzpicture}[anchorbase,scale=0.8]
\draw[mor] (0,0) rectangle (1,0.5) node[black,pos=0.5]{$d$};
\draw[mor] (2,0) rectangle (3,0.5) node[black,pos=0.5]{$x$};
\draw[mor] (1,1.5) rectangle (2,2) node[black,pos=0.5]{$x$};
\draw[mor] (1,-1.5) rectangle (2,-1) node[black,pos=0.5]{$x$};
\draw[usual] (0.5,0.5) to (1.5,1) to (1.5,1.5);
\draw[usual] (2.5,0.5) to (1.5,1);
\draw[usual] (0.5,0) to (1.5,-0.5) to (1.5,-1);
\draw[usual] (2.5,0) to (1.5,-0.5);
\rotcircleCCW{1.5}{1}{1}
\begin{scope}[yscale=-1]
\varrotcircleCCWforflip{1.5}{0.5}{1}
\end{scope}
\node[scale=0.7] at (1.5,0.25) {$\rho^{a}$};
\end{tikzpicture}
=
\ep_x \lambda_x \cdot
\begin{tikzpicture}[anchorbase,scale=1]
\draw[mor] (0,0) rectangle (1,0.5) node[black,pos=0.5]{$x$};
\end{tikzpicture},
\qquad
\begin{tikzpicture}[anchorbase,scale=0.8]
\draw[mor] (0,0) rectangle (1,0.5) node[black,pos=0.5]{$x$};
\draw[mor] (2,0) rectangle (3,0.5) node[black,pos=0.5]{$x^{-1}$};
\draw[mor] (1,1.5) rectangle (2,2) node[black,pos=0.5]{$d$};
\draw[mor] (1,-1.5) rectangle (2,-1) node[black,pos=0.5]{$d$};
\draw[usual] (0.5,0.5) to (1.5,1) to (1.5,1.5);
\draw[usual] (2.5,0.5) to (1.5,1);
\draw[usual] (0.5,0) to (1.5,-0.5) to (1.5,-1);
\draw[usual] (2.5,0) to (1.5,-0.5);
\rotcircleCCW{1.5}{1}{1}
\begin{scope}[yscale=-1]
\varrotcircleCCWforflip{1.5}{0.5}{1}
\end{scope}
\node[scale=0.7] at (1.5,0.25) {$\rho^{a}$};
\end{tikzpicture}
=
\ep_x \mu_x \cdot
\begin{tikzpicture}[anchorbase,scale=1]
\draw[mor] (0,0) rectangle (1,0.5) node[black,pos=0.5]{$d$};
\end{tikzpicture}
.
\end{gather}

Using \eqref{signflippinginclx} we can change orientation and replace this calculation with
\begin{gather}
\begin{tikzpicture}[anchorbase,scale=0.8]
\draw[mor] (0,0) rectangle (1,0.5) node[black,pos=0.5]{$d$};
\draw[mor] (2,0) rectangle (3,0.5) node[black,pos=0.5]{$x$};
\draw[mor] (1,1.5) rectangle (2,2) node[black,pos=0.5]{$x$};
\draw[mor] (1,-1.5) rectangle (2,-1) node[black,pos=0.5]{$x$};
\draw[usual] (0.5,0.5) to (1.5,1) to (1.5,1.5);
\draw[usual] (2.5,0.5) to (1.5,1);
\draw[usual] (0.5,0) to (1.5,-0.5) to (1.5,-1);
\draw[usual] (2.5,0) to (1.5,-0.5);
\varrotcircleCW{1.5}{1}{1}
\begin{scope}[yscale=-1]
\varrotcircleCCWforflip{1.5}{0.5}{1}
\end{scope}
\node[scale=0.7] at (1.5,0.25) {$\rho^{a}$};
\end{tikzpicture}
=
\lambda_x \cdot
\begin{tikzpicture}[anchorbase,scale=1]
\draw[mor] (0,0) rectangle (1,0.5) node[black,pos=0.5]{$x$};
\end{tikzpicture},
\qquad
\begin{tikzpicture}[anchorbase,scale=0.8]
\draw[mor] (0,0) rectangle (1,0.5) node[black,pos=0.5]{$x$};
\draw[mor] (2,0) rectangle (3,0.5) node[black,pos=0.5]{$x^{-1}$};
\draw[mor] (1,1.5) rectangle (2,2) node[black,pos=0.5]{$d$};
\draw[mor] (1,-1.5) rectangle (2,-1) node[black,pos=0.5]{$d$};
\draw[usual] (0.5,0.5) to (1.5,1) to (1.5,1.5);
\draw[usual] (2.5,0.5) to (1.5,1);
\draw[usual] (0.5,0) to (1.5,-0.5) to (1.5,-1);
\draw[usual] (2.5,0) to (1.5,-0.5);
\varrotcircleCW{1.5}{1}{1}
\begin{scope}[yscale=-1]
\varrotcircleCCWforflip{1.5}{0.5}{1}
\end{scope}
\node[scale=0.7] at (1.5,0.25) {$\rho^{a}$};
\end{tikzpicture}
=
\mu_x \cdot
\begin{tikzpicture}[anchorbase,scale=1]
\draw[mor] (0,0) rectangle (1,0.5) node[black,pos=0.5]{$d$};
\end{tikzpicture}
.
\end{gather}

\begin{Example} Clearly $\mu_d = \lambda_d$, so the categorical dimension of $A_d$ is $1$. This is as it should be, since the categorical dimension of the monoidal unit is always $1$.
\end{Example}

Note that all diagonal cells in type $A$ have size one and consist only of their Duflo involution, so there is no need to compute anything. However, because it is clarifying, we give some examples where we compute $\lambda_d$ in type $A$.

\begin{Example} As in \autoref{S4middle}, let $W = S_4$ and consider the middle cell. Let $d_1 = s_1 s_3$ and $d_2 = s_2 s_1 s_3 s_2$. Note that $B_{s_1 s_3}$ and $B_{s_2 s_1 s_3 s_2}$ are both Bott--Samelson bimodules, so no idempotents are necessary. We can choose maps $\incl_{d_1}$ and $\incl_{d_2}$ as follows.
\begin{equation}\incl_{d_1}=
\begin{tikzpicture}[anchorbase,scale=1]
\begin{scope}[yscale=-1]
\draw[soergelone] (0.5,0.5) to (1.5,1) to (1.5,1.5);
\draw[soergelone] (2.5,0.5) to (1.5,1);
\begin{scope}[shift={(0.25,0)}]
\draw[soergelthree] (0.5,0.5) to (1.5,1) to (1.5,1.5);
\draw[soergelthree] (2.5,0.5) to (1.5,1);
\end{scope}
\end{scope}
\end{tikzpicture},\quad
\incl_{d_{2}}=
\begin{tikzpicture}[anchorbase,scale=1]
\begin{scope}[yscale=-1]
\draw[soergelone] (0.5,0.5) to (1.5,1) to (1.5,1.5);
\draw[soergelone] (2.5,0.5) to (1.5,1);
\begin{scope}[shift={(0.25,0)}]
\draw[soergelthree] (0.5,0.5) to (1.5,1) to (1.5,1.5);
\draw[soergelthree] (2.5,0.5) to (1.5,1);
\end{scope}
\draw[soergeltwo] (1,0.5) to (1.625,0.75) to (2.25,0.5);
\draw[soergeltwo] (0.25,0.5) to (1.25,1) to (1.25,1.5);
\draw[soergeltwo] (3,0.5) to (2,1) to (2,1.5);
\end{scope}
\end{tikzpicture}\;.  \end{equation}
In particular, $\ep_{d_1} = 1 = \ep_{d_2}$. For any $f \in R$ we have
\begin{equation}
\begin{tikzpicture}[anchorbase,scale=1]
\draw[soergelone] (0.5,0) to[out=270,in=180] (1,-0.5) to[out=0,in=270] (1.5,0) to[out=90,in=0] (1,0.5) to[out=180,in=90] (0.5,0);
\begin{scope}[shift={(0.25,0)}]
\draw[soergelthree] (0.5,0) to[out=270,in=180] (1,-0.5) to[out=0,in=270] (1.5,0) to[out=90,in=0] (1,0.5) to[out=180,in=90] (0.5,0);
\end{scope}
\draw[soergelone] (1,-1) to (1,-0.5);
\draw[soergelone] (1,0.5) to (1,1);
\draw[soergelthree] (1.25,-1) to (1.25,-0.5);
\draw[soergelthree] (1.25,0.5) to (1.25,1);
\node at (1.125,0) {$f$};
\end{tikzpicture}
=
\partial_1\partial_3(f)
\begin{tikzpicture}[anchorbase,scale=1]
\diagrammaticmorphism{0}{0}{10}{1,3}\;.
\end{tikzpicture}
\end{equation}
\begin{equation} \label{itsd2thistime}
\begin{tikzpicture}[anchorbase,scale=1]
\draw[soergelone] (0.5,0) to[out=270,in=180] (1,-0.5) to[out=0,in=270] (1.5,0) to[out=90,in=0] (1,0.5) to[out=180,in=90] (0.5,0);
\begin{scope}[shift={(0.25,0)}]
\draw[soergelthree] (0.5,0) to[out=270,in=180] (1,-0.5) to[out=0,in=270] (1.5,0) to[out=90,in=0] (1,0.5) to[out=180,in=90] (0.5,0);
\end{scope}
\draw[soergelone] (1,-1) to (1,-0.5);
\draw[soergelone] (1,0.5) to (1,1);
\draw[soergelthree] (1.25,-1) to (1.25,-0.5);
\draw[soergelthree] (1.25,0.5) to (1.25,1);
\node at (1.125,0) {$f$};
\draw[soergeltwo] (1.125,0) circle (0.25) ;
\draw[soergeltwo] (0.75,-1) to (0.75, -0.75) to (0.25,-0.25) to (0.25,0.25) to (0.75,0.75) to (0.75,1);
\draw[soergeltwo] (1.5,-1) to (1.5, -0.75) to (2,-0.25) to (2,0.25) to (1.5,0.75) to (1.5,1);
\end{tikzpicture}
=
\partial_2\partial_1\partial_3(\alpha_2\partial_2(f))
\begin{tikzpicture}[anchorbase,scale=1]
\draw[soergeltwo,markedtwo=1] (-0.2,0) to (-0.2,0.25);
\draw[soergeltwo,markedtwo=1] (-0.2,1) to (-0.2,0.75);
\diagrammaticmorphism{0}{0}{10}{1,3,2}
\end{tikzpicture}
+
s_2\partial_1\partial_3(\alpha_2\partial_2(f))
\begin{tikzpicture}[anchorbase,scale=1]
\diagrammaticmorphism{0}{0}{10}{2,1,3,2}
\end{tikzpicture}	\;.
\end{equation}
When $f$ has degree $4 = 2a$, like $\rho^a$, the first coefficient $\partial_2 \partial_1 \partial_3(\alpha_2 \partial_2(f))$ on the right-hand side of \eqref{itsd2thistime} vanishes. We knew this had to happen, since the result is a degree zero endomorphism of $B_{d_2}$, so it must be a multiple of the identity, whereas the first diagram lives in lower terms (for the categorical filtration $\catstuff{I}_{< d_2}$, not the monoidal filtration $\catstuff{I}_{< \mathcal{J}}$).

When $\rho$ is a regular standard element we have
\begin{equation} \partial_1 \partial_3(\rho^2) = 2, \quad s_2 \partial_1 \partial_3(\alpha_2 \partial_2(\rho^2)) = -6. \end{equation}The proof is complete.
\end{Example}

\begin{Remark} The fact that $\mu_x/\lambda_x$ is independent of the choice of $\rho$ is obvious when you know that the result is the categorical dimension, but is highly non-obvious from the computational perspective. Consider the setting of \autoref{Ex:F4}, where indeed $d$ is the distinguished involution in a diagonal cell which also contains $w$. Though we will not prove it until later in this section, we claim that the ratio which is $\partial_d(\ptr_{d,(2,3,4,1)}(\rho^4))$ divided by $\partial_d(\rho^4)$ is the categorical dimension of $w$. In this case, the independence of $\rho$ followed from the surprising fact that $\partial_d$ and $\partial_d\circ \ptr_{d,(2,3,4,1)}$ were colinear on polynomials in this degree. In general, independence of $\rho$ has the same flavor as in this example. \end{Remark}

\subsection{Reduction to the partial trace}\label{subsection:simplifications}

\begin{Notation} For the rest of this paper we work only with clockwise-oriented trivalent vertices $\incl_d$ and $\incl_x$, and we omit the circle henceforth. \end{Notation}

\begin{Notation} 
Fix a regular dominant element $\rho$. In this section we fix a diagonal cell $\mathcal{H}$ with Duflo involution $d$, inside a two-sided cell $\mathcal{J}$, and we fix another
element $w \in \mathcal{H}$. We make the following simplying assumption:
\begin{equation} \label{assumeme}\text{There exists } t \in W \text{ such that } w = dt \text{ and } \len(w) = \len(d) + \len(t). \end{equation}
Equivalently, there is some reduced expression $\un{w}$ for $w$ obtained by concatenating reduced expressions for $d$ and $t$, i.e. $\un{w} = \un{d} \un{t}$. Equivalently, $w$ is greater than $d$ in the weak right Bruhat order.

We let $r = \len(t)$. If $\un{t} = (t_1, t_2, \ldots, t_r)$, then let $\un{t}_{\le k} = (t_1, t_2, \ldots, t_k)$ for $k \le r$. Let $x_k$ be the element with reduced expression $\un{d} \un{t}_{\le k}$. \end{Notation}

\begin{Lemma} The elements $x_k$ with reduced expression $d \un{t}_{\le k}$ are all in the same right cell as $d$. \end{Lemma}

\begin{proof} Clearly $B_{x_k}$ is a summand of $B_{x_{k-1}} \ot B_{t_k}$, so the right cells of $x_k$ form a (weakly) decreasing sequence. Since $x_0 = d$ and $x_r = w$ live in the same right cell, then all $x_k$ live in the same right cell. \end{proof}

\begin{Example} \label{longesteltworks} When $d = w_I$ is the longest element of a finite parabolic subgroup generated by $I \subset S$, then $B_d \cong R \ot_{R^I} R(\ell(d))$. Any element $w \le_R d$ also has $I$ in its left descent set (there are many proofs of this fact, using ordinary KL combinatorics, or using singular Soergel bimodules). Any element with $I$ in its left descent set has a reduced expression beginning with $w_I$. Thus, the condition \eqref{assumeme} holds for any $w$ in the same diagonal cell as $w_I$. \end{Example}

\begin{Example} If $d$ is less than $w$ in the weak right Bruhat order and they live in the same diagonal cell $\mathcal{H}$, then $w_0 d$ is greater than $w_0 w$ in the weak right Bruhat order and they live in the same diagonal cell $\mathcal{H}^{\prime}$. Thus, for every nontrivial example of \eqref{assumeme}, the $w_0$-dual is a counterexample. However, it is possible that there is a different diagonal cell in the same two-sided cell as $\mathcal{H}^{\prime}$ for which \eqref{assumeme} holds for all $w$. \end{Example} 

What is special about \eqref{assumeme} is that one can construct the morphism $\incl_w$ directly from $\incl_d$, as follows. Let us fix a map $\incl_d \colon B_d \ot B_d \to B_d$ of degree $-a$ for which $\lambda_d$ (defined in \autoref{thm:whatisdim}) is nonzero. Using $\incl_d$, for each $0 \le k \le r$ we can construct a morphism $\incl_k \colon B_d \to B_{x_k} B_{x_k^{-1}}$ of degree $-a$, namely 
\begin{equation} \incl_k :=  \begin{tikzpicture}[anchorbase,scale=1,yscale=-1]
\draw[mor] (0,0) rectangle (1.1,0.5) node[black,pos=0.5]{\raisebox{-0.06cm}{$x_k$}};
\draw[mor] (1.7,0) rectangle (2.8,0.5) node[black,pos=0.5]{\raisebox{-0.06cm}{$x_k^{-1}$}};
\draw[usual] (1.8,0.5)to[out=90,in=0] (1.4,0.9) to[out=180,in=90] (1,0.5);
\node[scale=0.8] at (1.4,1.5) {$t_{\leq k}$};
\draw[usual] (2,0.5)to[out=90,in=0] (1.4,1.1) to[out=180,in=90] (0.8,0.5);
\draw[usual] (2.2,0.5)to[out=90,in=0] (1.4,1.3) to[out=180,in=90] (0.6,0.5);
\draw[mor] (0,0.5) rectangle (0.4,1) node[black,pos=0.5]{\raisebox{-0.06cm}{$d$}};
\draw[mor] (2.4,0.5) rectangle (2.8,1) node[black,pos=0.5]{\raisebox{-0.06cm}{$d$}};
\draw[usual] (0.2,1) to (1.4,2) to (1.4,2.5);
\draw[usual] (2.6,1) to (1.4,2) to (1.4,2.5);
\draw[mor] (1.2,2.5) rectangle (1.6,3) node[black,pos=0.5]{\raisebox{-0.06cm}{$d$}};
\end{tikzpicture}. \end{equation} 
The calculations done below will justify that $(\incl_k,\incl_k)_{\rho}$ is nonzero, and hence $\incl_k$ is a valid choice of trivalent vertex. We let $\incl_w := \incl_r$.

\begin{Remark} \label{rem:epwepdspecialcase} In \autoref{Lem:sameep} we noted that $\ep_w = \ep_d$, though we postponed the proof in general. We note here that, when \eqref{assumeme} holds, the statement that $\ep_w = \ep_d$ follows easily from the construction of $\incl_w$ above, because the bottom half of the diagram is invariant under the horizontal flip. \end{Remark}

\begin{Definition}\label{D:partialtrace} Let $f \in R$ have degree $2a$ (we only need the case $f = \rho^a$), and let $x$ be in the same right cell as $d$. Write $\lambda_k(f)$ for the scalar appearing in
\begin{equation}
\begin{tikzpicture}[anchorbase,scale=1]
\draw[mor] (1,-1) rectangle (2,-0.5) node[black,pos=0.5]{$x_k$};
\draw[mor] (0,0) rectangle (1,0.5) node[black,pos=0.5]{$d$};
\draw[mor] (2,0) rectangle (3,0.5) node[black,pos=0.5]{$x_k$};
\draw[mor] (1,1) rectangle (2,1.5) node[black,pos=0.5]{$x_k$};
\node[scale=0.7] at (1.5,0.25) {$f$};
\draw[usual] (0.75,0.5) to (1.5,0.75) to (1.5,1);
\draw[usual] (2.25,0.5) to (1.5,0.75);
\draw[usual] (0.75,0) to (1.5,-0.25) to (1.5,-0.5);
\draw[usual] (2.25,0) to (1.5,-0.25);
\end{tikzpicture} \; = \; \lambda_k(f) \cdot \id_{x_k}.
\end{equation}
Note that $\lambda_k(f) \in \R$.
\end{Definition}

\begin{Theorem} \label{thm:lambdadlambdaw} Under the assumption \eqref{assumeme} we have $\lambda_d(f) = \lambda_w(f)$ for all $f \in R$ of degree $2a$. \end{Theorem}

\begin{proof} We claim that $\lambda_d(f) = \lambda_k(f)$ for all such $f$, and for all $0 \le k \le r$. The result is evident for $k=0$.

First we note that this computation can be safely performed in $\catstuff{S}_{\not< \mathcal{J}}$. This is because $\catstuff{I}_{< \mathcal{J}}$ intersects the space of degree zero endomorphisms of $B_{x_k}$ trivially. In particular, we can ignore any morphisms which factor through lower cells as we do the computation.

We now realize $B_{x_{k}}$ as the image of an idempotent in $B_{x_{k-1}} \ot B_{t_k}$. This idempotent is explicitly described in \autoref{Thm:inductiveclasp}. Thus, we have
\begin{equation}
\begin{tikzpicture}[anchorbase,scale=1]
\draw[mor] (0,0) rectangle (1,0.5) node[black,pos=0.5]{$d$};
\draw[mor] (2,0) rectangle (3,0.5) node[black,pos=0.5]{$x_k$};
\draw[mor] (1,1.5) rectangle (2,2) node[black,pos=0.5]{$x_k$};
\draw[usual] (0.5,0.5) to (1.5,1) to (1.5,1.5);
\draw[usual] (2.5,0.5) to (1.5,1);
\begin{scope}[yscale=-1,shift={(0,-0.5)}]
\draw[mor] (1,1.5) rectangle (2,2) node[black,pos=0.5]{$x_k$};
\draw[usual] (0.5,0.5) to (1.5,1) to (1.5,1.5);
\draw[usual] (2.5,0.5) to (1.5,1);
\end{scope}
\node[scale=0.8] at (1.5,0.25) {$f$};
\end{tikzpicture} 
=
\begin{tikzpicture}[anchorbase,scale=1]
\draw[mor] (0,0) rectangle (1,0.5) node[black,pos=0.5]{$d$};
\draw[mor] (2,0) rectangle (3,0.5) node[black,pos=0.5]{$x_{k-1}$};
\draw[mor] (1,1.5) rectangle (2,2) node[black,pos=0.5]{$x_{k-1}$};
\draw[mor] (1,2) rectangle (2.2,2.5) node[black,pos=0.5]{$x_{k}$};
\draw[usual] (2.1, 2) to (3.2,1) to (3.2,0);
\draw[usual] (0.5,0.5) to (1.5,1) to (1.5,1.5);
\draw[usual] (2.5,0.5) to (1.5,1);
\begin{scope}[yscale=-1,shift={(0,0.5)}]
\draw[mor] (0,0) rectangle (1,0.5) node[black,pos=0.5]{$d$};
\draw[mor] (2,0) rectangle (3,0.5) node[black,pos=0.5]{$x_{k-1}$};
\draw[mor] (1,1.5) rectangle (2,2) node[black,pos=0.5]{$x_{k-1}$};
\draw[mor] (1,2) rectangle (2.2,2.5) node[black,pos=0.5]{$x_{k}$};
\draw[usual] (0.5,0.5) to (1.5,1) to (1.5,1.5);
\draw[usual] (2.5,0.5) to (1.5,1);
\draw[usual] (2.1, 2) to (3.2,1) to (3.2,-0.25);
\end{scope}
\node[scale=0.8] at (1.5,-0.25) {$f$};
\draw[usual](0.5,0) to (0.5,-0.5);
\draw[mor] (2,-0.5) rectangle (3.4,0) node[black,pos=0.5]{$x_k$};
\end{tikzpicture}
=
\begin{tikzpicture}[anchorbase,scale=1]
\draw[mor] (0,0) rectangle (1,0.5) node[black,pos=0.5]{$d$};
\draw[mor] (2,0) rectangle (3,0.5) node[black,pos=0.5]{$x_{k-1}$};
\draw[mor] (1,1.5) rectangle (2,2) node[black,pos=0.5]{$x_{k-1}$};
\draw[mor] (1,2) rectangle (2.2,2.5) node[black,pos=0.5]{$x_{k}$};
\draw[usual] (2.1, 2) to (3.2,1) to (3.2,0);
\draw[usual] (0.5,0.5) to (1.5,1) to (1.5,1.5);
\draw[usual] (2.5,0.5) to (1.5,1);
\begin{scope}[yscale=-1,shift={(0,-0.5)}]
\draw[mor] (1,1.5) rectangle (2,2) node[black,pos=0.5]{$x_{k-1}$};
\draw[mor] (1,2) rectangle (2.2,2.5) node[black,pos=0.5]{$x_{k}$};
\draw[usual] (0.5,0.5) to (1.5,1) to (1.5,1.5);
\draw[usual] (2.5,0.5) to (1.5,1);
\draw[usual] (2.1, 2) to (3.2,1) to (3.2,-0.25);
\end{scope}
\node[scale=0.8] at (1.5,0.25) {$f$};
\end{tikzpicture}
+
\sum
\begin{tikzpicture}[anchorbase,scale=1]
\draw[mor] (0,-0.5) rectangle (1,0) node[black,pos=0.5]{$d$};
\draw[mor] (2,0) rectangle (3,0.5) node[black,pos=0.5]{$x_{k-1}$};
\draw[mor] (1,1.5) rectangle (2,2) node[black,pos=0.5]{$x_{k-1}$};
\draw[mor] (1,2) rectangle (2.2,2.5) node[black,pos=0.5]{$x_{k}$};
\draw[usual] (2.1, 2) to (3.2,1) to (3.2,0);
\draw[usual] (0.5,0) to (1.5,1) to (1.5,1.5);
\draw[usual] (2.5,0.5) to (1.5,1);
\begin{scope}[yscale=-1,shift={(0,0.5)}]
\draw[mor] (2,0) rectangle (3,0.5) node[black,pos=0.5]{$x_{k-1}$};
\draw[mor] (1,1.5) rectangle (2,2) node[black,pos=0.5]{$x_{k-1}$};
\draw[mor] (1,2) rectangle (2.2,2.5) node[black,pos=0.5]{$x_{k}$};
\draw[usual] (0.5,0) to (1.5,1) to (1.5,1.5);
\draw[usual] (2.5,0.5) to (1.5,1);
\draw[usual] (2.1, 2) to (3.2,1) to (3.2,-0.25);
\end{scope}
\node[scale=0.8] at (1.5,-0.25) {$f$};
\draw[mor] (2,-0.5) rectangle (3.4,0) node[black,pos=0.5]{$q_j$};
\end{tikzpicture}\;.
\end{equation}
The idempotent $q_j$ is described in \eqref{defqy}, and it is a composition of degree zero maps $B_{x_{k-1}} \ot B_{t_k} \to B_y \to B_{x_{k-1}} \ot B_{t_k}$, where $y \muless x_{k-1}$ and $yt_k < y$. In particular, the diagram containing $q_j$ factors through a morphism $B_{x_k} \to B_d \ot B_y$ of degree $-a$.

We claim that $\Hom_{\not< \mathcal{J}}(B_{x_k},B_d B_y)=0$. Since $B_y$ is a summand of $B_{x_{k-1}} \ot B_{t_k}$, we see that $y$ is either in $\mathcal{J}$ or in a strictly lower
two-sided cell. In the latter case the statement is obvious. If $y \in \mathcal{J}$ then the dimension of this morphism space is $\gamma_{d,y}^{x_k}$, which is zero unless $y = x_k$.
But $y < x_{k-1} < x_k$, proving the claim.

Thus, all terms with idempotents $q_j$ die, and we have
\begin{equation} 
\begin{tikzpicture}[anchorbase,scale=1]
\draw[mor] (0,0) rectangle (1,0.5) node[black,pos=0.5]{$d$};
\draw[mor] (2,0) rectangle (3,0.5) node[black,pos=0.5]{$x_k$};
\draw[mor] (1,1.5) rectangle (2,2) node[black,pos=0.5]{$x_k$};
\draw[usual] (0.5,0.5) to (1.5,1) to (1.5,1.5);
\draw[usual] (2.5,0.5) to (1.5,1);
\begin{scope}[yscale=-1,shift={(0,-0.5)}]
\draw[mor] (1,1.5) rectangle (2,2) node[black,pos=0.5]{$x_k$};
\draw[usual] (0.5,0.5) to (1.5,1) to (1.5,1.5);
\draw[usual] (2.5,0.5) to (1.5,1);
\end{scope}
\node[scale=0.8] at (1.5,0.25) {$f$};
\end{tikzpicture} 
=
\begin{tikzpicture}[anchorbase,scale=1]
\draw[mor] (0,0) rectangle (1,0.5) node[black,pos=0.5]{$d$};
\draw[mor] (2,0) rectangle (3,0.5) node[black,pos=0.5]{$x_{k-1}$};
\draw[mor] (1,1.5) rectangle (2,2) node[black,pos=0.5]{$x_{k-1}$};
\draw[mor] (1,2) rectangle (2.2,2.5) node[black,pos=0.5]{$x_{k}$};
\draw[usual] (2.1, 2) to (3.2,1) to (3.2,0);
\draw[usual] (0.5,0.5) to (1.5,1) to (1.5,1.5);
\draw[usual] (2.5,0.5) to (1.5,1);
\begin{scope}[yscale=-1,shift={(0,-0.5)}]
\draw[mor] (1,1.5) rectangle (2,2) node[black,pos=0.5]{$x_{k-1}$};
\draw[mor] (1,2) rectangle (2.2,2.5) node[black,pos=0.5]{$x_{k}$};
\draw[usual] (0.5,0.5) to (1.5,1) to (1.5,1.5);
\draw[usual] (2.5,0.5) to (1.5,1);
\draw[usual] (2.1, 2) to (3.2,1) to (3.2,-0.25);
\end{scope}
\node[scale=0.8] at (1.5,0.25) {$f$};
\end{tikzpicture}
. \end{equation}
This easily yields the inductive step $\lambda_{k-1}(f) = \lambda_k(f)$.
\end{proof}

\begin{Theorem} \label{thm:itsapartialtracedummy} Under the assumption \eqref{assumeme}, for any element $f \in R$ of degree $2a$ we have
\begin{equation}
\begin{tikzpicture}[anchorbase,scale=1]
\draw[mor] (1,-1) rectangle (2,-0.5) node[black,pos=0.5]{$d$};
\draw[mor] (0,0) rectangle (1,0.5) node[black,pos=0.5]{$w$};
\draw[mor] (2,0) rectangle (3,0.5) node[black,pos=0.5]{$w^{-1}$};
\draw[mor] (1,1) rectangle (2,1.5) node[black,pos=0.5]{$d$};
\node[scale=0.7] at (1.5,0.25) {$f$};
\draw[usual] (0.75,0.5) to (1.5,0.75) to (1.5,1);
\draw[usual] (2.25,0.5) to (1.5,0.75);
\draw[usual] (0.75,0) to (1.5,-0.25) to (1.5,-0.5);
\draw[usual] (2.25,0) to (1.5,-0.25);
\end{tikzpicture}
\; = \;
\begin{tikzpicture}[anchorbase,scale=1]
\draw[mor] (1,-1) rectangle (2,-0.5) node[black,pos=0.5]{$d$};
\draw[mor] (0,0) rectangle (1,0.5) node[black,pos=0.5]{$d$};
\draw[mor] (2,0) rectangle (3,0.5) node[black,pos=0.5]{$d$};
\draw[mor] (1,1) rectangle (2,1.5) node[black,pos=0.5]{$d$};
\node[scale=0.7] at (1.5,0.25) {$\ptr_{d,\un{t}}(f)$};
\draw[usual] (0.75,0.5) to (1.5,0.75) to (1.5,1);
\draw[usual] (2.25,0.5) to (1.5,0.75);
\draw[usual] (0.75,0) to (1.5,-0.25) to (1.5,-0.5);
\draw[usual] (2.25,0) to (1.5,-0.25);
\end{tikzpicture}
.
\end{equation}
\end{Theorem}

\begin{proof} The proof of this theorem is very similar to the previous theorem. We prove it by proving the analogous statement for all $x_k$, i.e. that
\begin{equation} \label{analyzeme2}
\begin{tikzpicture}[anchorbase,scale=1]
\draw[mor] (1,-1) rectangle (2,-0.5) node[black,pos=0.5]{$d$};
\draw[mor] (0,0) rectangle (1,0.5) node[black,pos=0.5]{$x_k$};
\draw[mor] (2,0) rectangle (3,0.5) node[black,pos=0.5]{$x_k^{-1}$};
\draw[mor] (1,1) rectangle (2,1.5) node[black,pos=0.5]{$d$};
\node[scale=0.7] at (1.5,0.25) {$f$};
\draw[usual] (0.75,0.5) to (1.5,0.75) to (1.5,1);
\draw[usual] (2.25,0.5) to (1.5,0.75);
\draw[usual] (0.75,0) to (1.5,-0.25) to (1.5,-0.5);
\draw[usual] (2.25,0) to (1.5,-0.25);
\end{tikzpicture}
\; = \;
\begin{tikzpicture}[anchorbase,scale=1]
\draw[mor] (1,-1) rectangle (2,-0.5) node[black,pos=0.5]{$d$};
\draw[mor] (0,0) rectangle (1,0.5) node[black,pos=0.5]{$d$};
\draw[mor] (2,0) rectangle (3,0.5) node[black,pos=0.5]{$d$};
\draw[mor] (1,1) rectangle (2,1.5) node[black,pos=0.5]{$d$};
\node[scale=0.6] at (1.5,0.25) {$\ptr_{d,\un{t}_{\le k}}(f)$};
\draw[usual] (0.75,0.5) to (1.5,0.75) to (1.5,1);
\draw[usual] (2.25,0.5) to (1.5,0.75);
\draw[usual] (0.75,0) to (1.5,-0.25) to (1.5,-0.5);
\draw[usual] (2.25,0) to (1.5,-0.25);
\end{tikzpicture}
.
\end{equation}
The base case is once again obvious.

Let us rewrite the right-hand side of \eqref{analyzeme2}, using the horizontal flip of \autoref{Thm:inductiveclasp} to expand $e_{x_k^{-1}}$ as

\begin{equation}
\begin{tikzpicture}[anchorbase,scale=1]
\draw[mor] (0,0) rectangle (1,0.5) node[black,pos=0.5]{$x_k$};
\draw[mor] (2,0) rectangle (3,0.5) node[black,pos=0.5]{$x_k^{-1}$};
\draw[mor] (1,1.5) rectangle (2,2) node[black,pos=0.5]{$d$};
\draw[usual] (0.5,0.5) to (1.5,1) to (1.5,1.5);
\draw[usual] (2.5,0.5) to (1.5,1);
\begin{scope}[yscale=-1,shift={(0,-0.5)}]
\draw[mor] (1,1.5) rectangle (2,2) node[black,pos=0.5]{$d$};
\draw[usual] (0.5,0.5) to (1.5,1) to (1.5,1.5);
\draw[usual] (2.5,0.5) to (1.5,1);
\end{scope}
\node[scale=0.8] at (1.5,0.25) {$f$};
\end{tikzpicture} 
=
\begin{tikzpicture}[anchorbase,scale=1]
\draw[mor] (0,0) rectangle (1.1,0.5) node[black,pos=0.5]{\raisebox{-0.06cm}{$x_k$}};
\draw[mor] (1.7,0) rectangle (2.8,0.5) node[black,pos=0.5]{\raisebox{-0.06cm}{$x_k^{-1}$}};
\draw[usual] (1.8,0.5)to[out=90,in=0] (1.4,0.9) to[out=180,in=90] (1,0.5);
\draw[usual] (2,0.5)to[out=90,in=0] (1.4,1.1) to[out=180,in=90] (0.8,0.5);
\draw[usual] (2.2,0.5)to[out=90,in=0] (1.4,1.3) to[out=180,in=90] (0.6,0.5);
\draw[mor] (0,0.5) rectangle (0.4,1) node[black,pos=0.5]{\raisebox{-0.06cm}{$d$}};
\draw[mor] (2.4,0.5) rectangle (2.8,1) node[black,pos=0.5]{\raisebox{-0.06cm}{$d$}};
\draw[usual] (0.2,1) to (1.4,2) to (1.4,2.5);
\draw[usual] (2.6,1) to (1.4,2) to (1.4,2.5);
\draw[mor] (1.2,2.5) rectangle (1.6,3) node[black,pos=0.5]{\raisebox{-0.06cm}{$d$}};
\begin{scope}[yscale=-1,shift={(0,-0.5)}]
\draw[usual] (1.8,0.5)to[out=90,in=0] (1.4,0.9) to[out=180,in=90] (1,0.5);
\draw[usual] (2,0.5)to[out=90,in=0] (1.4,1.1) to[out=180,in=90] (0.8,0.5);
\draw[usual] (2.2,0.5)to[out=90,in=0] (1.4,1.3) to[out=180,in=90] (0.6,0.5);
\draw[mor] (0,0.5) rectangle (0.4,1) node[black,pos=0.5]{\raisebox{-0.06cm}{$d$}};
\draw[mor] (2.4,0.5) rectangle (2.8,1) node[black,pos=0.5]{\raisebox{-0.06cm}{$d$}};
\draw[usual] (0.2,1) to (1.4,2) to (1.4,2.5);
\draw[usual] (2.6,1) to (1.4,2) to (1.4,2.5);
\draw[mor] (1.2,2.5) rectangle (1.6,3) node[black,pos=0.5]{\raisebox{-0.06cm}{$d$}};
\end{scope}
\node[scale=0.8] at (1.4,0.25) {$f$};
\end{tikzpicture} 
=
\begin{tikzpicture}[anchorbase,scale=1]
\draw[mor] (0,0) rectangle (1.1,0.5) node[black,pos=0.5]{\raisebox{-0.06cm}{$x_k$}};
\draw[mor] (1.9,0) rectangle (2.8,0.5) node[black,pos=0.5]{\raisebox{-0.06cm}{$x_{k-1}^{-1}$}};
\draw[usual] (1.8,0) to (1.8,0.5)to[out=90,in=0] (1.4,0.9) to[out=180,in=90] (1,0.5);
\draw[usual] (2,0.5)to[out=90,in=0] (1.4,1.1) to[out=180,in=90] (0.8,0.5);
\draw[usual] (2.2,0.5)to[out=90,in=0] (1.4,1.3) to[out=180,in=90] (0.6,0.5);
\draw[mor] (0,0.5) rectangle (0.4,1) node[black,pos=0.5]{\raisebox{-0.06cm}{$d$}};
\draw[mor] (2.4,0.5) rectangle (2.8,1) node[black,pos=0.5]{\raisebox{-0.06cm}{$d$}};
\draw[usual] (0.2,1) to (1.4,2) to (1.4,2.5);
\draw[usual] (2.6,1) to (1.4,2) to (1.4,2.5);
\draw[mor] (1.2,2.5) rectangle (1.6,3) node[black,pos=0.5]{\raisebox{-0.06cm}{$d$}};
\begin{scope}[yscale=-1,shift={(0,-0.5)}]
\draw[usual] (1.8,0.5)to[out=90,in=0] (1.4,0.9) to[out=180,in=90] (1,0.5);
\draw[usual] (2,0.5)to[out=90,in=0] (1.4,1.1) to[out=180,in=90] (0.8,0.5);
\draw[usual] (2.2,0.5)to[out=90,in=0] (1.4,1.3) to[out=180,in=90] (0.6,0.5);
\draw[mor] (0,0.5) rectangle (0.4,1) node[black,pos=0.5]{\raisebox{-0.06cm}{$d$}};
\draw[mor] (2.4,0.5) rectangle (2.8,1) node[black,pos=0.5]{\raisebox{-0.06cm}{$d$}};
\draw[usual] (0.2,1) to (1.4,2) to (1.4,2.5);
\draw[usual] (2.6,1) to (1.4,2) to (1.4,2.5);
\draw[mor] (1.2,2.5) rectangle (1.6,3) node[black,pos=0.5]{\raisebox{-0.06cm}{$d$}};
\end{scope}
\node[scale=0.8] at (1.4,0.25) {$f$};
\end{tikzpicture} 
+
\sum
\begin{tikzpicture}[anchorbase,scale=1]
\draw[mor] (0,0) rectangle (1.1,0.5) node[black,pos=0.5]{\raisebox{-0.06cm}{$x_k$}};
\draw[mor] (1.7,0) rectangle (2.8,0.5) node[black,pos=0.5]{\raisebox{-0.06cm}{\reflectbox{$q_j$}}};
\draw[usual] (1.8,0.5)to[out=90,in=0] (1.4,0.9) to[out=180,in=90] (1,0.5);
\draw[usual] (2,0.5)to[out=90,in=0] (1.4,1.1) to[out=180,in=90] (0.8,0.5);
\draw[usual] (2.2,0.5)to[out=90,in=0] (1.4,1.3) to[out=180,in=90] (0.6,0.5);
\draw[mor] (0,0.5) rectangle (0.4,1) node[black,pos=0.5]{\raisebox{-0.06cm}{$d$}};
\draw[mor] (2.4,0.5) rectangle (2.8,1) node[black,pos=0.5]{\raisebox{-0.06cm}{$d$}};
\draw[usual] (0.2,1) to (1.4,2) to (1.4,2.5);
\draw[usual] (2.6,1) to (1.4,2) to (1.4,2.5);
\draw[mor] (1.2,2.5) rectangle (1.6,3) node[black,pos=0.5]{\raisebox{-0.06cm}{$d$}};
\begin{scope}[yscale=-1,shift={(0,-0.5)}]
\draw[usual] (1.8,0.5)to[out=90,in=0] (1.4,0.9) to[out=180,in=90] (1,0.5);
\draw[usual] (2,0.5)to[out=90,in=0] (1.4,1.1) to[out=180,in=90] (0.8,0.5);
\draw[usual] (2.2,0.5)to[out=90,in=0] (1.4,1.3) to[out=180,in=90] (0.6,0.5);
\draw[mor] (0,0.5) rectangle (0.4,1) node[black,pos=0.5]{\raisebox{-0.06cm}{$d$}};
\draw[mor] (2.4,0.5) rectangle (2.8,1) node[black,pos=0.5]{\raisebox{-0.06cm}{$d$}};
\draw[usual] (0.2,1) to (1.4,2) to (1.4,2.5);
\draw[usual] (2.6,1) to (1.4,2) to (1.4,2.5);
\draw[mor] (1.2,2.5) rectangle (1.6,3) node[black,pos=0.5]{\raisebox{-0.06cm}{$d$}};
\end{scope}
\node[scale=0.8] at (1.4,0.25) {$f$};
\end{tikzpicture} \; .
\end{equation}

Here we have written \reflectbox{$q_j$} for the horizontal reflection of $q_j$. Again, any term with \reflectbox{$q_j$} factors through a morphism $B_d \to B_{x_k} \ot B_y$ of degree $-a$, where $y \muless x_{k-1}^{-1}$ and $t_k y < y$. By similar arguments to the previous proof, either $y$ is in strictly lower terms, or $\gamma_{x_k,y}^d = 0$, and either way we deduce that this term contributes zero to the final answer. Meanwhile, the first term on the right-hand side involves a partial trace of $f$, so

\begin{equation} 	\begin{tikzpicture}[anchorbase,scale=1]
\draw[mor] (0,0) rectangle (1,0.5) node[black,pos=0.5]{$x_k$};
\draw[mor] (2,0) rectangle (3,0.5) node[black,pos=0.5]{$x_k^{-1}$};
\draw[mor] (1,1.5) rectangle (2,2) node[black,pos=0.5]{$d$};
\draw[usual] (0.5,0.5) to (1.5,1) to (1.5,1.5);
\draw[usual] (2.5,0.5) to (1.5,1);
\begin{scope}[yscale=-1,shift={(0,-0.5)}]
\draw[mor] (1,1.5) rectangle (2,2) node[black,pos=0.5]{$d$};
\draw[usual] (0.5,0.5) to (1.5,1) to (1.5,1.5);
\draw[usual] (2.5,0.5) to (1.5,1);
\end{scope}
\node[scale=0.8] at (1.5,0.25) {$f$};
\end{tikzpicture}  \; = \; 
\begin{tikzpicture}[anchorbase,scale=1]
\draw[mor] (0,0) rectangle (1.1,0.5) node[black,pos=0.5]{\raisebox{-0.06cm}{$x_{k-1}$}};
\draw[mor] (1.7,0) rectangle (2.8,0.5) node[black,pos=0.5]{\raisebox{-0.06cm}{$x_{k-1}^{-1}$}};
\draw[usual] (1.8,0.5)to[out=90,in=0] (1.4,0.9) to[out=180,in=90] (1,0.5);
\draw[usual] (2,0.5)to[out=90,in=0] (1.4,1.1) to[out=180,in=90] (0.8,0.5);
\draw[usual] (2.2,0.5)to[out=90,in=0] (1.4,1.3) to[out=180,in=90] (0.6,0.5);
\draw[mor] (0,0.5) rectangle (0.4,1) node[black,pos=0.5]{\raisebox{-0.06cm}{$d$}};
\draw[mor] (2.4,0.5) rectangle (2.8,1) node[black,pos=0.5]{\raisebox{-0.06cm}{$d$}};
\draw[usual] (0.2,1) to (1.4,2) to (1.4,2.5);
\draw[usual] (2.6,1) to (1.4,2) to (1.4,2.5);
\draw[mor] (1.2,2.5) rectangle (1.6,3) node[black,pos=0.5]{\raisebox{-0.06cm}{$d$}};
\begin{scope}[yscale=-1,shift={(0,-0.5)}]
\draw[usual] (1.8,0.5)to[out=90,in=0] (1.4,0.9) to[out=180,in=90] (1,0.5);
\draw[usual] (2,0.5)to[out=90,in=0] (1.4,1.1) to[out=180,in=90] (0.8,0.5);
\draw[usual] (2.2,0.5)to[out=90,in=0] (1.4,1.3) to[out=180,in=90] (0.6,0.5);
\draw[mor] (0,0.5) rectangle (0.4,1) node[black,pos=0.5]{\raisebox{-0.06cm}{$d$}};
\draw[mor] (2.4,0.5) rectangle (2.8,1) node[black,pos=0.5]{\raisebox{-0.06cm}{$d$}};
\draw[usual] (0.2,1) to (1.4,2) to (1.4,2.5);
\draw[usual] (2.6,1) to (1.4,2) to (1.4,2.5);
\draw[mor] (1.2,2.5) rectangle (1.6,3) node[black,pos=0.5]{\raisebox{-0.06cm}{$d$}};
\end{scope}
\node[scale=0.8] at (1.4,0.25) {$\mathrm{tr}(f)$};
\end{tikzpicture}
\; . \end{equation}
This easily yields the inductive step.
\end{proof}

\begin{Corollary} \label{Cor:catdimunderassumption} Under the assumption \eqref{assumeme}, categorical dimension of $A_w$ is equal to the ratio
\begin{equation} \lambda_d(\ptr_{d,\un{t}}(\rho^a)) / \lambda_d(\rho^a). \end{equation} \end{Corollary}

\subsection{Parabolic Cells}\label{subsection:simplifications2}

For this section we assume $d = w_I$ is the longest element of a parabolic subgroup attached to $I \subset S$, and that $d$ lives in two-sided cell $\mathcal{J}$ and diagonal cell $\mathcal{H}$. We call $\mathcal{H}$ a \emph{parabolic cell}. As in \autoref{longesteltworks}, \eqref{assumeme} holds for any element $w \in \mathcal{H}$, so our discussion will continue the discussion of the previous section. One other special feature is that $\len(d) = a(d)$, so that we can hope to normalize $\incl_d$ so that it preserves the one-tensor.

Throughout this section let us fix a reduced expression $\un{d}$ for $d$. Another special feature of longest elements is that every simple reflection appearing in $\un{d}$ is
also in the right (or left) descent set of $d$. This allows us to construct $\incl_d$ using $(d,i)$-trivalent vertices for various $i \in I$.

\begin{Lemma} When $d = w_I$ for $I \subset S$, one can construct a morphism $\incl_d \colon B_d \to B_d \ot B_d$ of degree $-a = - \len(d)$ using $(d,i)$-trivalent vertices for various $i \in I$, as in the following diagram. There are two ways to do this for each reduced expression $\un{d} = (i_1, \ldots, i_a)$ of $d$, and they are all equal and nonzero.

\begin{equation}
\begin{tikzpicture}[anchorbase,scale=1,yscale=-1]
\draw[mor] (0,0) rectangle (1,0.5) node[black,pos=0.5]{$d$};
\draw[mor] (0,1) rectangle (1,1.5) node[black,pos=0.5]{$d$};
\draw[mor] (0,1.5) rectangle (1,2) node[black,pos=0.5]{$d$};
\node[scale=0.8] at (0.5,2.25) {$\vdots$};
\draw[mor] (0,2.5) rectangle (1,3) node[black,pos=0.5]{$d$};
\draw[mor] (0,3.5) rectangle (1,4) node[black,pos=0.5]{$d$};
\draw[mor] (1.5,0) rectangle (2.5,0.5) node[black,pos=0.5]{$d$};
\draw[usual] (0.5,0.5) to (0.5,1);
\draw[usual] (0.5,3) to (0.5,3.5);
\draw[usual] (1.6,0.5) to (1.6,1.25) to (1,1.25);
\node[scale=0.7] at (1.25,1.6) {$i_2$};
\node[scale=0.7] at (1.25,1.1) {$i_1$};
\node[scale=0.7] at (1.25,2.6) {$i_a$};
\draw[usual] (1.8,0.5) to (1.8,1.75) to (1,1.75);
\node[scale=0.7] at (2.1,0.75) {$\hdots$};
\draw[usual] (2.4,0.5) to (2.4,2.75) to (1,2.75);
\end{tikzpicture}
=
\begin{tikzpicture}[anchorbase,scale=1,yscale=-1]
\begin{scope}[xscale=-1]
\draw[mor] (0,0) rectangle (1,0.5) node[black,pos=0.5]{$d$};
\draw[mor] (0,1) rectangle (1,1.5) node[black,pos=0.5]{$d$};
\draw[mor] (0,1.5) rectangle (1,2) node[black,pos=0.5]{$d$};
\node[scale=0.8] at (0.5,2.25) {$\vdots$};
\draw[mor] (0,2.5) rectangle (1,3) node[black,pos=0.5]{$d$};
\draw[mor] (0,3.5) rectangle (1,4) node[black,pos=0.5]{$d$};
\draw[mor] (1.5,0) rectangle (2.5,0.5) node[black,pos=0.5]{$d$};
\draw[usual] (0.5,0.5) to (0.5,1);
\draw[usual] (0.5,3) to (0.5,3.5);
\draw[usual] (1.6,0.5) to (1.6,1.25) to (1,1.25);
\draw[usual] (1.8,0.5) to (1.8,1.75) to (1,1.75);
\node[scale=0.7] at (2.1,0.75) {$\hdots$};
\draw[usual] (2.4,0.5) to (2.4,2.75) to (1,2.75);
\node[scale=0.7] at (1.25,1.6) {$i_{a-1}$};
\node[scale=0.7] at (1.25,1.1) {$i_a$};
\node[scale=0.7] at (1.25,2.6) {$i_1$};
\end{scope}
\end{tikzpicture}\;.
\end{equation}
\end{Lemma}

\begin{proof} Clearly the diagrams above are well-defined morphisms of the appropriate degree. Their equality follows because both morphisms preserves the one-tensor. This also implies that the morphisms are nonzero. \end{proof}

\begin{Corollary} We have $\ep_d = 1$. \end{Corollary}

\begin{proof} The map $\incl_d$ defined above is invariant under the horizontal flip. \end{proof}

Now we simplify these trivalent vertices, to ease the calculations to come.

\begin{Lemma} \label{Lem:removed} As morphisms from $B_d \ot B_{\un{d}} \to B_d$, the following are equal.

\begin{equation}
\begin{tikzpicture}[anchorbase,scale=1]
\draw[mor] (0,0) rectangle (1,0.5) node[black,pos=0.5]{$d$};
\draw[mor] (0,1) rectangle (1,1.5) node[black,pos=0.5]{$d$};
\draw[mor] (0,1.5) rectangle (1,2) node[black,pos=0.5]{$d$};
\node[scale=0.8] at (0.5,2.25) {$\vdots$};
\draw[mor] (0,2.5) rectangle (1,3) node[black,pos=0.5]{$d$};
\draw[mor] (0,3.5) rectangle (1,4) node[black,pos=0.5]{$d$};
\draw[mor] (1.5,0) rectangle (2.5,0.5) node[black,pos=0.5]{$d$};
\draw[usual] (0.5,0.5) to (0.5,1);
\draw[usual] (0.5,3) to (0.5,3.5);
\draw[usual] (1.6,0.5) to (1.6,1.25) to (1,1.25);
\draw[usual] (1.8,0.5) to (1.8,1.75) to (1,1.75);
\node[scale=0.7] at (2.1,0.75) {$\hdots$};
\draw[usual] (2.4,0.5) to (2.4,2.75) to (1,2.75);
\end{tikzpicture}
=
\begin{tikzpicture}[anchorbase,scale=1]
\draw[mor] (0,0) rectangle (1,0.5) node[black,pos=0.5]{$d$};
\draw[mor] (0,1) rectangle (1,1.5) node[black,pos=0.5]{$d$};
\draw[mor] (0,1.5) rectangle (1,2) node[black,pos=0.5]{$d$};
\node[scale=0.8] at (0.5,2.25) {$\vdots$};
\draw[mor] (0,2.5) rectangle (1,3) node[black,pos=0.5]{$d$};
\draw[mor] (0,3.5) rectangle (1,4) node[black,pos=0.5]{$d$};
\draw[usual] (0.5,0.5) to (0.5,1);
\draw[usual] (0.5,3) to (0.5,3.5);
\draw[usual] (1.6,0) to (1.6,1.25) to (1,1.25);
\draw[usual] (1.8,0) to (1.8,1.75) to (1,1.75);
\node[scale=0.7] at (2.1,0.25) {$\hdots$};
\draw[usual] (2.4,0) to (2.4,2.75) to (1,2.75);
\end{tikzpicture}\;.
\end{equation}
\end{Lemma}

The proof is similar to those of \autoref{thm:lambdadlambdaw} and \autoref{thm:itsapartialtracedummy}, but replacing arguments involving $\catstuff{I}_{< \mathcal{J}}$ with a separate low degree vanishing argument.

\begin{proof} First we make some observations about longest elements in the Hecke algebra. The element $b_{w_I}$ is in the center of the Hecke algebra of $W_I$, and spans a one-dimensional ideal. For each $i \in I$ we have the well-known statement
\begin{equation} b_i b_{w_I} = [2] b_{w_I}, \end{equation}
where $[2] = v + v^{-1}$, which is categorified by
\begin{equation} B_i B_{w_I} \cong B_{w_I} B_i \cong B_{w_i}^{\oplus [2]}. \end{equation}
Similarly, if $\un{x}$ is any expression in $I$ of length $k$ then
\begin{equation} B_{\un{x}} B_{w_I} \cong B_{w_I} B_{\un{x}} \cong B_{w_I}^{\oplus [2]^k}. \end{equation}
The minimal degree of a morphism in $\Hom(B_{w_I} B_{\un{x}}, B_{w_I})$ is $-k$.

Meanwhile, for the standard generator $h_i$ (with $i \in I$) of the Hecke algebra we have
\begin{equation} h_i b_{w_I} = v^{-1} b_{w_I}. \end{equation}
As a consequence, for any $y \in W_I$ we have
\begin{equation} h_y b_{w_I} = b_{w_I} h_y = v^{-\len(y)} b_{w_I}. \end{equation}
Using the degree bounds on Kazhdan--Lusztig polynomials, we deduce that
\begin{equation} b_y b_{w_I} = b_{w_I} b_y = \mu_y b_{w_I} \end{equation}
where $\mu_y$ has exponents bounded above by $v^{\len(y)}$ and below by $v^{-\len(y)}$. Categorified we have
\begin{equation} B_{w_I} B_y \cong B_{w_I}^{\oplus \mu_y}, \end{equation}
and the minimal degree of a morphism in $\Hom(B_{w_I} B_y, B_{w_I})$ is $-\len(y)$.

One can combine these statements: the minimal degree of a morphism in $\Hom(B_{w_I} B_y B_{\un{x}}, B_{w_I})$ is $-(\len(y) + \len(\un{x}))$. That any morphism of even more negative degree must be zero we call \emph{low degree vanishing}.

Now we choose a reduced expression $\un{d} = (t_1, t_2, \ldots, t_a)$ for $d = w_I$, and we let $x_k$ be expressed by $(t_1, \ldots, t_k)$. We begin to expand the idempotent $e_d$ using \autoref{Thm:inductiveclasp}. Each diagram with $q_j$ factors through a degree $-a$ map $B_d \ot B_y \to B_d$, where $B_y$ is a non-top summand of $B_{x_{a-1}} \ot B_{t_k}$.  Any such map is zero by low degree vanishing, since $y < x_{a-1}$ so $y \in W_I$ and $\len(y) < \len(w_I) = a$. So we can replace $e_d$ with $e_{x_{a-1}} \ot \id_{t_a}$.

Now continue, rewriting $e_{x_{a-1}}$ using \autoref{Thm:inductiveclasp}. This time, each diagram with $q_j$ factors through a degree $-a$ map $B_d \ot B_y \ot B_{t_k} \to B_d$, where $B_y$ is a non-top summand of $B_{x_{a-2}} \ot B_{t_{a-1}}$. Thus, $y \in W_I$ and $\len(y) < \len(x_{a-2}) = a-2$. By low degree vanishing, this diagram is zero. Thus, we can replace $e_{x_{a-1}} \ot \id_{t_k}$ with $e_{x_{a-2}} \ot \id_{t_{k-1}} \ot \id_{t_k}$.

Repeating, we eventually replace $e_d$ with the identity map of the Bott--Samelson object $B_{\un{d}}$. \end{proof}

Recall that there is a Demazure operator $\partial_x$ for any $x \in W$, obtained by composing the simple Demazure operators $\partial_i$ along a reduced expression for $x$.

\begin{Theorem} \label{Thm:lambdaispartial} We have $\lambda_d(f) = \partial_d(f)$. \end{Theorem}

\begin{proof} Let $i$ be the first color in our chosen reduced expression for $d$.
\begin{equation}
\begin{tikzpicture}[anchorbase,scale=0.8]
\draw[mor] (0,0) rectangle (1,0.5) node[black,pos=0.5]{$d$};
\draw[mor] (2,0) rectangle (3,0.5) node[black,pos=0.5]{$d$};
\draw[mor] (1,1.5) rectangle (2,2) node[black,pos=0.5]{$d$};
\draw[usual] (0.5,0.5) to (1.5,1) to (1.5,1.5);
\draw[usual] (2.5,0.5) to (1.5,1);
\begin{scope}[yscale=-1,shift={(0,-0.5)}]
\draw[mor] (1,1.5) rectangle (2,2) node[black,pos=0.5]{$d$};
\draw[usual] (0.5,0.5) to (1.5,1) to (1.5,1.5);
\draw[usual] (2.5,0.5) to (1.5,1);
\end{scope}
\node[scale=0.8] at (1.5,0.25) {$f$};
\end{tikzpicture} 
=
\begin{tikzpicture}[anchorbase,scale=0.8]
\draw[mor] (0,0) rectangle (1,0.5) node[black,pos=0.5]{$d$};
\draw[mor] (0,1) rectangle (1,1.5) node[black,pos=0.5]{$d$};
\draw[mor] (0,1.5) rectangle (1,2) node[black,pos=0.5]{$d$};
\node[scale=0.8] at (0.5,2.25) {$\vdots$};
\draw[mor] (0,2.5) rectangle (1,3) node[black,pos=0.5]{$d$};
\draw[mor] (0,3.5) rectangle (1,4) node[black,pos=0.5]{$d$};
\draw[mor] (1.5,0) rectangle (2.5,0.5) node[black,pos=0.5]{$d$};
\draw[usual] (0.5,0.5) to (0.5,1);
\draw[usual] (0.5,3) to (0.5,3.5);
\draw[usual] (1.6,0.5) to (1.6,1.25) to (1,1.25);
\draw[usual] (1.8,0.5) to (1.8,1.75) to (1,1.75);
\node[scale=0.7] at (2.1,0.75) {$\hdots$};
\draw[usual] (2.4,0.5) to (2.4,2.75) to (1,2.75);
\begin{scope}[yscale=-1,shift={(0,-0.5)}]
\draw[mor] (0,1) rectangle (1,1.5) node[black,pos=0.5]{$d$};
\draw[mor] (0,1.5) rectangle (1,2) node[black,pos=0.5]{$d$};
\node[scale=0.8] at (0.5,2.25) {$\vdots$};
\draw[mor] (0,2.5) rectangle (1,3) node[black,pos=0.5]{$d$};
\draw[mor] (0,3.5) rectangle (1,4) node[black,pos=0.5]{$d$};
\draw[usual] (0.5,0.5) to (0.5,1);
\draw[usual] (0.5,3) to (0.5,3.5);
\draw[usual] (1.6,0.5) to (1.6,1.25) to (1,1.25);
\draw[usual] (1.8,0.5) to (1.8,1.75) to (1,1.75);
\node[scale=0.7] at (2.1,0.75) {$\hdots$};
\draw[usual] (2.4,0.5) to (2.4,2.75) to (1,2.75);
\end{scope}
\node[scale=0.8] at (1.25,0.25) {$f$};
\end{tikzpicture}
=
\begin{tikzpicture}[anchorbase,scale=0.8]
\draw[mor] (0,0) rectangle (1,0.5) node[black,pos=0.5]{$d$};
\draw[mor] (0,0.5) rectangle (1,1) node[black,pos=0.5]{$d$};
\draw[mor] (0,1.5) rectangle (1,2) node[black,pos=0.5]{$d$};
\node[scale=0.8] at (0.5,2.25) {$\vdots$};
\draw[mor] (0,2.5) rectangle (1,3) node[black,pos=0.5]{$d$};
\draw[mor] (0,3.5) rectangle (1,4) node[black,pos=0.5]{$d$};
\draw[usual] (0.5,1) to (0.5,1.5);
\draw[usual] (0.5,3) to (0.5,3.5);
\draw[usual] (1.6,0) to (1.6,0.75) to (1,0.75);
\draw[usual] (1.8,0) to (1.8,1.75) to (1,1.75);
\node[scale=0.7] at (2.1,0.75) {$\hdots$};
\draw[usual] (2.4,0) to (2.4,2.75) to (1,2.75);
\begin{scope}[yscale=-1,shift={(0,-0.5)}]
\draw[mor] (0,0.5) rectangle (1,1) node[black,pos=0.5]{$d$};
\draw[mor] (0,1.5) rectangle (1,2) node[black,pos=0.5]{$d$};
\node[scale=0.8] at (0.5,2.25) {$\vdots$};
\draw[mor] (0,2.5) rectangle (1,3) node[black,pos=0.5]{$d$};
\draw[mor] (0,3.5) rectangle (1,4) node[black,pos=0.5]{$d$};
\draw[usual] (0.5,1) to (0.5,1.5);
\draw[usual] (0.5,3) to (0.5,3.5);
\draw[usual] (1.6,0.5) to (1.6,0.75) to (1,0.75);
\draw[usual] (1.8,0.5) to (1.8,1.75) to (1,1.75);
\node[scale=0.7] at (2.1,0.75) {$\hdots$};
\draw[usual] (2.4,0.5) to (2.4,2.75) to (1,2.75);
\end{scope}
\node[scale=0.8] at (1.25,0.25) {$f$};
\end{tikzpicture}
=
\begin{tikzpicture}[anchorbase,scale=0.8]
\draw[mor] (0,0) rectangle (1,0.5) node[black,pos=0.5]{$d$};
\draw[mor] (0,0.5) rectangle (1,1) node[black,pos=0.5]{$d$};
\draw[mor] (0,1.5) rectangle (1,2) node[black,pos=0.5]{$d$};
\node[scale=0.8] at (0.5,2.25) {$\vdots$};
\draw[mor] (0,2.5) rectangle (1,3) node[black,pos=0.5]{$d$};
\draw[mor] (0,3.5) rectangle (1,4) node[black,pos=0.5]{$d$};
\draw[usual] (0.5,1) to (0.5,1.5);
\draw[usual] (0.5,3) to (0.5,3.5);
\draw[usual,marked=0] (1.3,0.75) to (1,0.75);
\draw[usual] (1.8,0) to (1.8,1.75) to (1,1.75);
\node[scale=0.7] at (2.1,0.75) {$\hdots$};
\draw[usual] (2.4,0) to (2.4,2.75) to (1,2.75);
\begin{scope}[yscale=-1,shift={(0,-0.5)}]
\draw[mor] (0,0.5) rectangle (1,1) node[black,pos=0.5]{$d$};
\draw[mor] (0,1.5) rectangle (1,2) node[black,pos=0.5]{$d$};
\node[scale=0.8] at (0.5,2.25) {$\vdots$};
\draw[mor] (0,2.5) rectangle (1,3) node[black,pos=0.5]{$d$};
\draw[mor] (0,3.5) rectangle (1,4) node[black,pos=0.5]{$d$};
\draw[usual] (0.5,1) to (0.5,1.5);
\draw[usual] (0.5,3) to (0.5,3.5);
\draw[usual,marked=0] (1.3,0.75) to (1,0.75);
\draw[usual] (1.8,0.5) to (1.8,1.75) to (1,1.75);
\node[scale=0.7] at (2.1,0.75) {$\hdots$};
\draw[usual] (2.4,0.5) to (2.4,2.75) to (1,2.75);
\end{scope}
\node[scale=0.8] at (1.4,0.25) {$\partial_i(f)$};
\end{tikzpicture}
=\ldots
=
\begin{tikzpicture}[anchorbase,scale=0.8]
\draw[mor] (0,0) rectangle (1,0.5) node[black,pos=0.5]{$d$};
\node[scale=0.8] at (1.5,0.25) {$\partial_d(f)$};
\end{tikzpicture}\;.
\end{equation} The second equality used \autoref{Lem:removed}.  Following that, we use polynomial forcing to break the innermost loop, and then absorb the dots using \eqref{triunit}; a similar way to perform this same calculation is to follow the last few equalities in the proof of \autoref{lemma:AtoB}. Thus, the innermost loop is removed, at the cost of replacing $f$ with $\partial_i(f)$. Repeating, we apply Demazure operators to $f$ along a reduced expression of $d$, whence the result. \end{proof}

Thus, for diagonal cells containing longest elements of parabolic subgroups, combining this theorem with \autoref{Cor:catdimunderassumption}, we can compute categorical dimensions
whenever we can compute partial traces. If the branching graph from $d$ to $w$ is linear, we can use \autoref{thm:linearpartialtracerecursionbetter}, and completely automate this
computation.

\begin{Example} This process is what justifies the claims made in \autoref{Ex:F4}. \end{Example}

\begin{Example} We continue \autoref{example:smoothH3}. We have $\partial_d = \partial_2 \partial_3 \partial_2$ and one can compute that $\partial_d(\rho^3) = 6$ for the standard regular element. If one sets $f_{3,0} = \rho^3$, then one can compute that $\partial_d(f_{0,0}) = -6$. Hence $\dim(\obstuff{A}_{232123})=-1$.
\end{Example}

\begin{Example} We now give some examples in type $H_4$, where $m_{12} = 5$ and $m_{23} = m_{34} = 3$. We assume that $[2]_{1,2} = [2]_{2,1} = \phi$ is the golden ratio. Let $\un{d} = 31$, a reduced expression for the longest element of a parabolic subgroup of type $A_1 \times A_1$. Let $\un{t} = 2143$ and $\underline{x} = \un{d}\un{t}$, a reduced expression for $x$ which is in the same diagonal cell as $d$. The graph $\something(\underline{x_{2}})$ is linear:
\begin{gather}\label{eq:h4a2-pgraph}
\Gamma(\underline{d}2143)=
\begin{tikzcd}[ampersand replacement=\&,row sep=scriptsize,column sep=scriptsize,arrows={shorten >=-0.5ex,shorten <=-0.5ex},labels={inner sep=0.05ex},arrow style=tikz]
\emptyset\ar[r,yshift=0.1cm,soergelthree]
\&
3\ar[r,yshift=0.1cm,soergelone]
\&
31=\underline{d}\ar[r,yshift=0.1cm,soergeltwo]
\&
\underline{d}2\ar[r,yshift=0.1cm,soergelone]
\ar[l,yshift=-0.1cm,soergelone]
\& 
\underline{d}21\ar[r,yshift=0.1cm,soergelfour]
\& 
\underline{d}214\ar[r,yshift=0.1cm,soergelthree]
\& 
\underline{d}2143
\end{tikzcd}.
\end{gather}
The only nontrivial LIF is $\partial_1(\alpha_2) = - \phi$. One can compute partial traces algorithmically, and one can verify that
\begin{equation} \ptr_{d,\un{t}}(f) = \alpha_2 \partial_2(\alpha_1 \partial_1(\alpha_4 \partial_4(\alpha_3 \partial_3(f))))+ \frac{1}{\phi} \alpha_2 \partial_1(\alpha_4 \partial_4(\alpha_3 \partial_3(f))). \end{equation}
Plugging in $f = \rho^2$ we get
\begin{equation} \partial_{31}(\ptr_{d,\un{t}}(\rho^2)) = 2 \phi, \qquad \partial_{31}(\rho^2) = 2, \end{equation}
from which we deduce that $\dim(\obstuff{A}_{x})=\phi$. \end{Example}

\begin{Example} Still in type $H_4$ as above, let $\un{d} = 232432$, a reduced expression for the longest element of a parabolic subgroup of type $A_3$. Let $\un{t} = 1234$ and $\underline{x}_2 = \un{d}\un{t}$, a reduced expression for $x_2$ which is one of fourteen elements in the same diagonal cell as $d$. The graph $\something(\underline{x_{2}})$ is linear:
\begin{gather*}
\begin{tikzcd}[ampersand replacement=\&,row sep=scriptsize,column sep=scriptsize,arrows={shorten >=-0.5ex,shorten <=-0.5ex},labels={inner sep=0.05ex},arrow style=tikz]
\emptyset\ar[r,yshift=0.1cm,soergeltwo]
\&
2 \ar[r,yshift=0.1cm,soergelthree]
\& 
23\ar[r,yshift=0.1cm,soergeltwo]\ar[l,yshift=-0.1cm,soergeltwo]
\&
232 \ar[r,yshift=0.1cm,soergelfour]
\&
2324 \ar[r,yshift=0.1cm,soergelthree] \ar[l,yshift=-0.1cm,soergelthree]
\&
23243 \ar[r,yshift=0.1cm,soergeltwo] \ar[l,yshift=-0.1cm,soergeltwo]
\&
\underline{d} \ar[r,yshift=0.1cm,soergelone]
\&
\underline{d}1 \ar[r,yshift=0.1cm,soergeltwo]\ar[l,yshift=-0.1cm,soergeltwo]
\&
\underline{d}12 \ar[r,yshift=0.1cm,soergelthree] \ar[l,yshift=-0.1cm,soergelthree]
\&
\underline{d}123 \ar[r,yshift=0.1cm,soergelfour] \ar[l,yshift=-0.1cm,soergelfour]
\&
\underline{x}_{2}
\end{tikzcd}
.
\end{gather*}
One can now compute local intersection forms and partial traces algorithmically, which we did with the aid of a computer. We also used a computer for computations with Demazure operators. We end up with 
\begin{equation} \partial_d(\ptr_{d,\un{t}}(\rho^6)) = 1440\phi, \qquad \partial_d(\rho^6) = 720, \end{equation}
from which we deduce that $\dim(\obstuff{A}_{x_2})=2 \phi$.
\end{Example}



\newcommand{\etalchar}[1]{$^{#1}$}

\end{document}